\documentclass[smallextended]{svjour3}       

\smartqed  
\usepackage{graphicx}
\usepackage{enumerate}
\usepackage[OT1]{fontenc}
\usepackage[usenames]{color}
\usepackage[colorlinks,
linkcolor=red,
anchorcolor=blue,
citecolor=blue
]{hyperref}
\usepackage{mathrsfs}
\usepackage[square,numbers,sort]{natbib}

\usepackage{hyperref}
\usepackage[protrusion=true, expansion=true]{microtype}
\usepackage{float}
\usepackage{subfigure}
\usepackage{amsfonts,amsmath,amssymb,url,xspace}
\usepackage{tikz}
\usepackage{verbatim}
\usetikzlibrary{arrows,shapes}
\usepackage{mathtools}
\usepackage{authblk}
\usepackage[bottom]{footmisc}
\usepackage{enumitem}
\usepackage{leftidx,mathtools, abraces}
\usepackage{lscape}
\usepackage{mathabx}
\usepackage{bm}
\usepackage{soul}

\usepackage{algorithm}
\usepackage{algorithmicx,algpseudocode}

\usepackage{chngcntr}
\counterwithout{equation}{section}

\allowdisplaybreaks[1]

\newcommand{\blue}{\color{blue}}

\newcommand{\mR}{\mathbb{R}}

\newcommand{\mE}{\mathbb{E}}

\newcommand{\E}{\mathcal{E}}
\newcommand{\mL}{\mathcal{L}}

\newcommand{\mM}{\mathcal{M}}
\newcommand{\mX}{\mathcal{X}}
\newcommand{\mF}{\mathcal{F}}
\newcommand{\mA}{\mathcal{A}}
\newcommand{\mB}{\mathcal{B}}
\newcommand{\mT}{\mathcal{T}}

\newcommand{\be}{\bm{e}}

\def\I{\mathcal{I}}

\def\P{{\mathcal P}}

\def\T{{\mathcal T}}
\def\H{{\mathcal H}}
\def\0{{\boldsymbol 0}}
\def\ba{{\boldsymbol{a}}}
\def\bb{{\boldsymbol{b}}}
\def\N{{\mathcal N}}

\newcommand{\mI}{{\mathcal{I}}}

\newcommand{\bw}{\bm{w}}
\newcommand{\bq}{\bm{q}}
\newcommand{\bmu}{\boldsymbol{\mu}}
\def\bb{{\boldsymbol{b}}}
\def\ba{{\boldsymbol{a}}}
\def\by{{\boldsymbol{y}}}
\def\bz{{\boldsymbol{z}}}
\def\bv{{\boldsymbol{v}}}
\def\b1{{\boldsymbol{1}}}
\def\bx{{\boldsymbol{x}}}
\def\bu{{\boldsymbol{u}}}

\def\bl{{\boldsymbol{l}}}
\def\bz{{\boldsymbol{z}}}
\def\blambda{{\boldsymbol{\lambda}}}
\def\tlambda{{\boldsymbol{\lambda}}^\star}

\def\tx{{\boldsymbol{x}}^\star}

\def\bardelta{\bar{\delta}}

\def\tDelta{{\tilde{\Delta}}}
\def\hDelta{{\hat{\Delta}}}
\newcommand{\diag}{{\rm diag}}

\def\barf{{\bar{f}}}

\def\bnabla{{\bar{\nabla}}}

\def\barDelta{{\bar{\Delta}}}
\def\tDelta{{\tilde{\Delta}}}

\def\barepsilon{{\bar{\epsilon}}}

\def\baralpha{{\bar{\alpha}}}
\def\barL{{\bar{\mL}}}

\def\tmu{{\bmu^\star}}

\def\barQ{{\bar{Q}}}
\def\barnu{{\bar{\nu}}}
\def\hatDelta{{\hat{\Delta}}}
\def\hatH{{\hat{H}}}
\def\cDelta{{\widecheck{\Delta}}}

\def\bbnabla{{\bar{\bnabla}}}
\def\barR{{\bar{R}}}
\def\tepsilon{{\tilde{\epsilon}}}
\def\tnu{{\tilde{\nu}}}
\def\bart{{\bar{t}}}

\mathtoolsset{showonlyrefs=true}

\usepackage{xargs}
\usepackage[colorinlistoftodos,prependcaption,textsize=tiny]{todonotes}
\newcommandx{\unsure}[2][1=]{\todo[linecolor=red,backgroundcolor=red!25,bordercolor=red,#1]{#2}}
\newcommandx{\change}[2][1=]{\todo[linecolor=blue,backgroundcolor=blue!25,bordercolor=blue,#1]{#2}}
\newcommandx{\info}[2][1=]{\todo[linecolor=OliveGreen,backgroundcolor=OliveGreen!25,bordercolor=OliveGreen,#1]{#2}}
\newcommandx{\improvement}[2][1=]{\todo[linecolor=Plum,backgroundcolor=Plum!25,bordercolor=Plum,#1]{#2}}

\allowdisplaybreaks[1]

\newcommand{\rbr}[1]{\left(#1\right)}
\newcommand{\sbr}[1]{\left[#1\right]}
\newcommand{\cbr}[1]{\left\{#1\right\}}
\newcommand{\nbr}[1]{\left\|#1\right\|}
\newcommand{\abr}[1]{\left|#1\right|}

\ifx\QED\undefined
\def\QED{~\rule[-1pt]{5pt}{5pt}\par\medskip}
\fi

\begin{document}

\title{Inequality Constrained Stochastic Nonlinear Optimization via Active-Set Sequential Quadratic Programming}

\titlerunning{Inequality Constrained Stochastic Optimization via Active-Set SQP}

\author{Sen Na \and Mihai Anitescu \and Mladen Kolar}

\institute{Sen Na (corresponding author)\at
Department of Statistics, University of California, Berkeley \\
International Computer Science Institute\\
\email{senna@berkeley.edu}
\and
Mihai Anitescu \at
Mathematics and Computer Science Division, Argonne National Laboratory\\
\email{anitescu@mcs.anl.gov}
\and
Mladen Kolar \at
Booth School of Business, The University of Chicago\\
\email{mladen.kolar@chicagobooth.edu}
}

\authorrunning{Na et al.} 
\date{Received: date / Accepted: date}
\maketitle

\begin{abstract}

We study nonlinear optimization problems with a stochastic objective and deterministic equality and inequality constraints, which emerge in numerous applications including finance, manufacturing, power systems and, recently, deep neural networks. We propose an active-set stochastic sequential quadratic programming (StoSQP) algorithm that utilizes a differentiable exact augmented Lagrangian as the merit function. The algorithm adaptively selects the penalty parameters of the augmented Lagrangian, and performs a stochastic line search to decide the stepsize. The global convergence is established: for any initialization, the KKT residuals converge to zero \textit{almost surely}. Our algorithm and analysis further develop the prior work of Na et al. \cite{Na2022adaptive}. Specifically, we allow nonlinear inequality constraints \textit{without} requiring the strict complementary condition; refine some of designs in \cite{Na2022adaptive} such as the feasibility error condition and the monotonically increasing sample size; strengthen the global convergence guarantee; and improve the sample complexity on the objective Hessian. We demonstrate the performance of the designed algorithm on a subset of nonlinear problems collected in CUTEst test set and on constrained logistic regression problems.


\end{abstract}

\section{Introduction}\label{sec:1}

We study stochastic nonlinear optimization problems with deterministic equality and inequality constraints:
\begin{align}\label{pro:1}
\min_{\bx\in \mR^d}\;\; & f(\bx) = \mE[F(\bx; \xi)], \nonumber\\
\text{s.t. }\; & c(\bx) = \0,\\
& g(\bx) \leq \0, \nonumber
\end{align}
where $f: \mR^d\rightarrow \mR$ is an expected objective, $c: \mR^d\rightarrow \mR^{m}$ are deterministic~equality constraints, $g: \mR^d\rightarrow \mR^{r}$ are deterministic inequality constraints,~$\xi\sim \P$ is a random variable following the distribution $\P$, and $F(\cdot\;; \xi):\mR^d\rightarrow\mR$ is~a realized objective. In stochastic optimization regime, the direct evaluation of $f$ and its derivatives is not accessible. Instead, it is assumed that one can generate independent and identically distributed samples $\{\xi_i\}_{i}$ from $\P$, and estimate $f$ and its derivatives based on the realizations $\{F(\cdot \;; \xi_i )\}_{i}$.

Problem \eqref{pro:1} widely appears in a variety of industrial applications including finance, transportation, manufacturing, and power systems \citep{Birge1997State, Silvapulle2004Constrained}. It includes constrained empirical risk minimization (ERM) as a special case, where $\P$ can be regarded as a uniform distribution over $n$ data points $\{\xi_i = (\by_i, \bz_i)\}_{i=1}^n$, with $(\by_i, \bz_i)$ being the feature-outcome pairs. Thus, the objective has a finite-sum form as
\begin{equation*}
f(\bx) = \frac{1}{n}\sum_{i=1}^{n}F(\bx; \xi_i) = \frac{1}{n}\sum_{i=1}^{n}F(\bx; \by_i, \bz_i).
\end{equation*}
The goal of \eqref{pro:1} is to find the optimal parameter $\tx$ that fits the data best.~One of the most~common choices of $F$ is the negative log-likelihood of the underlying distribution of $(\by_i, \bz_i)$. In this case, the optimizer $\tx$ is called the maximum likelihood estimator (MLE). Constraints on parameters are also common in practice, which are used to encode prior model knowledge or to restrict model complexity. For example, \cite{Liew1976Two, Liew1976Inequality} studied inequality constrained least-squares problems, where inequality constraints maintain structural consistency such as non-negativity of the elasticities. \cite{Phillips1991constrained, Onuk2015Constrained} studied statistical properties~of~constrained MLE, where constraints characterize the parameters space of interest. More recently, a growing literature on training constrained neural networks has been reported \citep{Goh2018Constrained, Chen2018Constraint, Livieris2019adaptive, Livieris2019improved}, where constraints are imposed to avoid weights either vanishing or exploding, and objectives are in the above finite-sum form.

This paper aims to develop a numerical procedure to solve \eqref{pro:1} with a~global convergence guarantee. When the objective $f$ is deterministic, numerous nonlinear optimization methods with well-understood convergence results are applicable, such as exact penalty methods, augmented Lagrangian methods,~sequential quadratic programming (SQP) methods, and interior-point methods \citep{Nocedal2006Numerical}. However, methods to solve constrained \textit{stochastic} nonlinear problems~with satisfactory convergence guarantees have been developed only recently. In particular, with only equality constraints, \cite{Berahas2021Sequential} designed a very first stochastic SQP (StoSQP) scheme using an $\ell_1$-penalized merit function, and showed that for any initialization, the KKT residuals $\{R_t\}_t$ converge in two different regimes, determined by a prespecified deterministic stepsize-related sequence $\{\alpha_t\}_t$:
\begin{enumerate}[label=(\alph*),topsep=0pt]	\setlength\itemsep{0.0em}
\item (constant sequence) if $\alpha_t = \alpha$ for some small $\alpha>0$, then $ \sum_{i=0}^{t-1}\mE[R_i^2]/t \leq \Upsilon/(\alpha t) + \Upsilon\alpha$ for some $\Upsilon>0$;

\item (decaying sequence) if $\alpha_t$ satisfies $\sum_{t=0}^{\infty}\alpha_t = \infty$ and $\sum_{t=0}^{\infty}\alpha_t^2 <\infty$, then $\liminf_{t\rightarrow \infty}\mE[R_t^2] = 0$.
\end{enumerate}
Both convergence regimes are well known for unconstrained stochastic problems where $R_t = \|\nabla f(\bx_t)\|$ (see \cite{Bottou2018Optimization} for a recent review), while \cite{Berahas2021Sequential} generalized~the results to equality constrained problems. Within the algorithm of \cite{Berahas2021Sequential}, the authors designed a stepsize selection scheme (based on the prespecified deterministic sequence) to bring some sort of adaptivity into the algorithm. However, it turns out that the prespecified sequence, which can be aggressive or conservative, still highly affects the performance. To address the adaptivity~issue,~\cite{Na2022adaptive} proposed an alternative StoSQP, which exploits a differentiable exact augmented Lagrangian merit function, and enables a stochastic line search~procedure to adaptively select~the stepsize. Under a different setup (where the model~is~precisely~estimated with high probability), \cite{Na2022adaptive} proved a different guarantee: for any~initialization, $\liminf_{t\rightarrow \infty} R_t = 0$ \textit{almost~surely}. Subsequently, a series of extensions have been reported. \cite{Berahas2021Stochastic} designed a StoSQP scheme to~deal~with~rank-deficient constraints. \cite{Curtis2021Inexact} designed a StoSQP that exploits inexact Newton directions. \cite{Berahas2022Accelerating} designed an accelerated StoSQP via variance reduction for finite-sum problems. \cite{Berahas2022Adaptive} further developed \cite{Berahas2021Sequential} to achieve adaptive sampling. \cite{Curtis2021Worst} established the worst-case~iteration complexity of StoSQP, and \cite{Na2022Asymptotic} established the asymptotic local rate of StoSQP and performed statistical inference. In addition, \cite{Oztoprak2021Constrained} investigated a deterministic SQP where the objective and constraints are evaluated with noise. However, all aforementioned literature does not include inequality constraints.

Our paper develops this line of research by designing a StoSQP method that works with nonlinear inequality constraints. In order to do so, we have to overcome a number of intrinsic difficulties that arise in dealing with inequality constraints, which were already noted in classical nonlinear optimization literature \citep{Bertsekas1982Constrained, Nocedal2006Numerical}. Our work is built upon \cite{Na2022adaptive}, where we exploited an augmented Lagrangian merit function under the SQP framework. We~\mbox{enhance}~some~of~\mbox{designs} in \cite{Na2022adaptive} (e.g., the feasibility error condition, the \mbox{increasing} batch size, and~the~complexity of Hessian sampling; more on these later), and the analysis of this~paper is more involved. To \mbox{generalize} \cite{Na2022adaptive}, we address the following two subtleties.
\begin{enumerate}[label=(\alph*),topsep=0pt]	\setlength\itemsep{0.0em}
\item With inequalities, SQP subproblems are inequality constrained (nonconvex) quadratic programs (IQPs), which themselves are difficult to solve in most cases. Some SQP literature (e.g., \cite{Boggs1995Sequential}) supposes to apply a QP~solver to solve IQPs exactly, however, a practical scheme should embed a finite number of inner loop iterations of active-set methods or interior-point methods~into~the main SQP loop, to solve IQPs approximately. Then, the inner loop may lead to an~approximation error for search direction in each iteration, which complicates the analysis.

\item When applied to deterministic objectives with inequalities, the SQP search~direction is a descent direction of the augmented Lagrangian only in~a~neighborhood of a KKT point \citep[Propositions 8.3, 8.4]{Pillo2002Augmented}. This is in contrast~to~\mbox{equality} constrained problems, where the descent property of the SQP direction holds globally, provided the penalty parameters of~the augmented Lagrangian are suitably chosen. Such a difference is indeed brought by inequality constraints: to make the (active-set) SQP direction informative, the estimated active set has to be close to the optimal active set (see Lemma \ref{lem:6} for details). Thus, simply changing the merit function in \cite{Na2022adaptive} does not work for~Problem~\eqref{pro:1}.
\end{enumerate}
The existing literature on inequality constrained SQP has addressed (a) and~(b) via various tools for deterministic objectives, while we provide new insights~into stochastic objectives. To resolve (a), we design an active-set StoSQP scheme, where given the current iterate, we first identify an active set which includes all inequality constraints that are likely to be equalities. We then obtain the~search direction by solving a SQP subproblem, where we include all inequality constraints in the identified active set but regard them as equalities. In this case, the subproblem is an equality constrained QP (EQP), and can be solved exactly provided the matrix factorization is within the computational~budget. To resolve (b), we provide a safeguarding direction to the scheme. In each~step, we check if the SQP subproblem is solvable and generates a descent direction of the augmented Lagrangian merit function. If yes, we maintain the SQP direction as it typically enjoys a~fast local rate; if no, we switch to the safeguarding direction (e.g., one gradient/Newton step of the augmented~Lagrangian), along which~the iterates still decrease the augmented Lagrangian although the convergence may not be as effective as that of SQP.

Furthermore, to design a scheme that adaptively selects the penalty~parameters and stepsizes for Problem \eqref{pro:1}, additional challenges have to be~resolved. In particular, we know that there are \textit{unknown deterministic} thresholds for penalty parameters to ensure one-to-one correspondence between a stationary point of the merit function and a KKT point of Problem \eqref{pro:1}. However, due to the scheme stochasticity, the stabilized penalty parameters are random. We are unsure~if the stabilized values are above (or below, depending on the context) the thresholds or not. Thus, we cannot directly conclude that the iterates converge to a KKT point, even if we ensure a sufficient decrease on the merit function in each step, and enforce the iterates to converge~to~one~of~its~stationary~points.

The above difficulty has been resolved for the $\ell_1$-penalized merit function in \cite{Berahas2021Sequential}, where the authors imposed a probability condition on the noise (satisfied by symmetric noise; see \cite[Proposition 3.16]{Berahas2021Sequential}). \cite{Na2022adaptive} resolved this difficulty for the augmented Lagrangian merit function by modifying the SQP scheme when selecting the penalty parameters. In particular, \cite{Na2022adaptive} required the feasibility~error to be bounded by the gradient magnitude of the augmented Lagrangian in \textit{each} step, and generated monotonically increasing samples to estimate the gradient. Although that analysis does not require noise conditions, adjusting the penalty parameters to enforce the feasibility error condition may not~be~necessary for the iterates that are far from stationarity. Also, generating increasing samples is not satisfactory since the sample size should be adaptively chosen based on the iterates. In this paper, we refine the techniques of \cite{Na2022adaptive} and generalize them to inequality constraints. We weaken the feasibility error condition by using a (large) multiplier to rescale the augmented Lagrangian gradient, and more significantly, enforcing it \textit{only when the magnitude of the rescaled augmented Lagrangian gradient is smaller than the estimated KKT residual}. In other words, the feasibility error condition is imposed only when we have a stronger evidence that the iterate is approaching to a stationary point than approaching to a KKT point. Such a relaxation matches the motivation~of~the~feasibility error condition, i.e., bridging the gap between stationary points and KKT~points.~We~also~get~rid of the increasing sample size requirement by adaptively controlling the absolute deviation of the augmented Lagrangian gradient \textit{for the new iterates only} (i.e. the previous step is a \textit{successful} step; see Section \ref{sec:4}). Following \cite{Na2022adaptive}, we perform a stochastic line search procedure. However, instead of using the same~sample~set to estimate the gradient $\nabla f$ and~Hessian $\nabla^2 f$ as in \cite{Na2022adaptive}, we sharpen~the~analysis and realize that the needed samples for $\nabla^2 f$ are significantly less than $\nabla f$.

With all above extensions from \cite{Na2022adaptive}, we finally prove that the KKT residual $R_t$ satisfies $\lim_{t\rightarrow \infty} R_t = 0$ \textit{almost surely} for any initialization. Such a result is stronger than \cite[Theorem 4.10]{Paquette2020Stochastic} for unconstrained problems and \cite[Theorem 4]{Na2022adaptive} for equality constrained~problems, which only showed the ``liminf" type~of convergence. Our result also differs from the (liminf) convergence~of~the~expected KKT residual $\mE[R_t^2]$ established in \cite{Berahas2021Sequential, Berahas2021Stochastic, Berahas2022Adaptive, Berahas2022Accelerating, Curtis2021Inexact} (under a different~setup).

\vskip3pt
\noindent{\bf Related work.} A number of methods have been proposed to optimize stochastic objectives without constraints, varying from first-order methods to second-order methods \citep{Bottou2018Optimization}. For all methods, adaptively choosing the stepsize is particularly important for practical deployment. A line of literature selects the stepsize by adaptively controlling the batch size and embedding natural (stochastic) line search into the schemes \citep{Friedlander2012Hybrid, Byrd2012Sample, Krejic2013Line, De2017Automated, Bollapragada2018Adaptive}. Although empirical experiments suggest the validity of stochastic line search, a rigorous analysis is~missing.~Until recently, researchers revisited unconstrained stochastic optimization via~the lens of classical nonlinear optimization methods, and were able to show promising convergence guarantees. In particular, \cite{Bandeira2014Convergence, Chen2017Stochastic, Gratton2017Complexity, Blanchet2019Convergence, Sun2022trust} studied stochastic trust-region methods, and \cite{Cartis2017Global, Serafino2020LSOS, Paquette2020Stochastic, Berahas2021Global} studied stochastic line search methods. Moreover, \cite{Berahas2021Sequential, Na2022adaptive, Berahas2021Stochastic, Curtis2021Inexact, Berahas2022Accelerating, Berahas2022Adaptive} designed a variety of StoSQP schemes to solve equality constrained stochastic problems. Our paper contributes to this line of works by proposing~an active-set StoSQP scheme to handle inequality~constraints.

There are numerous methods for solving deterministic problems with nonlinear constraints, varying from exact penalty methods, augmented Lagrangian methods, interior-point methods, and sequential quadratic programming (SQP) methods \citep{Nocedal2006Numerical}. Our paper is based~on~SQP,~which is a very effective (or at least competitive) approach for small or large problems. When inequality constraints are present, SQP can be classified into IQP and EQP approaches.~The~former solves inequality constrained subproblems; the latter, to which our method~belongs, solves equality constrained subproblems. A clear advantage of EQP~over IQP is that the subproblems are less expensive to solve, especially when the quadratic matrix is indefinite. See \cite[Chapter 18.2]{Nocedal2006Numerical} for a comparison. Within SQP schemes, an exact penalty function is used as the merit function to monitor the progress of the iterates towards a KKT point. The $\ell_1$-penalized merit function, $f(\bx) + \mu\rbr{\|c(\bx)\|_1 + \|\max\{g(\bx),\0\}\|_1}$, is always a plausible choice because of its simplicity. However, a disadvantage of such non-differentiable merit functions is their impedance of fast local rates. A nontrivial local modification of SQP has to be employed to relieve such an issue \cite{Boggs1995Sequential}. As a resolution,~multiple differentiable merit functions have been proposed \citep{Bertsekas1982Constrained}. We exploit an augmented Lagrangian merit function, which was first proposed for equality constrained problems by \cite{Pillo1979New, Pillo1980method}, and then extended to inequality constrained problems by \cite{Pillo1982new, Pillo1985Continuously}. \cite{Pillo2002Augmented} further improved this series of works by designing a new~augmented Lagrangian, and established the exact property under weaker conditions. Although not crucial for that exact property analysis, \cite{Pillo2002Augmented} did not include equality constraints. In this paper, we enhance the augmented~Lagrangian in \cite{Pillo2002Augmented} by containing both equality and inequality constraints; and study the case where the objective is stochastic. When inequality constraints are suppressed, our algorithm and analysis naturally reduce to \cite{Na2022adaptive} (with refinements). We should mention that differentiable merit functions are often more expensive to evaluate, and their benefits are mostly revealed for local rates (see \cite[Figure 1]{Na2021Global} for a comparison between the augmented Lagrangian and $\ell_1$ merit functions on an optimal control problem). Thus, with only established global analysis, we do not aim to claim the benefits of the augmented Lagrangian over the popular~$\ell_1$ merit function. On the other hand, the augmented Lagrangian is a very common alternative of~\mbox{non-differentiable penalty functions, which has been widely utilized}  for inequality constrained problems~and~achieved~promising~performance~\citep{Zavala2014Scalable, Pillo2005Convergence, Pillo2008truncated, Pillo2011primal, Pillo2011exact}. Also, our global analysis is the first step towards understanding~the~local rate of~StoSQP when differentiable merit functions~are~employed.

\vskip 2pt
\noindent{\bf Structure of the paper.} We introduce the exploited augmented Lagrangian merit function and active-set SQP subproblems in Section \ref{sec:2}. We propose our StoSQP scheme and analyze it in Section \ref{sec:4}. The experiments and conclusions are in Sections \ref{sec:5} and~\ref{sec:6}. Due to the space limit, we defer all proofs to Appendix.

\vskip 1pt
\noindent{\bf Notation.} 
We use $\|\cdot\|$ to denote the $\ell_2$ norm for vectors and spectrum~norm for matrices. For two scalars $a$ and $b$, $a\wedge b = \min\{a, b\}$ and $a\vee b = \max\{a, b\}$.~For two vectors $\ba$ and $\bb$ with the same dimension, $\min\{\ba, \bb\}$ and $\max\{\ba, \bb\}$ are~vectors by taking entrywise minimum and maximum, respectively. For $\ba\in \mR^r$, $\diag(\ba) \in \mR^{r\times r}$ is a diagonal matrix whose diagonal entries are specified by $\ba$ sequentially. $I$ denotes the identity matrix whose dimension is clear from the context. For a set $\mA \subseteq\{1,2,\ldots, r\}$ and a vector $\ba\in \mR^r$ (or a matrix $A\in \mR^{r\times d}$), $\ba_{\mA} \in \mR^{|\mA|}$ (or $A_{\mA}\in \mR^{|\mA|\times d}$) is a sub-vector (or a sub-matrix) including only~the indices in $\mA$; $\Pi_{\mA}(\cdot): \mR^r\rightarrow \mR^r$ (or $\mR^{r\times d} \rightarrow \mR^{r\times d}$) is a projection operator with $[\Pi_{\mA}(\ba)]_i = \ba_i$ if $i\in \mA$ and $[\Pi_{\mA}(\ba)]_i = 0$ if $i\notin \mA$ (for $A\in\mR^{r\times d}$, $\Pi_{\mA}(A)$ is applied column-wise); $\mA^c = \{1,2,\ldots, r\}\backslash\mA$.~Finally, we reserve the notation for the Jacobian~matrices of constraints: $J(\bx) = \nabla^T c(\bx) = (\nabla c_1(\bx), \ldots, \nabla c_m(\bx))^T \in \mR^{m\times d}$ and $G(\bx) = \nabla^T g(\bx) = (\nabla g_1(\bx), \ldots, \nabla g_r(\bx))^T \in \mR^{r\times d}$.

\section{Preliminaries}\label{sec:2}
\vskip-11.3pt
Throughout this section, we suppose $f, c, g$ are twice continuously differentiable (i.e., $f,g,c\in C^2$). The Lagrangian function of Problem \eqref{pro:1} is
\begin{equation*}
\mL(\bx, \bmu, \blambda) = f(\bx) + \bmu^Tc(\bx) + \blambda^Tg(\bx).
\end{equation*}
We denote by
\begin{equation}\label{equ:feasible}
\Omega = \{\bx\in \mR^d: c(\bx) = \0, g(\bx) \leq \0 \}
\end{equation}
the feasible set and
\begin{equation}\label{equ:I}
\mI(\bx) = \{i:  1\leq i\leq r, g_i(\bx) = \0\}
\end{equation}
the active set. We aim to find a KKT point $(\tx, \tmu, \tlambda)$ of \eqref{pro:1} satisfying
\begin{equation}\label{equ:KKT:cond}
\nabla_{\bx}\mL(\tx, \tmu, \tlambda) = \0,\; c(\tx) = \0,\; g(\tx)\leq \0, \; \tlambda \geq \0, \; (\tlambda)^Tg(\tx) = 0.
\end{equation}
When a constraint qualification holds, existing a dual pair $(\tmu, \tlambda)$ to satisfy~\eqref{equ:KKT:cond} is a first-order necessary condition for $\tx$ being a local solution of \eqref{pro:1}.~In~most cases, it is difficult to have an initial iterate that satisfies all inequality constraints, and enforce inequality constraints to hold as the iteration proceeds. This motivates us to consider a perturbed set. For $\nu>0$, we let
\begin{equation}\label{equ:domain}
\Omega \subsetneq\mT_{\nu} \coloneqq \cbr{\bx\in \mR^d: a(\bx) \leq \nu/2 } \quad \text{where }\; a(\bx) = \sum_{i=1}^{r} \max\{g_i(\bx), 0\}^3.
\end{equation}
Here, the perturbation radius $\nu/2$ is not essential and can be replaced by $\nu/\kappa$ for any $\kappa>1$. Also, the cubic power in $a(\bx)$ can be replaced by any power $s$ with $s>2$, which ensures that $a(\bx)\in C^2$ provided $g_i(\bx)\in C^2$, $\forall i$. We~also~define a scaling function
\begin{equation}\label{equ:q}
q_{\nu}(\bx, \blambda) = \frac{a_{\nu}(\bx)}{1+\|\blambda\|^2}\quad \text{ with }\; a_{\nu}(\bx) = \nu -a(\bx),
\end{equation}
where $a_{\nu}(\bx)$ measures the distance of $a(\bx)$ to the boundary $\nu$, and $q_{\nu}(\bx, \blambda)$ rescales $a_{\nu}(\bx)$ by penalizing $\blambda$ that has a large magnitude. In the definitions of \eqref{equ:domain} and \eqref{equ:q}, $\nu>0$ is a parameter to be chosen: given the current primal iterate $\bx_t$, we choose $\nu = \nu_t$ large enough so that $\bx_t \in \mT_{\nu}$. Note that while it is difficult to have $\bx_t\in \Omega$, it is easy to choose $\nu$ to have $\bx_t\in \mT_{\nu}$. We also note~that
\begin{equation*}
\frac{\nu}{2(1+\|\blambda\|^2)}\leq q_{\nu}(\bx, \blambda) \leq \nu \;\; \forall (\bx, \blambda) \in \mT_{\nu}\times \mR^r, \;\;  \text{ and }\;\; q_{\nu}(\bx, \blambda)\rightarrow 0 \; \text{ as } \|\blambda\|\rightarrow \infty.
\end{equation*}
With \eqref{equ:q} and a parameter $\epsilon>0$, we define a function to measure the dual~feasibility of inequality constraints:
\begin{multline}\label{equ:def:w}
\bw_{\epsilon, \nu}(\bx, \blambda)  \coloneqq g(\bx) - \bb_{\epsilon, \nu}(\bx, \blambda) \\
 \coloneqq g(\bx) - \min\{\0, g(\bx) + \epsilon q_{\nu}(\bx, \blambda)\blambda\} = \max\{g(\bx), - \epsilon q_{\nu}(\bx, \blambda)\blambda\}.
\end{multline}
The following lemma justifies the reasonability of the definition \eqref{equ:def:w}. The proof is immediate and omitted.

\begin{lemma}\label{lem:1}
Let $\epsilon, \nu>0$. For any $(\bx, \blambda) \in \mT_{\nu}\times \mR^r$, $\bw_{\epsilon, \nu}(\bx, \blambda) = \0 \Leftrightarrow g(\bx) \leq~\0,\\ \blambda\geq \0, \blambda^Tg(\bx) = 0$.
\end{lemma}

An implication of Lemma \ref{lem:1} is that, when the iteration sequence converges to a KKT point, $\bw_{\epsilon, \nu}(\bx, \blambda)$ converges to 0, i.e., $g(\bx) = \bb_{\epsilon, \nu}(\bx,\blambda)$. This motivates us to define the following augmented Lagrangian function:
\begin{multline}\label{equ:aug:Lagrange}
\mL_{\epsilon, \nu, \eta}(\bx, \bmu, \blambda) = \mL(\bx, \bmu, \blambda) + \frac{1}{2\epsilon}\|c(\bx)\|^2 + \frac{1}{2\epsilon q_{\nu}(\bx, \blambda)}\rbr{\|g(\bx)\|^2 - \|\bb_{\epsilon, \nu}(\bx, \blambda)\|^2}\\
+ \frac{\eta}{2}\nbr{\begin{pmatrix}
J(\bx) \nabla_{\bx}\mL(\bx, \bmu, \blambda)\\
G(\bx) \nabla_{\bx}\mL(\bx, \bmu, \blambda) + \diag^2(g(\bx))\blambda
\end{pmatrix}}^2,
\end{multline}
where $\eta>0$ is a prespecified parameter, which can be any positive number throughout the paper. The augmented Lagrangian \eqref{equ:aug:Lagrange} generalizes the one in \cite{Pillo2002Augmented} by including equality constraints and introducing $\eta$ to enhance flexibility ($\eta=2$ in \cite{Pillo2002Augmented}). Without inequalities, \eqref{equ:aug:Lagrange} reduces to the augmented Lagrangian studied in \cite{Na2022adaptive}. The penalty in \eqref{equ:aug:Lagrange} consists of two parts. The first part characterizes the feasibility error and consists of $\|c(\bx)\|^2$ and $\|g(\bx)\|^2 - \|\bb_{\epsilon, \nu}(\bx, \blambda)\|^2$. The latter term is rescaled by $1/q_{\nu}(\bx, \blambda)$ to penalize $\blambda$ with a large magnitude.~In~fact, if $\|\blambda\|\rightarrow \infty$, then $q_{\nu}(\bx, \blambda)\blambda \rightarrow \0$ so that $b_{\epsilon, \nu}(\bx, \blambda) \rightarrow \min\{\0, g(\bx)\}$ (cf. \eqref{equ:def:w}). Thus, the penalty term $(\|g(\bx)\|^2 - \|b_{\epsilon}(\bx, \blambda)\|^2)/q_{\nu}(\bx,\blambda)\rightarrow \infty$, which is impossible when the iterates decrease $\mL_{\epsilon, \nu, \eta}$. The second part characterizes the optimality error and does not depend on the parameters $\epsilon$ and $\nu$. 
We mention that there are alternative forms of the augmented Lagrangian, some of~which~transform nonlinear inequalities using (squared) slack variables \citep{Bertsekas1982Constrained, Zavala2014Scalable}. In that case, additional variables are involved and the strict complementarity condition is often needed to ensure the equivalence between the original and transformed~problems~\citep{Fukuda2017note}.

The exact property of \eqref{equ:aug:Lagrange} can be studied similarly as in \cite{Pillo2002Augmented}, however this~is incremental and not crucial for our analysis. We will only use (a stochastic version of) \eqref{equ:aug:Lagrange} to monitor the progress of the iterates. By direct calculation, we obtain the~gradient $\nabla\mL_{\epsilon, \nu, \eta}$. We first suppress the evaluation~point~for~conciseness, and define the following matrices
\begin{align}\label{equ:def:Qmatrices}
Q_{11} = & (\nabla_{\bx}^2\mL)J^T, \;\; Q_{12} = \sum_{i=1}^{m}(\nabla^2c_i)( \nabla_{\bx}\mL) \be_{i, m}^T, \;\; Q_{1} = Q_{11} + Q_{12} \in \mR^{d\times m}, \nonumber\\
Q_{21} = & (\nabla_{\bx}^2\mL)G^T,\;\; Q_{22} = \sum_{i=1}^{r}(\nabla^2g_i)( \nabla_{\bx}\mL) \be_{i, r}^T, \;\; Q_{23} = 2G^T\diag(g)\diag(\blambda), \\
Q_{2} = &\sum_{i=1}^{3}Q_{2i}\in \mR^{d\times r},\;\; M =  \left(\begin{smallmatrix}
M_{11} & M_{12}\\
M_{21} & M_{22}
\end{smallmatrix}\right) = \left(\begin{smallmatrix}
JJ^T & JG^T\\
GJ^T & GG^T + \diag^2(g)
\end{smallmatrix}\right) \in \mR^{(m+r)\times (m+r)}, \nonumber
\end{align}
where $\be_{i, m}\in \mR^{m}$ is the $i$-th canonical basis of $\mR^m$ (similar for $\be_{i, r}\in \mR^{r}$). Then, 
\begin{multline}\label{equ:aug:der}
\begin{pmatrix}
\nabla_{\bx}\mL_{\epsilon, \nu, \eta}\\ \nabla_{\bmu}\mL_{\epsilon, \nu, \eta}\\
\nabla_{\blambda}\mL_{\epsilon, \nu, \eta}
\end{pmatrix} = \begin{pmatrix}
I & \frac{1}{\epsilon}J^T & \frac{1}{\epsilon q_{\nu}}G^T\\
 & I \\
 & & I
\end{pmatrix}\begin{pmatrix}
\nabla_{\bx}\mL\\
c\\
\bw_{\epsilon, \nu}
\end{pmatrix} +  \begin{pmatrix}
\frac{3\|\bw_{\epsilon, \nu}\|^2}{2\epsilon q_{\nu}a_{\nu}}G^T\bl\\
\0\\
\frac{\|\bw_{\epsilon, \nu}\|^2}{\epsilon a_{\nu}}\blambda
\end{pmatrix}  \\ + \eta\begin{pmatrix}
Q_1 & Q_2\\
M_{11} & M_{12}\\
M_{21} & M_{22}
\end{pmatrix}\begin{pmatrix}
J\nabla_{\bx}\mL\\
G\nabla_{\bx}\mL + \diag^2(g)\blambda
\end{pmatrix},
\end{multline}
where $\bl = \bl(\bx) = \diag(\max\{g(\bx), \0\})\max\{g(\bx),\0\}$. Clearly, the evaluation~of $\nabla\mL_{\epsilon, \nu, \eta}$ requires $\nabla f$ and $\nabla^2 f$, which have to be replaced by their stochastic counterparts $\bnabla f$ and $\bnabla^2 f$ for Problem \eqref{pro:1}. Based on \eqref{equ:aug:der}, we note that, if~the feasibility error vanishes, then $\nabla\mL_{\epsilon, \nu, \eta} = \0$ implies the KKT conditions \eqref{equ:KKT:cond}~hold for any $\epsilon,\nu,\eta>0$. We summarize this observation in the next lemma. The~result holds without any constraint qualifications.

\begin{lemma}\label{lem:2}
Let $\epsilon, \nu, \eta >0$ and let $(\tx, \tmu, \tlambda) \in \mT_{\nu} \times \mR^{m}\times \mR^{r}$ be a primal-dual triple. If $\|c(\tx)\| = \|\bw_{\epsilon, \nu}(\tx, \tlambda)\| = \|\nabla\mL_{\epsilon, \nu, \eta}(\tx, \tmu, \tlambda)\| = 0$, then~$(\tx, \tmu, \tlambda)$ satisfies \eqref{equ:KKT:cond} and, hence, is a KKT point of Problem \eqref{pro:1}.
\end{lemma}

\begin{proof}
See Appendix \ref{pf:lem:2}.
\end{proof}

In the next subsection, we introduce an active-set SQP direction that~is~motivated by the augmented Lagrangian \eqref{equ:aug:Lagrange}.

\subsection{An active-set SQP direction via EQP}

Let $\epsilon, \nu, \eta>0$ be fixed parameters. Suppose we have the $t$-th iterate $(\bx_t,\bmu_t,\blambda_t)\in \mT_{\nu}\times \mR^m\times \mR^r$, let us denote $J_t = J(\bx_t)$, $G_t = G(\bx_t)$ (similar for $\nabla f_t, c_t, g_t$,~$q_{\nu}^t$ etc.) to be the quantities evaluated at the $t$-th iterate. We generally use index~$t$ as subscript, except for the quantities (e.g., $q_{\nu}^t$) that depend on~$\epsilon$,~$\nu$,~or~$\eta$,~which have been used as subscript. For an active set $\mA\subseteq \{1,\ldots, r\}$, we denote~$\blambda_{t_a} = (\blambda_t)_{\mA}$,~$\blambda_{t_c} = (\blambda_t)_{\mA^c}$ (similar for $g_{t_a}$, $g_{t_c}$, $G_{t_a}$, $G_{t_c}$ etc.) to be the sub-vectors (or~sub-matrices), and denote $\Pi_a(\cdot) = \Pi_{\mA}(\cdot)$, $\Pi_c(\cdot) = \Pi_{\mA^c}(\cdot)$ for shorthand.

With the $t$-th iterate $(\bx_t,\bmu_t,\blambda_t)$ and the above notation, we first define~the identified active set as
\begin{equation}\label{equ:active}
\mA_{\epsilon, \nu}^t\coloneqq  \mA_{\epsilon, \nu}(\bx_t, \blambda_t) \coloneqq \{i: 1 \leq i\leq r, \; (g_t)_i \geq -\epsilon q_{\nu}^t\cdot(\blambda_t)_i\}.
\end{equation}
We then solve the following coupled linear system
\begin{subequations}\label{equ:SQP:direction}
\begin{align}
\overbrace{\begin{pmatrix}
B_t & J_t^T & G_{t_a}^T\\
J_t\\
G_{t_a}
\end{pmatrix} }^{K_{t_a}}\begin{pmatrix}
\Delta\bx_t\\
\tDelta\bmu_t\\
\tDelta\blambda_{t_a}
\end{pmatrix} & =  -\begin{pmatrix}
\nabla_{\bx}\mL_t - G_{t_c}^T\blambda_{t_c}\\
c_t\\
g_{t_a}
\end{pmatrix},   \label{equ:SQP:direction:1}\\
\underbrace{\begin{pmatrix}
J_tJ_t^T & J_tG_t^T\\
G_tJ_t^T & G_tG_t^T + \diag^2(g_t)
\end{pmatrix}  }_{M_t}\begin{pmatrix}
\Delta\bmu_t\\
\Delta\blambda_t
\end{pmatrix}& \\
&\hskip-3.7cm =  - \bigg\{  \begin{pmatrix}
J_t\nabla_{\bx}\mL_t\\
G_t\nabla_{\bx}\mL_t +\Pi_c(\diag^2(g_t)\blambda_t)
\end{pmatrix} + \begin{pmatrix}
Q_{1,t}^T\\
Q_{2,t}^T
\end{pmatrix} \Delta\bx_t \bigg\}, \label{equ:SQP:direction:2}
\end{align}
\end{subequations}
for some $B_t$ that approximates the Hessian $\nabla_{\bx}^2\mL_t$. Our active-set SQP direction is then $\Delta_t\coloneqq (\Delta\bx_t, \Delta\bmu_t, \Delta\blambda_t)$. Finally, we update the iterate as
\begin{equation}\label{equ:update}
\begin{pmatrix}
\bx_{t+1}\\
\bmu_{t+1}\\
\blambda_{t+1}
\end{pmatrix} = \begin{pmatrix}
\bx_{t}\\
\bmu_{t}\\
\blambda_{t}
\end{pmatrix} + \alpha_t\begin{pmatrix}
\Delta\bx_t\\
\Delta\bmu_t\\
\Delta\blambda_t
\end{pmatrix}
\end{equation}
with $\alpha_t$ chosen to ensure a certain sufficient decrease on the merit function \eqref{equ:aug:Lagrange}.

The definition of active set was introduced in \cite[(8.5)]{Pillo2002Augmented} and has been utilized, e.g., in \cite{Pillo2008truncated}. Intuitively, for the $i$-th inequality constraint, if $g_i^\star = (g(\bx^\star))_i = 0$ and $\blambda_i^\star>0$, then $i$ will be identified when $(\bx_t, \blambda_t)$ is close to $(\tx,\tlambda)$; if $g_i^\star<0$ and $\tlambda_i = 0$, then $i$ will~not be identified. The stepsize $\alpha_t$ is usually chosen by line search. In Section \ref{sec:4}, we will design a stochastic line search~scheme to select $\alpha_t$ adaptively. Compared to fully stochastic SQP schemes \cite{Berahas2021Sequential, Berahas2021Stochastic, Curtis2021Inexact}, we need a more precise model estimation. We explain the SQP direction~\eqref{equ:SQP:direction}~in~the~next~remark.

\begin{remark}\label{rem:1}

Our dual direction $(\Delta\bmu_t, \Delta\blambda_t)$ differs from the usual SQP direction introduced, for example,~in \cite[(8.9)]{Pillo2002Augmented}. In particular, the system \eqref{equ:SQP:direction:1} is nothing but the KKT conditions of EQP:
\begin{align}\label{pro:2}
\min_{\Delta\bx_t}\;\; & \frac{1}{2}(\Delta\bx_t)^T B_t\Delta\bx_t + (\nabla f_t)^T\Delta\bx_t, \nonumber\\
\text{s.t.}\;\; & c_t + J_t\Delta\bx_t = \0,\\
& g_{t_a} + G_{t_a}\Delta\bx_t = \0. \nonumber
\end{align}
Thus, $(\Delta\bx_t, \bmu_t + \tDelta\bmu_t, \blambda_{t_a} + \tDelta\blambda_{t_a})$ solved from \eqref{equ:SQP:direction:1} is also the primal-dual~solution of~the above EQP. However, instead of using $(\tDelta\bmu_t, \tDelta\blambda_{t_a}, -\blambda_{t_c})$,~we~solve the dual direction $(\Delta\bmu_t, \Delta\blambda_t)$ for both active and inactive constraints from~\eqref{equ:SQP:direction:2}. As~$B_t$ converges to $\nabla_{\bx}^2\mL_t$ and $(\bx_t,\bmu_t,\blambda_t)$ converges to a KKT point $(\tx, \tmu, \tlambda)$, it is fairly easy to see that $(\Delta\bmu_t, \Delta\blambda_t)$ converges to $(\tDelta\bmu_t, \tDelta\blambda_t)$ (where we~denote $\tDelta\blambda_{t_c} = -\blambda_{t_c}$) in a higher order by noting that
\begin{multline*}
\begin{pmatrix}
	J_tJ_t^T & J_tG_t^T\\
	G_tJ_t^T & G_tG_t^T + \diag^2(g_t)
\end{pmatrix}\begin{pmatrix}
\tDelta\bmu_t\\
\tDelta\blambda_t
\end{pmatrix} \stackrel{\eqref{equ:SQP:direction:1}}{=} \begin{pmatrix}
\0\\
\Pi_a(\diag^2(g_t)\tDelta\blambda_t)
\end{pmatrix} \\ -\cbr{\begin{pmatrix}
J_t\nabla_{\bx}\mL_t\\
G_t\nabla_{\bx}\mL_t +\Pi_c(\diag^2(g_t)\blambda_t)
\end{pmatrix} + \begin{pmatrix}
J_t\\
G_t
\end{pmatrix}B_t\Delta\bx_t}.
\end{multline*}
Thus, the fast \textit{local} rate of the SQP direction $(\Delta\bx_t, \tDelta\bmu_t, \tDelta\blambda_t)$ is preserved~by~$\Delta_t$. However, it turns out that the adjustment of $\Delta_t$ is crucial for the merit function \eqref{equ:aug:Lagrange} when $B_t$ is far from $\nabla_{\bx}^2\mL_t$. A similar, coupled SQP system is employed~for equality constrained problems \citep{Lucidi1990Recursive, Na2022adaptive},~while~we extend to inequality constraints here. In fact, \cite[Proposition 8.2]{Pillo2002Augmented} showed~that $(\Delta\bx_t, \tDelta\bmu_t, \tDelta\blambda_t)$ is a descent direction of $\mL_{\epsilon, \nu, \eta}^t$ if $(\bx_t, \bmu_t, \blambda_t)$ is near a KKT point and $B_t = \nabla^2_{\bx}\mL_t$. However, $B_t = \nabla^2_{\bx}\mL_t$ (i.e., no Hessian modification) is restrictive even for a deterministic line search, and that descent result does not hold if~$B_t \neq \nabla^2_{\bx}\mL_t$. In contrast, as shown in Lemma \ref{lem:6}, $\Delta_t$ is a descent direction even if $B_t$ is not close to $\nabla_{\bx}^2\mL_t$.

\end{remark}

\subsection{The descent property of $\Delta_t$}

In this subsection, we present a descent property of $\Delta_t$. We focus on the term $(\nabla\mL_{\epsilon, \nu, \eta}^t)^T\Delta_t$. Different from SQP for equality constrained problems,~$\Delta_t$~may~not be a descent direction of $\mL_{\epsilon, \nu, \eta}^t$ for some points even if $\epsilon$ is chosen~small enough. To see it clearly, we suppress the iteration index, denote $g_a = g_{t_a}$ (similar for $\blambda_a$, $\blambda_c$ etc.), and divide $\nabla\mL_{\epsilon, \nu, \eta} $ (cf. \eqref{equ:aug:der}) into two terms: a dominating term that depends on $(g_a, \blambda_c)$ \textit{linearly}, and a higher-order term that depends on $(g_a, \blambda_c)$~at least \textit{quadratically}. In particular, we write $\nabla\mL_{\epsilon, \nu, \eta} = \nabla\mL_{\epsilon, \nu, \eta}^{(1)} + \nabla\mL_{\epsilon, \nu, \eta}^{(2)}$ where 
\begin{align}\label{equ:new:aug:der}
\begin{pmatrix}
\nabla_{\bx}\mL_{\epsilon, \nu, \eta}^{(1)}\\
\nabla_{\bmu}\mL_{\epsilon, \nu, \eta}^{(1)}\\
\nabla_{\blambda}\mL_{\epsilon, \nu, \eta}^{(1)}
\end{pmatrix} & = \begin{pmatrix}
I & \frac{1}{\epsilon}J^T & \frac{1}{\epsilon q_{\nu}}G^T\\
& I & \\
& & I
\end{pmatrix}\begin{pmatrix}
\nabla_{\bx}\mL\\
c\\
\bw_{\epsilon, \nu}
\end{pmatrix} + \eta\begin{pmatrix}
Q_1 & Q_2\\
M_{11} & M_{12}\\
M_{21} & M_{22}
\end{pmatrix}\left(\begin{smallmatrix}
J\nabla_{\bx} \mL\\
G\nabla_{\bx}\mL + \Pi_c(\diag^2(g)\blambda)
\end{smallmatrix}\right), \nonumber\\
\begin{pmatrix}
\nabla_{\bx}\mL_{\epsilon, \nu, \eta}^{(2)}\\
\nabla_{\bmu}\mL_{\epsilon, \nu, \eta}^{(2)}\\
\nabla_{\blambda}\mL_{\epsilon, \nu, \eta}^{(2)}
\end{pmatrix} & = \begin{pmatrix}
\frac{3\|\bw_{\epsilon, \nu}\|^2}{2\epsilon q_{\nu}a_{\nu}}G^T\bl\\
\0\\
\frac{\|\bw_{\epsilon, \nu}\|^2}{\epsilon a_{\nu}}\blambda
\end{pmatrix} + \eta \begin{pmatrix}
Q_{2,a}\\
M_{12,a}\\
M_{22,a}
\end{pmatrix}\diag^2(g_a)\blambda_a.
\end{align}

Loosely speaking (see Lemma \ref{lem:6} for a rigorous result), $(\nabla\mL_{\epsilon, \nu, \eta}^{(1)})^T\Delta$ provides a sufficient decrease provided the penalty parameters are suitably chosen, while $(\nabla\mL_{\epsilon, \nu, \eta}^{(2)})^T\Delta$ has no such guarantee in general. Since $\nabla\mL_{\epsilon, \nu, \eta}^{(2)}$ depends~on~$(g_a, \blambda_c)$ \textit{quadratically}, to ensure $\nabla\mL_{\epsilon, \nu, \eta}^T\Delta<0$, we require $\|g_a\|\vee\|\blambda_c\|$ to~be~small~enough to let the linear term $(\nabla\mL_{\epsilon, \nu, \eta}^{(1)})^T\Delta$ dominate. This essentially requires~the iterate to be close to a KKT point, since $\|g_a\| = \|\blambda_c\| = 0$ at a KKT point. With this discussion in mind, if the iterate is far from a KKT point, $\Delta$ may not be a descent direction of $\mL_{\epsilon, \nu, \eta}$. In fact, for an iterate that is far from~a~KKT~point,~the KKT matrix $K_a$ (and its component $G_a$) is likely to be singular due to the imprecisely identified active set. Thus, Newton system \eqref{equ:SQP:direction} is not solvable~at~this iterate at all, let alone it generates a descent direction. Without inequalities, the quadratic term $\nabla\mL_{\epsilon, \nu, \eta}^{(2)}$ disappears and our analysis reduces to the one in \cite{Na2022adaptive}. We realize that the existence of $\nabla\mL_{\epsilon, \nu, \eta}^{(2)}$ results in a very different augmented Lagrangian to the one in \cite{Na2022adaptive}; and brings difficulties in designing a global~algorithm to deal with inequality constraints.

We point out that requiring a local iterate is not an artifact of the~proof~technique. Such a requirement is imposed for different search directions in related literature. For example, \cite{Pillo2002Augmented} showed that the SQP direction obtained by either EQP or IQP is a descent direction of $\mL_{\epsilon, \nu, \eta}$ in a \textit{neighborhood} of a KKT point (cf. Propositions 8.2 and 8.4). That work also required $B_t = \nabla_{\bx}^2\mL_t$, which~we relax by considering a coupled Newton system. Subsequently, \cite{Pillo2008truncated, Pillo2011exact} studied truncated Newton directions, whose descent properties hold only \textit{locally} as well (cf. \cite[Proposition 3.7]{Pillo2008truncated}, \cite[Proposition 10]{Pillo2011exact}).

Now, we introduce two assumptions and formalize the descent property. 

\begin{assumption}[LICQ]\label{ass:1}
We assume at $\tx$ that $(J^T(\tx)\;\; G^T_{\I(\tx)}(\tx))$ has full column rank, where $\I(\tx)$ is the active inequality set defined in \eqref{equ:I}.
\end{assumption}

\begin{assumption}\label{ass:2}
For $\bz \in \{\bz\in\mR^d: J_t\bz = \0, G_{t_a}\bz = \0\}$, we have $\bz^TB_t\bz \geq \gamma_{B}\|\bz\|^2$ and $\|B_t\| \leq \Upsilon_{B}$ for constants $\Upsilon_{B}\geq 1\geq \gamma_{B}>0$.
\end{assumption}

The above condition on $B_t$ is standard in nonlinear optimization literature \citep{Bertsekas1982Constrained}. In fact, $B_t = I$ with $\gamma_{B} = \Upsilon_{B} = 1$ is sufficient for the analysis in this paper. The condition $\Upsilon_{B}\geq 1\geq \gamma_{B}>0$ (similar for other constants defined later) is inessential, which is only for simplifying the presentation. Without such a requirement, our analyses hold by replacing $\gamma_{B}$ with $\gamma_{B}\wedge 1$ and $\Upsilon_{B}$~with~$\Upsilon_{B}\vee 1$.

\begin{lemma}\label{lem:6}
Let $\nu, \eta>0$ and suppose Assumptions \ref{ass:1} and \ref{ass:2} hold. There exist a constant $\Upsilon>0$ depending on $\Upsilon_{B}$ but not on $(\nu,\eta,\gamma_{B})$, and a compact~set~$\mX_{\epsilon,\nu}\times \mM\times\Lambda_{\epsilon,\nu}$ around $(\tx, \tmu, \tlambda)$ depending on $(\epsilon, \nu)$ but not on $\eta$,\footnote{Here, we mean $\mX_{\epsilon,\nu}$ and $\Lambda_{\epsilon,\nu}$ only \textit{directly depend} on $\epsilon,\nu$ but not $\eta$, which are~in~contrast to neighborhoods $\mX_{\epsilon,\nu,\eta}$ and $\Lambda_{\epsilon,\nu,\eta}$. However, since the threshold of $\epsilon$, $\gamma_{B}^2(\gamma_{B}\wedge \eta)/\cbr{(1\vee \nu)\Upsilon}$, is also determined by $\eta$, the final local neighborhoods $\mX_{\epsilon,\nu}$ and $\Lambda_{\epsilon,\nu}$ with $\epsilon$ below the threshold also \textit{indirectly depend} on $\eta$. Recall that $\eta$ can be any positive constant throughout the paper.} such that if $(\bx_t, \bmu_t, \blambda_t) \in \mX_{\epsilon,\nu}\times\mM\times \Lambda_{\epsilon,\nu}$ with $\epsilon$ satisfying  $\epsilon\leq \gamma_{B}^2(\gamma_{B}\wedge \eta)/\cbr{(1\vee \nu)\Upsilon}$,~then
\begin{equation*}
(\nabla\mL_{\epsilon, \nu, \eta}^{t\;(1)})\Delta_t \leq -\frac{\gamma_{B}\wedge\eta}{2} \nbr{\left(\begin{smallmatrix}
\Delta\bx_t\\
J_t\nabla_{\bx}\mL_t\\
G_t\nabla_{\bx}\mL_t+ \Pi_c(\diag^2(g_t)\blambda_t)
\end{smallmatrix} \right)}^2.
\end{equation*}
Furthermore, there exists a compact subset $\mX_{\epsilon,\nu,\eta}\times \mM\times\Lambda_{\epsilon,\nu,\eta}\subseteq \mX_{\epsilon,\nu}\times\mM\times \Lambda_{\epsilon,\nu}$ depending additionally on $\eta$,  such that if $(\bx_t, \bmu_t, \blambda_t) \in \mX_{\epsilon,\nu,\eta}\times\mM\times \Lambda_{\epsilon,\nu,\eta}$, then
\begin{equation*}
(\nabla\mL_{\epsilon, \nu, \eta}^{t\;(2)})\Delta_t \leq \frac{\gamma_{B}\wedge\eta}{4} \nbr{\left(\begin{smallmatrix}
\Delta\bx_t\\
J_t\nabla_{\bx}\mL_t\\
G_t\nabla_{\bx}\mL_t+ \Pi_c(\diag^2(g_t)\blambda_t)
\end{smallmatrix} \right)}^2.
\end{equation*}

\end{lemma}

\begin{proof}
See Appendix \ref{pf:lem:6}.
\end{proof}

Similar arguments for other directions can be found in \cite[Proposition~3.5]{Pillo2008truncated} and \cite[Proposition 9]{Pillo2011exact}. By the proof of Lemma \ref{lem:6}, we know that as long~as~$M_t$ and $(J_t^T\;\; G_{t_a}^T)$ in the SQP system \eqref{equ:SQP:direction} have full (column) rank, $(\nabla\mL_{\epsilon, \nu, \eta}^{t\; (1)})^T\Delta_t$ ensures a~sufficient decrease provided $\epsilon$ is small enough. However, from \eqref{npequ:1} in the proof, we also see that $(\nabla\mL_{\epsilon, \nu, \eta}^{t\; (2)})^T\Delta_t$ is only bounded by
\begin{equation*}
(\nabla\mL_{\epsilon, \nu, \eta}^{t\;(2)})\Delta_t \leq \Upsilon'\rbr{\frac{1\vee \nu}{\epsilon(1\wedge\nu^2)}\vee \eta}(\|g_{t_a}\| + \|\blambda_{t_c}\|)
\nbr{\left(\begin{smallmatrix}
\Delta\bx_t\\
J_t\nabla_{\bx}\mL_t\\
G_t\nabla_{\bx}\mL_t+ \Pi_c(\diag^2(g_t)\blambda_t)
\end{smallmatrix} \right)}^2,
\end{equation*}
where $\Upsilon'>0$ is a constant independent of $(\epsilon, \nu, \eta)$. Thus, to ensure $(\nabla\mL_{\epsilon, \nu, \eta}^t)^T\Delta_t$ to be negative, we have to restrict to a neighborhood, in which $\|g_{t_a}\|\vee \|\blambda_{t_c}\|$~is small enough so that $\Upsilon'(\frac{1\vee\nu}{\epsilon(1\wedge\nu^2)}\vee \eta)(\|g_{t_a}\| + \|\blambda_{t_c}\|) \leq (\gamma_{B}\wedge \eta)/4$. This~requirement is achievable near a KKT pair $(\tx, \tlambda)$, where the active set is correctly identified (implying that $\|g_{t_a}\|\leq \|(g_t)_{\I(\tx)}\|$ and $\|\blambda_{t_c}\|\leq \|(\blambda_t)_{\{i: 1\leq i\leq r, \tlambda_i=0\}}\|$); and the radius of the neighborhood clearly depends on $(\epsilon, \nu, \eta)$.

In the next section, we exploit the introduced augmented Lagrangian merit function \eqref{equ:aug:Lagrange} and the active-set SQP direction \eqref{equ:SQP:direction} to design a StoSQP~scheme~for Problem \eqref{pro:1}. We will adaptively choose proper $\epsilon$ and $\nu$ (recall that $\eta>0$ can be any positive number in this paper), incorporate stochastic line search to select the stepsize, and globalize the scheme by utilizing a safeguarding direction (e.g., Newton or steepest descent step) of the merit function $\mL_{\epsilon,\nu,\eta}$. If the system \eqref{equ:SQP:direction} is not solvable, or is solvable but does not generate a descent direction, we search along the alternative direction~to~decrease~the~merit~function. However, since $\Delta_t$ usually enjoys a fast local rate (see \cite[Proposition 8.3]{Pillo2002Augmented} for a local analysis of $(\Delta\bx_t, \tDelta\bmu_t, \tDelta\blambda_t)$ and Remark \ref{rem:1}), we prefer to preserve $\Delta_t$ as much as possible.

\section{An Adaptive Active-Set StoSQP Scheme}\label{sec:4}

We design an adaptive scheme for Problem \eqref{pro:1} that embeds stochastic~line~search, originally designed and analyzed for unconstrained problems in \cite{Cartis2017Global, Paquette2020Stochastic}, into~an active-set StoSQP. There are two challenges to design adaptive schemes for constrained problems. First, the merit function has penalty parameters that~are random and adaptively specified; while for unconstrained problems one simply uses the objective function in line search. To show the global convergence, it is crucial that the stochastic penalty parameters are stabilized \textit{almost surely}.~Thus, for each run, after few iterations we always target a stabilized merit function. Otherwise, if each iteration decreases a different merit function, the decreases across iterations may not accumulate. Second, since the stabilized parameters are random, they may not be below \textit{unknown deterministic} thresholds. Such a condition is critical to ensure the equivalence between the stationary points of the merit function and the KKT points of Problem \eqref{pro:1}. Thus, even~if~we~converge to a stationary point of the (stabilized) merit function, it is not necessarily true that the stationary point is a KKT point of Problem \eqref{pro:1}.

With only equality constraints, \cite{Berahas2021Sequential, Na2022adaptive} addressed the first challenge under a boundedness condition, and our paper follows the same type of analysis.~Similar boundedness condition is also required for deterministic analyses to have~the penalty parameters stabilized \cite[Chapter 4.3.3]{Bertsekas1982Constrained}. \cite{Berahas2021Sequential} resolved the second challenge by introducing a noise condition (satisfied by symmetric noise), while~\cite{Na2022adaptive} resolved it by adjusting the SQP scheme when selecting the penalty parameters. As introduced in Section \ref{sec:1}, the technique of \cite{Na2022adaptive} has multiple flaws: (i) it~requires generating increasing samples to estimate the gradient of the~\mbox{augmented}~Lagrangian (cf. \cite[Step 1]{Na2022adaptive}); (ii) it imposes a feasibility error condition for each step (cf. \cite[(19)]{Na2022adaptive}). In this paper, we refine the technique of \cite{Na2022adaptive} and enable inequality constraints. As revealed by Section \ref{sec:2}, the present~analysis~of~inequality constraints is much more involved; and more importantly, our ``lim"~convergence	guarantee strengthens the existing ``liminf" convergence of the stochastic line search in \cite{Paquette2020Stochastic, Na2022adaptive}. In what follows, we use $\bar{(\cdot)}$ to denote random quantities, except for the iterate $(\bx_t, \bmu_t, \blambda_t)$. For example, $\baralpha_t$ denotes a random stepsize.

\subsection{The proposed scheme}

Let $\eta, \alpha_{max}, \kappa_{grad}, \chi_{grad}, \chi_{f},\chi_{err}>0; \rho>1$; $\gamma_{B}\in(0,1]$; $\beta, p_{grad}, p_f\in (0, 1)$; $\kappa_f\in(0, \beta/(4\alpha_{max})]$ be fixed tuning parameters. Given quantities $(\bx_t, \bmu_t, \blambda_t, \barnu_t, \\\barepsilon_t, \baralpha_t,\bardelta_t)$ at the $t$-th iteration with $\bx_t\in \mT_{\barnu_t}$, we perform the following five steps to derive quantities at the $(t+1)$-th iteration.

\vskip3pt
\noindent\textbf{Step 1: Estimate objective derivatives.}
We generate a batch of independent samples $\xi_1^t$ to estimate the gradient $\nabla f_t$ and Hessian $\nabla^2f_t$. The estimators $\bnabla f_t$ and $\bnabla^2f_t$ may not be computed with the same amount of samples, since they have different sample complexities. For example, we can compute $\bnabla f_t$ using $\xi_1^t$ while compute $\bnabla^2 f_t$ using a fraction of $\xi_1^t$ (more on this in Section \ref{sec:4.3}).~With $\bnabla f_t$, $\bnabla^2f_t$, we then compute $\bnabla_{\bx}\mL_t$, $\barQ_{1,t}$, and $\barQ_{2,t}$ used in the system \eqref{equ:SQP:direction}.

We require the batch size $|\xi_1^t|$ to be large enough to make the gradient error of the merit function small. In particular, we define
\begin{equation*}
\barDelta(\nabla\mL_\eta^t) \coloneqq \|\bnabla\mL_{\barepsilon_t,\barnu_t,\eta}^t - \nabla\mL_{\barepsilon_t,\barnu_t,\eta}^t\|. 
\end{equation*}
\hskip-3.5pt A simple observation from \eqref{equ:aug:der} is that $\barDelta(\nabla\mL_\eta^t)$ is independent of $\barepsilon_t$ (and $\barnu_t$), which will be selected later (Step 2). We require $|\xi_1^t|$ to satisfy two conditions: 

\noindent\textbf{(a)} the event $\E_1^t$,
\begin{equation}\label{event:E1}
\E_1^t = \bigg\{\barDelta(\nabla\mL_\eta^t)
\leq \kappa_{grad} \baralpha_t \underbrace{\nbr{\left(\begin{smallmatrix}
\bnabla_{\bx}\mL_t\\
c_t\\
\max\{g_t, -\blambda_t\}
\end{smallmatrix}\right) } }_{\barR_t} \bigg\},
\end{equation}
satisfies
\begin{equation}\label{cond:event:E1}
P_{\xi_1^t}\rbr{\E_1^t} \geq 1 - p_{grad};
\end{equation}
\noindent\textbf{(b)} if $t-1$ is a \textit{successful} step (see Step 5 for the meaning), then
\begin{equation}\label{cond:grad:var}
\mE_{\xi_1^t}[ \barDelta(\nabla\mL_\eta^t)] \leq \chi_{grad}\cdot (\bardelta_t/\baralpha_t)^{1/2}.
\end{equation}

The sample complexities to ensure \eqref{cond:event:E1} and \eqref{cond:grad:var} will be discussed in~Section \ref{sec:4.3}. Compared to \cite{Na2022adaptive}, we do not let $|\xi_1^t|$ increase~monotonically, while we~impose an expectation condition \eqref{cond:grad:var} when we arrive at a new iterate. By our analysis, it is easy to see that \eqref{cond:grad:var} can also be replaced by requiring the \textit{subsequence} $\{|\xi_1^t|: t-1 \text{ is a successful step}\}$ to increase to the infinity (e.g., increase by~at~least~one each time), which is still weaker than \cite{Na2022adaptive}. The right hand side of \eqref{cond:grad:var} will be clear when we utilize $\bardelta_t$ later in Step 5 (cf. \eqref{cond:decrease:trust}). We use $P_{\xi_1^t}(\cdot)$ and $\mE_{\xi_1^t}[\cdot]$ to
denote the probability and expectation that are evaluated over the randomness of sampling $\xi_1^t$ only, while other random quantities are conditioned on, such as $(\bx_t, \bmu_t, \blambda_t)$ and $\baralpha_t$. More precisely, we mean $P_{\xi_1^t}(\E_1^t) = P(\E_1^t\mid \mF_{t-1})$ (similar~for $\mE_{\xi_1^t}[\cdot]$) where the $\sigma$-algebra $\mF_{t-1}$ is defined in \eqref{sigma:alg} below.

\vskip 3pt
\noindent\textbf{Step 2: Set parameter $\barepsilon_t$.}
With current $\barnu_t$, we decrease $\barepsilon_t \leftarrow \barepsilon_t/\rho$ until $\barepsilon_t$ is small enough to satisfy the following two conditions simultaneously:

\noindent\textbf{(a)} the feasibility error is proportionally bounded by the gradient of the merit function, whenever the iterate is closer to a stationary point than a KKT point: \begin{equation}\label{cond:bound:fes:error}
\nbr{(c_t, \bw_{\barepsilon_t, \barnu_t}^t
)}\leq \chi_{err}\|\bnabla \mL_{\barepsilon_t, \barnu_t, \eta}^t \| \quad\quad \text{ if }\; \chi_{err}\|\bnabla \mL_{\barepsilon_t, \barnu_t, \eta}^t\| \leq \barR_t ;
\end{equation}
(we use the same multiplier $\chi_{err}$ only for simplifying the notation.)

\noindent\textbf{(b)} if the SQP system \eqref{equ:SQP:direction} with $\bnabla_{\bx}\mL_t$, $\barQ_{1,t}$, and $\barQ_{2,t}$ is solvable, then we~obtain $\barDelta_t = (\barDelta\bx_t, \barDelta\bmu_t, \barDelta\blambda_t)$ and require
\begin{equation}\label{cond:decrease:1}
(\bnabla\mL_{\barepsilon_t, \barnu_t, \eta}^{t\;(1)})^T\barDelta_t \leq -\frac{(\gamma_{B}\wedge \eta)}{2} \nbr{\left(\begin{smallmatrix}
\barDelta\bx_t\\
J_t\bnabla_{\bx}\mL_t\\
G_t\bnabla_{\bx}\mL_t+ \Pi_c(\diag^2(g_t)\blambda_t)
\end{smallmatrix}\right) }^2.
\end{equation}

We prove in Lemma \ref{lem:7} and Lemma~\ref{lem:8} that both \eqref{cond:bound:fes:error} and \eqref{cond:decrease:1} can be satisfied for sufficiently small $\barepsilon_t$. In fact, Lemma \ref{lem:6} has already established \eqref{cond:decrease:1} for the deterministic case. Even though $\barDelta_t$ is not always used as the search direction, we still enforce \eqref{cond:decrease:1} to hold for $(\bnabla\mL_{\barepsilon_t, \barnu_t, \eta}^{t\;(1)})^T\barDelta_t$. The reason for this is to avoid ruling out $\barDelta_t$ just because $\barepsilon_t$ is not small enough, which~would~result~in~a~positive dominating term $(\bnabla\mL_{\barepsilon_t, \barnu_t, \eta}^{t\;(1)})^T\barDelta_t$. If \eqref{equ:SQP:direction} is not solvable (e.g., the active set is imprecisely identified so that $K_{t_a}$ is singular), then \eqref{cond:decrease:1} is not needed.

The condition \eqref{cond:bound:fes:error} is the key to ensure that the stationary point of the merit function that we converge to is a KKT point of \eqref{pro:1}. Motivated by Lemma \ref{lem:2}, we know that ``the stationarity of the merit function plus vanishing feasibility error" implies vanishing KKT residual. \eqref{cond:bound:fes:error} states that the feasibility error is roughly controlled by the gradient of the merit function. \eqref{cond:bound:fes:error} relaxes \cite[(19)]{Na2022adaptive} from two aspects. First, \cite{Na2022adaptive} had no multiplier while we allow any (large) multiplier $\chi_{err}$. Second, \cite{Na2022adaptive} enforced \eqref{cond:bound:fes:error} for each step, while we enforce it only when~we~observe a stronger evidence that the scheme is approaching to a stationary point than to a KKT point. The above relaxations are driven by the intention of imposing the condition. When adjusting~$\barepsilon_t$,~if~$\|\bnabla\mL_{\barepsilon_t,\barnu_t,\eta}^t\|$~first~\mbox{exceeds}~$\barR_t$~\mbox{before}~$\|(c_t, \bw_{\barepsilon_t, \barnu_t}^t)\|$ (which easily happens for a large $\barnu_t$), then one can immediately stop the~adjustment of $\barepsilon_t$. Compared to \cite{Na2022adaptive} where the SQP system is supposed to be always solvable, \eqref{cond:bound:fes:error} has extra usefulness: when $\barDelta_t$ is not available, \eqref{cond:bound:fes:error} ensures that the safeguarding direction can be computed using the samples in Step 1. Such a desire is not easily achieved, and further relaxations of \eqref{cond:bound:fes:error} can be designed if we generate new samples for the safeguarding direction (in Step~3).~The~subtlety lies in the fact that no penalty parameters are involved when we generate $\xi_1^t$ in Step 1, while \eqref{cond:bound:fes:error} builds a connection between $\xi_1^t$~and~the penalty parameters.~It implies that the set $\xi_1^t$ satisfying \eqref{cond:event:E1} and \eqref{cond:grad:var} also satisfies the corresponding conditions for the safeguarding direction.

\vskip 3pt
\noindent\textbf{Step 3: Decide the search direction.}
We may obtain a stochastic SQP direction $\barDelta_t$ from Step~2. However, if \eqref{equ:SQP:direction} is not solvable, or it is solvable~but~$\barDelta_t$ is not a sufficient descent direction because
\vskip-4pt\begin{equation}\label{cond:decrease:2}
(\bnabla\mL_{\barepsilon_t, \barnu_t, \eta}^{t\;(2)})^T\barDelta_t > \frac{(\gamma_{B}\wedge \eta)}{4} \nbr{\left(\begin{smallmatrix}
\barDelta\bx_t\\
J_t\bnabla_{\bx}\mL_t\\
G_t\bnabla_{\bx}\mL_t+ \Pi_c(\diag^2(g_t)\blambda_t)
\end{smallmatrix}\right) }^2,
\end{equation}
then an alternative safeguarding direction $\hatDelta_t$ must be employed to ensure~the~decrease of the merit function. In that case, we follow \cite{Pillo2008truncated, Pillo2011exact} and regard $\mL_{\barepsilon_t,\barnu_t,\eta}$ as a \mbox{penalized} objective. We require $\hatDelta_t$ to~satisfy
\begin{equation}\label{cond:alter:dir}
(\bnabla\mL_{\barepsilon_t, \barnu_t, \eta}^t)^T\hatDelta_t \leq - 1/\chi_{u}\cdot \|\bnabla\mL_{\barepsilon_t, \barnu_t, \eta}^t\|^2\;\; \text{ and }\;\; \|\hatDelta_t\|\leq \chi_{u}\cdot\nbr{\bnabla\mL_{\barepsilon_t, \barnu_t, \eta}^t}
\end{equation}
for a constant $\chi_{u}\geq 1$. Similar to \eqref{cond:bound:fes:error}, we use the same constant $\chi_{u}$ for the~two multipliers to simplify the notation. When using two different constants $\chi_{1,u}$ and $\chi_{2,u}$, we can always set $\chi_{u} = 1/\chi_{1,u}\vee \chi_{2,u}$ to let \eqref{cond:alter:dir} hold. The condition 
\eqref{cond:alter:dir} is standard in the literature  \cite[(60a,b)]{Pillo2008truncated}  \cite[(52a,b)]{Pillo2011exact}. One~example~that satisfies \eqref{cond:alter:dir} and is computationally cheap is~the steepest descent direction~$\hatDelta_t = - \bnabla\mL_{\barepsilon_t, \barnu_t, \eta}^t$ with $\chi_{u} = 1$. Such a direction can be computed (almost) without any extra cost since the two components of $\bnabla\mL_{\barepsilon_t, \barnu_t, \eta}^t$, $\bnabla\mL_{\barepsilon_t, \barnu_t, \eta}^{t\; (1)}$ and $\bnabla\mL_{\barepsilon_t, \barnu_t, \eta}^{t\; (2)}$, have been computed when checking \eqref{cond:decrease:1} and \eqref{cond:decrease:2}. Another example that is more computationally expensive is the regularized Newton step $\hatH_t\hatDelta_t = -\bnabla\mL_{\barepsilon_t, \barnu_t, \eta}^t$, where $\hatH_t$ captures second-order information of $\mL_{\barepsilon_t,\barnu_t,\eta}^t$ and satisfies~$1/\chi_{u}I \preceq \hatH_t\preceq \chi_{u}I$. In particular, $\hatH_t$ can be obtained by regularizing the (generalized) 
Hessian matrix $H_t$, which is provided and discussed in \cite{Pillo2002Augmented, Pillo2008truncated}, and has the~form\footnote{See (6.1)-(6.3) in \cite{Pillo2002Augmented} for a similar expression to \eqref{equ:def:H}. Our $H_t$ generalizes that definition by including equality constraints and approximating the Hessian $\nabla_{\bx}^2\mL_t$ by $B_t$.}
\begin{align}\label{equ:def:H}
&H_{t,\bx\bx} =  B_t + \eta B_t\cbr{J_t^TJ_t + G_t^TG_t}B_t + \frac{1}{\barepsilon_t}J_t^TJ_t + \frac{1}{\barepsilon_t q_{\barnu_t}^t}G_{t_a}^TG_{t_a}, \nonumber\\
&H_{t,(\bmu, \blambda)\bx} = \left(\begin{smallmatrix}
J_t\\
\Pi_a(G_t)
\end{smallmatrix}\right) + \eta \left(\begin{smallmatrix}
J_tJ_t^T & J_tG_t^T\\
G_tJ_t^T & G_tG_t^T + \diag^2(\Pi_c(g_t))
\end{smallmatrix}\right)\left(\begin{smallmatrix}
J_t\\
G_t
\end{smallmatrix}\right)B_t,\\
&H_{t,(\bmu,\blambda)(\bmu, \blambda)} =  \left(\begin{smallmatrix}
\0 & \0 \\
\0 & -\barepsilon_t q_{\barnu_t}^t\diag(\Pi_c(\b1))
\end{smallmatrix}\right) + \eta \left(\begin{smallmatrix}
J_tJ_t^T & J_tG_t^T\\
G_tJ_t^T & G_tG_t^T + \diag^2(\Pi_c(g_t))
\end{smallmatrix}\right)^2. \nonumber
\end{align}
Here, $\b1 = (1,\ldots, 1)\in\mR^r$ is the all one vector. Other~examples that improve upon the regularized Newton step include the choices in \cite{Pillo2011primal, Fasano2009nonmonotone}, where a~truncated conjugate gradient~method is applied to an \textit{indefinite} Newton system \cite[Proposition 3.3, (14)]{Pillo2011primal}. We will numerically implement the regularized Newton and the steepest descent steps in Section \ref{sec:5}.

\vskip 3pt
\noindent\textbf{Step 4: Estimate the merit function.}
Let $\cDelta_t$ denote the adopted search direction; thus $\cDelta_t = \barDelta_t$ from Step 2 or $\cDelta_t = \hDelta_t$ from Step 3. We aim to~perform stochastic line search by checking the Armijo condition \eqref{cond:armijo} at the trial~point
\begin{equation*}
\bx_{s_t} = \bx_t + \baralpha_t\cDelta\bx_t, \quad\quad \bmu_{s_t} = \bmu_t + \baralpha_t\cDelta\bmu_t, \quad\quad \blambda_{s_t} = \blambda_t + \baralpha_t\cDelta\blambda_t.
\end{equation*}
We estimate the merit function in this~step and perform line search in Step 5.

First, we check if the trial primal point $\bx_{s_t}$ is in $\mT_{\barnu_t}$. In particular, if~$\bx_{s_t} \notin \mT_{\barnu_t}$, that is $a_{s_t} = a(\bx_{s_t}) > \barnu_t/2$ (cf. \eqref{equ:domain}), then we stop the current iteration and~reject the trial point by letting $(\bx_{t+1}, \bmu_{t+1}, \blambda_{t+1}) = (\bx_t, \bmu_t, \blambda_t)$, $\barepsilon_{t+1} = \barepsilon_t$, $\baralpha_{t+1} = \baralpha_t$, and $\bardelta_{t+1} = \bardelta_t$. We also increase $\barnu_t$ by letting
\begin{equation}\label{cond:nu:t}
\barnu_{t+1} = \rho^j\barnu_t\quad\text{ with }\quad j = \lceil\log(2a_{s_t}/\barnu_t)/\log \rho\rceil,
\end{equation}
where $\lceil y\rceil$ denotes the least integer that exceeds $y$. The definition~of $j\geq 1$~in~\eqref{cond:nu:t} ensures $\bx_{s_t} \in \mT_{\barnu_{t+1}}$. However, $j=1$ works as well, since $\bx_{t+1}=\bx_t \in \mT_{\barnu_t} \subseteq \mT_{\barnu_{t+1}}$, as required for performing the next iteration. In the case of $\bx_{s_t} \notin \mT_{\barnu_t}$, particularly if $a_{s_t} \geq \barnu_t$, evaluating the merit function $\mL_{\barepsilon_t, \barnu_t, \eta}^{s_t}$ is not informative since the penalty term in $\mL_{\barepsilon_t, \barnu_t, \eta}^{s_t}$ may be rescaled by a negative multiplier.~Thus, we increase $\barnu_t$ and rerun the iteration at the current point.

Otherwise $\bx_{s_t} \in \mT_{\barnu_t}$, then we generate a batch of independent samples $\xi_2^t$, that are independent from $\xi_1^t$ as well, and estimate $f_t, f_{s_t}, \nabla f_t, \nabla f_{s_t}$. Similar to Step 1, the~estimators $\barf_t, \barf_{s_t}$ and $\bbnabla f_t, \bbnabla f_{s_t}$ may not be computed with~the~same amount of samples. For example, $\barf_t$ and $\barf_{s_t}$ can be computed using $\xi_2^t$ while~$\bbnabla f_t$ and $\bbnabla f_{s_t}$ can be computed using a fraction of $\xi_2^t$. The sample complexities are discussed in Section \ref{sec:4.3}. Here, we distinguish $\bbnabla f_t$ from $\bnabla f_t$ in Step 1. While both of them are estimates of $\nabla f_t$, the former is computed based on $\xi_2^t$ and the latter is computed based on $\xi_1^t$. Using $\barf_t,\barf_{s_t},\bbnabla f_t,\bbnabla f_{s_t}$, we compute~$\barL_{\barepsilon_t, \barnu_t, \eta}^t$~and $\barL_{\barepsilon_t, \barnu_t, \eta}^{s_t}$ according to \eqref{equ:aug:Lagrange}.

We require $|\xi_2^t|$ is large enough such that the event $\E_2^t$,
\begin{equation}\label{event:E2}
\text{\footnotesize $\E_2^t = \cbr{\abr{\barL_{\barepsilon_t, \barnu_t, \eta}^t - \mL_{\barepsilon_t, \barnu_t, \eta}^t} \vee \abr{\barL_{\barepsilon_t, \barnu_t, \eta}^{s_t} - \mL_{\barepsilon_t, \barnu_t, \eta}^{s_t}} \leq -\kappa_f\baralpha_t^2(\bnabla\mL_{\barepsilon_t, \barnu_t, \eta}^t)^T\cDelta_t }$,}
\end{equation}
satisfies
\begin{equation}\label{cond:event:E2}
P_{\xi_2^t}\rbr{\E_2^t} \geq 1 - p_f
\end{equation}
and
\begin{equation}\label{cond:merit:var}
\mE_{\xi_2^t}[|\barL_{\barepsilon_t, \barnu_t, \eta}^t - \mL_{\barepsilon_t, \barnu_t, \eta}^t|^2]\vee \mE_{\xi_2^t}[|\barL_{\barepsilon_t, \barnu_t, \eta}^{s_t} - \mL_{\barepsilon_t, \barnu_t, \eta}^{s_t}|^2] \leq \chi_{f}\cdot \bardelta_t^2.
\end{equation}
Similar to \eqref{cond:event:E1} and \eqref{cond:grad:var}, $P_{\xi_2^t}(\cdot)$ and $\mE_{\xi_2^t}[\cdot]$ denote that the~randomness is taken over sampling $\xi_2^t$ only, while other random quantities are conditioned on. That is, $P_{\xi_2^t}(\E_2^t) = P(\E_2^t\mid \mF_{t-0.5})$ (similar for $\mE_{\xi_2^t}[\cdot]$) where the $\sigma$-algebra~$\mF_{t-0.5} = \mF_{t-1}\cup \sigma(\xi_1^t)$ is defined in \eqref{sigma:alg} below.

\vskip 3pt
\noindent\textbf{Step 5: Perform line search.} With the merit function estimates, we check the Armijo condition next.

\noindent\textbf{(a)} If the Armijo condition holds,
\begin{equation}\label{cond:armijo}
\barL_{\barepsilon_t, \barnu_t, \eta}^{s_t} \leq \barL_{\barepsilon_t, \barnu_t, \eta}^t + \beta\baralpha_t (\bnabla\mL_{\barepsilon_t, \barnu_t, \eta}^t)^T\cDelta_t,
\end{equation}
then the trial point is accepted by letting $(\bx_{t+1}, \bmu_{t+1}, \blambda_{t+1}) = (\bx_{s_t}, \bmu_{s_t}, \blambda_{s_t})$~and the stepsize is increased by $\baralpha_{t+1} = \rho\baralpha_t\wedge \alpha_{max}$. Furthermore, we check if the decrease of the merit function is reliable. In particular, if
\begin{equation}\label{cond:decrease:trust}
-\beta\baralpha_t (\bnabla\mL_{\barepsilon_t, \barnu_t, \eta}^t)^T\cDelta_t \geq \bardelta_t,
\end{equation}
then we increase $\bardelta_t$ by $\bardelta_{t+1} = \rho\bardelta_t$; otherwise, we decrease $\bardelta_t$ by $\bardelta_{t+1} = \bardelta_t/\rho$.

\noindent\textbf{(b)} If the Armijo condition \eqref{cond:armijo} does not hold, then the trial point is rejected~by letting $(\bx_{t+1}, \bmu_{t+1}, \blambda_{t+1}) = (\bx_{t}, \bmu_{t},\blambda_{t})$, $\baralpha_{t+1} = \baralpha_t/\rho$ and $\bardelta_{t+1} = \bardelta_t/\rho$.

Finally, for both cases \textbf{(a)} and \textbf{(b)}, we let $\barepsilon_{t+1} = \barepsilon_t$, $\barnu_{t+1} = \barnu_t$ and repeat~the procedure from Step 1. From \eqref{cond:decrease:trust}, we can see that $\bardelta_t$ (roughly) has the order $\baralpha_t\|\bnabla\mL_{\barepsilon_t,\barnu_t,\eta}^t\|^2$, which justifies the definition of the right hand side of \eqref{cond:grad:var}.

The proposed scheme is summarized in Algorithm \ref{alg:ASto}. We define three~types~of iterations for line search. If the Armijo condition \eqref{cond:armijo} holds, we call the~iteration a \textit{\textbf{successful step}}, otherwise we call it an \textit{\textbf{unsuccessful step}}. For a successful step, if the sufficient decrease in \eqref{cond:decrease:trust} is satisfied, we call it a \textit{\textbf{reliable step}}, otherwise we call it an \textit{\textbf{unreliable step}}. Same notion is used in \cite{Cartis2017Global, Paquette2020Stochastic, Na2022adaptive}.

To end this section, let us introduce the filtration induced by the randomness of the algorithm. Given a random sample sequence $\{\xi_1^t,\xi_2^t\}_{t=0}^\infty$,\footnote{We note that $\xi_2^t$ may not be generated if Lines 13 and 14 of Algorithm \ref{alg:ASto} are performed. However, for simplicity we suppose a sample $\xi_2^t$ is still generated in this case, although no quantity is determined by this sample.} we let $\mF_t = \sigma(\{\xi_1^j, \xi_2^j\}_{j=0}^t)$, $t\geq 0$, be the $\sigma$-algebra generated by all the samples till $t$; $\mF_{t-0.5} = \sigma(\{\xi_1^j, \xi_2^j\}_{j=0}^{t-1}\cup \xi_1^t)$, $t\geq 0$, be the $\sigma$-algebra generated by all the~samples till $t-1$ and the sample $\xi_1^t$; and $\mF_{-1}$ be the trivial $\sigma$-algebra generated by the initial iterate (which is deterministic). Throughout the presentation, we let $\barepsilon_t$ be the quantity obtained after Step 2; that is, $\barepsilon_t$ satisfies \eqref{cond:bound:fes:error} and \eqref{cond:decrease:1}.~With~this setup, it is easy to see that
\begin{equation}\label{sigma:alg}
\begin{aligned}
\sigma(\bx_t, \bmu_t, \blambda_t)\cup\sigma(\barnu_t)\cup \sigma(\baralpha_t)\cup\sigma(\bardelta_t) \subseteq & \mF_{t-1},\\
\sigma(\bx_{s_t}, \bmu_{s_t}, \blambda_{s_t})\cup\sigma(\barDelta_t, \hDelta_t, \cDelta_t)\cup \sigma(\barepsilon_t) \subseteq & \mF_{t-0.5}.
\end{aligned}
\end{equation}
We analyze Algorithm \ref{alg:ASto} in the next subsection.

\begin{algorithm}[!tb]
\caption{An Adaptive Active-Set StoSQP with Augmented Lagrangian}\label{alg:ASto}
\begin{algorithmic}[1]
\State \textbf{Input:} $(\bx_0, \bmu_0, \blambda_0)$, $\baralpha_0 = \alpha_{max}>0$, $\barepsilon_0$, $\bardelta_0$, $\eta$,  $\kappa_{grad}$, $\chi_{grad}$, $\chi_{f}$, $\chi_{err}>0$, $\rho>1$,~$\gamma_{B}\in(0, 1]$, $\beta$, $p_{grad}$, $p_f\in(0, 1)$, $\kappa_f\in(0, \beta/(4\alpha_{max})]$, $\barnu_0 = 2\sum_{i=1}^{r}\max\{(g_0)_i,0\}^3+1$;
		
\For{$t = 0,1,2\ldots$}
		
\State Generate $\xi_1^t$ so that \textbf{(a)} \eqref{cond:event:E1} holds; \textbf{(b)} \eqref{cond:grad:var} holds if $t-1$ is a successful~step; compute $\bnabla_{\bx}\mL_t$, $\barQ_{1,t}$, $\barQ_{2,t}$ as in \eqref{equ:def:Qmatrices}; \Comment{{\blue Step 1: estimate derivatives}}
\item[]

\While{\{\eqref{cond:bound:fes:error} does not hold\} OR \{\eqref{equ:SQP:direction} is solvable AND \eqref{cond:decrease:1} does not hold\}}
\State $\barepsilon_t \leftarrow \barepsilon_t/\rho$; \Comment{{\blue Step 2: set $\barepsilon_t$}}
\EndWhile
\item[]

\If{\{\eqref{equ:SQP:direction} is not solvable\} OR \{\eqref{equ:SQP:direction} is solvable AND   \eqref{cond:decrease:2} holds\}}
\State Obtain a backup direction $\hDelta_t$ and let $\cDelta_t = \hDelta_t$; \Comment{{\blue Step 3: decide $\cDelta_t$}}
\Else
\State $\cDelta_t = \barDelta_t$;
\EndIf
\item[]

\If{$\bx_{s_t} \notin \mT_{\barnu_t}$} \Comment{{\blue Step 4: estimate merit function}}
\State $(\bx_{t+1}, \bmu_{t+1}, \blambda_{t+1}) = (\bx_{t}, \bmu_{t}, \blambda_{t})$, $\baralpha_{t+1} = \baralpha_t$, $\bardelta_{t+1} = \bardelta_t$, $\barepsilon_{t+1} = \barepsilon_t$;
\State $\barnu_{t+1} = \rho^j\barnu_t$ with $j = \lceil\log(2a_{s_t}/\barnu_t)/\log \rho\rceil$;
\Else
\State Generate $\xi_2^t$ and compute $\barL_{\barepsilon_t, \barnu_t, \eta}^t$, $\barL_{\barepsilon_t, \barnu_t, \eta}^{s_t}$ so that \eqref{cond:event:E2} and \eqref{cond:merit:var} hold;
\If{$\barL_{\barepsilon_t, \barnu_t, \eta}^{s_t} \leq \barL_{\barepsilon_t, \barnu_t, \eta}^t + \beta\baralpha_t (\bnabla\mL_{\barepsilon_t, \barnu_t, \eta}^t)^T\cDelta_t$} \Comment{{\blue Step 5: line search}}
\State $(\bx_{t+1}, \bmu_{t+1}, \blambda_{t+1}) = (\bx_{s_t}, \bmu_{s_t}, \blambda_{s_t})$,  $\baralpha_{t+1} = \rho\baralpha_t\wedge \alpha_{max}$; \Comment{successful step}
\If{$-\beta\baralpha_t (\bnabla\mL_{\barepsilon_t, \barnu_t, \eta}^t)^T\cDelta_t \geq \bardelta_t$}\Comment{reliable step}
\State $\bardelta_{t+1} = \rho\bardelta_t$;
\Else\Comment{unreliable step}
\State $\bardelta_{t+1} = \bardelta_t/\rho$;
\EndIf
\Else \Comment{unsuccessful step}
\State $(\bx_{t+1}, \bmu_{t+1}, \blambda_{t+1}) = (\bx_{t}, \bmu_{t}, \blambda_{t})$, $\baralpha_{t+1} = \baralpha_t/\rho$, $\bardelta_{t+1} = \bardelta_t/\rho$;
\EndIf
\State $\barepsilon_{t+1} = \barepsilon_t$, $\barnu_{t+1} = \barnu_t$;
\EndIf
\EndFor
\end{algorithmic}
\end{algorithm}

\subsection{Assumptions and stability of parameters}\label{sec:4.2}

We study the stability of the parameter sequence $\{\barepsilon_t, \barnu_t\}_t$. We will show that,~for each run of the algorithm, the sequence is stabilized after a finite number of iterations. Thus, Lines 5 and 14 of Algorithm \ref{alg:ASto} will not be performed when the iteration index $t$ is large enough. We begin by introducing the assumptions.

\begin{assumption}[Regularity condition]\label{ass:3}
We assume the iterate $\{(\bx_t, \bmu_t, \blambda_t)\}$ and trial point $\{(\bx_{s_t}, \bmu_{s_t}, \blambda_{s_t})\}$ are contained in a convex compact region~$\mX\times \mM\\ \times\Lambda$. Further, if $\bx_{s_t}\in \mT_{\barnu_t}$, then the segment $\{\zeta\bx_t + (1-\zeta)\bx_{s_t}: \zeta\in(0, 1)\}\subseteq \mT_{\theta\barnu_t}$ for some $\theta \in [1, 2)$. We also assume the functions $f, g, c$ are thrice continuously differentiable over $\mX$, and realizations $|F(\bx, \xi)|$, $\|\nabla F(\bx, \xi)\|$, $\|\nabla^2 F(\bx, \xi)\|$ are uniformly bounded over $\bx\in \mX$ and $\xi\sim \P$.
\end{assumption}

\begin{assumption}[Constraint qualification]\label{ass:4}
For any $\bx \in \Omega$, we assume that $(J^T(\bx)\;\; G^T_{\I(\bx)}(\bx))$~has~full column rank, where $\Omega$ is the feasible set in \eqref{equ:feasible} and $\I(\bx)$ is the active set in~\eqref{equ:I}. For any $\bx\in \mX\backslash \Omega$, we assume the linear system
\begin{equation}\label{equ:LS}
\begin{aligned}
c_i(\bx) + \nabla^T c_i(\bx)\bz = \0, \quad\quad i: c_i(\bx) \neq 0,\\
g_i(\bx) + \nabla^Tg_i(\bx)\bz\leq \0, \quad\quad i: g_i(\bx)> 0,
\end{aligned}
\end{equation}
has a solution for $\bz\in\mR^d$. 
\end{assumption}

The boundedness condition on realizations in Assumption \ref{ass:3} is widely used in StoSQP analysis to have a well-behaved stochastic penalty parameter sequence \citep{Berahas2021Sequential, Na2022adaptive, Berahas2021Stochastic, Curtis2021Inexact}. The third derivatives of $f,g,c$ are only required in the analysis and not needed in the implementation. They are required since the existence of the (generalized) Hessian of the augmented Lagrangian needs the~third~derivatives. See, for example, \cite[Section 6]{Pillo2002Augmented} for the same requirement. For deterministic schemes, the compactness condition on the iterates is typical for the augmented Lagrangian and SQP analyses \citep[Chapter~4]{Bertsekas1982Constrained} \cite[Chapter 18]{Nocedal2006Numerical}. Some literature relaxed it by assuming all quantities (e.g., the objective gradient and constraints Jacobian, etc.) are uniformly upper bounded with a lower bounded objective~(so as the merit function). However, either condition is rather restrictive for StoSQP due to the underlying randomness of the~scheme.~That said,~given the StoSQP iterates presumably contract to a deterministic feasible set, we believe that an unbounded iteration sequence is rare in general. Furthermore, \mbox{compared to fully} \mbox{stochastic schemes in \citep{Berahas2021Sequential, Berahas2021Stochastic, Curtis2021Inexact}, we generate a batch of samples} to have a more precise estimation of the true model in each iteration; thus,~our~stochastic~iterates have a higher chance to closely track the underlying deterministic iterates.

The convexity of $\mM\times \Lambda$ can be removed by defining a closed convex hull $\overline{\text{conv}(\mM)} \times \overline{\text{conv}(\mM)}$.~However, the convexity of the set for the primal iterates is essential to enable a valid Taylor expansion. See \cite[Proposition 2.2 and~Section 4]{Pillo2011primal} \cite[Proposition 2.4 and (14)]{Pillo2005Convergence} and references therein for the same requirement for doing line search with \eqref{equ:aug:Lagrange} and applying its Taylor expansion.

In particular, by the design of Algorithm~\ref{alg:ASto}, we have $\bx_t \in \mT_{\barnu_t}$ for any $t$,~while the trial step $\bx_{s_t}$ may be outside $\mT_{\barnu_t}$. If $\bx_{s_t}\notin \mT_{\barnu_t}$, we enlarge $\barnu_t$ (Line 14) and rerun the iteration from the beginning. Assumption \ref{ass:3} states that if it turns out that $\bx_{s_t}\in \mT_{\barnu_t}$, then the whole segment $\zeta\bx_t+(1-\zeta)\bx_{s_t}$, which may not completely lie in $\mT_{\barnu_t}$ as $\mT_{\barnu_t}$ may be nonconvex, is supposed to lie in a larger space $\mT_{\theta\barnu_t}$ with $\theta \in[1,2)$. Since $\mL_{\barepsilon_t, \barnu_t, \eta}$ is SC$^1$ in $\mT_{2\barnu_t}^\circ \times \mR^m\times\mR^r$~and~$\mT_{\theta\barnu_t} \subseteq \mT_{2\barnu_t}^\circ$, where $\mT_{2\barnu_t}^\circ$ denotes the interior of $\mT_{2\barnu_t}$, the second-order Taylor expansion at $(\bx_t, \bmu_t, \blambda_t)$ is allowed \citep{Pillo2002Augmented}. Note that the range of $\theta$ is inessential. If we replace $\nu/2$ in \eqref{equ:domain} by $\nu/\kappa$ for any $\kappa>1$, then we would allow the existence of $\theta$ in $[1, \kappa)$.~In other words, $\theta$ can be as large as any $\kappa$. In fact, the condition on the segment always holds when the input $\alpha_{max}$,~the~upper bound of $\baralpha_t$ (cf. Line 18), is suitably upper bounded. Specifically, supposing $\sup_{\mX}\|\nabla a(\bx)\|\vee \sup_t\|\cDelta\bx_t\| \leq \Upsilon$ (ensured by compactness of iterates), for any $\theta> 1$ and $\zeta\in(0, 1)$, as long as $\alpha_{max} \leq (\theta-1)\barnu_0/(2\Upsilon^2)$, we have $\zeta\bx_t+(1-\zeta)\bx_{s_t} \in \mT_{\theta\barnu_t}$ by noting that
\begin{multline*}
a(\zeta\bx_t+(1-\zeta)\bx_{s_t}) = a(\bx_t + \baralpha_t(1-\zeta)\cDelta\bx_t)
\leq  a(\bx_t) + \baralpha_t(1-\zeta)\Upsilon^2\\ \leq\frac{\barnu_t}{2} + \alpha_{max}\Upsilon^2 \leq \frac{\barnu_t}{2} + \frac{(\theta-1)\barnu_0}{2} \leq \frac{\barnu_t}{2} + \frac{(\theta-1)\barnu_t}{2} = \frac{\theta\barnu_t}{2}.
\end{multline*}
Clearly, the condition on the segment is not required if $\mT_{\nu}$ in \eqref{equ:domain} is a convex set, which is the case, for example, if we have linear inequality constraints $\bx\leq \0$; or more generally, each $g_i(\cdot)$ is a convex function. We further investigate the effect of the range of $\theta$ by varying $\kappa$ ($\kappa = 2$ by default; cf. \eqref{equ:domain}) in the experiments.

By the compactness condition and noting that $\barnu_t$ is increased by at least a factor of $\rho$ each time in \eqref{cond:nu:t}, we immediately know that $\barnu_t$ stabilizes when $t$ is large. Moreover, if we let
\begin{equation}\label{bound:nu}
\tnu = \rho^{\tilde{j}}\barnu_0\quad \text{ with } \quad \tilde{j} = \lceil \log(2\max_{\mX}a(\bx)/\barnu_0)/\log \rho\rceil,
\end{equation}
then $\barnu_t \leq \tnu$, $t\geq 0$, almost surely. We will show a similar result for $\barepsilon_t$.

Assumption \ref{ass:4} imposes the constraint qualifications. In particular, for feasible points $\Omega$, we assume the linear independence constraint qualification~(LICQ), which is a standard condition to ensure the existence and uniqueness of the Lagrangian multiplier \citep{Nocedal2006Numerical}. For infeasible points $\mX\backslash \Omega$, we assume that the solution set of the linear system \eqref{equ:LS} is nonempty. The condition \eqref{equ:LS} restricts the behavior of the constraint functions outside the feasible set, which, together with~the compactness condition, implies $\Omega \neq \emptyset$ (cf. \cite[Proposition 2.5]{Lucidi1992New}). In fact, the condition \eqref{equ:LS} weakens the generalized Mangasarian-Fromovitz constraint qualification (MFCQ) \cite[Definition 2.5]{Xu2014Smoothing}; and relates to the weak MFCQ, which is proposed for problems with only inequalities in \cite[Definition 1]{Lucidi1992New} and adopted in \cite[Assumption A3]{Pillo2002Augmented} and \cite[Assumption 3.2]{Pillo2008truncated}. However, \cite{Lucidi1992New} requires the weak MFCQ to hold for feasible points in addition to LICQ;~while~\cite{Pillo2002Augmented, Pillo2008truncated} and this paper remove such a condition. The condition \eqref{equ:LS} simplifies and generalizes the weak MFCQ in \cite{Lucidi1992New, Pillo2002Augmented, Pillo2008truncated} by including equality constraints. We note that the weak MFCQ is slightly weaker than \eqref{equ:LS}. By the Gordan's theorem \citep{Goldman19574.}, \eqref{equ:LS} implies that $\{c_i\cdot \nabla c_i\}_{i: c_i\neq 0}\cup \{\nabla g_i\}_{i: g_i>0}$ are positively linearly independent:
\begin{equation*}
\sum_{i: c_i\neq 0} a_ic_i\nabla c_i + \sum_{i: g_i> 0}b_i \nabla g_i \neq \0,
\end{equation*}
for any coefficients $a_i, b_i\geq 0$ and $\sum_i a_i^2+b_i^2 >0$. In contrast, the weak~MFCQ only requires that the above linear combination is nonzero for a particular~set~of coefficients. However, we adopt the simplified but a bit stronger condition only because \eqref{equ:LS} has a cleaner form and a clearer connection to SQP subproblems. The coefficients of the weak MFCQ in \cite{Lucidi1992New,Pillo2002Augmented, Pillo2008truncated} are relatively hard~to~interpret. Instead of regarding the constraint qualification as the essence of~constraints, those coefficients depend on particular choice of the merit function, although that assumption statement is sharper. That said, \eqref{equ:LS} is still weaker than other literature on the augmented Lagrangian \citep{Pillo1982new, Pillo1986exact, Lucidi1988New}; and weaker than~what~is widely assumed in SQP analysis \citep{Boggs1995Sequential}, where the IQP system, $c_i + \nabla^Tc_i\bz = \0$, $1\leq i\leq m$, $g_i + \nabla^Tg_i\bz \leq \0$, $1\leq i\leq r$, is supposed to have a solution.~Moreover, we do not require the strict complementary condition, which is~often~imposed for the merit functions that apply (squared) slack variables to transform nonlinear inequality~constraints  \citep[A2]{Zavala2014Scalable}, \citep[Proposition 3.8]{Fukuda2017note}.

The first lemma shows that $\eqref{cond:bound:fes:error}$ is satisfied for a sufficiently small $\barepsilon_t$.~Although \eqref{cond:bound:fes:error} is inspired by \cite[(19)]{Na2022adaptive} for equalities, the proof is quite different from that paper (cf. Lemma 4 there).

\begin{lemma}\label{lem:7}
Under Assumptions \ref{ass:3} and \ref{ass:4}, there exists a deterministic threshold $\tilde{\epsilon}_1>0$ such that \eqref{cond:bound:fes:error} holds for any $\barepsilon_t \leq \tepsilon_1$.
\end{lemma}
\vskip-3pt
\begin{proof}
See Appendix \ref{pf:lem:7}.
\end{proof}

The second lemma shows that \eqref{cond:decrease:1} is satisfied for small $\barepsilon_t$. The analysis is similar to Lemma \ref{lem:6}. We need the following condition on the SQP system \eqref{equ:SQP:direction}.

\begin{assumption}\label{ass:5}
We assume that, whenever \eqref{equ:SQP:direction} is solvable, $(J_t^T\; G_{t_a}^T)$ has full column rank, and there exist positive constants $\Upsilon_{B}\geq 1\geq  \gamma_{B}\vee \gamma_{H}$ such~that
\begin{equation*}
B_t\preceq \Upsilon_{B} I, \quad\quad M_t\succeq\gamma_{H} I, \quad\quad \begin{pmatrix}
J_t\\
G_{t_a}
\end{pmatrix}\begin{pmatrix}
J_t^T & G_{t_a}^T
\end{pmatrix}\succeq \gamma_{H} I,
\end{equation*}
and $\bz^TB_t\bz \geq \gamma_{B}\|\bz\|^2$, $\forall \bz\in\{\bz\in\mR^d:  J_t\bz = \0, G_{t_a}\bz = \0\}$.
\end{assumption}

Assumption \ref{ass:5} summarizes Assumptions \ref{ass:1} and \ref{ass:2}. As shown in Lemma \ref{lem:6}, the conditions on $M_t$ and $(J_t^T\; G_{t_a}^T)$ hold locally. For the presented global~analysis, the Hessian approximation $B_t$ is easy to construct to satisfy the condition, e.g., $B_t = I$; however, such a choice is not proper for fast local rates. In practice, given a lower bound $\gamma_{B}>0$, $B_t$ is constructed by doing a regularization on a subsampled Hessian (e.g., for finite-sum objectives) or a sketched~Hessian~(e.g.,~for regression objectives), which can preserve certain~second-order information and~be obtained with less expense. With Assumption \ref{ass:5}, we have~\mbox{the following result}.

\begin{lemma}\label{lem:8}
Under Assumptions \ref{ass:3} and \ref{ass:5}, there exists a deterministic threshold $\tilde{\epsilon}_2>0$ such that \eqref{cond:decrease:1} holds for any $\barepsilon_t \leq \tepsilon_2$.
\end{lemma}
\vskip-3pt
\begin{proof}
See Appendix \ref{pf:lem:8}.
\end{proof}

We summarize \eqref{bound:nu}, Lemmas \ref{lem:7} and \ref{lem:8} in the next theorem.

\begin{theorem}\label{thm:2a}
Under Assumptions \ref{ass:3}, \ref{ass:4}, and \ref{ass:5}, there exist deterministic thresholds $\tnu$, $\tepsilon>0$ such that $\{\barnu_t, \barepsilon_t\}_t$ generated by Algorithm \ref{alg:ASto} satisfy $\barnu_t \leq \tnu$, $\barepsilon_t \geq \tepsilon$. Moreover, almost surely, there exists an iteration threshold $\bar{t}<\infty$, such that $\barepsilon_t = \barepsilon_{\bar{t}}$, $\barnu_t = \barnu_{\bar{t}}$, $ t\geq \bar{t}$.
\end{theorem}
\vskip-3pt
\begin{proof}
The existence of $\tnu$ is showed in \eqref{bound:nu}. By Lemmas \ref{lem:7} and \ref{lem:8}, and defining $\tepsilon = (\tepsilon_1\wedge \tepsilon_2)/\rho$, we show the existence of $\tepsilon$. The existence of the iteration threshold $\bart$ is ensured by noting that $\{\barnu_t, 1/\barepsilon_t\}_t$ are bounded from above; and each update increases the parameters by at least a factor of $\rho>1$.
\end{proof}

We mention that the iteration threshold $\bart$ is random for stochastic schemes and it changes between different runs. However, it always exists. The following analysis supposes $t$ is large enough such that $t\geq \bart$ and $\barepsilon_t, \barnu_t$ have stabilized. We condition our analysis on the $\sigma$-algebra $\mF_{\bart}$, which means that we only~consider the randomness of the generated samples after $\bart+1$ iterations and, by \eqref{sigma:alg}, the parameters $\barepsilon_{\bart}, \barnu_{\bart}$ are fixed. We should point out that, although it is standard~to focus only on the tail of the iteration sequence to show the global convergence (even for the deterministic case \cite[Theorem 18.3]{Nocedal2006Numerical}), an important aspect that is missed by such an analysis is the non-asymptotic guarantees. In particular, we know the scheme~changes~the~merit~\mbox{parameters}~for~at~most~$\log(\tnu\barepsilon_0/(\barnu_0\tepsilon))/\log(\rho)$ times; however, how many iterations it spans for all the changes is not answered by our analysis. Establishing a bound on $\bart$ in expectation or high probability sense would help us further understand the efficiency of the scheme. However, since any characterization of $\bart$ is difficult even for deterministic schemes, we leave such a study to the future. \mbox{Another missing aspect is the iteration complexity}, where we are interested in the number of iterations to attain an $\epsilon$-first- or~second-order stationary point (we abuse $\epsilon$ notation here to refer to the accuracy level). The iteration complexity is recently studied for two StoSQP schemes under~very particular setups \citep{Curtis2021Worst, Berahas2022Adaptive}; none of the existing works allow either stochastic~line search or inequality constraints. We leave the iteration complexity of our scheme to the future as well.

\subsection{Convergence analysis}

We conduct the global convergence analysis for Algorithm \ref{alg:ASto}. We prove that $\lim_{t\rightarrow \infty} R_t = 0$ \textit{almost surely}, where $R_t = \|(\nabla_{\bx}\mL_t, c_t, \max\{g_t,-\blambda_t\})\|$ is~the~KKT residual. We suppose the line search conditions \eqref{cond:event:E1}, \eqref{cond:grad:var}, \eqref{cond:event:E2}, \eqref{cond:merit:var} hold. We will discuss the sample complexities that ensure these generic conditions in~Section \ref{sec:4.3}. It is fairly easy to see that all conditions hold for large batch sizes.

Our proof structure closely follows \citep{Na2022adaptive}. The analyses are more involved in Lemmas \ref{lem:10}, \ref{lem:12}, \ref{lem:13}, \ref{lem:14} and Theorem \ref{thm:3}, which account for the differences between equality and inequality constraints, and account for our relaxations~of~the~feasibility error condition and the increasing sample size requirement of \cite{Na2022adaptive}. The analysis in Theorem \ref{thm:5} is new, which strengthens the ``liminf" \mbox{convergence}~in~\cite{Na2022adaptive}. The analyses are slightly adjusted in Theorem \ref{thm:4}, and the same in Lemma \ref{lem:11} and Theorem \ref{thm:2}. The adopted potential function (or Lyapunov function) is
\begin{equation}\label{equ:potential}
\Theta_{\barepsilon_{\bart}, \barnu_{\bart}, \eta, \omega}^t = \omega\mL_{\barepsilon_{\bart}, \barnu_{\bart}, \eta}^t + \frac{1-\omega}{2}\baralpha_t\|\nabla \mL_{\barepsilon_{\bart}, \barnu_{\bart}, \eta}^t\|^2 + \frac{1-\omega}{2}\bardelta_t, \quad t\geq \bart + 1,
\end{equation}
where $\omega \in(0, 1)$ is a coefficient to be specified later. We note that using $\mL_{\barepsilon_{\bart}, \barnu_{\bart}, \eta}^t$ by itself (i.e., $\omega=1$) to monitor the iteration progress is not suitable for the stochastic setting; it is possible that $\mL_{\barepsilon_{\bart}, \barnu_{\bart}, \eta}^t$ increases while $\barL_{\barepsilon_{\bart}, \barnu_{\bart}, \eta}^t$ decreases. In contrast, $\Theta_{\barepsilon_{\bart}, \barnu_{\bart}, \eta, \omega}^t$ linearly combines different components and has a composite measure of the progress. For example, the decrease of $\Theta_{\barepsilon_{\bart}, \barnu_{\bart}, \eta, \omega}^t$ may come from $\bardelta_t$ (Lines 22 and 25 of Algorithm \ref{alg:ASto}).

Since parameters $\barepsilon_{\bar{t}}, \barnu_{\bart}, \eta$ in $\mL_{\barepsilon_{\bart}, \barnu_{\bart}, \eta}$ are fixed (conditional on $\mF_{\bart}$), we denote $\Theta_{\omega}^t = \Theta_{\barepsilon_{\bart}, \barnu_{\bart},\eta, \omega}^t$ for notational simplicity. In the presentation of theoretical~results, we only track the parameters $(\beta, \alpha_{max}, \kappa_{grad}, \kappa_{f}, p_{grad}, p_{f}, \chi_{grad}, \chi_{f})$~that relate to the line search conditions. In particular, we use $C_1, C_2\ldots$ and $\Upsilon_1, \Upsilon_2\ldots$ to denote deterministic constants that are independent from these parameters, but may depend on $(\gamma_{B}, \gamma_{H}, \Upsilon_B, \chi_{u}, \chi_{err},\rho, \eta, \barepsilon_0, \barnu_0)$, and thus~\mbox{depend}~on~the~deterministic thresholds $\tepsilon$ and $\tnu$. 
Recall that $(\gamma_{B}, \gamma_{H}, \Upsilon_B, \chi_{u})$ come from Assumption \ref{ass:5} and \eqref{cond:alter:dir}, while $(\chi_{err}, \rho, \eta, \barepsilon_0, \barnu_0)$ are any algorithm inputs.

The first lemma presents a preliminary result.

\begin{lemma}\label{lem:9}
Under Assumptions \ref{ass:3}, \ref{ass:4}, \ref{ass:5}, the following results hold deterministically conditional on $\mF_{t-1}$.
\begin{enumerate}[label=(\alph*),topsep=0pt]
\setlength\itemsep{0.0em}
\item There exists $C_1>0$ such that the following two inequalities hold for any iteration $t \geq 0$ ((a2) also holds for $s_t$), any parameters $\epsilon,\nu$, and any~generated sample set $\xi$:

\noindent \textbf{(a1)} $\nbr{\bnabla\mL_{\epsilon, \nu, \eta}^t - \nabla\mL_{\epsilon, \nu, \eta}^t } \leq C_1\cbr{\nbr{\bnabla f_t - \nabla f_t} \vee (\barR_t\wedge 1)}\cdot \nbr{\bnabla^2 f_t - \nabla^2 f_t} $;

\noindent\textbf{(a2)} $\abr{\barL_{\epsilon, \nu, \eta}^t - \mL_{\epsilon, \nu, \eta}^t } \leq C_1\{|\barf_t - f_t| \vee [(\barR_t\vee \|\bnabla f_t-\nabla f_t\|)\wedge 1]\cdot \nbr{\bnabla f_t - \nabla f_t} \}$.

\item There exists $C_2>0$ such that for any $t\geq 0$ and set $\xi$,
\begin{equation*}
\nbr{\bnabla_{\bx}\mL_t} \leq C_2\cbr{\|\bnabla\mL_{\barepsilon_t, \barnu_t, \eta}^t\| + \nbr{(
c_t,\;\bw_{\barepsilon_t, \barnu_t}^t) } }.
\end{equation*}

\item There exists $C_3>0$ such that for any $t\geq 0$ and set $\xi$, if \eqref{equ:SQP:direction} is solvable,~then 
\begin{equation*}
\nbr{\bnabla\mL_{\barepsilon_t, \barnu_t, \eta}^t} \leq C_3\nbr{\left(\begin{smallmatrix}
\barDelta\bx_t\\
J_t\bnabla_{\bx}\mL_t\\
G_t\bnabla_{\bx}\mL_t + \Pi_c(\diag^2(g_t)\blambda_t)
\end{smallmatrix} \right)}.
\end{equation*}
\end{enumerate}
\end{lemma}

\vskip-3pt
\begin{proof}
See Appendix \ref{pf:lem:9}.
\end{proof}

The results in Lemma \ref{lem:9} hold deterministically conditional on $\mF_{t-1}$, because the samples $\xi$ for computing $\bnabla\mL_{\barepsilon_t, \barnu_t, \eta}^t$, $\bnabla_{\bx}\mL_t$ are supposed to be also~given~by~the statement. The following result suggests that if both the gradient $\nabla\mL_{\barepsilon_{\bart}, \barnu_{\bart}, \eta}^t$ and the function evaluations $\mL_{\barepsilon_{\bart}, \barnu_{\bart}, \eta}^t$, $\mL_{\barepsilon_{\bart}, \barnu_{\bart}, \eta}^{s_t}$ are precisely estimated, in the sense that the event $\E_1^t\cap \E_2^t$ happens (cf. \eqref{event:E1}, \eqref{event:E2}), then there is a uniform~lower bound on $\baralpha_t$ to make the Armijo condition hold.

\begin{lemma}\label{lem:10}
For $t\geq \bart + 1$, suppose $\E_1^t\cap \E_2^t$ happens. There exists $\Upsilon_1>0$ such that the $t$-th step satisfies the Armijo condition \eqref{cond:armijo} (i.e., is a successful step)~if
\begin{equation*}
\baralpha_t \leq \frac{1-\beta}{\Upsilon_1(\kappa_{grad}+\kappa_{f} + 1)}.
\end{equation*}
\end{lemma}

\vskip-3pt
\begin{proof}
See Appendix \ref{pf:lem:10}.
\end{proof}

The next result suggests that, if only the function evaluations $\mL_{\barepsilon_{\bart}, \barnu_{\bart}, \eta}^t$,~$\mL_{\barepsilon_{\bart}, \barnu_{\bart}, \eta}^{s_t}$ are precisely estimated, in the sense that the event $\E_2^t$ happens, then a sufficient decrease of $\barL_{\barepsilon_{\bart}, \barnu_{\bart}, \eta}^t$ implies a sufficient decrease of $\mL_{\barepsilon_{\bart}, \barnu_{\bart}, \eta}^t$. The proof directly follows \cite[Lemma 6]{Na2022adaptive}, and thus is omitted.

\begin{lemma}\label{lem:11}
For $t\geq \bart + 1$, suppose $\E_2^t$ happens. If the $t$-th step satisfies the Armijo condition \eqref{cond:armijo}, then
\begin{equation*}
\mL_{\barepsilon_{\bart}, \barnu_{\bart}, \eta}^{s_t} \leq \mL_{\barepsilon_{\bart}, \barnu_{\bart}, \eta}^t + \frac{\baralpha_t\beta}{2}(\bnabla\mL_{\barepsilon_{\bart}, \barnu_{\bart}, \eta}^t)^T\cDelta_t.
\end{equation*}
\end{lemma}

Based on Lemmas \ref{lem:10} and \ref{lem:11}, we now establish an error recursion for the potential function $\Theta_{\omega}^t$ in \eqref{equ:potential}. Our analysis is separated into three cases according to the events: $\E_1^t\cap\E_2^t$, $(\E_1^t)^c\cap \E_2^t$ and $(\E_2^t)^c$. We will show that $\Theta_{\omega}^t$ decreases in the case of $\E_1^t\cap\E_2^t$, while may increase in the other two cases. Fortunately, by letting $p_{grad}$ and $p_{f}$ be small, $\Theta_{\omega}^t$ always decreases in expectation.

We first show in Lemma \ref{lem:12} that $\Theta_{\omega}^t$ decreases when $\E_1^t\cap\E_2^t$ happens. We~note that the decrease of $\Theta_{\omega}^t$ exceeds $\baralpha_t\|\nabla\mL_{\barepsilon_{\bart}, \barnu_{\bart}, \eta}^t\|^2$ by $\bardelta_t$ (up to a multiplier).

\begin{lemma}\label{lem:12}
For $t\geq \bart+1$, suppose $\E_1^t\cap\E_2^t$ happens. There exists $\Upsilon_2>0$, such that if $\omega$ satisfies
\begin{equation}\label{cond:omega}
\frac{1-\omega}{\omega} \leq \frac{\beta}{\Upsilon_2(\kappa_{grad}\alpha_{max}+\alpha_{max} + 1)^2} \wedge \frac{1}{18(\rho-1)},
\end{equation}
then
\begin{equation*}
\Theta_{\omega}^{t+1} - \Theta_{\omega}^t \leq -\frac{1}{2}\rbr{1-\omega}\rbr{1 - \frac{1}{\rho}} \rbr{\baralpha_t\nbr{\nabla\mL_{\barepsilon_{\bart}, \barnu_{\bart}, \eta}^t}^2 + \bardelta_t}.
\end{equation*}
\end{lemma}

\vskip-3pt
\begin{proof}
See Appendix \ref{pf:lem:12}.
\end{proof}

We then show in Lemma \ref{lem:13} that $\Theta_{\omega}^t$ may increase, if $\nabla\mL_{\barepsilon_{\bart}, \barnu_{\bart}, \eta}^t$ is not precisely estimated (i.e., $(\E_1^t)^c$ happens) but $\mL_{\barepsilon_{\bart}, \barnu_{\bart}, \eta}^t$, $\mL_{\barepsilon_{\bart}, \barnu_{\bart}, \eta}^{s_t}$ are precisely estimated (i.e., $\E_2^t$ happens). The increase is proportional to $\baralpha_t\|\nabla\mL_{\barepsilon_{\bart}, \barnu_{\bart}, \eta}^t\|^2$.

\begin{lemma}\label{lem:13}
For $t\geq \bart+1$, suppose $(\E_1^t)^c\cap \E_2^t$ happens. Under \eqref{cond:omega}, we have
\begin{equation*}
\Theta_{\omega}^{t+1} - \Theta_{\omega}^t \leq \rho(1-\omega)\baralpha_t\nbr{\nabla\mL_{\barepsilon_{\bart}, \barnu_{\bart}, \eta}^t}^2.
\end{equation*}
\end{lemma}

\vskip-3pt
\begin{proof}
See Appendix \ref{pf:lem:13}.
\end{proof}

We finally show in Lemma \ref{lem:14} that $\Theta_{\omega}^t$ increases and the increase can exceed $\baralpha_t\|\nabla\mL_{\barepsilon_{\bart}, \barnu_{\bart}, \eta}^t\|^2$, if $\mL_{\barepsilon_{\bart}, \barnu_{\bart}, \eta}^t$, $\mL_{\barepsilon_{\bart}, \barnu_{\bart}, \eta}^{s_t}$ are not precisely estimated. In this case, the exceeding terms have to be controlled by making use of the condition \eqref{cond:merit:var}.

\begin{lemma}\label{lem:14}
For $t\geq \bart + 1$, suppose $(\E_2^t)^c$ happens. Under \eqref{cond:omega}, we have
\begin{equation*}
\Theta_{\omega}^{t+1} - \Theta_{\omega}^t \leq \rho(1-\omega)\baralpha_t\nbr{\nabla\mL_{\barepsilon_{\bart}, \barnu_{\bart}, \eta}^t}^2 + \omega\cbr{\abr{\barL_{\barepsilon_{\bart}, \barnu_{\bart}, \eta}^{s_t} - \mL_{\barepsilon_{\bart}, \barnu_{\bart}, \eta}^{s_t}} + \abr{\barL_{\barepsilon_{\bart}, \barnu_{\bart}, \eta}^t - \mL_{\barepsilon_{\bart}, \barnu_{\bart}, \eta}^t}}.
\end{equation*}
\end{lemma}

\vskip-3pt
\begin{proof}
See Appendix \ref{pf:lem:14}.
\end{proof}

Combining Lemmas \ref{lem:12}, \ref{lem:13}, \ref{lem:14}, we derive the one-step error recursion of $\Theta_{\omega}^t$. The proof directly follows that of \cite[Theorem 2]{Na2022adaptive} and is omitted.

\begin{theorem}[One-step error recursion]\label{thm:2}
For $t\geq \bart+1$, suppose $\omega$ satisfies \eqref{cond:omega} and $p_{grad}$ and $p_{f}$ satisfy
\begin{equation}\label{cond:p}
\frac{p_{grad} + \sqrt{{ (1\vee \chi_{f})}\cdot p_{f}}}{(1-p_{grad})(1-p_{f})} \leq \frac{\rho-1}{8\rho}\cbr{\frac{1}{\rho}\wedge \frac{1-\omega}{\omega}}.
\end{equation}
Then
\begin{equation*}
\mE\sbr{\Theta_{\omega}^{t+1} - \Theta_{\omega}^t\mid \mF_{t-1}} 
\leq -\frac{1}{4}(1-p_{grad})(1-p_{f})(1-\omega)\rbr{1-\frac{1}{\rho}}\rbr{\bardelta_t + \baralpha_t\nbr{\nabla\mL_{\barepsilon_{\bart}, \barnu_{\bart}, \eta}^t}^2}.
\end{equation*}
\end{theorem}

With Theorem \ref{thm:2}, we derive the convergence of $\baralpha_tR_t^2$ in the next theorem, where $R_t = \|(\nabla_{\bx}\mL_t, c_t, \max\{g_t,-\blambda_t\})\|$ is the KKT residual.

\begin{theorem}\label{thm:3}
Under the conditions of Theorem \ref{thm:2}, $\lim\limits_{t\rightarrow \infty}\baralpha_t R_t^2 = 0$ almost surely.
\end{theorem}

\vskip-3pt
\begin{proof}
See Appendix \ref{pf:thm:3}.
\end{proof}

Then, we show that the ``liminf" of the KKT residuals converges to zero.

\begin{theorem}[``liminf" convergence]\label{thm:4}
Consider Algorithm \ref{alg:ASto} under \mbox{Assumptions} \ref{ass:3}, \ref{ass:4}, \ref{ass:5}. Suppose $\omega$ satisfies \eqref{cond:omega} and $p_{grad}, p_{f}$ satisfy \eqref{cond:p}. Then, almost surely, we have that $\liminf_{t\rightarrow \infty} R_t = 0$.
\end{theorem}

\vskip-4.6pt
\begin{proof}
See Appendix \ref{pf:thm:4}.
\end{proof}

Finally, we strengthen the statement in Theorem \ref{thm:4} and complete the global convergence analysis of Algorithm \ref{alg:ASto}.

\begin{theorem}[Global convergence]\label{thm:5}
Under the same conditions of Theorem \ref{thm:4}, we have that
\begin{equation*}
\lim\limits_{t\rightarrow\infty}R_t = 0, \quad \text{almost surely}.
\end{equation*}
\end{theorem}

\vskip-4.6pt
\begin{proof}
See Appendix \ref{pf:thm:5}.
\end{proof}

Our analysis generalizes the results of \citep{Na2022adaptive} to inequality constrained problems. The ``lim" convergence guarantee in Theorem \ref{thm:5} strengthens~the~existing ``liminf" convergence guarantee of stochastic line search for both unconstrained problems \cite[Theorem 4.10]{Paquette2020Stochastic} and equality constrained problems \cite[\mbox{Theorem}~4]{Na2022adaptive}. Theorem \ref{thm:5} also differs from the results in \cite{Berahas2021Sequential, Berahas2021Stochastic, Curtis2021Inexact}, where the authors showed the (liminf) convergence of the \textit{expected} KKT residual under a fully stochastic setup. Compared to \cite{Berahas2021Sequential, Berahas2021Stochastic, Curtis2021Inexact}, our scheme does not tune a deterministic sequence that controls the stepsizes and determines the convergence behavior (i.e., converging to a KKT point or only its neighborhood). Our scheme tunes two probability parameters $p_{grad}, p_f$. Seeing from \eqref{cond:omega} and \eqref{cond:p}, the upper bound conditions on $p_{grad}, p_f$ depend on the inputs $(\rho, \beta, \kappa_{grad}, \alpha_{max})$ and a universal constant $\Upsilon_2$. Estimating $\Upsilon_2$ is often difficult in practice; however, $p_{grad}, p_f$ affect the algorithm's performance only via the generated batch sizes, and the batch sizes depend on $p_{grad}, p_f$ only via the logarithmic factors (see \eqref{cond:xi1} and \eqref{cond:xi2} later). Thus, the algorithm is robust to $p_{grad}, p_f$. We will also empirically test the robustness to parameters for Algorithm \ref{alg:ASto} in Section \ref{sec:5}. In addition, \eqref{cond:omega} and \eqref{cond:p} suggest that the larger the parameters $(\rho, 1/\beta, \kappa_{grad}, \alpha_{max})$ we use, the smaller the probabilities $p_{grad}, p_f$ have to be. Such a dependence is consistent with the general intuition: the algorithm performs more aggressive updates with less restrictive Armijo condition when $(\rho, 1/\beta, \kappa_{grad}, \alpha_{max})$ are large; thus, a more precise model estimation in each iteration is desired in this case.

\subsection{Discussion on sample complexities}\label{sec:4.3}

As introduced in Section \ref{sec:1}, the stochastic line search is performed by generating a batch of samples in~each iteration to have a precise model estimation, which is standard in the literature \citep{Friedlander2012Hybrid, Byrd2012Sample, Krejic2013Line, De2017Automated, Bollapragada2018Adaptive, Paquette2020Stochastic, Cartis2017Global}. The batch sizes are adaptively controlled based on the iteration progress. We now discuss the batch sizes $|\xi_1^t|$ and $|\xi_2^t|$ to ensure the generic conditions \eqref{cond:event:E1}, \eqref{cond:grad:var}, \eqref{cond:event:E2},~\eqref{cond:merit:var}~of~Algorithm~\ref{alg:ASto}. We show that, if the KKT residual $R_t$ does not~vanish, all the~conditions are satisfied by properly choosing $|\xi_1^t|$ and $|\xi_2^t|$.

\vskip3pt
\noindent\textbf{Sample complexity of $\xi_1^t$.} The samples $\xi_1^t$ are used to estimate $\nabla f_t$ and~$\nabla^2f_t$ in Step 1 of Algorithm \ref{alg:ASto}. The estimators $\bnabla f_t$ and $\bnabla^2 f_t$ can be computed with different amount of samples, and their samples may or may not be independent. Let us suppose $\bnabla f_t$ is computed by samples $\xi_1^t$, while $\bnabla^2 f_t$ is computed by a subset of samples $\tau_1^t\subseteq \xi_1^t$. The case where $\bnabla f_t$ and $\bnabla^2 f_t$ are computed by two disjoint subsets of $\xi_1^t$ can be studied following the same analysis. We define
\begin{equation*}
\bnabla f_t = \frac{1}{|\xi_1^t|}\sum_{\xi \in \xi_1^t} \nabla F(\bx_t; \xi),\quad\quad \bnabla ^2 f_t = \frac{1}{ |\tau_1^t|}\sum_{\xi\in \tau_1^t}\nabla^2 F(\bx_t;\xi).
\end{equation*}
By Lemma \ref{lem:9}(a1), we know that \eqref{cond:event:E1} holds if, with probability $1 - p_{grad}$,
\begin{equation}\label{equ:gradHes}
\hskip -1pt\|\bnabla f_t - \nabla f_t\| \leq O(\kappa_{grad}\baralpha_t\barR_t),\;\; \|\bnabla^2 f_t - \nabla^2 f_t\| \leq O(\kappa_{grad}\baralpha_t\barR_t/ (\barR_t\wedge 1)),
\end{equation}
where we suppress universal constants (such as the variance of a single sample) in $O(\cdot)$ notation. By matrix Bernstein inequality \citep[Theorem 7.7.1]{Tropp2015Introduction}, \eqref{equ:gradHes} is satisfied if
\begin{equation}\label{Step:1:a}
|\xi_1^t| \geq O\rbr{\frac{\log(d/p_{grad})}{\kappa_{grad}^2\baralpha_t^2\barR_t^2}}\quad \text{ and }\quad |\tau_1^t| \geq (\barR_t^2\wedge 1) \cdot |\xi_1^t|.
\end{equation}
Furthermore, we use the bound $\mE[\|\bnabla^2f_t - \nabla^2 f_t\|^2\mid \mF_{t-1}] \leq O(\log d/|\tau_1^t|)$~(cf.~\cite[(6.1.6)]{Tropp2015Introduction}) and know that \eqref{cond:grad:var} holds if
\begin{equation}\label{Step:1:b}
|\xi_1^t|\geq O\rbr{\frac{\baralpha_t\log d}{\chi_{grad}^2\bardelta_t}}\quad \text{ and }\quad |\tau_1^t| \geq (\barR_t^2\wedge 1)\cdot |\xi_1^t|.
\end{equation}
\hskip-3.6pt Combining \eqref{Step:1:a} and \eqref{Step:1:b} together, we know that the conditions \eqref{cond:event:E1} and~\eqref{cond:grad:var}~are satisfied if
\begin{equation}\label{cond:xi1}
|\xi_1^t| \geq O\rbr{\frac{\log(d/p_{grad})}{\kappa_{grad}^2\baralpha_t^2\barR_t^2 \wedge \chi_{grad}^2\bardelta_t/\baralpha_t} } , \quad |\tau_1^t| \geq (\barR_t^2\wedge 1)\cdot|\xi_1^t|.
\end{equation}
Since \eqref{cond:grad:var} is imposed only when $t-1$ is a successful step, the term $\chi_{grad}^2\bardelta_t/\baralpha_t$~on the denominator in \eqref{cond:xi1} can be removed when $t-1$ is an unsuccessful step. In contrast to \cite{Na2022adaptive}, where the gradient $\nabla f_t$ and Hessian $\nabla^2 f_t$ are \mbox{computed}~based~on the same set of samples, we sharpen the calculation and realize that the batch size $|\tau_1^t|$ for $\nabla^2 f_t$ can be significantly less than $|\xi_1^t|$ for $\nabla f_t$. When $\barR_t$ gets close to zero, the ratio $|\tau_1^t|/|\xi_1^t|$ will also decay to zero.

We mention that $\barR_t$ on the right hand side of the condition $|\xi_1^t|$ in $\eqref{cond:xi1}$ has to be computed by samples $\xi_1^t$. A practical algorithm can first specify $\xi_1^t$, then compute $\barR_t$, and finally check if \eqref{cond:xi1} holds. For example, a While~loop~can~be designed to gradually increase $|\xi_1^t|$ until \eqref{cond:xi1} holds (cf. \cite[Algorithm 4]{Na2022adaptive}).~Such~a While loop always terminates in finite time when $R_t>0$, because $\barR_t\rightarrow R_t$ as $|\xi_1^t|$ increases (by the law of large number) so that the right hand side of \eqref{cond:xi1} does not diverge.

\vskip3pt
\noindent\textbf{Sample complexity of $\xi_2^t$.} The samples $\xi_2^t$ are used to estimate $f_t, f_{s_t}, \nabla f_t, \nabla f_{s_t}$ in Step 4 of Algorithm \ref{alg:ASto}. Similar to the discussion above, the estimators $\barf_t, \barf_{s_t}$ and $\bbnabla f_t, \bbnabla f_{s_t}$ can be computed with different amount of samples, and their samples may or may not be independent. Let us suppose $\barf_t, \barf_{s_t}$ are computed by samples $\xi_2^t$, while $\bbnabla f_t, \bbnabla f_{s_t}$ are computed by a subset of samples $\tau_2^t\subseteq\xi_2^t$. We define (similar for $\barf_{s_t}, \bbnabla f_{s_t}$)
\begin{equation*}
\barf_t = \frac{1}{|\xi_2^t|}\sum_{\xi \in \xi_2^t} F(\bx_t;\xi),\quad\quad \bbnabla f_t = \frac{1}{ |\tau_2^t|}\sum_{\xi \in  \tau_2^t}\nabla F(\bx_t;\xi).
\end{equation*}
By Lemma \ref{lem:9}(a2), we know that \eqref{cond:event:E2} holds if, with probability $1-p_f$,
\begin{align}\label{equ:fgrad}
|\barf_t - f_t|\vee |\barf_{s_t} - f_{s_t}| & \leq O(-\kappa_f\baralpha_t^2(\bnabla \mL_{\barepsilon_t,\barnu_t,\eta}^t)^T\cDelta_t),\\
\|\bbnabla f_t - \nabla f_t\|\vee \|\bbnabla f_{s_t} - \nabla f_{s_t}\| & \leq \text{\footnotesize $O\cbr{\frac{-\kappa_f\baralpha_t^2(\bnabla \mL_{\barepsilon_t,\barnu_t,\eta}^t)^T\cDelta_t}{ \cbr{\bar{\barR}_t\vee \bar{\barR}_{s_t}\vee \{-\kappa_f\baralpha_t^2(\bnabla \mL_{\barepsilon_t,\barnu_t,\eta}^t)^T\cDelta_t\}^{1/2}}  \wedge 1  }}$ } ,
\end{align}
where $\bar{\barR}_t$ and $\bar{\barR}_{s_t}$ are computed by $\tau_2^t$ and we use the fact that, for scalars $a,b$, $(a\wedge 1)\vee(b\wedge 1) = (a\vee b)\wedge1$. By Bernstein inequality, \eqref{equ:fgrad} is satisfied if
\begin{equation}\label{Step:2:a}
\begin{aligned}
& |\xi_2^t|  \geq O\rbr{\frac{\log(d/p_f)}{\kappa_f^2\baralpha_t^4\{(\bnabla\mL_{\barepsilon_t, \barnu_t, \eta}^t)^T\cDelta_t\}^2}},\\
& |\tau_2^t| \geq 
\big( \{\bar{\barR}_t^2\vee \bar{\barR}_{s_t}^2\vee -\kappa_f\baralpha_t^2(\bnabla \mL_{\barepsilon_t,\barnu_t,\eta}^t)^T\cDelta_t \}\wedge 1 \big)\cdot |\xi_2^t| \\
& \hskip 0.6cm = \big(\{\bar{\barR}_t^2\vee \bar{\barR}_{s_t}^2\}\cdot|\xi_2^t| \vee \{\log(d/p_f)\cdot|\xi_2^t|\}^{1/2} \; \big) \wedge|\xi_2^t|. 
\end{aligned}
\end{equation}
Moreover, by $\mE[|\barf_t-f_t|\mid \mF_{t-0.5}]\leq O(1/|\xi_2^t|)$ and $\mE[\|\bbnabla f_t - \nabla f_t\|^4]\leq O(1/|\tau_2^t|^2)$, we can see that \eqref{cond:merit:var} holds if
\begin{equation}\label{Step:2:b}
|\xi_2^t|\geq O(1/(\chi_{f}\bardelta_t^2)),\quad |\tau_2^t| \geq \big( \{\bar{\barR}_t^2\vee \bar{\barR}_{s_t}^2\vee \chi_{f}\bardelta_t^2 \}\wedge 1 \big)\cdot |\xi_2^t|. 
\end{equation}
Combining \eqref{Step:2:a} and \eqref{Step:2:b} together, the conditions \eqref{cond:event:E2} and~\eqref{cond:merit:var}~are satisfied if
\begin{equation}\label{cond:xi2}
\begin{aligned}
|\xi_2^t| & \geq O\rbr{\frac{\log(d/p_f)}{\kappa_f^2\baralpha_t^4\{(\bnabla\mL_{\barepsilon_t, \barnu_t, \eta}^t)^T\cDelta_t\}^2\wedge\chi_{f}\bardelta_t^2}},\\
|\tau_2^t| & \geq \big(\{\bar{\barR}_t^2\vee \bar{\barR}_{s_t}^2\}\cdot|\xi_2^t| \vee \{\log(d/p_f)\cdot |\xi_2^t|\}^{1/2}\big) \wedge|\xi_2^t|.
\end{aligned}
\end{equation}
Similar to the complexity \eqref{cond:xi1}, \eqref{cond:xi2} suggests that the batch size $|\tau_2^t|$ for $\nabla f_t$, $\nabla f_{s_t}$ is significantly less than $|\xi_2^t|$ for $f_t$, $f_{s_t}$, with the ratio $|\tau_2^t|/|\xi_2^t|$~decaying to zero when $t$ increases. The denominator in \eqref{cond:xi2} is nonzero if $\barR_t\neq 0$ (which is always the case; otherwise, we should stop the iteration). In particular,~if~$\cDelta_t = \barDelta_t$, then
\begin{align*}
-\kappa_f\baralpha_t^2(\bnabla\mL_{\barepsilon_t, \barnu_t, \eta}^t)^T\cDelta_t
& \stackrel{\eqref{pequ:34}}{\geq} \frac{\kappa_{f}\baralpha_t^2(\gamma_{B}\wedge\eta)}{4}\nbr{\left(\begin{smallmatrix}
\barDelta\bx_t\\
J_t\bnabla_{\bx}\mL_t\\
G_t\bnabla_{\bx}\mL_t + \Pi_c(\diag^2(g_t)\blambda_t)
\end{smallmatrix} \right)}^2 \\
& \stackrel{\eqref{pequ:47}}{\geq }O(\kappa_{f}\baralpha_t^2\barR_t^2)>0;
\end{align*}
if $\cDelta_t=\hDelta_t$, then
\begin{equation*}
-\kappa_f\baralpha_t^2(\bnabla\mL_{\barepsilon_t, \barnu_t, \eta}^t)^T\cDelta_t \stackrel{\eqref{cond:alter:dir}}{\geq}\kappa_{f}\baralpha_t^2/\chi_{u}\nbr{\bnabla\mL_{\barepsilon_t, \barnu_t, \eta}^t}^2 \stackrel{\eqref{pequ:57}}{\geq} O(\kappa_{f}\baralpha_t^2\barR_t^2)>0.
\end{equation*}

\subsection{Discussion on computations and limitations}\label{sec:4.4}

We now briefly discuss the per-iteration computational cost of Algorithm~\ref{alg:ASto},~and present some limitations and extensions of the algorithm.

\vskip3pt
\noindent\textbf{Objective evaluations.} By Section \ref{sec:4.3} and the complexities in \eqref{cond:xi1} and \eqref{cond:xi2}, Algorithm~\ref{alg:ASto} generates $|\xi_1^t| + |\xi_2^t|$ samples in each iteration, and evaluates $2|\xi_2^t|$ function values, $|\xi_1^t| + 2|\tau_2^t|$ gradients, and $|\tau_1^t|$ Hessians for the objective.~To~see their orders from \eqref{cond:xi1} and \eqref{cond:xi2} clearly, let us suppose $\baralpha_t$ stabilizes at $\alpha_{max}$ (i.e., the steps are successful) and $\bardelta_t = O((\bnabla\mL_{\barepsilon_t, \barnu_t, \eta}^t)^T\cDelta_t)$ (see \eqref{cond:decrease:trust} for the reasonability). We also replace the stochastic quantities in \eqref{cond:xi1}, \eqref{cond:xi2} by deterministic counterparts and let $R_t\approx R_{s_t}$. Then, we can see that $|\xi_1^t| = |\tau_2^t| = O(1/R_t^2)$, $|\xi_2^t| = O(1/R_t^4)$, and $|\tau_1^t| = O(1)$. Thus, the objective evaluations are 
\begin{equation*}
\text{function values: } O(1/R_t^4), \quad \text{gradients: } O(1/R_t^2), \quad \text{Hessians: }  O(1).
\end{equation*}
We note that the evaluations for the function values and gradients are increasing as the iteration proceeds, and the function evaluations are square of the gradient evaluations. Under the same setup, our evaluation complexities for the functions and gradients are consistent with the unconstrained stochastic line search \cite[Section 2.3]{Paquette2020Stochastic} with $R_t$ replaced by $\|\nabla f_t\|$. Although the augmented Lagrangian merit function requires the Hessian evaluations,~the~\mbox{Hessian complexity}~is~significantly less than that of functions and gradients, and does not have to~increase during the iteration. Such an observation is missing in the prior work \citep{Na2022adaptive}.

\vskip3pt
\noindent\textbf{Constraint evaluations.} Since the constraints are deterministic, Algorithm~\ref{alg:ASto} has the same constraint evaluations as deterministic schemes. In particular,~the algorithm evaluates four function values (two for equalities and two for inequalities; and for each type of constraint, one for current point and one for trial point), four Jacobians, and two Hessians in each iteration.

\vskip3pt
\noindent\textbf{Computational cost.} Same as deterministic SQP schemes, solving Newton system dominates the computational cost. If we do not consider the potential sparse or block-diagonal structures that many problems have, solving the system \eqref{equ:SQP:direction} requires $O((d+m+|\text{active set}|)^3) + O((m+r)^3) = O(d^3+m^3+r^3)$ flops. Such computational cost is larger than solving a standard SQP system (see \cite[(8.9)]{Pillo2002Augmented}) by the extra term $O((m+r)^3)$. However, as explained in Remark~\ref{rem:1}, the analysis of standard SQP system relies on the exact Hessian, which is inaccessible in our stochastic setting. When the SQP direction is not employed, the backup direction can be obtained with $O(d+m+r)$ flops for the gradient step, $O((d+m+r)^3)$ flops for the regularized Newton step, and between for the truncated Newton step. Such computational cost is standard in the literature \cite{Pillo2008truncated, Pillo2011exact}, where a safeguarding direction satisfying \eqref{cond:alter:dir} is required to minimize the augmented Lagrangian. We should mention that, as the EQP scheme, the above computations are not very comparable with the IQP schemes. In that case, the SQP systems include inequality constraints and are more expensive~to solve, although less iterations may be performed. 

\vskip3pt
\noindent\textbf{Limitations of the design.} Algorithm \ref{alg:ASto} has few limitations. First, it solves the SQP systems exactly. In practice, one may apply conjugate gradient (CG) or minimum residual (MINRES) methods, or apply randomized iterative solvers to solve the systems inexactly. The inexact direction can reduce the computational cost significantly \citep{Curtis2021Inexact}. Second, our backup direction does not fully utilize the computations of the SQP direction. Although our analysis allows any backup direction satisfying \eqref{cond:alter:dir}, and utilizing Newton direction as a backup is standard in the literature \cite{Pillo2008truncated, Pillo2011exact}, a better choice is to directly modify the SQP direction. Then, we may derive a direction that has a faster convergence than the gradient direction, and less computations than the (regularized) Newton direction. We leave the refinements of these two limitations to the future.

\section{Numerical Experiments}\label{sec:5}

We implement the following two algorithms on 39 nonlinear problems collected in CUTEst test set \citep{Gould2014CUTEst}. We select the problems that have a non-constant~objective with less than 1000 free variables. We also require the problems to~have at least one inequality constraint, no infeasible constraints, no network constraints; and require the number of constraints to be less than the number of variables. The setup of each algorithm is as follows.
\begin{enumerate}[label=(\alph*),topsep=0pt]	\setlength\itemsep{0.0em}
\item \textbf{AdapNewton}: the adaptive scheme in Algorithm \ref{alg:ASto} with the safeguarding direction given by the regularized Newton step. We set the inputs as~$\baralpha_0 = \alpha_{max} = 1.5$, $\beta=0.3$, $\kappa_{f} = \beta/(4\alpha_{max}) = 0.05$, $ \kappa_{grad} = \chi_{grad} = \chi_{f} = \bardelta_0=1$, $\barepsilon_0=10^{-2}$, $\eta=10^{-4}$, $p_{grad}=p_f=0.1$, $\rho=2$. Here, we set $\alpha_{max}>1$ since a stochastic scheme can select a stepsize that is greater than one (cf. Figure \ref{fig:4}). $\beta$ is close to the middle of the interval $(0,0.5)$, which is a common~range for deterministic schemes. $(\barepsilon_0, \bardelta_0)$ are adaptively selected during the iteration, while we prefer a small initial $\barepsilon_0$ to run less adjustments on it. $\kappa_{f}$ is set as the allowed largest value $\beta/(4\alpha_{max})$ (cf. Algorithm \ref{alg:ASto}); however, the parameters $(\kappa_{grad},\kappa_{f},\chi_{grad},\chi_{f},p_{grad},p_f)$ all affect the batch sizes and play the same role as the constant $C$ that we study later. We let $\eta$ be small so that the~last penalty term of \eqref{equ:aug:Lagrange} is almost negligible, and the merit function \eqref{equ:aug:Lagrange} is~close~to a standard augmented Lagrangian function. 

\hskip0.5cm We also test the robustness of the algorithm to three parameters $C$, $\kappa$, $\chi_{err}$. Here, $C$ is the constant multiplier of the big ``$O$" notation in \eqref{cond:xi1} and \eqref{cond:xi2} (the variance $\sigma^2$ of a single sample is also absorbed in ``$O$", which we introduce later). $\kappa$ is a parameter of the set $\mT_\nu$ ($\kappa=2$ in \eqref{equ:domain}), and $\chi_{err}$ is a parameter of the feasibility error condition \eqref{cond:bound:fes:error}. Their default values are $C = \kappa = 2$ and $\chi_{err} = 1$, while we allow to vary them in wide ranges: $C,\kappa \\ \in \{2,2^3,2^6\}$ and $\chi_{err}\in\{1,10,10^2\}$. When we vary one parameter, the other two are set as default.

\item \textbf{AdapGD}: the adaptive scheme in Algorithm \ref{alg:ASto} with the safeguarding~direction given by the steepest descent step. The setup is the same as (b).
\end{enumerate}
For both algorithms, the initial iterate $(\bx_0, \bmu_0, \blambda_0)$ is specified by the CUTEst package. The package also provides the deterministic function, gradient and Hessian evaluation, $f_t, \nabla f_t, \nabla^2f_t$, in each iteration. We generate their stochastic counterparts by adding a Gaussian noise with variance $\sigma^2$. In particular, we let $\barf_t \sim \N(f_t, \sigma^2)$, $\bnabla f_t\sim \N(\nabla f_t, \sigma^2(I + \b1\b1^T))$, and $(\bnabla^2f_t)_{ij}\sim \N((\nabla f_t)_{ij}, \sigma^2)$. We try four levels of variance: $\sigma^2 \in\{ 10^{-8}, 10^{-4}, 10^{-2}, 10^{-1}\}$. Throughout the implementation, we let $B_t = I$ (cf. \eqref{equ:SQP:direction}, \eqref{equ:def:H}) and set the iteration budget to be $10^4$. The stopping criterion is
\begin{equation}\label{stop}
\baralpha_t\|\cDelta_t\| \leq 10^{-7}\quad \text{OR} \quad R_t \leq 10^{-5}\quad \text{OR}\quad t\geq 10^4.
\end{equation}
The former two cases suggest that the iteration converges within the budget. For each algorithm, each problem, and each setup, we average the results~of~all convergent runs among $5$ runs. Our code is available at \url{https://github.com/senna1128/Constrained-Stochastic-Optimization-Inequality}.

\begin{figure}[!htp]
\centering     
\subfigure[$C =  2$]{\label{KKTC1}\includegraphics[width=37mm]{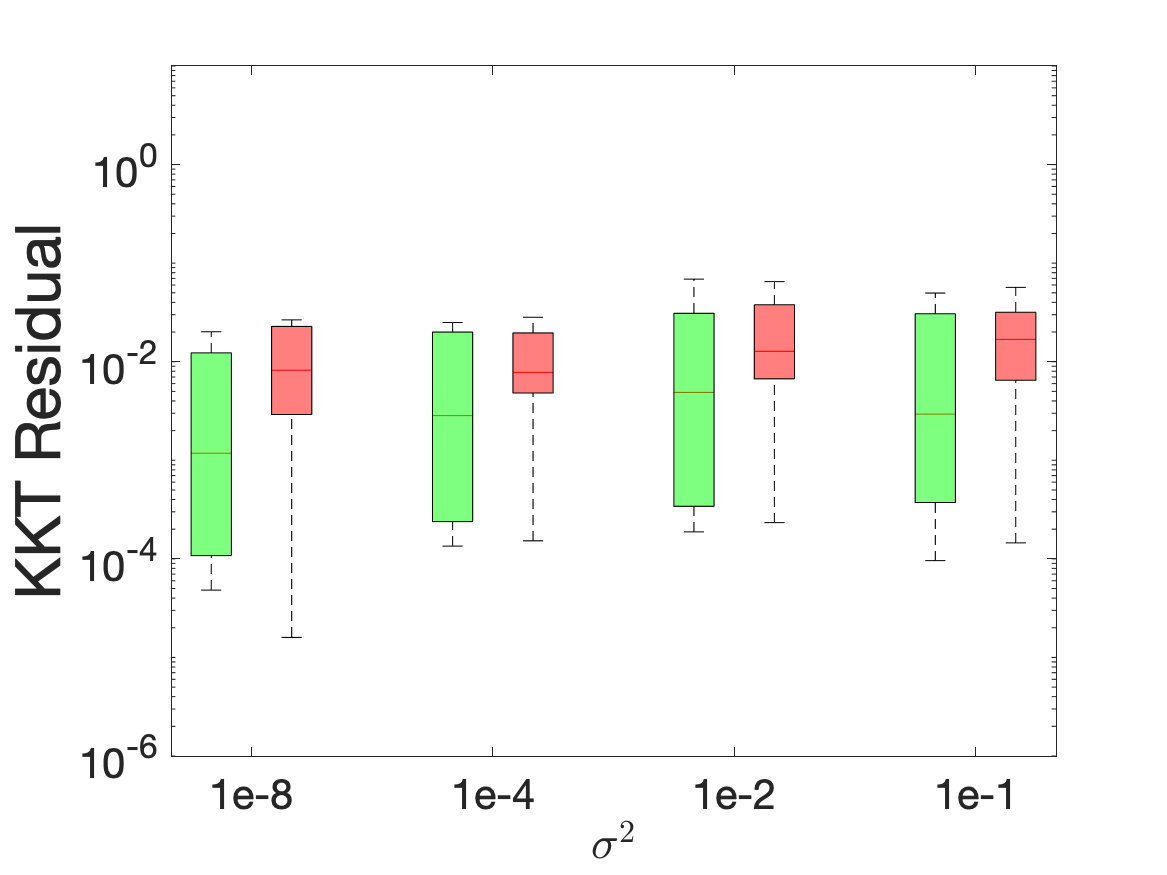}}
\subfigure[$C =  2^3$]{\label{KKTC2}\includegraphics[width=37mm]{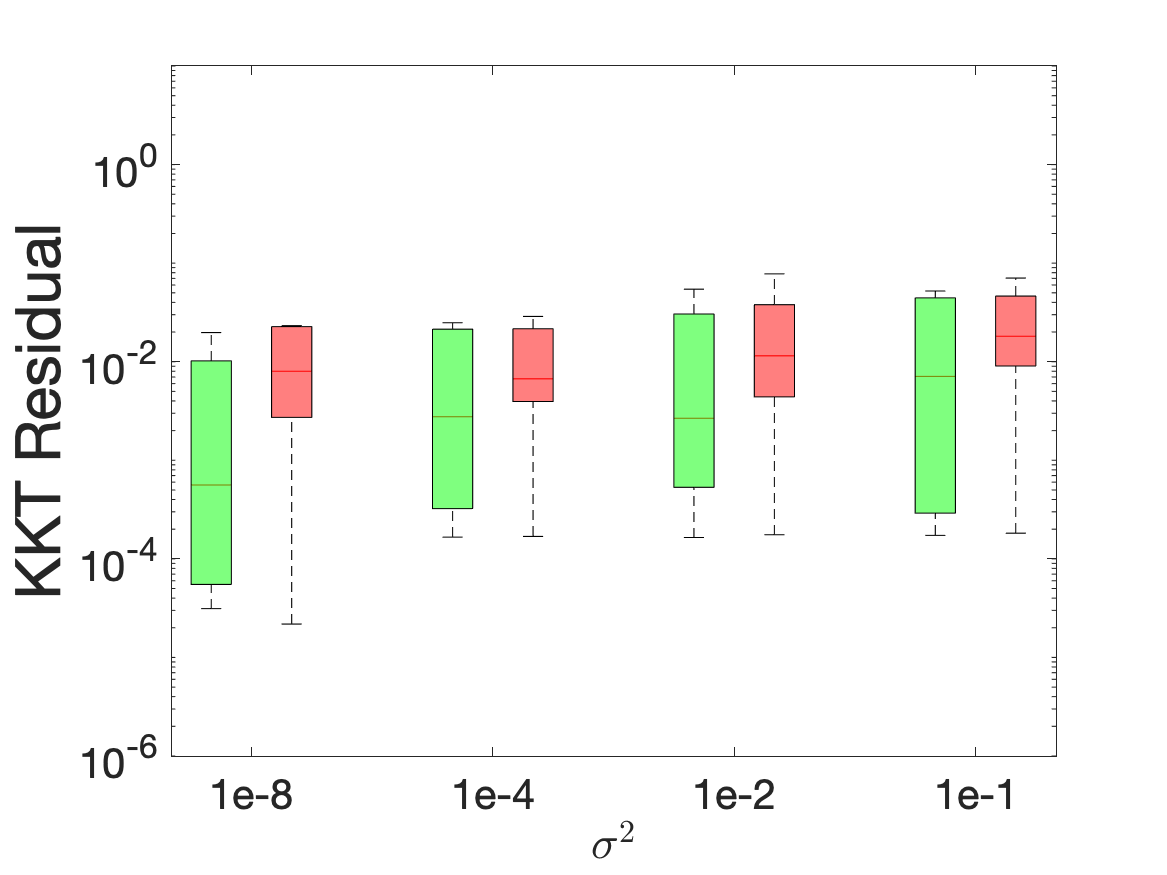}}
\subfigure[$C =  2^6$]{\label{KKTC3}\includegraphics[width=37mm]{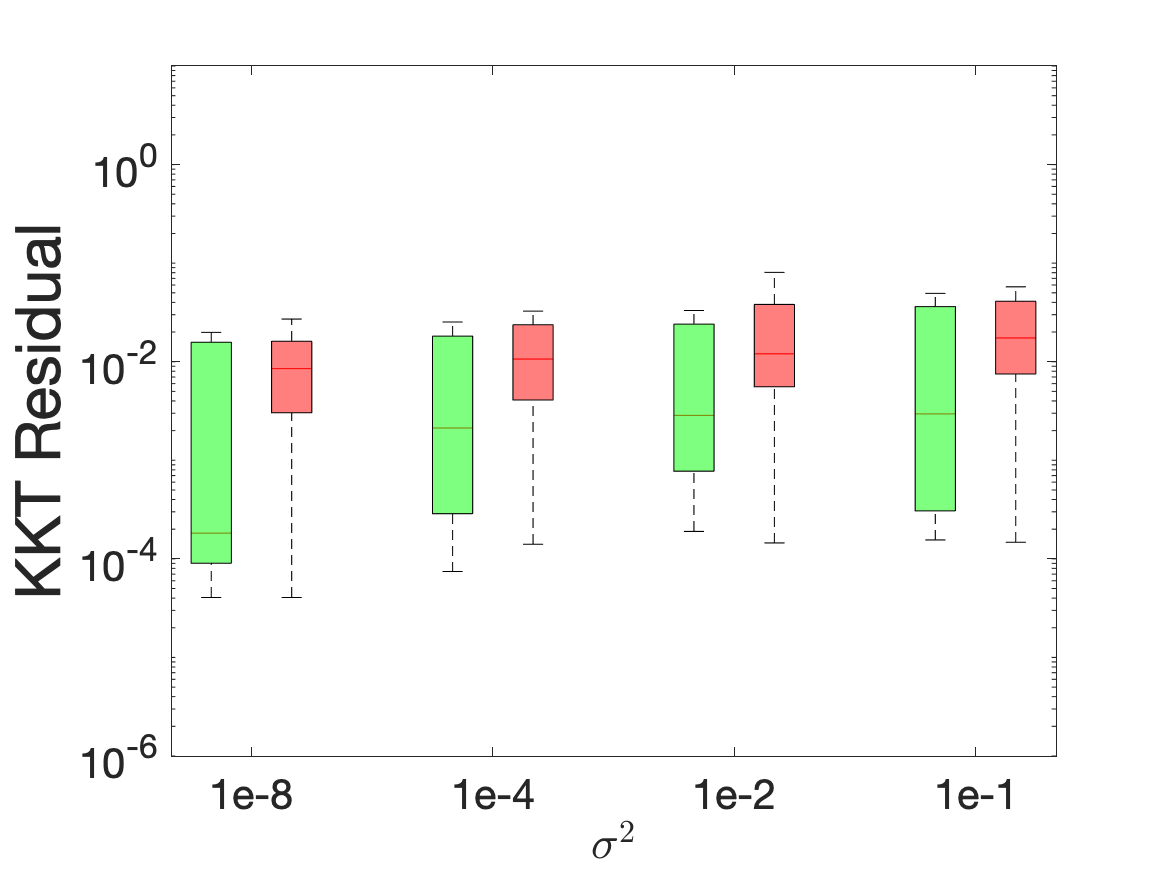}}

\subfigure[$\kappa =  2$]{\label{KKTK1}\includegraphics[width=37mm]{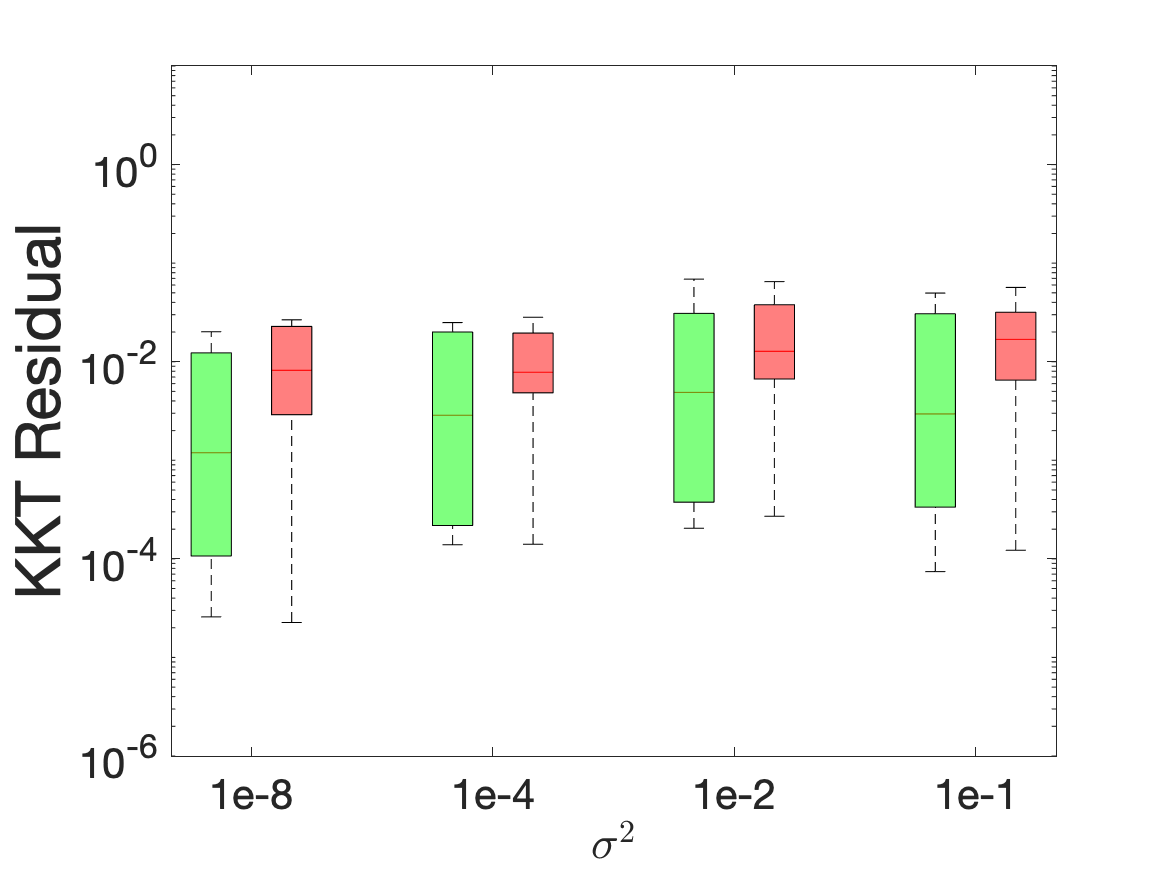}}
\subfigure[$\kappa =  2^3$]{\label{KKTK2}\includegraphics[width=37mm]{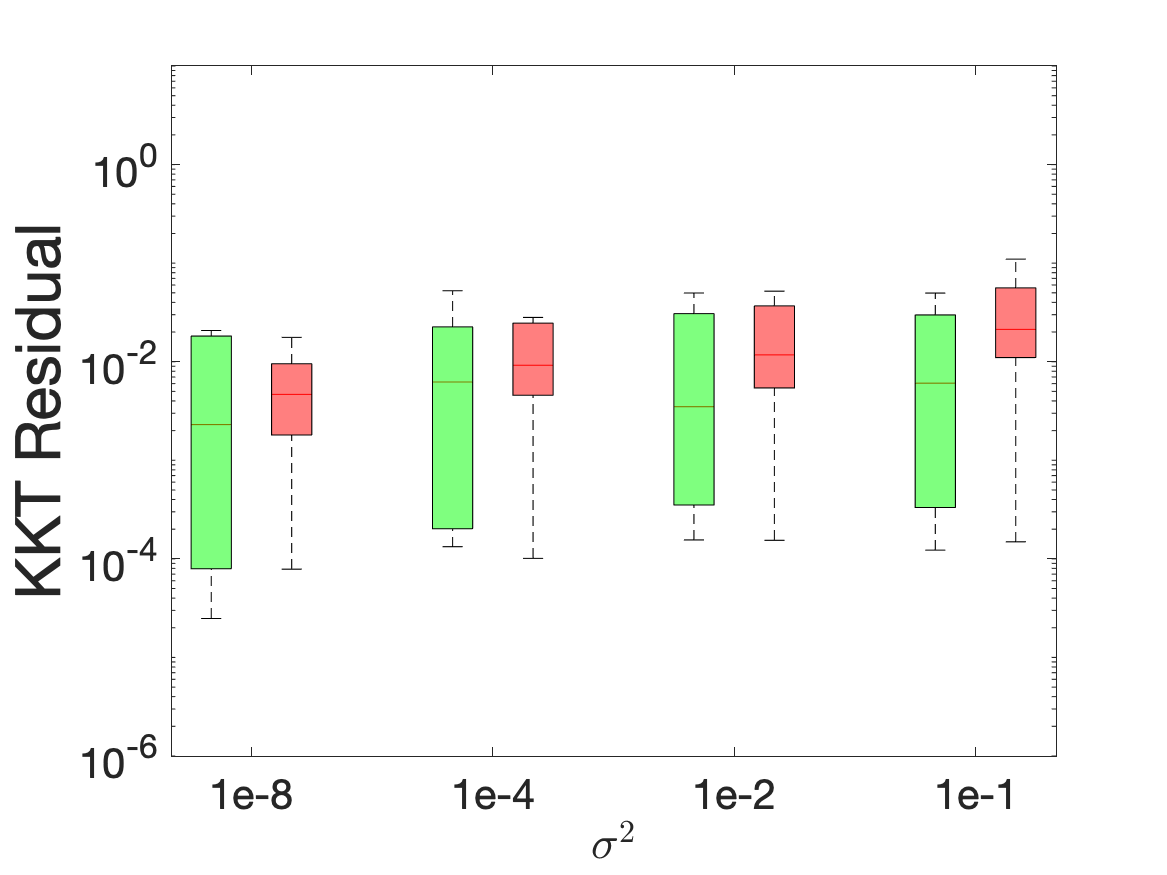}}
\subfigure[$\kappa =  2^6$]{\label{KKTK3}\includegraphics[width=37mm]{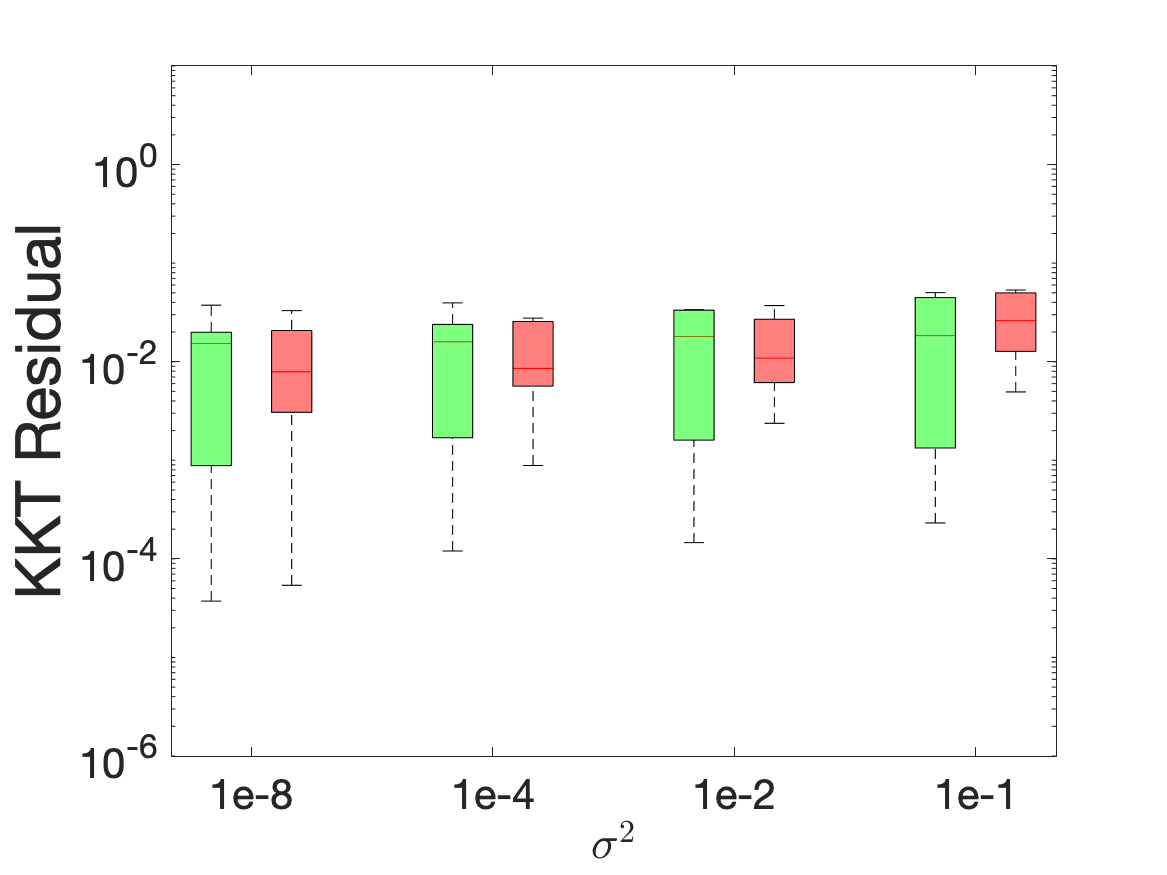}}

\subfigure[$\chi_{err} =  1$]{\label{KKTChi1}\includegraphics[width=37mm]{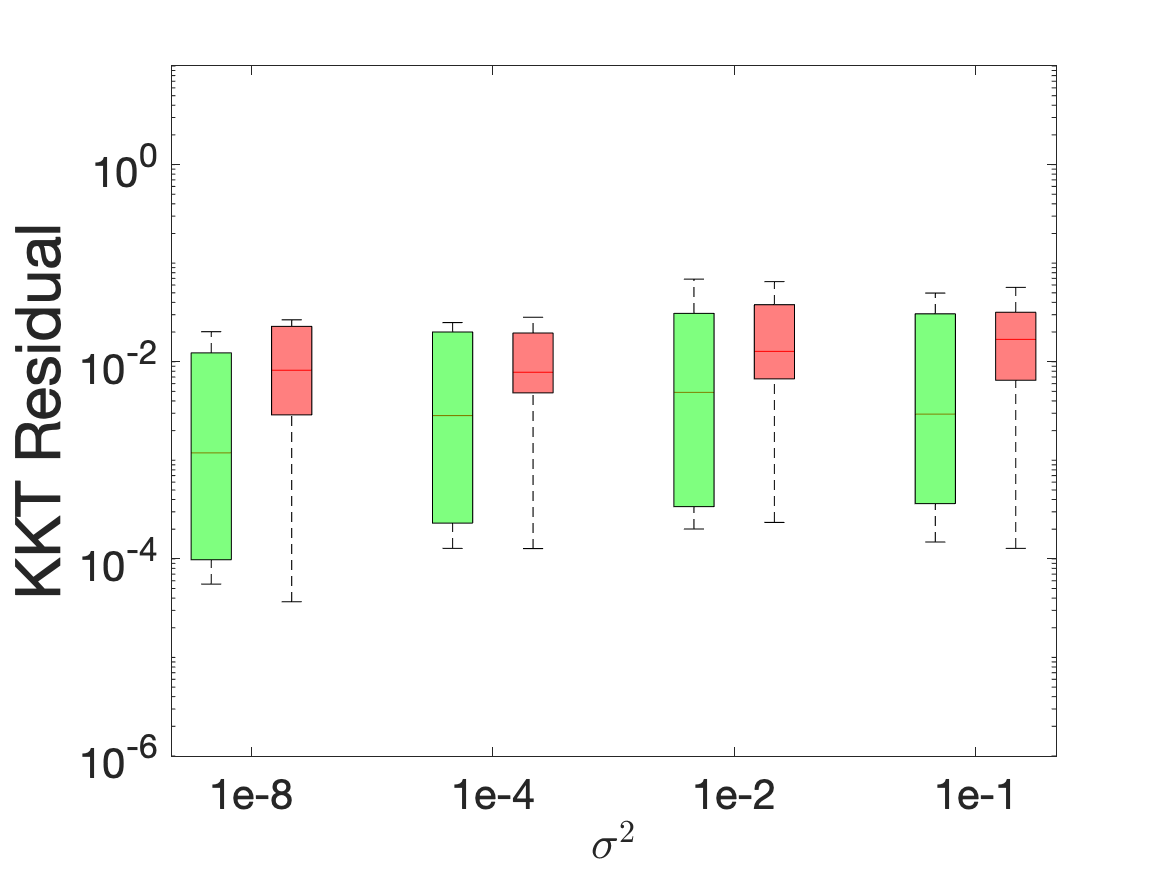}}
\subfigure[$\chi_{err} =  10$]{\label{KKTChi2}\includegraphics[width=37mm]{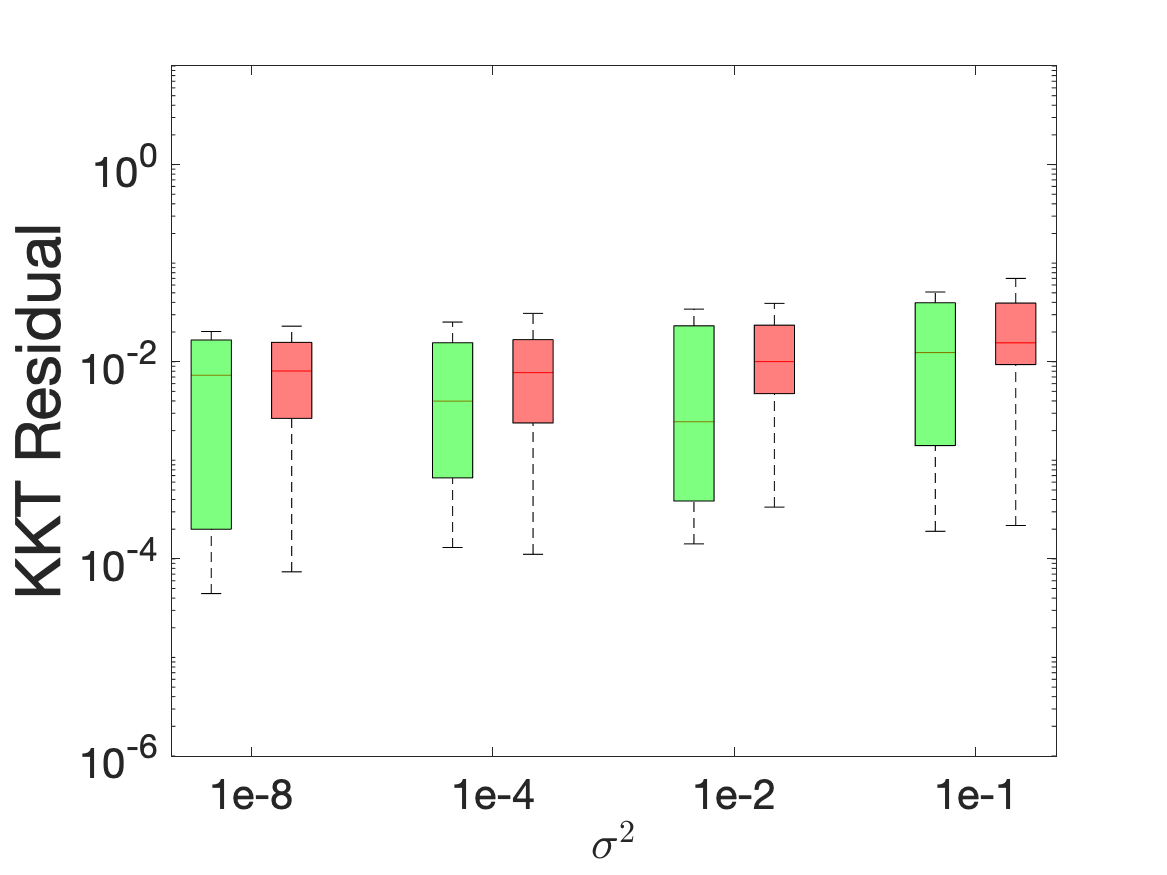}}
\subfigure[$\chi_{err} =  10^2$]{\label{KKTChi3}\includegraphics[width=37mm]{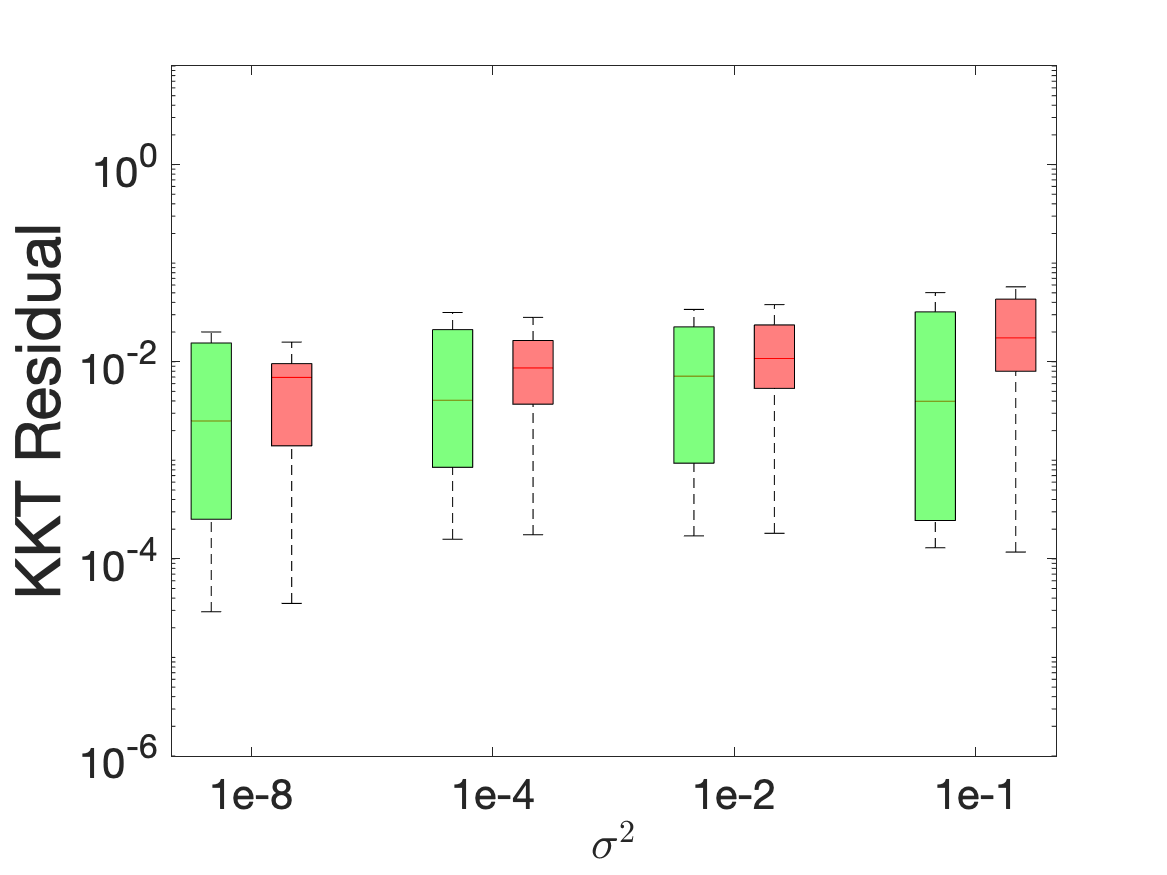}}

\includegraphics[width=0.3\textwidth]{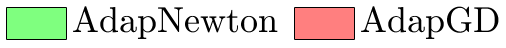}
\caption{KKT residual boxplots. Each panel corresponds to a setup of $(C, \kappa, \chi_{err})$. The default values are $C =\kappa= 2$ and $\chi_{err}=1$. When we vary one parameter, the other two are set as default. Thus, the three figures on the left column are the same.
}\label{fig:1}
\end{figure}

\vskip 4pt
\noindent\textbf{KKT residuals.} We draw the KKT residual boxplots for AdapNewton and AdapGD in Figure \ref{fig:1}. From the figure, we see that both algorithms are robust to tuning parameters $(C,\kappa,\chi_{err})$. For both algorithms, the median of the KKT residuals gradually increases as $\sigma^2$ increases, which is reasonable since~the~model estimation of each sample is more noisy when $\sigma^2$ is larger. However, the increase of the KKT residuals is mild since, regardless of $\sigma^2$, both methods generate enough samples in each iteration to enforce the model accuracy conditions (i.e., \eqref{cond:event:E1}, \eqref{cond:grad:var}, \eqref{cond:event:E2}, \eqref{cond:merit:var}). Figure \ref{fig:1} also suggests that AdapNewton outperforms AdapGD although the improvement is limited. In fact, the convergence on a few problems may be improved by utilizing the regularized Newton step as the backup of the SQP step; however, the SQP step will be employed eventually.

\begin{figure}[!htp]
\centering     
\subfigure[$C =  2$]{\label{SamC1}\includegraphics[width=37mm]{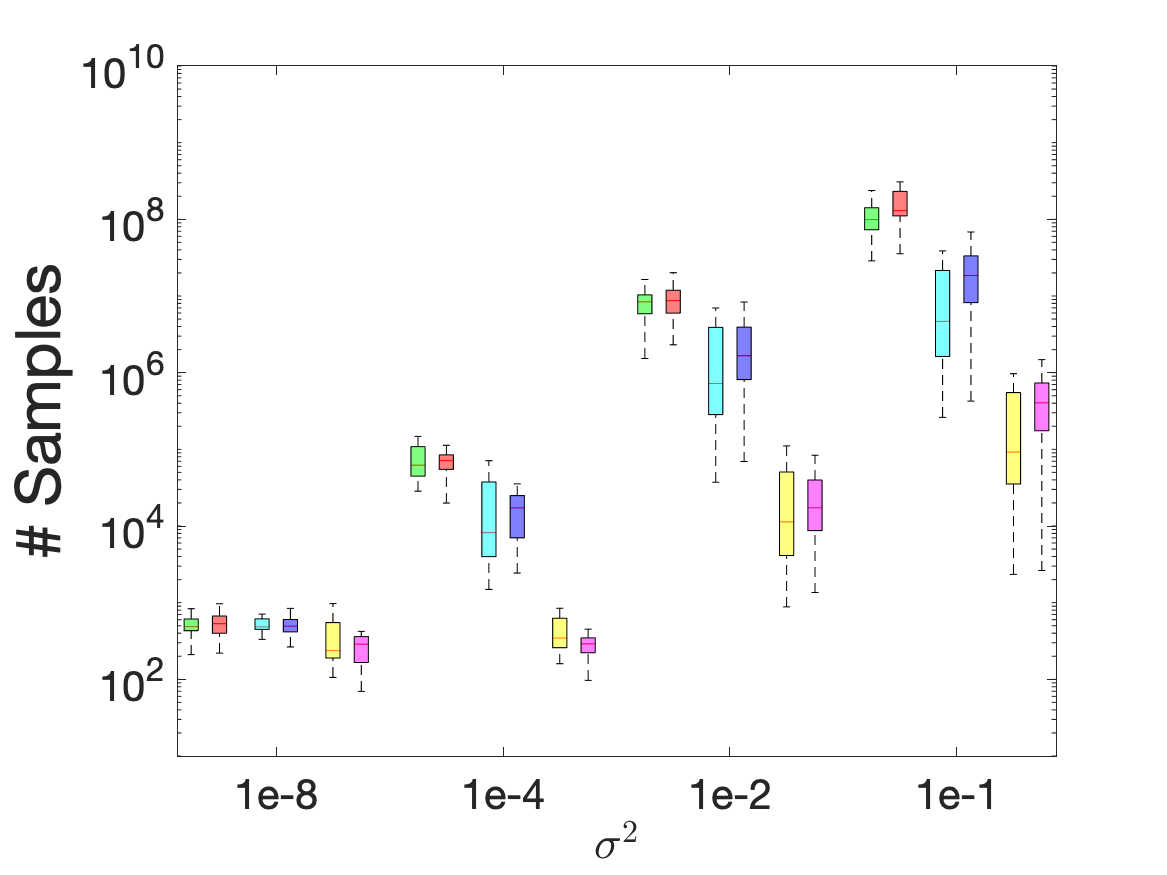}}
\subfigure[$C =  2^3$]{\label{SamC2}\includegraphics[width=37mm]{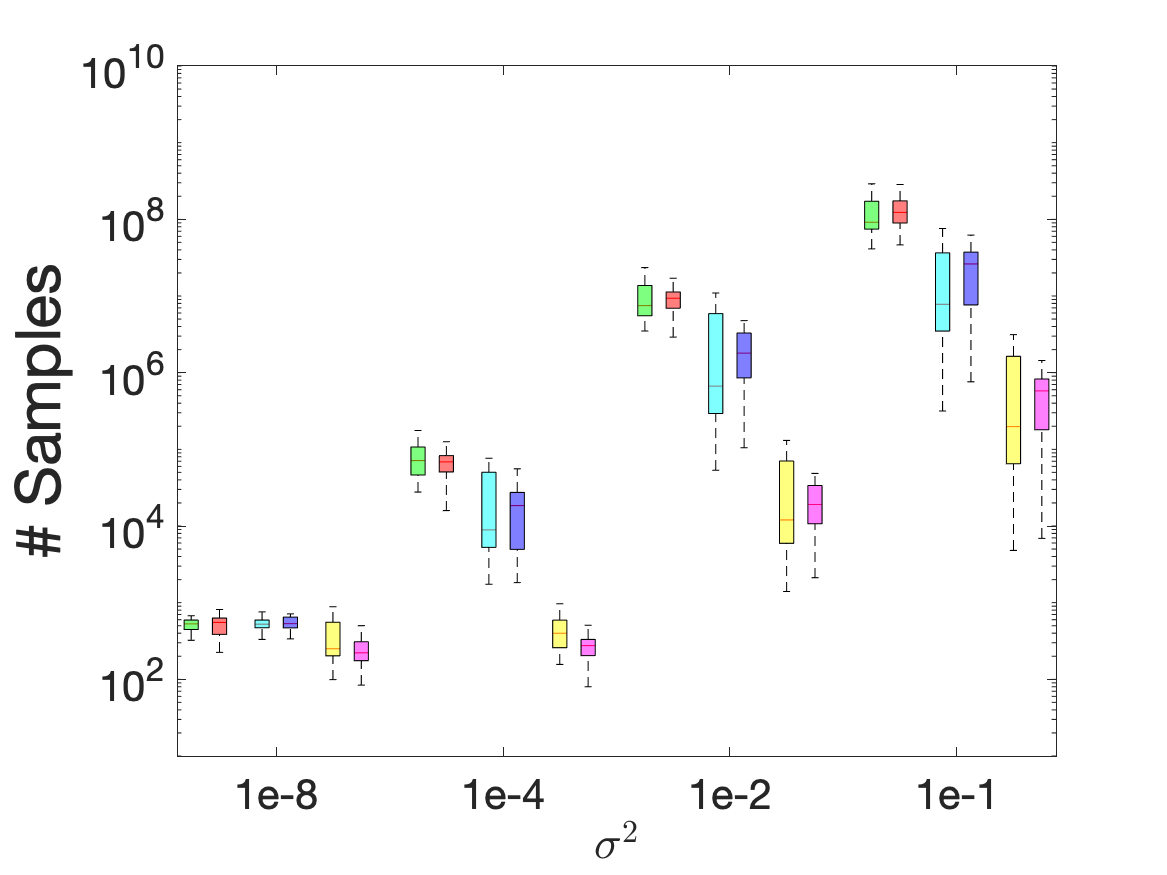}}
\subfigure[$C =  2^6$]{\label{SamC3}\includegraphics[width=37mm]{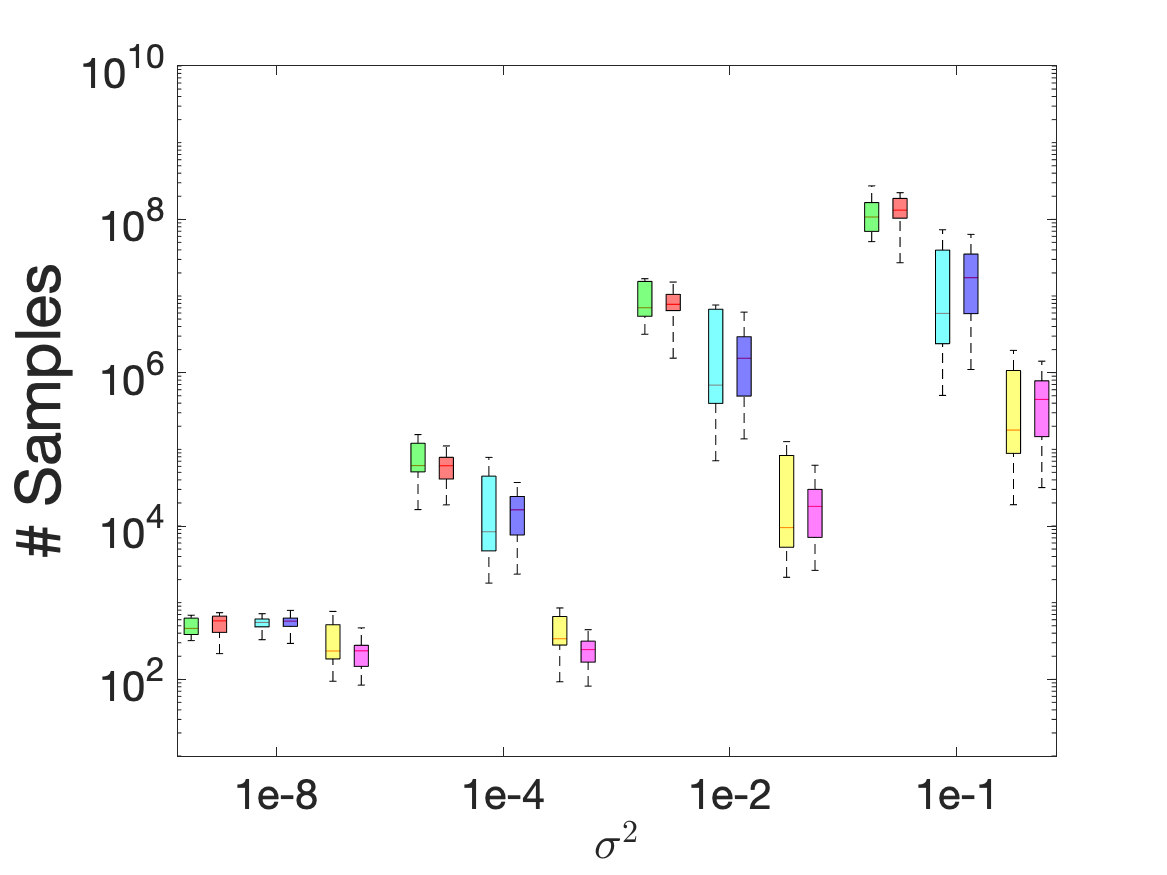}}
	
\subfigure[$\kappa =  2$]{\label{SamK1}\includegraphics[width=37mm]{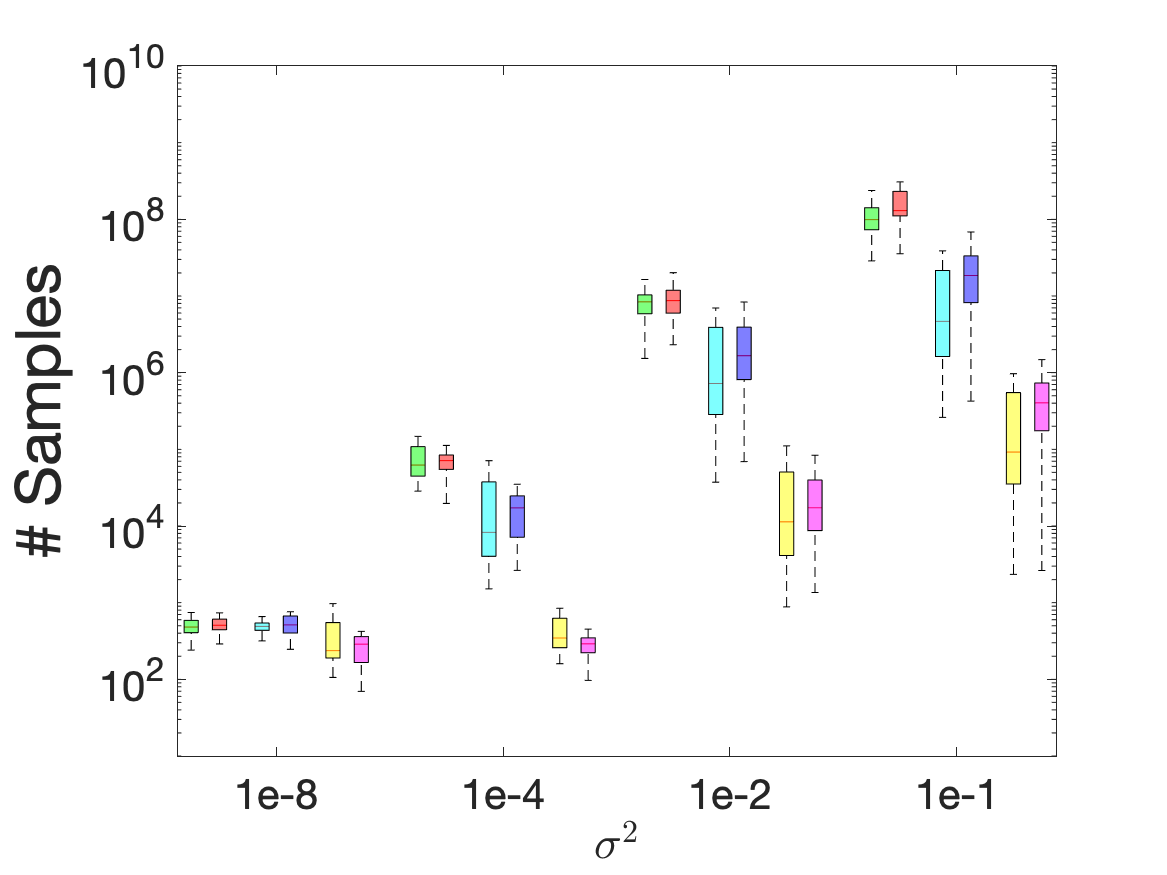}}
\subfigure[$\kappa =  2^3$]{\label{SamK2}\includegraphics[width=37mm]{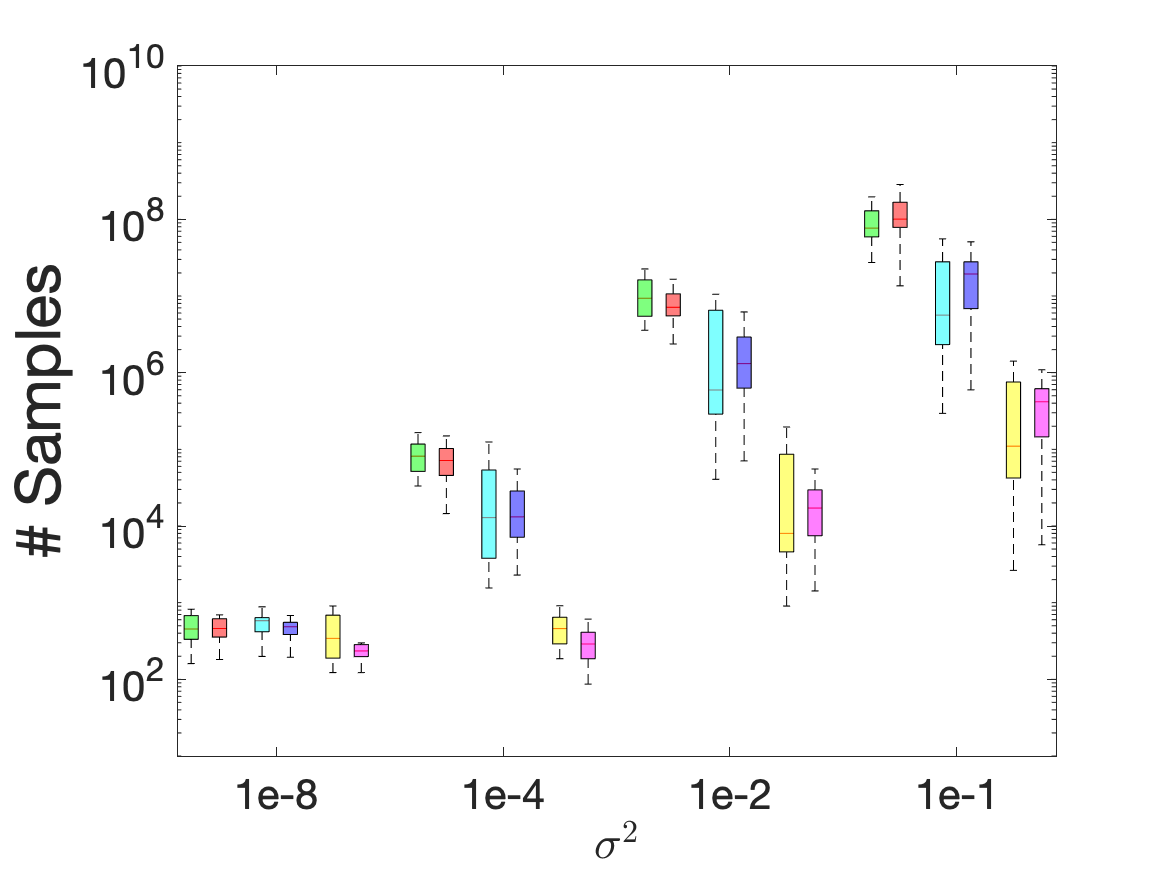}}
\subfigure[$\kappa =  2^6$]{\label{SamK3}\includegraphics[width=37mm]{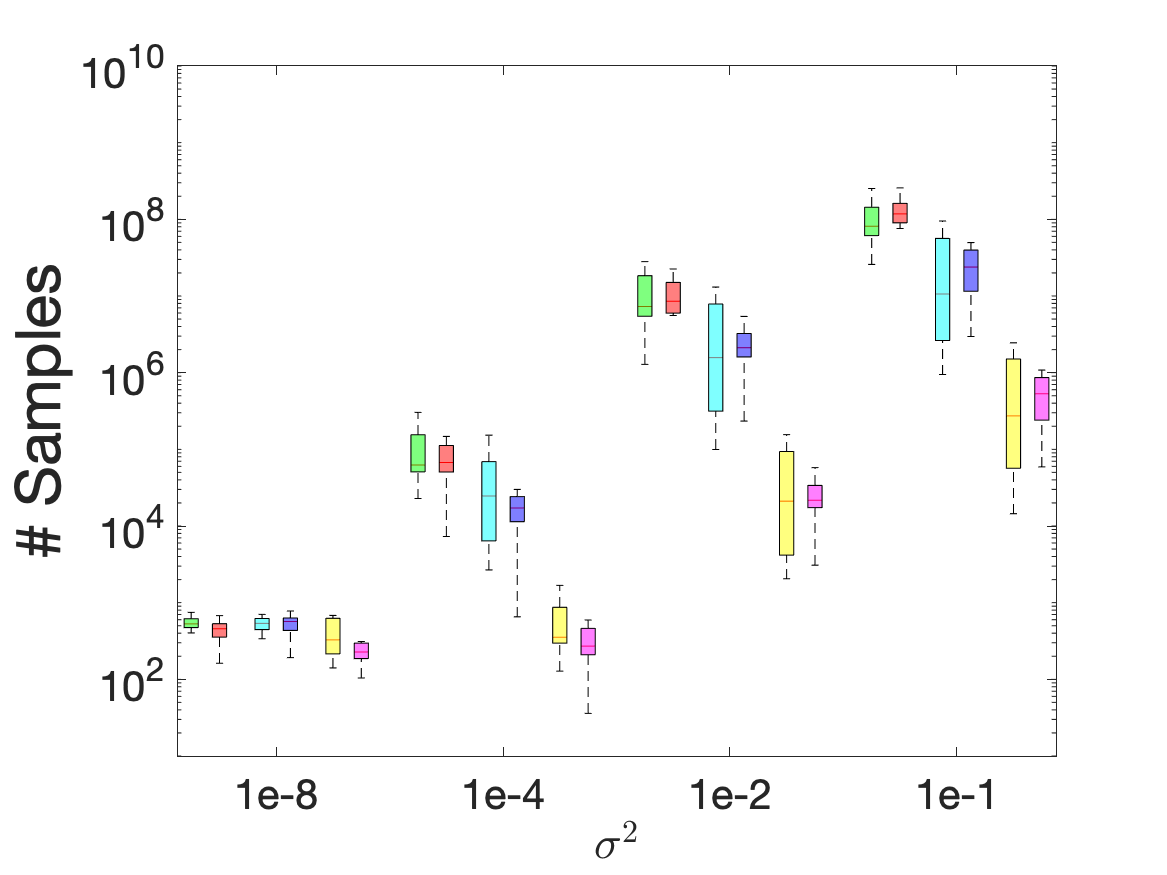}}
	
\subfigure[$\chi_{err} =  1$]{\label{SamChi1}\includegraphics[width=37mm]{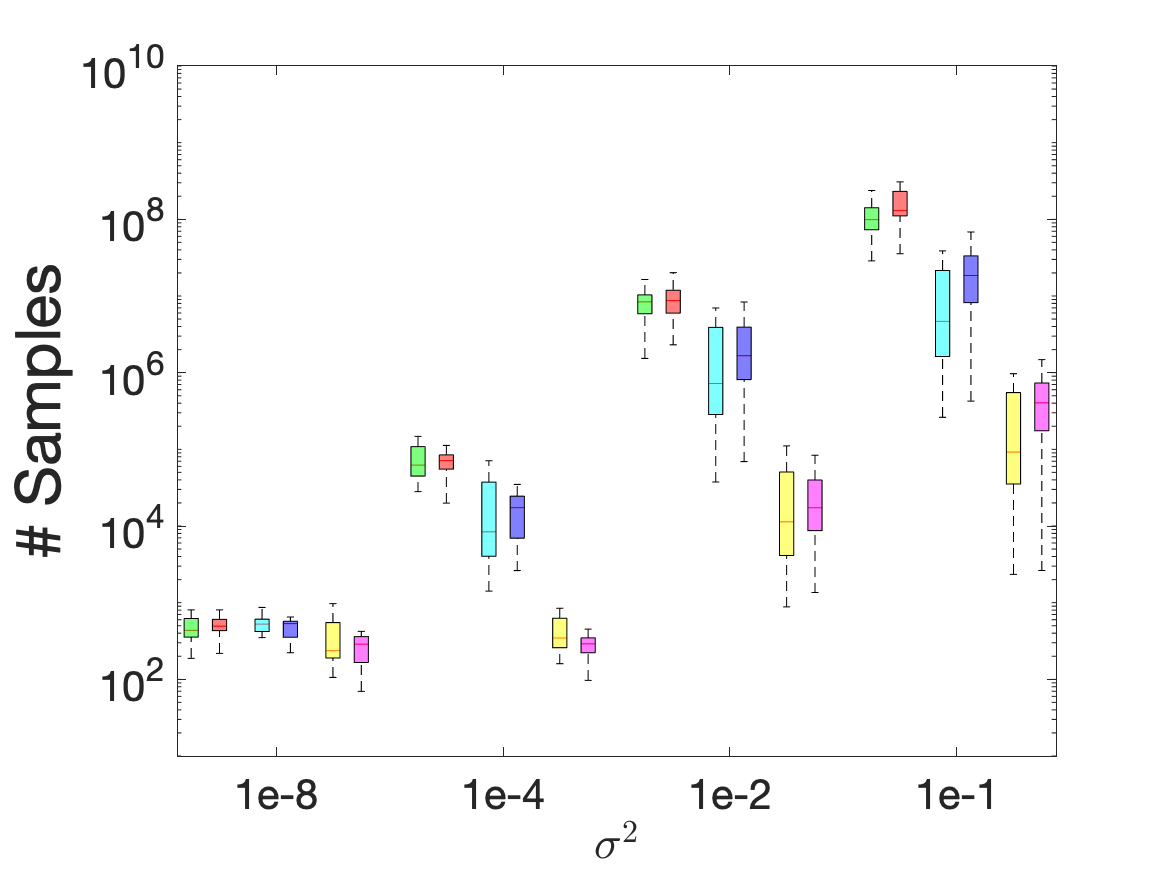}}
\subfigure[$\chi_{err} =  10$]{\label{SamChi2}\includegraphics[width=37mm]{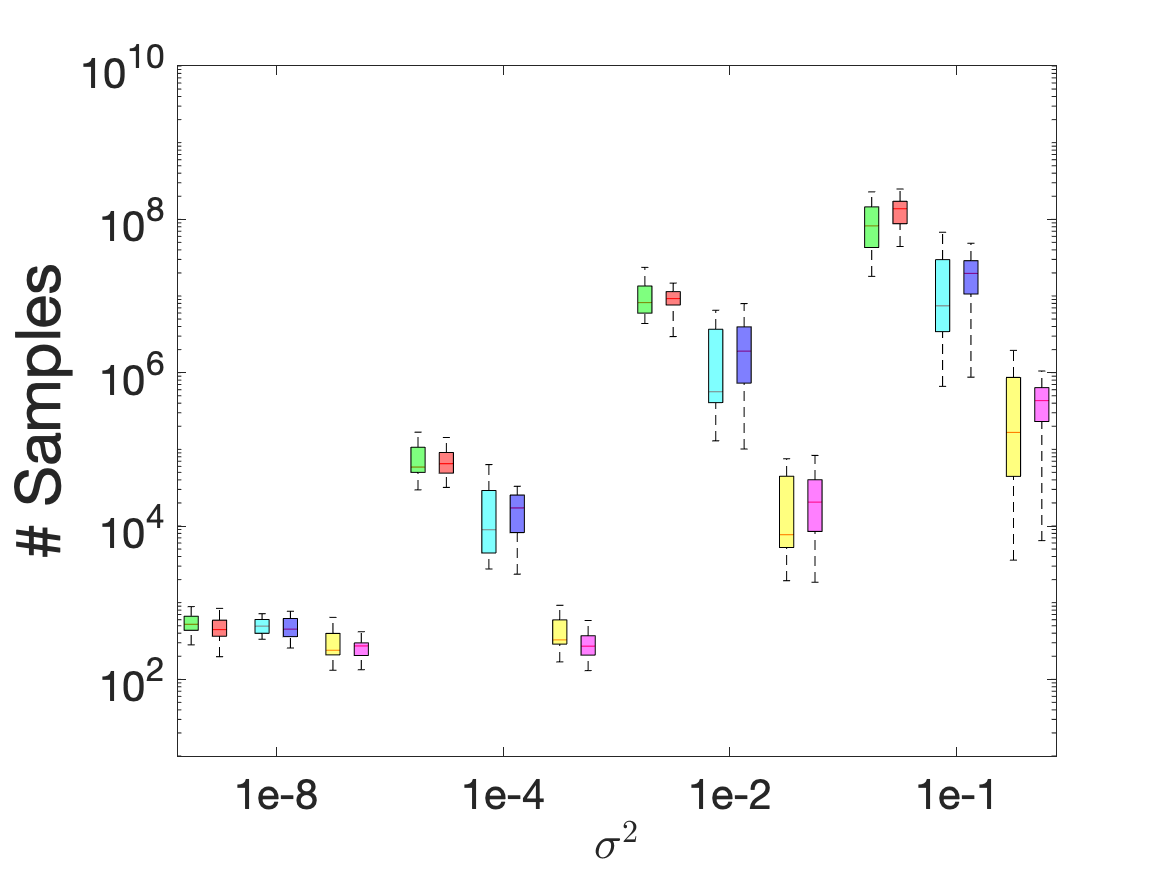}}
\subfigure[$\chi_{err} =  10^2$]{\label{SamChi3}\includegraphics[width=37mm]{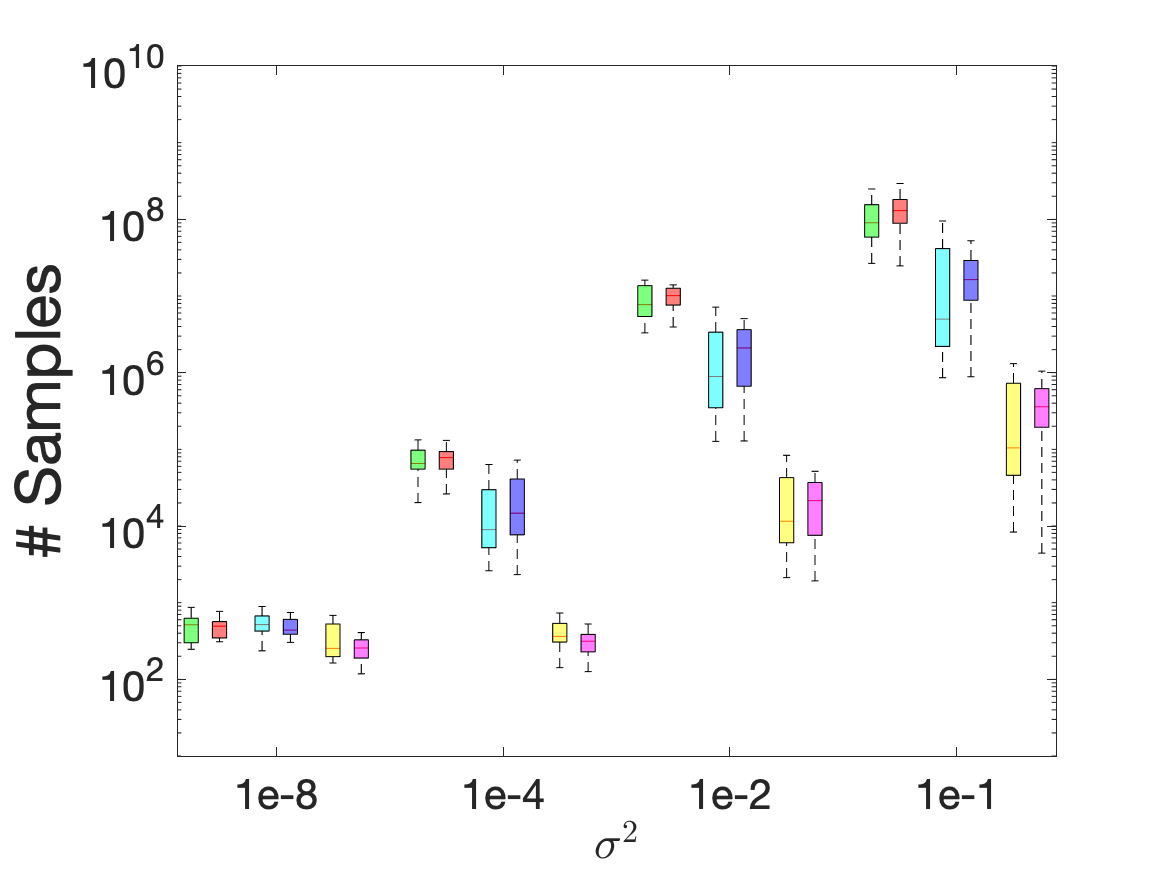}}
	
\includegraphics[width=\textwidth]{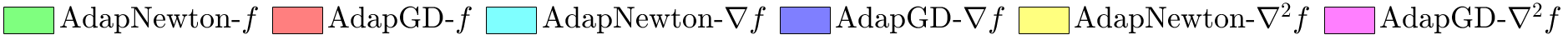}
\caption{Sample size boxplots. Each panel corresponds to a setup of $(C, \kappa, \chi_{err})$. The default values are $C =\kappa= 2$ and $\chi_{err}=1$. When we vary one parameter, the other two are set as default. Thus, the three figures on the left column are the same.}\label{fig:2}
\end{figure}

\vskip4pt
\noindent\textbf{Sample sizes.} We draw the sample size boxplots for AdapNewton and AdapGD in Figure \ref{fig:2}. From the figure, we see that both methods generate much less samples for estimating the objective Hessian compared to estimating the~objective value and gradient, between which the the objective gradient is estimated with less samples than the objective value. The sample size differences of the three quantities---objective value, gradient, Hessian---are clearer as $\sigma^2$ increases. For a fixed $\sigma^2$, the sample sizes of different setups of $(C, \kappa, \chi_{err})$ do not vary much. In fact, the parameters $\kappa$, $\chi_{err}$ do not directly affect the sample complexities. The parameter $C$ plays a similar role to $\sigma^2$ and affects the sample complexities via changing the multipliers in \eqref{cond:xi1} and \eqref{cond:xi2}. However, varying $C$ from 2 to 64 is marginal compared to varying $\sigma^2$ from $10^{-8}$ to $10^{-1}$. Thus, Figure \ref{fig:2} again illustrates the robustness of the designed adaptive algorithm. 

Moreover, as discussed in Sections \ref{sec:4.3} and \ref{sec:4.4}, the objective value, gradient, and Hessian have different sample complexities in each iteration, which depend on different powers of the reciprocal of the KKT residual $1/R_t$. When $\sigma^2=10^{-8}$, the small variance dominates the effect of $1/R_t$ so that all three quantities can be estimated with very few samples. When $\sigma^2=0.1$, the different dependencies of the sample sizes on $1/R_t$ are more evident. Overall, Figure \ref{fig:2} reveals the fact that different objective quantities can be estimated with different amount of samples. Such an aspect improves the prior work \cite{Na2022adaptive}, where the quantities with different sample complexities are estimated based on the same set of samples, and the effect of the variance $\sigma^2$ on the sample complexities is neglected.

\begin{figure}[!htp]
\centering     
\subfigure[AdapNewton]{\label{SamTraN}\includegraphics[width=54mm]{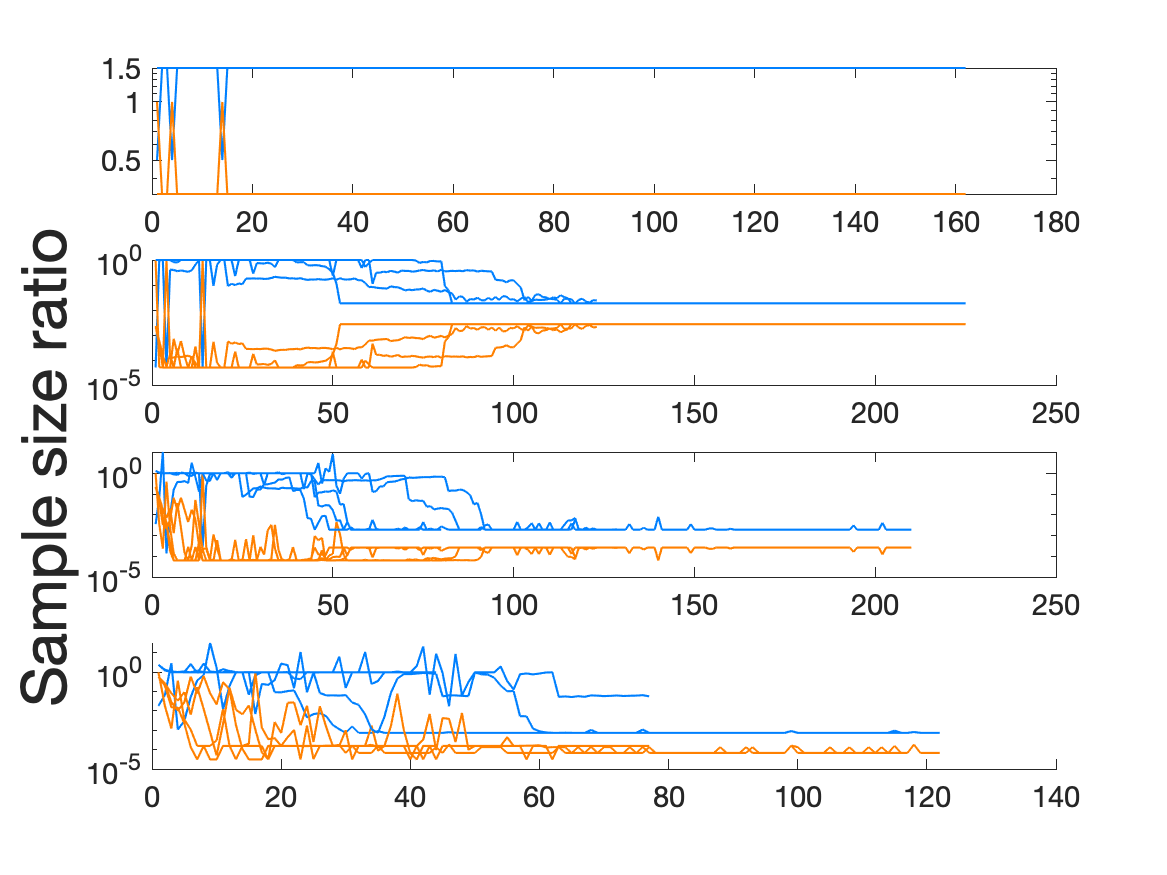}}
\subfigure[AdapGD]{\label{SamTraG}\includegraphics[width=54mm]{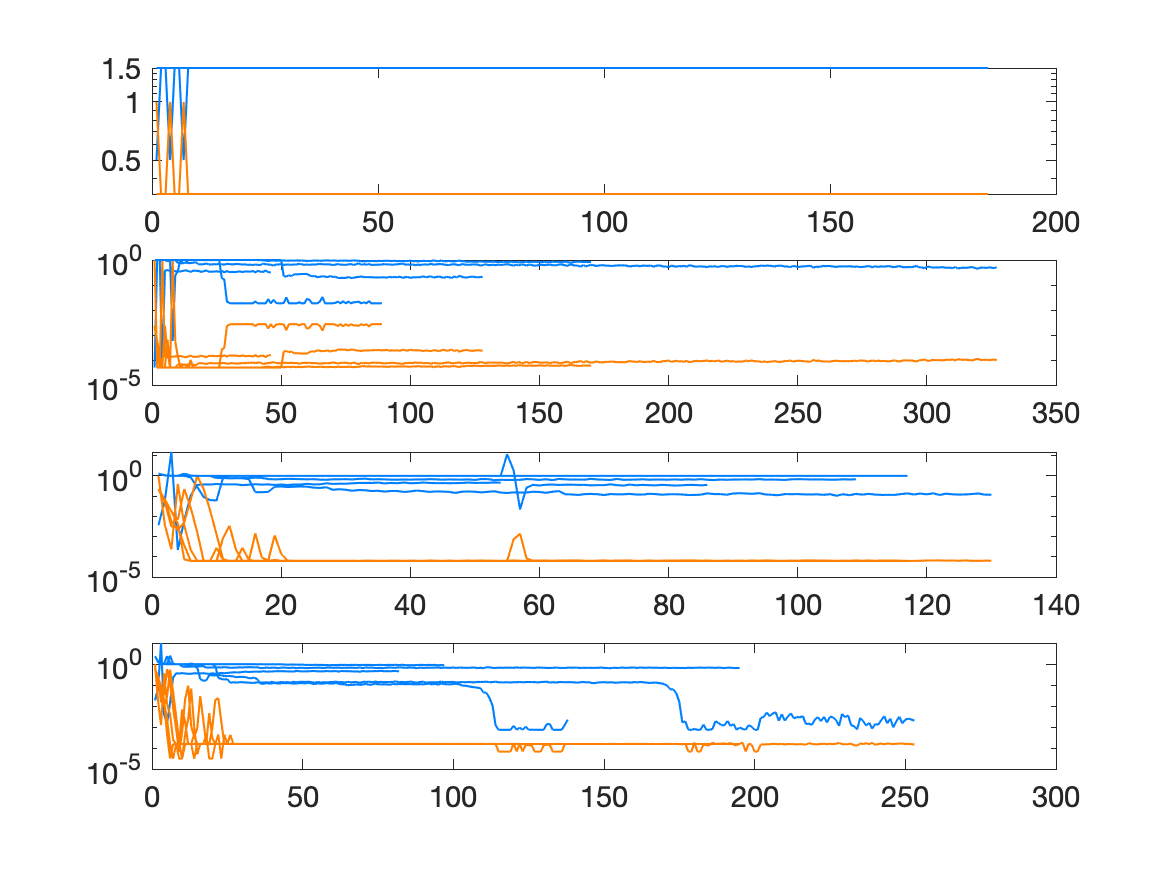}}
	
\includegraphics[width=0.4\textwidth]{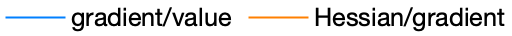}
\caption{Sample size ratio trajectories ($C=64$). Each plot has four rows, from top to bottom, corresponding to $\sigma^2 = 10^{-8}, 10^{-4}, 10^{-2}, 10^{-1}$. Each plot has ten lines with two colors. The five lines with the same color correspond to the five convergent problems.
}\label{fig:3}
\end{figure}

\begin{figure}[!tp]
\centering     
\subfigure[AdapNewton]{\label{StepTraN}\includegraphics[width=55mm]{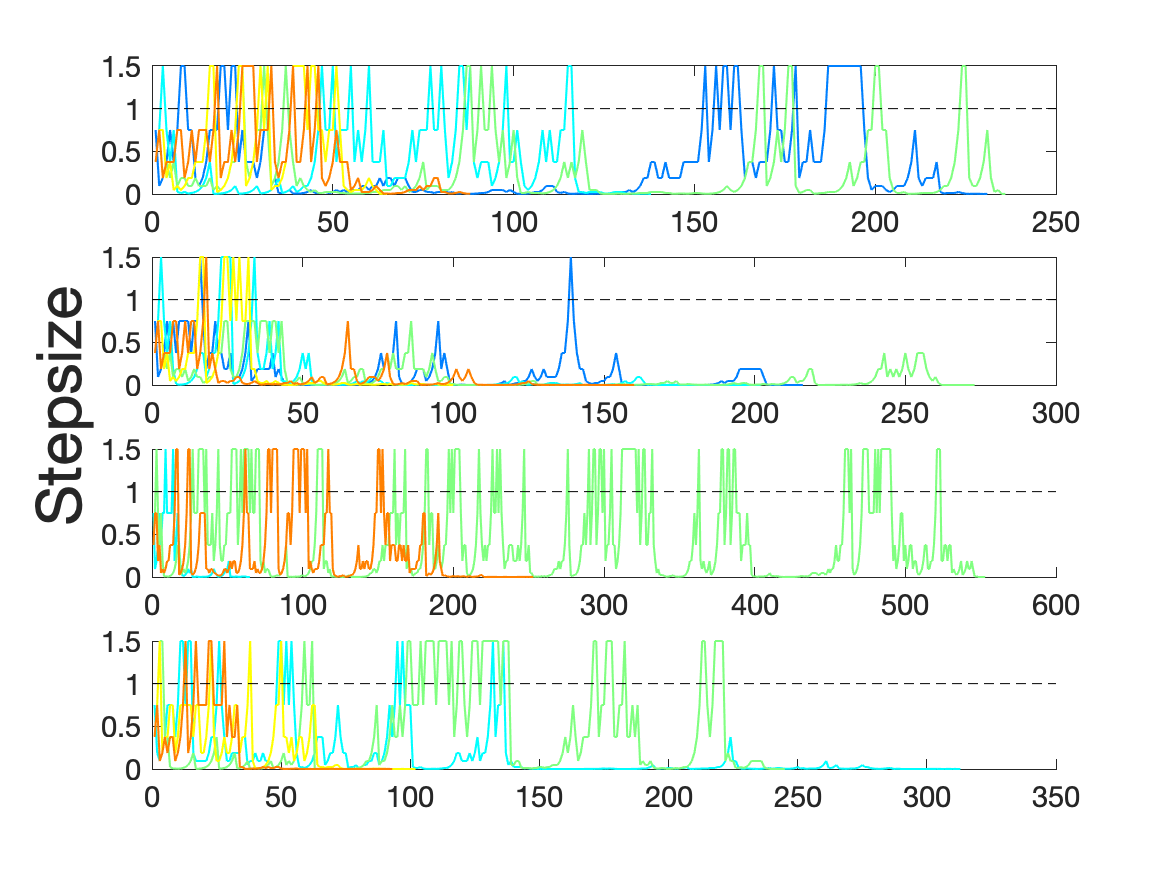}}
\subfigure[AdapGD]{\label{StepTraG}\includegraphics[width=55mm]{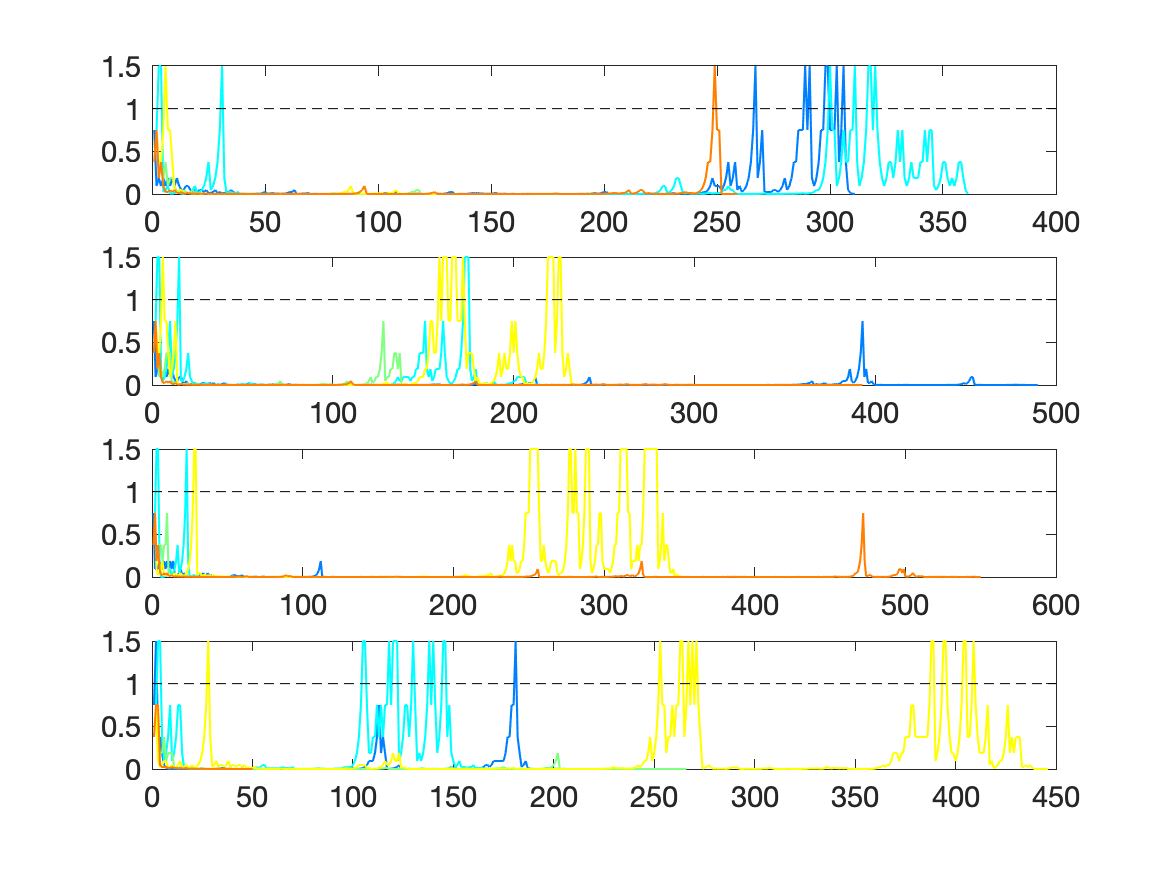}}
	
\caption{Stepsize trajectories. Each plot has four rows, from top to bottom, corresponding to $\sigma^2 = 10^{-8}, 10^{-4}, 10^{-2}, 10^{-1}$. Each plot has five lines, corresponding to the five problems.}\label{fig:4}
\end{figure}

In addition, we draw the trajectories of the sample size ratios. In particular, for both algorithms, we randomly pick 5 convergent problems and draw two ratio trajectories for each problem: one is the sample size of the gradient~over~the sample size of the value, and one is the sample size of the Hessian over the~sample size of the gradient. We take $C = 64$ as an example. The plot is shown in~Figure \ref{fig:3}. From the figure, we note that the sample size ratios tend to be stabilized at a small level, and the trend is more evident when $\sigma^2=0.1$. As we explained~for Figure \ref{fig:2} above, such an observation is consistent with our discussions in Section \ref{sec:4.3}, and illustrates the improvement of our analysis over \cite{Na2022adaptive} for performing the stochastic line search on the augmented Lagrangian merit function.

\vskip4pt
\noindent\textbf{Stepsize trajectories.} Figure \ref{fig:4} plots the stepsize trajectories that are selected by stochastic line search. We take the default setup as an example, i.e., $C=\kappa=2$, $\chi_{err}=1$. Similar to Figure \ref{fig:3}, for each level of $\sigma^2$, we randomly pick $5$ convergent problems to show the trajectories. Although there is no clear trend for the stepsize trajectories due to stochasticity, we clearly see for both~methods that the stepsize can increase significantly from a very small value and even exceed $1$. This exclusive property of the line search procedure ensures a fast convergence of the scheme, which is not enjoyed by many non-adaptive schemes where the stepsize often monotonically decays to zero.

We also examine some other aspects of the algorithm, such as the~proportion of the iterations with failed SQP steps, with unstabilized penalty parameters,~or with a triggered feasibility error condition \eqref{cond:bound:fes:error}. We also study the effect of a multiplicative noise, and implement the algorithm on an inequality  constrained logistic regression problem. Due to the space limit, these auxiliary experiments are provided in Appendix \ref{sec:aux:exp}.

\section{Conclusion}\label{sec:6}

This paper studied inequality constrained stochastic nonlinear optimization problems. We designed an active-set StoSQP algorithm that exploits the exact augmented Lagrangian merit function. The algorithm adaptively selects the penalty parameters of the augmented Lagrangian, and selects the stepsize~via stochastic line search. We proved that the KKT residuals converge to zero~almost surely, which generalizes and strengthens the result for unconstrained and equality constrained problems in \citep{Paquette2020Stochastic, Na2022adaptive} to enable wider applications.

The extension of this work includes studying more advanced StoSQP~schemes. As mentioned in Section \ref{sec:4.4}, the proposed StoSQP scheme has to solve the~SQP system exactly. We note that, recently, \cite{Curtis2021Inexact} designed a StoSQP scheme where~an inexact Newton direction is employed, and  \cite{Berahas2021Stochastic} designed a StoSQP scheme to relax LICQ condition. It is still open how to design related schemes to achieve relaxations with inequality constraints. In addition, some advanced SQP schemes solve inequality constrained problems by mixing IQP with EQP: one solves a \textit{convex} IQP to obtain an active set, and then solves an EQP to obtain the~search direction. See the ``SQP+" scheme in \cite{Morales2011sequential} for example.  Investigating this kind of mixed scheme with a stochastic objective is promising. Besides SQP,~there~are other classical methods for solving nonlinear problems that can be exploited to deal with stochastic objectives, such as the augmented Lagrangian methods and interior point methods. Different methods have different benefits and all of them deserve studying in the setup where the model can only be accessed with certain noise.

Finally, as mentioned in Section \ref{sec:4.2}, non-asymptotic analysis and iteration complexity of the proposed scheme are missing in our global analysis. Further, it is known for deterministic setting that differentiable merit functions can overcome the Maratos effect and facilitate a fast local rate, while non-smooth~merit functions (without advanced local modifications) cannot. This~raises~the questions: what is the local rate of the proposed StoSQP, and is the local rate better than the one using non-smooth merit functions? To answer these~questions, we need a better understanding on the local behavior of stochastic line search. Such a local study would complement the established global analysis, recognize the benefits of the differentiable merit functions, and bridge the understanding gap between stochastic SQP and deterministic SQP.

\section*{Acknowledgments}

We thank Associate Editor and two anonymous reviewers for instructive comments, which help us further enhance the algorithm design and presentation. This material was~\mbox{completed} in part with resources provided by the University of Chicago Research Computing Center. This material was based upon work supported by the U.S. Department of Energy, Office of Science, Office of Advanced Scientific Computing Research (ASCR) under Contract DE-AC02-06CH11347 and by NSF through award CNS-1545046.

\appendix
\numberwithin{equation}{section}
\numberwithin{theorem}{section}

\section{Proofs of Section \ref{sec:2}}

\subsection{Proof of Lemma \ref{lem:2}}\label{pf:lem:2}

Throughout the proof, we denote $g^\star = g(\bx^\star)$, $\bw_{\epsilon, \nu}^\star = \bw_{\epsilon, \nu}(\tx, \tlambda)$, $\nabla\mL^\star = \nabla\mL(\tx, \tmu, \tlambda)$ (similar for $c^\star$, $a_{\nu}^\star$, $q_{\nu}^\star$ etc.) to be the quantities evaluated at $(\tx, \tmu, \tlambda)\in\mT_{\nu}\times \mR^m\times\mR^r$.~Since $\bw_{\epsilon, \nu}^\star = \0$, we know from Lemma \ref{lem:1} that $g^\star\leq \0$, $\tlambda \geq \0$, $(\tlambda)^T g^\star = 0$. This implies that $\diag^2(g^\star)\tlambda = \0$. Furthermore, by $c^\star = \0$, $\bw_{\epsilon, \nu}^\star = \0$, $a_{\nu}^\star, \eta, \epsilon >0$, and $\nabla_{\bmu, \blambda}\mL_{\epsilon, \nu, \eta}^\star = \0$, we obtain from \eqref{equ:aug:der} that
\begin{equation}\label{pequ:1}
\begin{pmatrix}
M_{11}^\star & M_{12}^\star\\
M_{21}^\star & M_{22}^\star
\end{pmatrix}\begin{pmatrix}
J^\star\\
G^\star
\end{pmatrix}\nabla_{\bx}\mL^\star = \0.
\end{equation}
Recalling the definition of $M^\star$ in \eqref{equ:def:Qmatrices}, we multiply the matrix $\nabla_{\bx}^T\mL^\star((J^\star)^T\;\; (G^\star)^T)$ from the left and obtain
\begin{align*}
\0 \;\;& \stackrel{\mathclap{\eqref{pequ:1}}}{=}\;\; \nabla_{\bx}^T\mL^\star\begin{pmatrix}
(J^\star)^T & (G^\star)^T
\end{pmatrix}\begin{pmatrix}
J^\star (J^\star)^T & J^\star (G^\star)^T\\
G^\star (J^\star)^T & G^\star (G^\star)^T + \diag^2(g^\star)
\end{pmatrix}\begin{pmatrix}
J^\star\\
G^\star
\end{pmatrix}\nabla_{\bx}\mL^\star\\
& = \nbr{\rbr{(J^\star)^TJ^\star + (G^\star)^TG^\star}\nabla_{\bx}\mL^\star}^2 + \nbr{\diag(g^\star)G^\star\nabla_{\bx}\mL^\star}^2.
\end{align*}
This implies $\rbr{(J^\star)^TJ^\star + (G^\star)^TG^\star}\nabla_{\bx}\mL^\star = \0$. Multiplying $\nabla_{\bx}\mL^\star$ from the left, we have $J^\star\nabla_{\bx}\mL^\star = \0$ and $G^\star\nabla_{\bx}\mL^\star = \0$. Plugging into \eqref{equ:aug:der} and noting that $\nabla_{\bx}\mL_{\epsilon, \nu, \eta}^{\star} = \0$, $\bw_{\epsilon, \nu}^\star = \0$, $c^\star = \0$, $\diag^2(g^\star)\tlambda = \0$, and $q_{\nu}^\star, a_{\nu}^\star, \epsilon>0$, we obtain $\nabla_{\bx}\mL^\star = \0$. This shows $(\tx, \tmu, \tlambda)$ satisfies \eqref{equ:KKT:cond}, and we complete the proof.

\subsection{Proof of Lemma \ref{lem:6}}\label{pf:lem:6}

We require the following two preparation lemmas.

\begin{lemma}\label{lem:5}
Let $\I(\tx)$ be the active set defined in \eqref{equ:I}, and $\I^+(\tx, \tlambda) = \cbr{i\in \I(\tx): \tlambda_i>0}$. For any $\epsilon, \nu>0$, there exists a compact set $\mX_{\epsilon,\nu}\times \Lambda_{\epsilon,\nu}\ni(\tx, \tlambda)$ depending on $(\epsilon,\nu)$,~such~that
\begin{equation*}
\I^+(\tx, \tlambda) \subseteq \mA_{\epsilon, \nu}(\bx, \blambda) \subseteq \I(\tx), \quad\quad \forall (\bx, \blambda) \in \mX_{\epsilon,\nu}\times\Lambda_{\epsilon,\nu}.
\end{equation*}
\end{lemma}

\begin{proof}
See Appendix \ref{pf:lem:5}.
\end{proof}

\begin{lemma}\label{lem:aux:1}

Under Assumption \ref{ass:1}, there exist a compact set $X\ni\tx$ and a constant~$\gamma_{H}\in~(0,1]$ such that $M(\bx)\succeq\gamma_{H} I$ for any $\bx\in X$, where $M(\bx)$ is defined in \eqref{equ:def:Qmatrices}. Furthermore, for any $\epsilon,\nu>0$, there exists a compact set $\mX_{\epsilon,\nu}\times \Lambda_{\epsilon,\nu}\ni (\tx, \tlambda)$ depending on $(\epsilon, \nu)$, such~that
\begin{equation*}
\begin{pmatrix}
J(\bx)\\
G_{\mA_{\epsilon, \nu}(\bx, \blambda)}(\bx)
\end{pmatrix}\begin{pmatrix}
J(\bx)^T & G_{\mA_{\epsilon, \nu}(\bx, \blambda)}(\bx)^T
\end{pmatrix}\succeq\gamma_{H} I, \quad\quad \forall (\bx, \blambda) \in \mX_{\epsilon,\nu}\times\Lambda_{\epsilon,\nu}. 
\end{equation*}
\end{lemma}

\begin{proof}
See Appendix \ref{pf:lem:aux:1}.
\end{proof}

We now prove Lemma \ref{lem:6}. We suppress the evaluation point and the iteration index~$t$.~Let $\mX\times \mM\times \Lambda \subseteq \mT_{\nu}\times \mR^m\times \mR^r$ be any compact set around $(\tx,\tmu, \tlambda)$ (independent of $\epsilon, \nu, \eta$) and suppose $(\bx, \bmu, \blambda)\in \mX\times \mM\times \Lambda$. By Lemma \ref{lem:aux:1}, we know there exist a constant~$\gamma_{H}\in(0,1]$~and, for any $\epsilon, \nu>0$, a compact subset $\mX_{\epsilon,\nu}\times\Lambda_{\epsilon,\nu} \subseteq \mX\times\Lambda$ such that for any point in the subset,
\begin{equation}\label{bound:M}
M \succeq \gamma_{H} I \quad \text{ and }\quad \begin{pmatrix}
J\\
G_a
\end{pmatrix}\begin{pmatrix}
J^T & G_a^T
\end{pmatrix}\succeq\gamma_{H} I.
\end{equation}
Thus, by Assumption \ref{ass:2}, we know from \cite[Lemma 16.1]{Nocedal2006Numerical} that $K_a$ is invertible, and thus~\eqref{equ:SQP:direction}~is solvable. Furthermore, we can also show that (see \cite[Lemma 1]{Na2022adaptive} for a simple proof)
\begin{equation}\label{bound:Ka}
\|K_a^{-1}\| \leq 7\Upsilon_{B}^2/(\gamma_{B}\gamma_{H}).
\end{equation}
With the above two results, we conduct our analysis. Throughout the proof, we use $\Upsilon_1, \Upsilon_2\ldots$~to denote generic upper bounds of functions evaluated in the set $\mX\times \mM\times \Lambda$, which are~independent of $(\epsilon, \nu, \eta, \gamma_{B}, \gamma_{H})$. As they are upper bounds, without loss of generality, $\Upsilon_i \geq 1$,~$\forall i$.

We start from $(\nabla\mL_{\epsilon, \nu, \eta}^{(1)})^T\Delta$ and suppose $(\bx, \bmu, \blambda)\in \mX_{\epsilon,\nu}\times\mM\times\Lambda_{\epsilon,\nu}\subseteq\mX\times\mM\times\Lambda$, where $\mX_{\epsilon,\nu}$ and $\Lambda_{\epsilon,\nu}$ come from Lemma \ref{lem:aux:1}. We have
\begin{align}\label{pequ:10}
& (\nabla\mL_{\epsilon, \nu, \eta}^{(1)})^T\Delta \nonumber\\
& \stackrel{\mathclap{\eqref{equ:new:aug:der}}}{=} \Delta\bx^T\nabla_{\bx}\mL + \eta\Delta\bx^T\begin{pmatrix}
Q_1 & Q_2
\end{pmatrix}\left(\begin{smallmatrix}
J\nabla_{\bx}\mL\\
G\nabla_{\bx}\mL+ \Pi_c(\diag^2(g)\blambda)
\end{smallmatrix}\right) + \frac{1}{\epsilon}\Delta\bx^TJ^Tc  \nonumber\\
& \quad \quad + \frac{1}{\epsilon q_{\nu}}\Delta\bx^TG^T\bw_{\epsilon, \nu} + \left(\begin{smallmatrix}
\Delta\bmu\\
\Delta\blambda
\end{smallmatrix}\right)^T\left(\begin{smallmatrix}
c\\
\bw_{\epsilon, \nu}
\end{smallmatrix}\right)+ \eta\left(\begin{smallmatrix}
\Delta\bmu\\
\Delta\blambda
\end{smallmatrix}\right)^T\left(\begin{smallmatrix}
M_{11} & M_{12}\\
M_{21} & M_{22}
\end{smallmatrix}\right)\left(\begin{smallmatrix}
J\nabla_{\bx}\mL\\
G\nabla_{\bx}\mL+ \Pi_c(\diag^2(g)\blambda)
\end{smallmatrix} \right) \nonumber\\
& \stackrel{\mathclap{\eqref{equ:SQP:direction:2}}}{=} \Delta\bx^T\nabla_{\bx}\mL + \frac{1}{\epsilon}\Delta\bx^TJ^Tc + \frac{1}{\epsilon q_{\nu}}\Delta\bx^TG^T\bw_{\epsilon, \nu} + \left(\begin{smallmatrix}
\Delta\bmu\\
\Delta\blambda
\end{smallmatrix}\right)^T\left(\begin{smallmatrix}
c\\
\bw_{\epsilon, \nu}
\end{smallmatrix}\right) - \eta \nbr{\left(\begin{smallmatrix}
J\nabla_{\bx}\mL\\
G\nabla_{\bx}\mL+ \Pi_c(\diag^2(g)\blambda)
\end{smallmatrix} \right)}^2 \nonumber\\
&\stackrel{\mathclap{\substack{\eqref{equ:def:w}\\ \eqref{equ:active}}}}{=} \Delta\bx^T(\nabla_{\bx}\mL - G_c^T\blambda_c) + \frac{1}{\epsilon}\Delta\bx^TJ^Tc + \frac{1}{\epsilon q_\nu}\Delta\bx^TG_a^Tg_a + \left(\begin{smallmatrix}
c\\
g_a
\end{smallmatrix}\right)^T\left(\begin{smallmatrix}
\Delta\bmu\\
\Delta\blambda_a
\end{smallmatrix}\right) - \epsilon q_{\nu}\Delta\blambda_c^T\blambda_c  \nonumber\\
&\quad\quad - \eta \nbr{\left(\begin{smallmatrix}
J\nabla_{\bx}\mL\\
G\nabla_{\bx}\mL+ \Pi_c(\diag^2(g)\blambda)	\end{smallmatrix}\right) }^2 \nonumber\\
& \stackrel{\mathclap{\eqref{equ:SQP:direction:1}}}{=} - \Delta\bx^TB\Delta\bx + \left(\begin{smallmatrix}
c\\
g_a
\end{smallmatrix}\right)^T\left(\begin{smallmatrix}
\tDelta\bmu + \Delta\bmu\\
\tDelta\blambda_a + \Delta\blambda_a
\end{smallmatrix}\right) - \frac{1}{\epsilon}\|c\|^2 - \frac{1}{\epsilon q_{\nu}}\|g_a\|^2 - \epsilon q_{\nu}\Delta\blambda_c^T\blambda_c \nonumber\\
&\quad\quad - \eta \nbr{\left(\begin{smallmatrix}
J\nabla_{\bx}\mL\\
G\nabla_{\bx}\mL+ \Pi_c(\diag^2(g)\blambda)	\end{smallmatrix}\right) }^2.
\end{align}
Since $(\bx, \bmu, \blambda)\in\mX\times\mM\times\Lambda$, there exists $\Upsilon_1\geq 1$ such that $\nbr{(Q_1\; Q_2)} \leq \Upsilon_1$. Thus, we have
\begin{align}\label{pequ:11}
\nbr{\begin{pmatrix}
\Delta\bmu\\
\Delta\blambda
\end{pmatrix}}  & \;\; \stackrel{\mathclap{\eqref{equ:SQP:direction:2}}}{=} \;\; \nbr{M^{-1}\begin{pmatrix}
J\nabla_{\bx}\mL\\
G\nabla_{\bx}\mL+ \Pi_c(\diag^2(g)\blambda)	\end{pmatrix} + M^{-1}\begin{pmatrix}
Q_1^T\\
Q_2^T
\end{pmatrix}\Delta\bx } \nonumber\\
&\;\; \stackrel{\mathclap{\eqref{bound:M}} }{\leq} \;\; \frac{1}{\gamma_{H}}\nbr{\begin{pmatrix}
J\nabla_{\bx}\mL\\
G\nabla_{\bx}\mL+ \Pi_c(\diag^2(g)\blambda)	\end{pmatrix} } + \frac{\Upsilon_1}{\gamma_{H}} \nbr{\Delta\bx}\nonumber\\
& \;\;\leq   \frac{2\Upsilon_1}{\gamma_{H}} \nbr{\left(\begin{smallmatrix}
\Delta\bx\\
J\nabla_{\bx}\mL\\
G\nabla_{\bx}\mL+ \Pi_c(\diag^2(g)\blambda)
\end{smallmatrix}\right)  } \quad\quad (\text{since } 1\leq \Upsilon_1).
\end{align}
Moreover, we note that
\begin{align*}
&\cbr{ \left(\begin{smallmatrix}
J\\
G_a\\
G_c
\end{smallmatrix}\right)\left(\begin{smallmatrix}
J^T & G_a^T & G_c
\end{smallmatrix}\right) + \left(\begin{smallmatrix}
\0 & \0 & \0 \\
\0 & \diag^2(g_a) & \0\\
\0 & \0 & \diag^2(g_c)
\end{smallmatrix}\right)  } \left(\begin{smallmatrix}
\tDelta\bmu\\
\tDelta\blambda_a\\
-\blambda_c
\end{smallmatrix}\right)\\
& \stackrel{{\eqref{equ:SQP:direction:1}}}{=} - \left(\begin{smallmatrix}
J\\
G_a\\
G_c
\end{smallmatrix}\right) B\Delta\bx - \left(\begin{smallmatrix}
J\nabla_{\bx}\mL\\
G_a\nabla_{\bx}\mL\\
G_c\nabla_{\bx}\mL + \diag^2(g_c)\blambda_c
\end{smallmatrix}\right) + \left(\begin{smallmatrix}
\0 \\
\diag^2(g_a)\tDelta\blambda_a\\
\0
\end{smallmatrix}\right)\\
& \stackrel{\eqref{equ:SQP:direction:1}}{=}- \left(\begin{smallmatrix}
J\\
G_a\\
G_c
\end{smallmatrix}\right) B\Delta\bx - \left(\begin{smallmatrix}
J\nabla_{\bx}\mL\\
G_a\nabla_{\bx}\mL\\
G_c\nabla_{\bx}\mL + \diag^2(g_c)\blambda_c
\end{smallmatrix}\right) - \left(\begin{smallmatrix}
\0 \\
\diag(g_a)\diag(\tDelta\blambda_a)G_a\Delta\bx\\
\0
\end{smallmatrix}\right).
\end{align*}
By $(\bx, \bmu, \blambda)\in\mX\times\mM\times\Lambda$, there exist $\Upsilon_2, \Upsilon_3, \Upsilon_4\geq 1$ such that
\begin{equation*}
\nbr{ (J^T\; G^T) } \leq \Upsilon_2,\quad\quad \|\tDelta\blambda_a\| \stackrel{\eqref{equ:SQP:direction:1}}{\leq} \nbr{K_a^{-1}\left(\begin{smallmatrix}
\nabla_{\bx}\mL - G_c^T\blambda_c\\
c\\
g_a
\end{smallmatrix}\right)} \stackrel{\eqref{bound:Ka}}{\leq}\frac{\Upsilon_3}{\gamma_{H}\gamma_{B}},
\end{equation*}
and
\begin{equation*}
\|\diag(g_a)\diag(\tDelta\blambda_a)G_a\|\leq \frac{\Upsilon_4}{\gamma_{H}\gamma_{B}}.
\end{equation*}
Combining the above three displays,
\begin{multline}\label{pequ:12}
\nbr{\left(\begin{smallmatrix}
\tDelta\bmu\\
\tDelta\blambda_a\\
-\blambda_c
\end{smallmatrix}\right)}  \;\stackrel{\mathclap{\eqref{bound:M}}}{\leq}\frac{1}{\gamma_{H}} \rbr{\Upsilon_2\|B\Delta\bx\| + \frac{\Upsilon_4}{\gamma_{H}\gamma_{B}}\|\Delta\bx\| } + \frac{1}{\gamma_{H}}\nbr{\begin{pmatrix}
J\nabla_{\bx}\mL\\
G\nabla_{\bx}\mL+ \Pi_c(\diag^2(g)\blambda)
\end{pmatrix}} \\
\leq \frac{\Upsilon_2\Upsilon_B+\Upsilon_4 +1}{\gamma_{H}^2\gamma_{B}}\nbr{\left(\begin{smallmatrix}
\Delta\bx\\
J\nabla_{\bx}\mL\\
G\nabla_{\bx}\mL+ \Pi_c(\diag^2(g)\blambda)
\end{smallmatrix}\right) } \quad (\text{since } \gamma_{H}\vee\gamma_{B} \leq 1),
\end{multline}
where the second inequality also uses $\|B\|\leq \Upsilon_{B}$ by Assumption \ref{ass:2}. Combining \eqref{pequ:10}, \eqref{pequ:11}, \eqref{pequ:12}, and using $0<q_\nu\leq \nu$ and $\gamma_{H}\vee\gamma_{B} \leq 1$,
\begin{align}\label{pequ:13}
& (\nabla\mL_{\epsilon, \nu, \eta}^{(1)})^T\Delta \nonumber\\
& \stackrel{\mathclap{\eqref{pequ:10}}}{\leq}  -\Delta\bx^TB\Delta\bx + \nbr{\left(\begin{smallmatrix}
c\\
g_a
\end{smallmatrix}\right)}\cbr{\nbr{\left(\begin{smallmatrix}
\tDelta\bmu\\
\tDelta\blambda_a
\end{smallmatrix}\right)} + \nbr{\left(\begin{smallmatrix}
\Delta\bmu\\
\Delta\blambda_a
\end{smallmatrix}\right)}} - \frac{1}{\epsilon(1\vee \nu)}\nbr{\left(\begin{smallmatrix}
c\\
g_a
\end{smallmatrix}\right)}^2 + \epsilon\nu \|\Delta\blambda_c\|\|\blambda_c\| \nonumber\\
&\quad\quad  - \eta\nbr{\left(\begin{smallmatrix}
J\nabla_{\bx}\mL\\
G\nabla_{\bx}\mL+ \Pi_c(\diag^2(g)\blambda)	\end{smallmatrix}\right)}^2 \nonumber\\
& \stackrel{\mathclap{\substack{\substack{\eqref{pequ:11} \\ \eqref{pequ:12}}}}}{\leq} -\Delta\bx^TB\Delta\bx + \frac{2\Upsilon_1 + \Upsilon_2\Upsilon_B + \Upsilon_4 +1}{\gamma_{H}^2\gamma_{B}}\nbr{\left(\begin{smallmatrix}
c\\
g_a
\end{smallmatrix}\right)}\nbr{\left(\begin{smallmatrix}
\Delta\bx\\
J\nabla_{\bx}\mL\\
G\nabla_{\bx}\mL+ \Pi_c(\diag^2(g)\blambda)
\end{smallmatrix}\right) } - \frac{1}{\epsilon(1\vee \nu)}\nbr{\left(\begin{smallmatrix}
c\\
g_a
\end{smallmatrix}\right)}^2 \nonumber\\
&\quad\quad   + \epsilon\nu\cdot \frac{2\Upsilon_1(\Upsilon_2\Upsilon_B + \Upsilon_4 +1)}{\gamma_{H}^3\gamma_{B}} \nbr{\left(\begin{smallmatrix}
\Delta\bx\\
J\nabla_{\bx}\mL\\
G\nabla_{\bx}\mL+ \Pi_c(\diag^2(g)\blambda)
\end{smallmatrix} \right)}^2  - \eta\nbr{\left(\begin{smallmatrix}
J\nabla_{\bx}\mL\\
G\nabla_{\bx}\mL+ \Pi_c(\diag^2(g)\blambda)	\end{smallmatrix}\right)}^2 \nonumber\\
& \leq  -\Delta\bx^TB\Delta\bx + \frac{\Upsilon_5}{\gamma_{H}^2\gamma_{B}}\nbr{\left(\begin{smallmatrix}
c\\
g_a
\end{smallmatrix}\right)}\nbr{\left(\begin{smallmatrix}
\Delta\bx\\
J\nabla_{\bx}\mL\\
G\nabla_{\bx}\mL+ \Pi_c(\diag^2(g)\blambda)
\end{smallmatrix}\right) } - \frac{1}{\epsilon(1\vee \nu)}\nbr{\left(\begin{smallmatrix}
c\\
g_a
\end{smallmatrix}\right)}^2 \nonumber\\
& \quad\quad + \frac{\epsilon\nu\Upsilon_5}{\gamma_{H}^3\gamma_{B}} \nbr{\left(\begin{smallmatrix}
\Delta\bx\\
J\nabla_{\bx}\mL\\
G\nabla_{\bx}\mL+ \Pi_c(\diag^2(g)\blambda)
\end{smallmatrix}\right) }^2 - \eta\nbr{\left(\begin{smallmatrix}
J\nabla_{\bx}\mL\\
G\nabla_{\bx}\mL+ \Pi_c(\diag^2(g)\blambda)	\end{smallmatrix}\right)}^2,
\end{align}
where the last inequality holds by defining 
\begin{equation*}
\Upsilon_5 = 2\Upsilon_1 + \Upsilon_2\Upsilon_B + \Upsilon_4 +1 \vee 2\Upsilon_1(\Upsilon_2\Upsilon_B + \Upsilon_4 +1).
\end{equation*}
To deal with $\Delta\bx^TB\Delta\bx$ in \eqref{pequ:13}, we decompose $\Delta\bx$ as $\Delta\bx = \Delta\bu + \Delta\bv$ where $\Delta\bu\in~\text{Image}\cbr{(J^T\;\; G_a^T)}$ and $\Delta\bv\in \text{Ker}\cbr{(J^T\;\; G_a^T)^T}$. Note that
\begin{align}\label{pequ:14}
-\begin{pmatrix}
c\\
g_a
\end{pmatrix} = \begin{pmatrix}
J\\
G_a
\end{pmatrix}\Delta\bx = \begin{pmatrix}
J\\
G_a
\end{pmatrix}\Delta\bu & \Longrightarrow \Delta\bmu = -\begin{pmatrix}
J^T & G_a^T
\end{pmatrix}\cbr{\begin{pmatrix}
J\\
G_a
\end{pmatrix}\begin{pmatrix}
J^T & G_a^T
\end{pmatrix}}^{-1}\begin{pmatrix}
c\\
g_a
\end{pmatrix} \nonumber\\
& \stackrel{\mathclap{\eqref{bound:M}}}{\Longrightarrow}  \|\Delta\bmu\| \leq \frac{1}{\sqrt{\gamma_{H}}} \nbr{\begin{pmatrix}
c\\
g_a
\end{pmatrix}}.
\end{align}
Thus, by Assumption \ref{ass:2},
\begin{align}\label{pequ:31}
- & \Delta\bx^TB\Delta\bx \nonumber\\
& =  -\Delta\bv^TB\Delta\bv -2\Delta\bu^TB\Delta\bv - \Delta\bu^TB\Delta\bu
\leq -\gamma_{B} \|\Delta\bv\|^2 + 2\Upsilon_{B}\|\Delta\bv\|\|\Delta\bu\| + \Upsilon_{B}\|\Delta\bu\|^2 \nonumber\\
& \leq -\frac{3\gamma_{B}}{4}\|\Delta\bv\|^2 + (\Upsilon_{B} + \frac{4\Upsilon_{B}^2}{\gamma_{B}})\|\Delta\bu\|^2
=  -\frac{3\gamma_{B}}{4}\|\Delta\bx\|^2 + (\Upsilon_{B} + \frac{4\Upsilon_{B}^2}{\gamma_{B}} + \frac{3\gamma_{B}}{4})\|\Delta\bu\|^2 \nonumber\\
&\stackrel{\mathclap{\eqref{pequ:14}}}{\leq}  -\frac{3\gamma_{B}}{4}\|\Delta\bx\|^2 + (\Upsilon_{B} + \frac{4\Upsilon_{B}^2}{\gamma_{B}} + \frac{3\gamma_{B}}{4})\frac{1}{\gamma_{H}}\nbr{\begin{pmatrix}
c\\
g_a
\end{pmatrix}}^2 \nonumber\\
& \leq -\frac{3\gamma_{B}}{4}\|\Delta\bx\|^2 + \frac{\Upsilon_6}{\gamma_{H}\gamma_{B}}\nbr{\begin{pmatrix}
c\\
g_a
\end{pmatrix}}^2,
\end{align}
where the last inequality holds with $\Upsilon_6 = \Upsilon_B + 4\Upsilon_B^2 + 1$ by noting that $\gamma_{B}\leq 1$. Combining the above display with \eqref{pequ:13} and using the following Young's inequality,
\begin{multline*}
\frac{\Upsilon_5}{\gamma_{H}^2\gamma_{B}}\nbr{\left(\begin{smallmatrix}
c\\
g_a
\end{smallmatrix}\right)}\nbr{ \left(\begin{smallmatrix}
\Delta\bx\\
J\nabla_{\bx}\mL\\
G\nabla_{\bx}\mL+ \Pi_c(\diag^2(g)\blambda)
\end{smallmatrix}\right) } \\
\leq \rbr{\frac{\gamma_{B}}{8}\wedge \frac{\eta}{4}}\nbr{\left(\begin{smallmatrix}
\Delta\bx\\
J\nabla_{\bx}\mL\\
G\nabla_{\bx}\mL+ \Pi_c(\diag^2(g)\blambda)
\end{smallmatrix}\right) }^2 + \frac{2\Upsilon_5^2}{\gamma_{H}^4\gamma_{B}^2(\gamma_{B} \wedge \eta)}\nbr{\left(\begin{smallmatrix}
c\\
g_a
\end{smallmatrix}\right)}^2,
\end{multline*}
we have
\begin{align*}
(\nabla\mL_{\epsilon, \nu, \eta}^{(1)})^T\Delta \leq & -\frac{3\gamma_{B}}{4}\|\Delta\bx\|^2  + \cbr{\rbr{\frac{\gamma_{B}}{8}\wedge \frac{\eta}{4}} + \frac{\epsilon\nu \Upsilon_5}{\gamma_{H}^3\gamma_{B}}}\nbr{\left(\begin{smallmatrix}
\Delta\bx\\
J\nabla_{\bx}\mL\\
G\nabla_{\bx}\mL+ \Pi_c(\diag^2(g)\blambda)
\end{smallmatrix}\right) }^2\\
& \hskip-1.5cm + \cbr{\frac{\Upsilon_6}{\gamma_{H}\gamma_{B}} + \frac{2\Upsilon_5^2}{\gamma_{H}^4\gamma_{B}^2(\gamma_{B} \wedge \eta)} - \frac{1}{\epsilon(1\vee \nu)}}\nbr{\left(\begin{smallmatrix}
c\\
g_a
\end{smallmatrix}\right)}^2 -\eta\nbr{\left(\begin{smallmatrix}
J\nabla_{\bx}\mL\\
G\nabla_{\bx}\mL+ \Pi_c(\diag^2(g)\blambda)	\end{smallmatrix}\right)}^2\\
& \hskip-1.8cm \leq  -\cbr{\frac{\gamma_{B}\wedge \eta}{2} + \rbr{\frac{\gamma_{B}}{8} \wedge \frac{\eta}{4}} - \frac{\epsilon\nu \Upsilon_5}{\gamma_{H}^3\gamma_{B}} } \nbr{\left(\begin{smallmatrix}
\Delta\bx\\
J\nabla_{\bx}\mL\\
G\nabla_{\bx}\mL+ \Pi_c(\diag^2(g)\blambda)
\end{smallmatrix}\right) }^2\\
& \hskip-1.5cm - \cbr{\frac{1}{\epsilon(1\vee \nu)}  - \frac{\Upsilon_6}{\gamma_{H}\gamma_{B}} - \frac{2\Upsilon_5^2}{\gamma_{H}^4\gamma_{B}^2(\gamma_{B} \wedge \eta)} }\nbr{\left(\begin{smallmatrix}
c\\
g_a
\end{smallmatrix}\right)}^2.
\end{align*}
Therefore, as long as
\begin{subequations}\label{pequ:19}
\begin{align}
\frac{\gamma_{B}}{8}\wedge\frac{\eta}{4} \geq \frac{\epsilon\nu\Upsilon_5}{\gamma_{H}^3\gamma_{B}} \Longleftarrow& \frac{1}{\epsilon} \geq \frac{8\nu\Upsilon_5}{\gamma_{H}^3\gamma_{B}(\gamma_{B}\wedge\eta)}, \label{pequ:19:a}\\
\frac{1}{\epsilon(1\vee \nu)} - \frac{\Upsilon_6}{\gamma_{H}\gamma_{B}} -\frac{2\Upsilon_5^2}{\gamma_{H}^4\gamma_{B}^2(\gamma_{B} \wedge \eta)}  \geq 0 \Longleftarrow & \frac{1}{\epsilon} \geq \frac{(1\vee \nu )(2\Upsilon_5^2 + \Upsilon_6)}{\gamma_{H}^4\gamma_{B}^2(\gamma_{B} \wedge \eta)},\label{pequ:19:b}
\end{align}	
\end{subequations}
we have
\begin{equation*}
(\nabla\mL_{\epsilon, \nu, \eta}^{(1)})^T\Delta \leq - \frac{\gamma_{B} \wedge \eta}{2}\nbr{\left(\begin{smallmatrix}
\Delta\bx\\
J\nabla_{\bx}\mL\\
G\nabla_{\bx}\mL+ \Pi_c(\diag^2(g)\blambda)
\end{smallmatrix}\right) }^2.
\end{equation*}
Thus, letting $\Upsilon = \cbr{8\Upsilon_5 \vee (2\Upsilon_5^2 + \Upsilon_6)}/\gamma_{H}^4$ and noting that \eqref{pequ:19:a} is implied by \eqref{pequ:19:b}, we complete the first part of the statement.

We now prove the second part of the statement. By \eqref{equ:new:aug:der}, $(\bx,\bmu,\blambda)\in \mX\times\mM\times\Lambda$ (and hence \eqref{pequ:8}), and the fact that $a_\nu \geq \nu/2$, there exists $\Upsilon_7\geq 1$ such that
\begin{align*}
(\nabla\mL_{\epsilon, \nu, \eta}^{(2)})^T\Delta & \;\stackrel{\mathclap{\eqref{equ:new:aug:der}}}{=} \; \frac{3\|\bw_{\epsilon, \nu}\|^2}{2\epsilon q_{\nu}a_{\nu}} \Delta\bx^TG^T\bl + \eta \Delta\bx^TQ_{2,a}\diag^2(g_a)\blambda_a + \frac{\|\bw_{\epsilon, \nu}\|^2}{\epsilon a_{\nu}}\Delta\blambda^T \blambda \\
& \;\quad + \eta (\Delta\bmu^T \;\; \Delta\blambda^T)\begin{pmatrix}
M_{12,a}\\
M_{22,a}
\end{pmatrix}\diag^2(g_a)\blambda_a\\
&\; \leq \Upsilon_7\bigg\{\frac{1}{\epsilon \nu^2}\rbr{\|g_a\|^2 + \epsilon^2\nu^2\|\blambda_c\|^2 }\|\Delta\bx\|+ \eta \|g_a\|^2\|\Delta\bx\| \\
&\;\quad  + \frac{1}{\epsilon \nu}\rbr{\|g_a\|^2 + \epsilon^2\nu^2\|\blambda_c\|^2 }\|\Delta\blambda\|  + \eta \|g_a\|^2\|(\Delta\bmu, \Delta\blambda)\|\bigg\}.
\end{align*}
Since $\epsilon \leq 1$ by \eqref{pequ:19} (noting that $\Upsilon\geq 1 \geq \gamma_{H}\vee \gamma_{B}$), we simplify the above display by
\begin{align*}
(\nabla\mL_{\epsilon, \nu, \eta}^{(2)})^T\Delta & \leq  \Upsilon_7\cbr{\frac{1\vee \nu^2}{\epsilon\nu(1\wedge \nu)} (\|g_a\|^2 + \|\blambda_c\|^2)(\|\Delta\bx\| + \|\Delta\blambda\|) + \sqrt{2}\eta\|g_a\|^2 \|(\Delta\bx, \Delta\bmu, \Delta\blambda)\| }\\
&\leq  \sqrt{2}\Upsilon_7\cbr{\frac{1\vee \nu^2}{\epsilon\nu(1\wedge \nu)} (\|g_a\|^2 + \|\blambda_c\|^2)\|(\Delta\bx, \Delta\blambda)\|+ \eta\|g_a\|^2 \|(\Delta\bx, \Delta\bmu, \Delta\blambda)\|}\\
& \leq 2\sqrt{2}\Upsilon_7\rbr{\frac{1\vee \nu}{\epsilon (1\wedge \nu^2)}\vee \eta}(\|g_a\|^2 + \|\blambda_c\|^2) \|(\Delta\bx, \Delta\bmu, \Delta\blambda)\|.
\end{align*}
Noting that
\begin{align*}
\nbr{\left(\begin{smallmatrix}
\Delta\bx\\
\Delta\bmu\\
\Delta\blambda
\end{smallmatrix}\right)}& \leq \|\Delta\bx\| + \nbr{\begin{pmatrix}
\Delta\bmu\\
\Delta\blambda
\end{pmatrix} }\; \stackrel{\mathclap{\eqref{pequ:11}}}{\leq} \;\|\Delta\bx\| + \frac{2\Upsilon_1}{\gamma_{H}}\nbr{\left(\begin{smallmatrix}
\Delta\bx\\
J\nabla_{\bx}\mL\\
G\nabla_{\bx}\mL+ \Pi_c(\diag^2(g)\blambda)
\end{smallmatrix}\right) }\\
& \leq \frac{3\Upsilon_1}{\gamma_{H}}\nbr{\left(\begin{smallmatrix}
\Delta\bx\\
J\nabla_{\bx}\mL\\
G\nabla_{\bx}\mL+ \Pi_c(\diag^2(g)\blambda)
\end{smallmatrix}\right) } \quad\quad (\text{since }  \gamma_{H}\leq 1\leq \Upsilon_1),
\end{align*}
and
\begin{align*}
\nbr{\begin{pmatrix}
g_a\\
\blambda_c
\end{pmatrix}} & \leq \|g_a\| + \|\blambda_c\| \stackrel{\substack{\substack{\eqref{equ:SQP:direction:1}\\ \eqref{pequ:12}} }}{\leq} \Upsilon_2\|\Delta\bx\| + \frac{\Upsilon_2\Upsilon_{B} + \Upsilon_4 +1}{\gamma_{H}^2\gamma_{B}}\nbr{\left(\begin{smallmatrix}
\Delta\bx\\
J\nabla_{\bx}\mL\\
G\nabla_{\bx}\mL+ \Pi_c(\diag^2(g)\blambda)
\end{smallmatrix}\right) }\\
& \leq \frac{\Upsilon_2(\Upsilon_{B} +1)+ \Upsilon_4 +1}{\gamma_{H}^2\gamma_{B}}\nbr{\left(\begin{smallmatrix}
\Delta\bx\\
J\nabla_{\bx}\mL\\
G\nabla_{\bx}\mL+ \Pi_c(\diag^2(g)\blambda)
\end{smallmatrix}\right)} \quad\quad (\text{since }  \gamma_{H}\vee\gamma_{B}\leq 1),
\end{align*}
we define $\Upsilon_8 = 6\sqrt{2}\Upsilon_7\Upsilon_1(\Upsilon_2(\Upsilon_{B}+1) + \Upsilon_4 +1)$ and have
\begin{equation}\label{npequ:1}
(\nabla\mL_{\epsilon, \nu, \eta}^{(2)})^T\Delta  \leq 
\frac{\Upsilon_8}{\gamma_{H}^3\gamma_{B}}\rbr{\frac{1\vee \nu}{\epsilon (1\wedge \nu^2)}\vee \eta}(\|g_a\| + \|\blambda_c\|) \nbr{ \left(\begin{smallmatrix}
\Delta\bx\\
J\nabla_{\bx}\mL\\
G\nabla_{\bx}\mL+ \Pi_c(\diag^2(g)\blambda)
\end{smallmatrix}\right) }^2.
\end{equation}
By Lemma \ref{lem:5}, we can find a compact subset of $\mX_{\epsilon,\nu}\times \Lambda_{\epsilon,\nu}$ depending only on $(\epsilon,\nu)$ such that~$\mA_{\epsilon, \nu} \subseteq \I(\tx)$ and $\mA_{\epsilon, \nu}^c\subseteq \{\I^+(\tx, \tlambda)\}^c$; thus
\begin{equation*}
\|g_a\| \leq \|g_{\I(\tx)}\|\quad\quad \text{ and }\quad\quad \|\blambda_c\| \leq  \|\blambda_{(\I^+(\tx, \tlambda))^c}\|. 
\end{equation*}
Furthermore, we let $\mX_{\epsilon,\nu,\eta}\times \Lambda_{\epsilon,\nu,\eta} \subseteq \mX_{\epsilon,\nu}\times \Lambda_{\epsilon,\nu}$ be a compact subset depending additionally on $\eta$, such that 
\begin{align*}
\|g_{\I(\tx)}\| & \leq \frac{\gamma_{H}^3\gamma_{B}}{\Upsilon_8} \rbr{\frac{\epsilon (1\wedge\nu^2)}{1\vee \nu}\wedge \frac{1}{\eta}}\frac{\gamma_{B}\wedge \eta}{8},\\
\|\blambda_{(\I^+(\tx, \tlambda))^c}\| & \leq \frac{\gamma_{H}^3\gamma_{B}}{\Upsilon_8} \rbr{\frac{\epsilon (1\wedge\nu^2)}{1\vee \nu}\wedge \frac{1}{\eta}}\frac{\gamma_{B}\wedge \eta}{8}.
\end{align*}
Then, combining \eqref{npequ:1} with the above two displays leads to
\begin{equation*}
(\nabla\mL_{\epsilon, \nu, \eta}^{(2)})^T\Delta \leq \frac{\gamma_{B}\wedge \eta}{4}\nbr{ \left(\begin{smallmatrix}
\Delta\bx\\
J\nabla_{\bx}\mL\\
G\nabla_{\bx}\mL+ \Pi_c(\diag^2(g)\blambda)
\end{smallmatrix}\right) }^2.
\end{equation*}
This completes the proof.

\subsection{Proof of Lemma \ref{lem:5}}\label{pf:lem:5}

Let $\mX\times \Lambda\subseteq \mT_{\nu}\times \mR^r$ be any compact set around $(\tx,\tlambda)$. For any $(\bx, \blambda)\in \mX\times \Lambda$, we have
\begin{equation}\label{pequ:8}
q_{\nu}(\bx, \blambda) \stackrel{\eqref{equ:q}}{\geq} \frac{\nu}{2}\cdot \frac{1}{1+\max_{ \blambda\in \Lambda}\|\blambda\|^2}\eqqcolon \kappa_\nu.
\end{equation}
For any $i \in \I^+(\tx, \tlambda)$, we know $g_i^\star = 0$ and $\tlambda_i >0$. Thus, $g_i^\star + \epsilon \kappa_\nu \tlambda_i >0$. Consider~the~ball $\mB_i^{\bx} = \{\bx: \|\bx - \tx\| \leq r_i\}\cap\mX$ and $\mB_i^{\blambda} = \{\blambda: \|\blambda - \tlambda\|\leq r_i\}\cap\Lambda$. For a sufficiently small $r_i$ (depending on $\epsilon$ and $\nu$), we have $(\tx,\tlambda)\in\mB_i^{\bx}\times \mB_i^{\blambda} \subseteq\mX\times  \Lambda$ and, for any $(\bx, \blambda)\in \mB_i^{\bx}\times \mB_i^{\blambda}$,
\begin{equation*}
g_i(\bx) \geq -\epsilon\kappa_\nu\blambda_i \stackrel{\eqref{pequ:8}}{\geq} -\epsilon q_{\nu}(\bx, \blambda) \blambda_i.
\end{equation*}
The first inequality is due to the continuity of $g_i$. This implies $i \in \mA_{\epsilon, \nu}(\bx, \blambda)$. Therefore,~for~any $(\bx, \blambda)$ in the compact set  $\cap_{i\in \I^+(\tx, \tlambda)}\mB_i^{\bx}\times \mB_i^{\blambda}$, we have $\I^+(\tx, \tlambda)\subseteq\mA_{\epsilon, \nu}(\bx, \blambda)$. The argument $\mA_{\epsilon, \nu}(\bx, \blambda)\subseteq \I(\tx)$ can be proved in the same way.

\subsection{Proof of Lemma \ref{lem:aux:1}}\label{pf:lem:aux:1}

By Assumption \ref{ass:1}, there exists a compact set $X\ni \tx$ small enough such that $(J^T(\bx)\;G^T_{\I(\tx)}(\bx))$ has full column rank for all $\bx\in X$. Furthermore, for any $(\ba, \bb)\in\mR^{m+r}$, we note that
\begin{multline}\label{pequ:6}
0 = (\ba^T\;\; \bb^T)M(\bx)\begin{pmatrix}
\ba\\
\bb
\end{pmatrix} \Longrightarrow \bb_{\I^c(\tx)} = \0\\ \Longrightarrow \|J^T(\bx)\ba + G^T_{\I(\tx)}(\bx)\bb_{\I(\tx)}\| = 0\Longrightarrow (\ba, \bb) = \0,
\end{multline}
where the first implication is due to $\diag(g(\bx))\bb = 0$ and $\I^c(\tx) \subseteq \I^c(\bx)$ (since $X$ is small), and the second implication is due to $\|J^T(\bx)\ba + G^T(\bx)\bb\| = 0$. Therefore, $M(\bx)$ is invertible. Moreover, for any $\mA \subseteq \I(\tx)$, we have
\begin{equation}\label{pequ:7}
\sigma_{\min}\cbr{\left(\begin{smallmatrix}
J(\bx)\\
G_\mA(\bx)
\end{smallmatrix}\right)\left(\begin{smallmatrix}
J^T(\bx) & G^T_{\mA}(\bx)
\end{smallmatrix}\right) } \geq \sigma_{\min} \cbr{\left(\begin{smallmatrix}
J(\bx)\\
G_{\I(\tx)}(\bx)
\end{smallmatrix}\right)\left(\begin{smallmatrix}
J^T(\bx) & G^T_{\I(\tx)}(\bx)
\end{smallmatrix} \right) } >0,
\end{equation}
where $\sigma_{\min}(\cdot)$ denotes the least singular value of a matrix. By \eqref{pequ:6}, \eqref{pequ:7}, and the~compactness of $X$, we know that there exists $\gamma_{H}\in (0, 1]$ such that
\begin{equation}\label{pequ:9}
M(\bx) \succeq \gamma_{H} I, \quad \begin{pmatrix}
J(\bx)\\
G_\mA(\bx)
\end{pmatrix}\begin{pmatrix}
J^T(\bx) & G^T_{\mA}(\bx)
\end{pmatrix}\succeq \gamma_{H} I, \quad \forall \bx\in X \text{ and } \mA \subseteq\I(\tx).
\end{equation}
To show the second part of the statement, we apply Lemma \ref{lem:5}, and know that there exists a compact set $\mX_{\epsilon,\nu}\times \Lambda_{\epsilon,\nu}\subseteq X\times \mR^r$ such that $\mA(\bx,\blambda)\subseteq \I(\tx)$, $\forall (\bx, \blambda)\in \mX_{\epsilon,\nu}\times~\Lambda_{\epsilon,\nu}$.~Combining this fact with \eqref{pequ:9}, we complete the proof.

\section{Proofs of Section \ref{sec:4}}

\subsection{Proof of Lemma \ref{lem:7}}\label{pf:lem:7}

It suffices to show that there exists a threshold $\tilde{\epsilon}>0$ such that for any samples $\xi_1$, any~parameter $\nu\in[\barnu_0, \tilde{\nu}]$, where $\barnu_0$ is the fixed initial input of Algorithm \ref{alg:ASto} and $\tnu$ is defined in \eqref{bound:nu}, and any point $(\bx, \bmu, \blambda) \in \mX\times\mM\times \Lambda$ with $\bx \in \mT_{\nu}$, if $\epsilon \leq \tilde{\epsilon}$, then
\begin{equation*}
\nbr{\begin{pmatrix}
c(\bx), \bw_{\epsilon, \nu}(\bx, \blambda)
\end{pmatrix}} \leq \chi_{err} \cdot \nbr{\bnabla \mL_{\epsilon, \nu, \eta}(\bx, \bmu, \blambda)},
\end{equation*}
where $\bnabla\mL_{\epsilon, \nu, \eta}$ is computed using samples in $\xi_1$ and $\eta, \chi_{err}>0$ are any given positive~constants. Note that everything above is deterministic; that is, our analysis does not depend on a~specific iteration sequence $\{(\bx_t, \bmu_t, \blambda_t)\}_t$. Thus, the threshold $\tilde{\epsilon}$ is deterministic. Let~us~prove the above statement by contradiction. Without loss of generality, we suppose $\chi_{err}\leq 1$.

Suppose the statement is false, then there exist a sequence $\{\epsilon_j, \xi_1^j, \nu_j\}_j$ and an evaluation point sequence $\{(\bx_j, \bmu_j, \blambda_j)\}_j\in \mX\times \mM\times \Lambda$ such that $\nu_j\in[\barnu_0, \tnu]$, $\bx_j \in \mT_{\nu_j}$, $\epsilon_j\searrow 0$ and 
\begin{equation}\label{pequ:20}
\|\bnabla \mL_{\epsilon_j, \nu_j, \eta}^j \| < 1/\chi_{err}\cdot \|(c_j, \bw_{\epsilon_j, \nu_j}^j )\|, \quad\quad \forall j\geq 0,
\end{equation}
where $\bnabla \mL_{\epsilon_j, \nu_j, \eta}^j$ is computed using samples $\xi_1^j$, and $\eta$ and $\chi_{err}$ are fixed constants. By the compactness condition, we suppose $(\bx_j, \bmu_j, \blambda_j) \rightarrow (\tilde{\bx}, \tilde{\bmu}, \tilde{\blambda}) \in \mX\times \mM\times \Lambda$ and $\nu_j \rightarrow \nu$ as $j \rightarrow \infty$ (otherwise, we can consider a convergent subsequence, which must exist). Noting that $c_j = c(\bx_j)$ and $\bw_{\epsilon_j, \nu_j}^j = \max\{g(\bx_j), -\epsilon_jq_{\nu_j}(\bx_j, \blambda_j)\blambda_j\}$ are bounded due to the compactness of $(\bx_j, \bmu_j, \blambda_j)$ and the boundedness of $\nu_j$ and $\epsilon_j$, we have from \eqref{pequ:20} that
\begin{equation}\label{pequ:21}
\epsilon_j \|\bnabla_{\bx}\mL_{\epsilon_j, \nu_j, \eta}^j\| \rightarrow 0\quad \text{ as } \quad j\rightarrow \infty.
\end{equation}
Moreover, since $\bx_j \in \mT_{\nu_j}$, we have $\sum_{i=1}^r\max\{(g_j)_i, 0\}^3\leq \nu_j/2$. Taking limit $j\rightarrow\infty$ leads to $\tilde{\bx}\in \mT_{\nu}$. Furthermore, by \eqref{equ:aug:der}, \eqref{pequ:21}, and the convergence of $(\bx_j, \bmu_j, \blambda_j)$, we get
\begin{equation*}
J^T(\tilde{\bx})c(\tilde{\bx}) + \frac{1}{q_{\nu}(\tilde{\bx}, \tilde{\blambda})}G^T(\tilde{\bx})\max\{g(\tilde{\bx}), \0\} +  \frac{3\|\max\{g(\tilde{\bx}), \0\}\|^2}{2q_{\nu}(\tilde{\bx}, \tilde{\blambda})a_{\nu}(\tilde{\bx})}G(\tilde{\bx})^T\bl(\tilde{\bx}) = \0,
\end{equation*}
which is further simplified as
\begin{equation}\label{pequ:22}
\sum_{i: c_i(\tilde{\bx})\neq 0}c_i(\tilde{\bx})\nabla c_i(\tilde{\bx}) + \sum_{i: g_i(\tilde{\bx})> 0}\cbr{\frac{1}{q_{\nu}(\tilde{\bx}, \tilde{\blambda})} + \frac{3\|\max\{g(\tilde{\bx}), \0\}\|^2g_i(\tilde{\bx})}{2q_{\nu}(\tilde{\bx}, \tilde{\blambda})a_{\nu}(\tilde{\bx}) } }g_i(\tilde{\bx})\nabla g_i(\tilde{\bx}) = \0.
\end{equation}

Suppose $\tilde{\bx} \in \mX\backslash \Omega$ and let $\I_c(\tilde{\bx}) = \{i: 1\leq i\leq m, c_i(\tilde{\bx}) \neq0\}$, and $\I_g(\tilde{\bx}) = \{i: 1\leq i\leq r, g_i(\tilde{\bx}) >0\}$. By Assumption \ref{ass:4}, the set
\begin{equation*}
\cbr{\bz\in \mR^d: c_i(\tilde{\bx})\nabla^Tc_i(\tilde{\bx})\bz < 0, i \in \I_c(\tilde{\bx}) \text{ and } \nabla^Tg_i(\tilde{\bx})\bz < 0, i \in \I_g(\tilde{\bx})}
\end{equation*}	
is nonempty. By the Gordan's theorem \citep{Goldman19574.}, for any $a_i, b_i\geq 0$ such that
\begin{equation}\label{pequ:23}
\sum_{i \in \I_c(\tilde{\bx})} a_i c_i(\tilde{\bx})\nabla c_i(\tilde{\bx}) + \sum_{i \in \I_g(\tilde{\bx})}b_i\nabla g_i(\tilde{\bx}) = \0,
\end{equation}
we have $a_i = b_i = 0$. Comparing \eqref{pequ:23} with \eqref{pequ:22}, and noting that the coefficients of \eqref{pequ:22} are all positive (since $\tilde{\bx}\in \mT_{\nu}$), we immediately get the contradiction. Thus, $\tilde{\bx}\in \Omega$.

By Assumption \ref{ass:4} and following the same reasoning as \eqref{pequ:6}, $M(\tilde{\bx})$ is invertible and, particularly, is positive definite. Thus, $M_j$ is invertible for large enough $j$. Let us suppose $\|M_j^{-1}\|\leq \Upsilon_M$ for some $\Upsilon_M>0$. Further, by direct calculation, we have
\begin{equation}\label{pequ:24}
\diag(g_j)\blambda_j = \diag(\blambda_j)\bw_{\epsilon_j, \nu_j}^j - \frac{1}{\epsilon_jq_{\nu_j}^j}(\diag(g_j)  - \diag(\bw_{\epsilon_j, \nu_j}^j))\bw_{\epsilon_j, \nu_j}^j.
\end{equation}
Thus, we can obtain
\begin{align}\label{pequ:25}
\begin{pmatrix}
J_j\\
G_j
\end{pmatrix}\bnabla_{\bx}\mL_{\epsilon_j, \nu_j, \eta}^j \stackrel{\mathclap{\eqref{equ:aug:der}}}{=} &\;\; \begin{pmatrix}
J_j\\
G_j
\end{pmatrix}\bnabla_{\bx}\mL_j + \eta\begin{pmatrix}
J_j\\
G_j
\end{pmatrix}\begin{pmatrix}
Q_{1,j} & Q_{2,j}
\end{pmatrix}\left(\begin{smallmatrix}
J_j\bnabla_{\bx}\mL_j\\
G_j\bnabla_{\bx}\mL_j + \diag^2(g_j)\blambda_j
\end{smallmatrix}\right) \nonumber\\
& \quad + \frac{1}{\epsilon_j}\begin{pmatrix}
J_j\\
G_j
\end{pmatrix}\begin{pmatrix}
J_j^T & \frac{G_j^T}{q_{\nu_j}^j} + \frac{3G_j^T\bl_j(\bw_{\epsilon_j, \nu_j}^j)^T}{2q_{\nu_j}^ja_{\nu_j}^j}
\end{pmatrix}\begin{pmatrix}
c_j\\
\bw_{\epsilon_j, \nu_j}^j
\end{pmatrix} \nonumber\\
\stackrel{\mathclap{\eqref{pequ:24}}}{=} & \;\;\cbr{I + \eta\begin{pmatrix}
J_j\\
G_j
\end{pmatrix}\begin{pmatrix}
Q_{1,j} & Q_{2,j}
\end{pmatrix}} \left(\begin{smallmatrix}
J_j\bnabla_{\bx}\mL_j\\
G_j\bnabla_{\bx}\mL_j + \diag^2(g_j)\blambda_j
\end{smallmatrix}\right) \nonumber\\
& \quad + \frac{1}{\epsilon_j}\bigg\{\begin{pmatrix}
J_j\\
G_j
\end{pmatrix}\begin{pmatrix}
J_j^T & \frac{G_j^T}{q_{\nu_j}^j} + \frac{3G_j^T\bl_j(\bw_{\epsilon_j, \nu_j}^j)^T}{2q_{\nu_j}^ja_{\nu_j}^j}
\end{pmatrix} \nonumber\\
&\quad + \left(\begin{smallmatrix}
\0 & \0\\
\0 & \frac{\diag^2(g_j) - \diag(g_j)\diag(\bw_{\epsilon_j, \nu_j}^j)}{q_{\nu_j}^j} - \epsilon_j\diag(g_j)\diag(\blambda_j)
\end{smallmatrix}\right) \bigg\}\begin{pmatrix}
c_j\\
\bw_{\epsilon_j, \nu_j}^j
\end{pmatrix} \nonumber\\
\eqqcolon & \H_{1,j} \left(\begin{smallmatrix}
J_j\bnabla_{\bx}\mL_j\\
G_j\bnabla_{\bx}\mL_j + \diag^2(g_j)\blambda_j
\end{smallmatrix} \right)+ \frac{1}{\epsilon_j}\H_{2,j}\begin{pmatrix}
c_j\\
\bw_{\epsilon_j, \nu_j}^j
\end{pmatrix}.
\end{align}
Let us focus on $\H_{2,j}$. We know that
\begin{align*}
\H_{2,j} & =  \left(\begin{smallmatrix}
J_jJ_j^T &  J_jG_j^T/q_{\nu_j}^j\\
G_jJ_j^T & \cbr{G_jG_j^T + \diag^2(g_j)}/q_{\nu_j}^j
\end{smallmatrix}\right) + \underbrace{ \left(\begin{smallmatrix}
\0 & \frac{3}{2q_{\nu_j}^ja_{\nu_j}^j}J_jG_j^T\bl_j(\bw_{\epsilon_j, \nu_j}^j)^T\\
\0 & \substack{\frac{3}{2q_{\nu_j}^ja_{\nu_j}^j}G_jG_j^T\bl_j(\bw_{\epsilon_j, \nu_j}^j)^T \\ - \frac{\diag(g_j)\diag(\bw_{\epsilon_j, \nu_j}^j)}{q_{\nu_j}^j}  - \epsilon_j\diag(g_j)\diag(\blambda_j)}
\end{smallmatrix} \right)}_{\Delta\H_{2,j}}\\
& =  M_j\begin{pmatrix}
I & \0\\
\0 & \frac{1}{q_{\nu_j}^j}I
\end{pmatrix} + \Delta\H_{2,j}.
\end{align*}
Recalling that $\sigma_{\min}(\cdot)$ denotes the least singular value of a matrix, by the Weyl's~inequality,
\begin{equation*}
\sigma_{\min}(\H_{2,j}) \geq \sigma_{\min}\cbr{M_j\begin{pmatrix}
I & \0\\
\0 & \frac{1}{q_{\nu_j}^j}I
\end{pmatrix}} - \|\Delta\H_{2,j}\|
\geq  \frac{\sigma_{\min}(M_j)}{1\vee q_{\nu_j}^j} -\|\Delta\H_{2,j}\|.
\end{equation*} 
Since $\epsilon_j\rightarrow 0$ and $\bw_{\epsilon_j, \nu_j}^j \rightarrow 0$ as $j\rightarrow \infty$ (because $\tilde{\bx}\in \Omega$), we know $\Delta\H_{2,j}\rightarrow \0$. In addition, since $M_j\rightarrow M(\tilde{\bx})$ with $M(\tilde{\bx})$ being positive definite, and $q_{\nu_j}^j \leq \nu_j = \tnu$, we know for some constant $\varphi>0$ and sufficiently large $j$,
\begin{equation}\label{pequ:26}
\sigma_{\min}(\H_{2,j}) \geq \varphi.
\end{equation}
Now we bound the first term in \eqref{pequ:25}. By \eqref{equ:aug:der} and the invertibility of $M_j$, we know
\begin{align}\label{pequ:27}
\hskip-0.5cm \nbr{\left(\begin{smallmatrix}
J_j\bnabla_{\bx}\mL_j\\
G_j\bnabla_{\bx}\mL_j + \diag^2(g_j)\blambda_j
\end{smallmatrix} \right)} \;\stackrel{\mathclap{\eqref{equ:aug:der}}}{=} & \;\; \frac{1}{\eta}\nbr{M_j^{-1}\cbr{\left(\begin{smallmatrix}
\bnabla_{\bmu}\mL_{\epsilon_j, \nu_j, \eta}^j\\
\bnabla_{\blambda}\mL_{\epsilon_j, \nu_j, \eta}^j
\end{smallmatrix}\right) - \left(\begin{smallmatrix}
c_j\\
\bw_{\epsilon_j, \nu_j}^j + \frac{\|\bw_{\epsilon_j, \nu_j}^j\|^2}{\epsilon_ja_{\nu_j}^j}\blambda_j
\end{smallmatrix} \right)} } \nonumber\\
\leq & \frac{\Upsilon_M}{\eta}\cbr{
\left(\begin{smallmatrix}
\bnabla_{\bmu}\mL_{\epsilon_j, \nu_j, \eta}^j\\
\bnabla_{\blambda}\mL_{\epsilon_j, \nu_j, \eta}^j
\end{smallmatrix}\right) + \nbr{\left(\begin{smallmatrix}
c_j\\
\bw_{\epsilon_j, \nu_j}^j
\end{smallmatrix} \right)}  + \frac{\|\bw_{\epsilon_j, \nu_j}^j\|^2\|\blambda_j\|}{\epsilon_ja_{\nu_j}^j}
} \nonumber\\
\stackrel{\mathclap{\eqref{pequ:20}}}{\leq} &\;\; \frac{\Upsilon_M}{\eta}\cbr{\rbr{1 + \frac{1}{\chi_{err}}}\nbr{\left(\begin{smallmatrix}
c_j\\
\bw_{\epsilon_j, \nu_j}^j
\end{smallmatrix} \right)} + \frac{\|\bw_{\epsilon_j, \nu_j}^j\|^2\|\blambda_j\|}{\epsilon_ja_{\nu_j}^j} }\nonumber\\
\stackrel{\mathclap{\eqref{equ:q}}}{\leq} &\;\; \frac{2\Upsilon_M}{\chi_{err}\eta}\cbr{\nbr{\left(\begin{smallmatrix}
c_j\\
\bw_{\epsilon_j, \nu_j}^j
\end{smallmatrix} \right)} + \frac{\|\bw_{\epsilon_j, \nu_j}^j\|^2\|\blambda_j\|}{\epsilon_j\nu_j} } \quad (\text{also use } \chi_{err}\leq 1)\nonumber\\
\leq & \frac{2\Upsilon_M}{\chi_{err}\eta\epsilon_j}\cbr{ \epsilon_j+ \frac{\|\bw_{\epsilon_j, \nu_j}^j\|\|\blambda_j\|}{\nu_j}}\nbr{\left(\begin{smallmatrix}
c_j\\
\bw_{\epsilon_j, \nu_j}^j
\end{smallmatrix} \right)}.
\end{align}
Moreover, by the compactness condition, we have $\|\H_{1,j}\|\leq \Upsilon_1$ and $\|(J_j^T\; G_j^T)\|\leq \Upsilon_2$ for some constants $\Upsilon_1, \Upsilon_2>0$. Combining \eqref{pequ:26}, \eqref{pequ:27} with \eqref{pequ:25}, we have
\begin{align*}
\epsilon_j\Upsilon_2\nbr{\bnabla_{\bx}\mL_{\epsilon_j, \nu_j, \eta}^j} \geq & \epsilon_j\nbr{\begin{pmatrix}
J_j\\
G_j
\end{pmatrix}\bnabla_{\bx}\mL_{\epsilon_j, \nu_j, \eta}^j} \\
\stackrel{\mathclap{\eqref{pequ:25}}}{\geq} &\;\; \nbr{\H_{2,j}\left(\begin{smallmatrix}
c_j\\
\bw_{\epsilon_j, \nu_j}^j
\end{smallmatrix} \right)} - \epsilon_j\nbr{\H_{1,j}\left(\begin{smallmatrix}
J_j\bnabla_{\bx}\mL_j\\
G_j\bnabla_{\bx}\mL_j + \diag^2(g_j)\blambda_j
\end{smallmatrix}\right) }\\
\stackrel{\mathclap{\eqref{pequ:26}}}{\geq} &\;\; \varphi\cdot \nbr{\left(\begin{smallmatrix}
c_j\\
\bw_{\epsilon_j, \nu_j}^j
\end{smallmatrix} \right)} - \epsilon_j\Upsilon_1\nbr{ \left(\begin{smallmatrix}
J_j\bnabla_{\bx}\mL_j\\
G_j\bnabla_{\bx}\mL_j + \diag^2(g_j)\blambda_j
\end{smallmatrix}\right) }\\
\stackrel{\mathclap{\eqref{pequ:27}}}{\geq} &\;\; \cbr{\varphi - \frac{2\Upsilon_1\Upsilon_M}{\chi_{err}\eta}\rbr{\epsilon_j + \frac{\|\bw_{\epsilon_j, \nu_j}^j\|\|\blambda_j\|}{\nu_j}} }\nbr{\begin{pmatrix}
c_j\\
\bw_{\epsilon_j, \nu_j}^j
\end{pmatrix}}\\ 
\eqqcolon & (\varphi - \varphi_j)\|(c_j, \bw_{\epsilon_j, \nu_j}^j)\|.
\end{align*}
Noting that $\varphi_j \rightarrow 0$ as $j\rightarrow\infty$ (since $\bw_{\epsilon_j, \nu_j}^j\rightarrow 0$ and $\epsilon_j\rightarrow 0$), we obtain for large $j$ that
\begin{equation*}
\epsilon_j\Upsilon_2/\chi_{err}\cdot \|(c_j, \bw_{\epsilon_j, \nu_j}^j)\| \stackrel{\eqref{pequ:20}}{\geq} \epsilon_j\Upsilon_2\|\bnabla_{\bx}\mL_{\epsilon_j, \nu_j, \eta}^j\| \geq \varphi/2\cdot \|(c_j, \bw_{\epsilon_j, \nu_j}^j)\|,
\end{equation*}
which cannot hold because $\epsilon_j\searrow 0$. This is a contradiction, and thus we complete the proof.

\subsection{Proof of Lemma \ref{lem:8}}\label{pf:lem:8}

The proof closely follows the proof of Lemma \ref{lem:6} in Appendix \ref{pf:lem:6}. We suppress the iteration $t$ and assume $\xi_1^t$ is any sample set. Our analysis is independent of the sample set $\xi_1^t$ for computing $\bnabla\mL_{\barepsilon_t, \barnu_t, \eta}^t$, and we will see that the threshold is independent of $t$. Like Lemma~\ref{lem:6}, we use $\Upsilon_1, \Upsilon_2, \ldots$ to denote generic constants that are independent of $(\barepsilon_t, \barnu_t, \eta, \gamma_{B}, \gamma_{H})$,~whose existence is ensured by the compactness of the iterates.

Following the derivation of \eqref{pequ:10}, we have
\begin{multline}\label{pequ:28}
(\bnabla\mL_{\barepsilon, \barnu, \eta}^{(1)})^T\barDelta =  - \barDelta\bx^TB\barDelta\bx + \left(\begin{smallmatrix}
c\\
g_a
\end{smallmatrix}\right)^T\left(\begin{smallmatrix}
\bar{\tDelta}\bmu + \barDelta\bmu\\
\bar{\tDelta}\blambda_a + \barDelta\blambda_a
\end{smallmatrix}\right) - \frac{1}{\barepsilon}\|c\|^2 - \frac{1}{\barepsilon q_{\barnu}}\|g_a\|^2  \\ - \barepsilon q_{\barnu}\barDelta\blambda_c^T\blambda_c 
- \eta \nbr{\left(\begin{smallmatrix}
J\bnabla_{\bx}\mL\\
G\bnabla_{\bx}\mL+ \Pi_c(\diag^2(g)\blambda)	\end{smallmatrix} \right)}^2,
\end{multline}
where $(\bar{\tDelta}\bmu, \bar{\tDelta}\blambda_a)$ is the dual solution of \eqref{equ:SQP:direction:1} with $\nabla_{\bx}\mL$ being replaced by $\bnabla_{\bx}\mL$. Following the derivation of \eqref{pequ:11}, there exists $\Upsilon_1>0$ such that
\begin{equation}\label{pequ:29}
\nbr{\begin{pmatrix}
\barDelta\bmu\\
\barDelta\blambda
\end{pmatrix}} \leq \frac{\Upsilon_1}{\gamma_{H}} \nbr{\left(\begin{smallmatrix}
\barDelta\bx\\
J\bnabla_{\bx}\mL\\
G\bnabla_{\bx}\mL+ \Pi_c(\diag^2(g)\blambda)
\end{smallmatrix} \right)}.
\end{equation}
Following the derivation of \eqref{pequ:12}, there exists $\Upsilon_2>0$ such that
\begin{equation}\label{pequ:30}
\nbr{\left(\begin{smallmatrix}
\bar{\tDelta}\bmu\\
\bar{\tDelta}\blambda_a\\
-\blambda_c
\end{smallmatrix} \right)}
\leq \frac{\Upsilon_2}{\gamma_{H}^2\gamma_{B}}\nbr{ \left(\begin{smallmatrix}
\barDelta\bx\\
J\bnabla_{\bx}\mL\\
G\bnabla_{\bx}\mL+ \Pi_c(\diag^2(g)\blambda)
\end{smallmatrix} \right)}.
\end{equation}
Following the derivation of \eqref{pequ:13} by combining \eqref{pequ:28}, \eqref{pequ:29}, and \eqref{pequ:30}, and noting that $0<q_{\barnu}\leq \barnu \leq \tilde{\nu}$ where $\tnu$ is defined in \eqref{bound:nu}, there exists $\Upsilon_3>0$ such that
\begin{align}\label{npequ:2}
(\bnabla\mL_{\barepsilon, \barnu, \eta}^{(1)})^T\barDelta \leq & -\barDelta\bx^TB\barDelta\bx + \frac{\Upsilon_3}{\gamma_{H}^2\gamma_{B}}\nbr{\left(\begin{smallmatrix}
c\\
g_a
\end{smallmatrix} \right)}\nbr{\left(\begin{smallmatrix}
\barDelta\bx\\
J\bnabla_{\bx}\mL\\
G\bnabla_{\bx}\mL+ \Pi_c(\diag^2(g)\blambda)
\end{smallmatrix} \right)} - \frac{1}{\barepsilon(1\vee \tilde{\nu})}\nbr{\left(\begin{smallmatrix}
c\\
g_a
\end{smallmatrix} \right)}^2 \nonumber\\
&   + \frac{\barepsilon\tnu\Upsilon_3}{\gamma_{H}^3\gamma_{B}} \nbr{\left(\begin{smallmatrix}
\barDelta\bx\\
J\bnabla_{\bx}\mL\\
G\bnabla_{\bx}\mL+ \Pi_c(\diag^2(g)\blambda)
\end{smallmatrix} \right)}^2 - \eta\nbr{\left(\begin{smallmatrix}
J\bnabla_{\bx}\mL\\
G\bnabla_{\bx}\mL+ \Pi_c(\diag^2(g)\blambda)	\end{smallmatrix} \right)}^2.
\end{align}
Following the derivation of \eqref{pequ:31}, there exists $\Upsilon_4>0$ such that
\begin{equation*}
-\barDelta\bx^TB\barDelta\bx \leq  -\frac{3\gamma_{B}}{4}\|\barDelta\bx\|^2 + \frac{\Upsilon_4}{\gamma_{H}\gamma_{B}}\nbr{\left(\begin{smallmatrix}
c\\
g_a
\end{smallmatrix} \right)}^2.
\end{equation*}
Combining the above display with \eqref{npequ:2} and using the following Young's inequality
\begin{multline*}
\frac{\Upsilon_3}{\gamma_{H}^2\gamma_{B}}\nbr{\left(\begin{smallmatrix}
c\\
g_a
\end{smallmatrix} \right)}\nbr{\left(\begin{smallmatrix}
\barDelta\bx\\
J\bnabla_{\bx}\mL\\
G\bnabla_{\bx}\mL+ \Pi_c(\diag^2(g)\blambda)
\end{smallmatrix} \right)} \\
\leq \rbr{\frac{\gamma_{B}}{8}\wedge \frac{\eta}{4}}\nbr{\left(\begin{smallmatrix}
\barDelta\bx\\
J\bnabla_{\bx}\mL\\
G\bnabla_{\bx}\mL+ \Pi_c(\diag^2(g)\blambda)
\end{smallmatrix} \right)}^2 + \frac{2\Upsilon_3^2}{\gamma_{H}^4\gamma_{B}^2(\gamma_{B} \wedge \eta)}\nbr{\left(\begin{smallmatrix}
c\\
g_a
\end{smallmatrix} \right)}^2,
\end{multline*}
we have
\begin{align*}
&(\bnabla\mL_{\epsilon, \nu, \eta}^{(1)})^T\barDelta \leq  -\frac{3\gamma_{B}}{4}\|\barDelta\bx\|^2  + \cbr{\rbr{\frac{\gamma_{B}}{8}\wedge \frac{\eta}{4}} + \frac{\barepsilon\tnu \Upsilon_3}{\gamma_{H}^3\gamma_{B}}}\nbr{\left(\begin{smallmatrix}
\barDelta\bx\\
J\bnabla_{\bx}\mL\\
G\bnabla_{\bx}\mL+ \Pi_c(\diag^2(g)\blambda)
\end{smallmatrix} \right)}^2\\
& \quad + \cbr{\frac{\Upsilon_4}{\gamma_{H}\gamma_{B}} + \frac{2\Upsilon_3^2}{\gamma_{H}^4\gamma_{B}^2(\gamma_{B} \wedge \eta)} - \frac{1}{\barepsilon(1\vee \tnu)}}\nbr{\left(\begin{smallmatrix}
c\\
g_a
\end{smallmatrix} \right)}^2 -\eta\nbr{\left(\begin{smallmatrix}
J\bnabla_{\bx}\mL\\
G\bnabla_{\bx}\mL+ \Pi_c(\diag^2(g)\blambda)	\end{smallmatrix} \right)}^2\\
&\leq  -\cbr{\frac{\gamma_{B}\wedge \eta}{2} + \rbr{\frac{\gamma_{B}}{8} \wedge \frac{\eta}{4}} - \frac{\barepsilon\tnu \Upsilon_3}{\gamma_{H}^3\gamma_{B}} } \nbr{ \left(\begin{smallmatrix}
\barDelta\bx\\
J\bnabla_{\bx}\mL\\
G\bnabla_{\bx}\mL+ \Pi_c(\diag^2(g)\blambda)
\end{smallmatrix} \right)}^2\\
& \quad - \cbr{\frac{1}{\barepsilon(1\vee \tnu)}  - \frac{\Upsilon_4}{\gamma_{H}\gamma_{B}}- \frac{2\Upsilon_3^2}{\gamma_{H}^4\gamma_{B}^2(\gamma_{B} \wedge \eta)} }\nbr{\left(\begin{smallmatrix}
c\\
g_a
\end{smallmatrix} \right)}^2.
\end{align*}
Therefore, as long as
\begin{equation}\label{pequ:32}
\begin{aligned}
\frac{\gamma_{B}}{8}\wedge\frac{\eta}{4} \geq \frac{\barepsilon\tnu\Upsilon_3}{\gamma_{H}^3\gamma_{B}} \Longleftarrow& \frac{1}{\barepsilon} \geq \frac{8\tnu\Upsilon_3}{\gamma_{H}^3\gamma_{B}(\gamma_{B}\wedge\eta)}, \\
\frac{1}{\barepsilon(1\vee \tnu)} - \frac{\Upsilon_4}{\gamma_{H}\gamma_{B}} -\frac{2\Upsilon_3^2}{\gamma_{H}^4\gamma_{B}^2(\gamma_{B} \wedge \eta)} \geq 0 \Longleftarrow & \frac{1}{\barepsilon} \geq \frac{(1\vee \tnu )(2\Upsilon_3^2 + \Upsilon_4)}{\gamma_{H}^4\gamma_{B}^2(\gamma_{B} \wedge \eta)},
\end{aligned}	
\end{equation}
we have
\begin{equation*}
(\bnabla\mL_{\epsilon, \nu, \eta}^{(1)})^T\barDelta \leq - \frac{\gamma_{B} \wedge \eta}{2}\nbr{\left(\begin{smallmatrix}
\barDelta\bx\\
J\bnabla_{\bx}\mL\\
G\bnabla_{\bx}\mL+ \Pi_c(\diag^2(g)\blambda)
\end{smallmatrix} \right)}^2.
\end{equation*}
Thus, we can define
\begin{equation*}
\tepsilon_2 \coloneqq \frac{\gamma_{H}^4\gamma_{B}^2(\gamma_{B}\wedge\eta)}{(2\Upsilon_3^2 + 8\Upsilon_3 + \Upsilon_4 )(\tnu\vee 1)},
\end{equation*}
which implies \eqref{pequ:32} and completes the proof.

\subsection{Proof of Lemma \ref{lem:9}}\label{pf:lem:9}

We let $C_1, C_2,\ldots$ be generic constants that are independent of $(\beta, \alpha_{max}, \kappa_{grad}, \kappa_{f}, p_{grad}, p_{f},\\ \chi_{grad},\chi_{f})$. These constants may not be consistent with the constants $C_1, C_2, C_3$ in the statement. However, the existence of $C_1,C_2, C_3$ in the statement follows directly from~our~proof.

\noindent\textbf{(a1).} By the definition of $\nabla\mL_{\epsilon, \nu, \eta}$ in \eqref{equ:aug:der}, all quantities depending on $\epsilon, \nu$ do not depend on the batch samples. We have
\begin{align*}
\bnabla\mL_{\epsilon,\nu,\eta}^t - \nabla\mL_{\epsilon,\nu,\eta}^t & \stackrel{\mathclap{\eqref{equ:aug:der}}}{=} \left(\begin{smallmatrix}
\bnabla_{\bx}\mL_t - \nabla_{\bx}\mL_t\\
\0\\
\0
\end{smallmatrix}\right) + \eta\left(\begin{smallmatrix}
\barQ_{1,t} & \barQ_{2,t}\\
M_{11,t} & M_{12,t}\\
M_{21,t} & M_{22,t}
\end{smallmatrix}\right)\left(\begin{smallmatrix}
J_t\bnabla_{\bx}\mL_t\\
G_t\bnabla_{\bx}\mL_t + \diag^2(g_t)\blambda_t
\end{smallmatrix}\right)\\
& \quad - \eta\left(\begin{smallmatrix}
Q_{1,t} & Q_{2,t}\\
M_{11,t} & M_{12,t}\\
M_{21,t} & M_{22,t}
\end{smallmatrix}\right)\left(\begin{smallmatrix}
J_t\nabla_{\bx}\mL_t\\
G_t\nabla_{\bx}\mL_t + \diag^2(g_t)\blambda_t
\end{smallmatrix}\right)\\
& = \left(\begin{smallmatrix}
\bnabla_{\bx}\mL_t - \nabla_{\bx}\mL_t\\
\0\\
\0
\end{smallmatrix}\right) +  \eta\left(\begin{smallmatrix}
Q_{1,t} & Q_{2,t}\\
M_{11,t} & M_{12,t}\\
M_{21,t} & M_{22,t}
\end{smallmatrix}\right)\left(\begin{smallmatrix}
J_t(\bnabla_{\bx}\mL_t - \nabla_{\bx}\mL_t)\\
G_t(\bnabla_{\bx}\mL_t - \nabla_{\bx}\mL_t)
\end{smallmatrix}\right)\\
&\quad + \eta\left(\begin{smallmatrix}
\barQ_{1,t} - Q_{1,t} & \barQ_{2,t} - Q_{2,t}\\
\0 & \0\\
\0 & \0
\end{smallmatrix}\right)\left(\begin{smallmatrix}
J_t\bnabla_{\bx}\mL_t\\
G_t\bnabla_{\bx}\mL_t + \diag^2(g_t)\blambda_t
\end{smallmatrix}\right). 
\end{align*}
By Assumption \ref{ass:3}, the definition \eqref{equ:def:Qmatrices}, and the facts that $\bnabla_{\bx}\mL_t - \nabla_{\bx}\mL_t = \bnabla f_t - \nabla f_t$~and~$\bnabla_{\bx}^2\mL_t - \nabla_{\bx}^2\mL_t = \bnabla^2 f_t - \nabla^2 f_t$, there exists $C_1>0$ (depending on $\eta$) such that
\begin{equation*}
\|\bnabla\mL_{\epsilon,\nu,\eta}^t - \nabla\mL_{\epsilon,\nu,\eta}^t\| \leq C_1\|\bnabla f_t - \nabla f_t\| + C_1\|\bnabla^2 f_t - \nabla^2 f_t\|(\|\bnabla_{\bx}\mL_t\| + \|\diag^2(g_t)\blambda_t\|).
\end{equation*}
Since $\|\diag^2(g_t)\blambda_t\| \leq C_2\|\max\{g_t, -\blambda_t\}\|$ for some constant $C_2>0$, we apply the definition~of $\barR_t$ in \eqref{event:E1} and the uniform boundedness of $\barR_t$, and know that the above inequality leads to the statement.

\noindent\textbf{(a2).} By the definition of $\mL_{\epsilon,\nu,\eta}$ in \eqref{equ:aug:Lagrange}, all quantities depending on $\epsilon,\nu$ do not depend on the batch samples. We have
\begin{equation*}
\barL_{\epsilon, \nu, \eta}^t - \mL_{\epsilon,\nu,\eta}^t \stackrel{\eqref{equ:aug:Lagrange}}{=} \barL_t - \mL_t + \frac{\eta}{2}\left(\begin{smallmatrix}
J_t(\bnabla_{\bx}\mL_t - \nabla_{\bx}\mL_t)\\
G_t(\bnabla_{\bx}\mL_t - \nabla_{\bx}\mL_t)
\end{smallmatrix}\right)^T\left(\begin{smallmatrix}
J_t(\bnabla_{\bx}\mL_t + \nabla_{\bx}\mL_t)\\
G_t(\bnabla_{\bx}\mL_t + \nabla_{\bx}\mL_t)+2\diag^2(g_t)\blambda_t
\end{smallmatrix}\right).
\end{equation*}
By Assumption \ref{ass:3} and the facts that $\barL_t - \mL_t = \barf_t-f_t$ and $\|\diag^2(g_t)\blambda_t\| \leq C_2\|\max\{g_t, -\blambda_t\}\|$, there exists $C_3>0$ (depending on $\eta$) such that
\begin{equation*}
|\barL_{\epsilon, \nu, \eta}^t - \mL_{\epsilon,\nu,\eta}^t| \leq C_3|\barf_t - f_t| + C_3\|\bnabla f_t - \nabla f_t\|(R_t + \barR_t).
\end{equation*}
Using $R_t\leq \barR_t + \|\bnabla f_t-\nabla f_t\| \leq 2(\barR_t \vee \|\bnabla f_t-\nabla f_t\|)$, we prove the statement.

\noindent\textbf{(b).} By \eqref{equ:aug:der} and Assumption \ref{ass:3}, there exists $C_4>0$ such that
\begin{multline*}
\|\bnabla_{\bx}\mL_t\| \leq  \|\bnabla_{\bx}\mL_{\barepsilon_t, \barnu_t, \eta}^t\| + C_4\cbr{\nbr{\begin{pmatrix}
J_t\bnabla_{\bx}\mL_t\\
G_t\bnabla_{\bx}\mL_t
\end{pmatrix}} + \|\diag^2(g_t)\blambda_t\|} \\
+ \frac{C_4}{\barepsilon_t(1\wedge q_{\barnu_t}^t)}\nbr{\begin{pmatrix}
c_t\\
\bw_{\barepsilon_t, \barnu_t}^t
\end{pmatrix}} + \frac{C_4}{\barepsilon_tq_{\barnu_t}^ta_{\barnu_t}^t}\|\bw_{\barepsilon_t, \barnu_t}^t\|^2.
\end{multline*}
By Theorem \ref{thm:2a}, we have
\begin{equation}\label{pequ:38}
\barepsilon_0\geq \barepsilon_t \geq \tepsilon, \quad \tnu\geq \barnu_t\geq q_{\barnu_t}^t \stackrel{\eqref{pequ:8}}{\geq} \kappa_{\barnu_t}\geq \kappa_{\barnu_0}, \quad \tnu\geq\barnu_t \geq a_{\barnu_t}^t \geq \frac{\barnu_t}{2}\geq \frac{\barnu_0}{2}.
\end{equation}
Thus, there exists $C_5>0$ such that
\begin{equation*}
\|\bnabla_{\bx}\mL_t\| \leq  \|\bnabla_{\bx}\mL_{\barepsilon_t, \barnu_t, \eta}^t\| + C_5\cbr{\nbr{\begin{pmatrix}
J_t\bnabla_{\bx}\mL_t\\
G_t\bnabla_{\bx}\mL_t
\end{pmatrix}} + \|\diag^2(g_t)\blambda_t\| + \nbr{\begin{pmatrix}
c_t\\
\bw_{\barepsilon_t, \barnu_t}^t
\end{pmatrix}} + \|\bw_{\barepsilon_t, \barnu_t}^t\|^2}.
\end{equation*}
Moreover, there exists $C_6>0$ such that
\begin{equation}\label{pequ:39}
\|\diag^2(g_t)\blambda_t\| \leq C_6 \nbr{\begin{pmatrix}
g_{t_a}\\
\blambda_{t_c}
\end{pmatrix}} \leq \frac{C_6}{\barepsilon_tq_{\barnu_t}^t \wedge 1}\nbr{\begin{pmatrix}
g_{t_a}\\
-\barepsilon_tq_{\barnu_t}^t\blambda_{t_c}
\end{pmatrix}} \stackrel{\eqref{pequ:38}}{\leq}\frac{C_6}{\tepsilon\kappa_{\barnu_0}\wedge 1}\|\bw_{\barepsilon_t, \barnu_t}^t\|,
\end{equation}
and
\begin{equation}\label{pequ:40}
\|\bw_{\barepsilon_t, \barnu_t}^t\| \stackrel{\text{Lem. } \ref{lem:3}}{\leq} C_6(\barepsilon_tq_{\barnu_t}^t \vee 1)\leq C_6(\barepsilon_0\tnu\vee 1).
\end{equation}
Combining the above three displays, there exists $C_7>0$ such that
\begin{equation}\label{pequ:37}
\|\bnabla_{\bx}\mL_t\| \leq  \|\bnabla_{\bx}\mL_{\barepsilon_t, \barnu_t, \eta}^t\| + C_7\cbr{\nbr{\begin{pmatrix}
J_t\bnabla_{\bx}\mL_t\\
G_t\bnabla_{\bx}\mL_t
\end{pmatrix}} + \nbr{\begin{pmatrix}
c_t\\
\bw_{\barepsilon_t, \barnu_t}^t
\end{pmatrix}} }. 
\end{equation}
We deal with the middle term. We know that
\begin{multline}\label{npequ:3}
\begin{pmatrix}
M_{11,t} & M_{12,t}\\
M_{21,t} & M_{22,t}
\end{pmatrix}\begin{pmatrix}
J_t\bnabla_{\bx}\mL_t\\
G_t\bnabla_{\bx}\mL_t
\end{pmatrix} \\\stackrel{\eqref{equ:aug:der}}{=} \frac{1}{\eta}\begin{pmatrix}
\bnabla_{\bmu}\mL_{\barepsilon_t, \barnu_t, \eta}^t\\
\bnabla_{\blambda}\mL_{\barepsilon_t, \barnu_t, \eta}^t
\end{pmatrix} - \frac{1}{\eta}\left(\begin{smallmatrix}
c_t\\
\bw_{\barepsilon_t, \barnu_t}^t + \frac{\|\bw_{\barepsilon_t, \barnu_t}^t\|^2}{\barepsilon_ta_{\barnu_t}^t}\blambda_t
\end{smallmatrix}\right) - \begin{pmatrix}
M_{12,t}\\
M_{22,t}
\end{pmatrix}\diag^2(g_t)\blambda_t.
\end{multline}
Multiplying $((J_t\bnabla_{\bx}\mL_t)^T\; (G_t\bnabla_{\bx}\mL_t)^T)$ on both sides, there exists $C_8>0$ such that
\begin{multline}\label{pequ:41}
\|J_t^TJ_t\bnabla_{\bx}\mL_t + G_t^TG_t\bnabla_{\bx}\mL_t\|^2 \leq  \begin{pmatrix}
J_t\bnabla_{\bx}\mL_t\\
G_t\bnabla_{\bx}\mL_t
\end{pmatrix}^T\begin{pmatrix}
M_{11,t} & M_{12,t}\\
M_{21,t} & M_{22,t}
\end{pmatrix}\begin{pmatrix}
J_t\bnabla_{\bx}\mL_t\\
G_t\bnabla_{\bx}\mL_t
\end{pmatrix}\\
\stackrel{\substack{\eqref{npequ:3}, \eqref{pequ:38}-\eqref{pequ:40}}}{\leq} C_8\|\bnabla_{\bx}\mL_t\|\cbr{\nbr{
\begin{pmatrix}
\bnabla_{\bmu}\mL_{\barepsilon_t, \barnu_t, \eta}^t\\
\bnabla_{\blambda}\mL_{\barepsilon_t, \barnu_t, \eta}^t
\end{pmatrix} } + \nbr{\begin{pmatrix}
c_t\\
\bw_{\barepsilon_t, \barnu_t}^t
\end{pmatrix}} }.
\end{multline}
Furthermore,
\begin{multline*}
\nbr{\begin{pmatrix}
J_t\bnabla_{\bx}\mL_t\\
G_t\bnabla_{\bx}\mL_t
\end{pmatrix}}^2\leq  \|\bnabla_{\bx}\mL_t\|\cdot \|J_t^TJ_t\bnabla_{\bx}\mL_t + G_t^TG_t\bnabla_{\bx}\mL_t\|\\
\stackrel{\eqref{pequ:41}}{\leq} \sqrt{C_8} \|\bnabla_{\bx}\mL_t\|^{\frac{3}{2}}\cbr{\nbr{
\begin{pmatrix}
\bnabla_{\bmu}\mL_{\barepsilon_t, \barnu_t, \eta}^t\\
\bnabla_{\blambda}\mL_{\barepsilon_t, \barnu_t, \eta}^t
\end{pmatrix} } + \nbr{\begin{pmatrix}
c_t\\
\bw_{\barepsilon_t, \barnu_t}^t
\end{pmatrix}} }^{\frac{1}{2}}.
\end{multline*}
Combining the above display with \eqref{pequ:37}, there exists $C_9>0$ such that
\begin{align*}
\|\bnabla_{\bx}\mL_t\|\leq &C_9\cbr{\nbr{\bnabla\mL_{\barepsilon_t, \barnu_t, \eta}^t} + \nbr{\begin{pmatrix}
c_t\\
\bw_{\barepsilon_t, \barnu_t}^t
\end{pmatrix}} } + C_9^{1/4}\|\bnabla_{\bx}\mL_t\|^{\frac{3}{4}}\cbr{\nbr{\bnabla\mL_{\barepsilon_t, \barnu_t, \eta}^t} + \nbr{\begin{pmatrix}
c_t\\
\bw_{\barepsilon_t, \barnu_t}^t
\end{pmatrix}} }^{\frac{1}{4}}\\
\leq & \frac{5C_9}{4}\cbr{\nbr{\bnabla\mL_{\barepsilon_t, \barnu_t, \eta}^t} + \nbr{\begin{pmatrix}
c_t\\
\bw_{\barepsilon_t, \barnu_t}^t
\end{pmatrix}} } + \frac{3}{4}\|\bnabla_{\bx}\mL_t\|,
\end{align*}
where the second inequality is due to Young's inequality $a^{3/4}b^{1/4} \leq 3a/4 + b/4$. Thus,
\begin{equation*}
\|\bnabla_{\bx}\mL_t\| \leq 5C_9\cbr{\nbr{\bnabla\mL_{\barepsilon_t, \barnu_t, \eta}^t} + \nbr{\begin{pmatrix}
c_t\\
\bw_{\barepsilon_t, \barnu_t}^t
\end{pmatrix}}}.
\end{equation*}
\noindent\textbf{(c).} By \eqref{equ:aug:der} and using \eqref{pequ:38}, \eqref{pequ:39} and \eqref{pequ:40}, there exists $C_{10}>0$ such that
\begin{align}\label{npequ:4}
&\nbr{\bnabla\mL_{\barepsilon_t, \barnu_t, \eta}^t}  \leq \|\bnabla_{\bx}\mL_t\| + C_{10}\nbr{\left(\begin{smallmatrix}
J_t\bnabla_{\bx}\mL_t\\
G_t\bnabla_{\bx}\mL_t + \Pi_c(\diag^2(g_t)\blambda_t)
\end{smallmatrix} \right)} + C_{10}\nbr{\left(\begin{smallmatrix}
c_t\\
\bw_{\barepsilon_t, \barnu_t}^t
\end{smallmatrix} \right)} \nonumber\\
& \leq \|\bnabla_{\bx}\mL_t\| + C_{10}\nbr{\left(\begin{smallmatrix}
J_t\bnabla_{\bx}\mL_t\\
G_t\bnabla_{\bx}\mL_t + \Pi_c(\diag^2(g_t)\blambda_t)
\end{smallmatrix} \right)} + C_{10}(\barepsilon_tq_{\barnu_t}^t\vee 1) \nbr{\left(\begin{smallmatrix}
c_t\\
g_{t_a}\\
\blambda_{t_c}
\end{smallmatrix} \right)} \nonumber\\
& \stackrel{\mathclap{\eqref{pequ:38}}}{\leq} \;\;\; \|\bnabla_{\bx}\mL_t\| + C_{10}\nbr{\left(\begin{smallmatrix}
J_t\bnabla_{\bx}\mL_t\\
G_t\bnabla_{\bx}\mL_t + \Pi_c(\diag^2(g_t)\blambda_t)
\end{smallmatrix} \right)} + C_{10}(\barepsilon_0\tnu\vee 1) \nbr{\left(\begin{smallmatrix}
c_t\\
g_{t_a}\\
\blambda_{t_c}
\end{smallmatrix} \right)}.
\end{align}
For $\bnabla_{\bx}\mL_t$, we have the following decomposition
\begin{align*}
\bnabla_{\bx}\mL_t = \underbrace{\cbr{I - (J_t^T\; G_{t_a}^T)\cbr{\begin{pmatrix}
J_t\\
G_{t_a}
\end{pmatrix}(J_t^T\; G_{t_a}^T) }^{-1}\begin{pmatrix}
J_t\\
G_{t_a}
\end{pmatrix} }}_{\P_{JG}^t}\bnabla_{\bx}\mL_t + (I - \P_{JG}^t)\bnabla_{\bx}\mL_t.
\end{align*}
By Assumptions \ref{ass:3} and \ref{ass:5}, we know $\|(I - \P_{JG}^t)\bnabla_{\bx}\mL_t\|\leq C_{11}\|(J_t\bnabla_{\bx}\mL_t, G_{t_a}\bnabla_{\bx}\mL_t)\|$ for some constant $C_{11}>0$. Furthermore, for some constant $C_{12}>0$, we also have
\begin{multline*}
\|\P_{JG}^t\bnabla_{\bx}\mL_t\| \stackrel{\eqref{equ:SQP:direction:1}}{=}  \nbr{\P_{JG}^t\cbr{B_t\barDelta\bx_t + J_t^T\bar{\tDelta}\bmu_t + G_{t_a}^T\bar{\tDelta}\blambda_{t_a} - G_{t_c}^T\blambda_{t_c}}}\\
\leq  \|\P_{JG}^tB_t\barDelta\bx_t\| + \|\P_{JG}^tG_{t_c}^T\blambda_{t_c}\|
\stackrel{\eqref{pequ:30}}{\leq}C_{12}\nbr{\left(\begin{smallmatrix}
\barDelta\bx_t\\
J_t\bnabla_{\bx}\mL_t\\
G_t\bnabla_{\bx}\mL_t+ \Pi_c(\diag^2(g_t)\blambda_t)
\end{smallmatrix}\right)}.
\end{multline*}
Combining the last two displays, we have
\begin{equation}\label{pequ:45}
\|\bnabla_{\bx}\mL_t\| \leq (C_{11}+ C_{12})\nbr{\left(\begin{smallmatrix}
\barDelta\bx_t\\
J_t\bnabla_{\bx}\mL_t\\
G_t\bnabla_{\bx}\mL_t+ \Pi_c(\diag^2(g_t)\blambda_t)
\end{smallmatrix} \right)}.
\end{equation}
Moreover, there exists $C_{13}>0$ such that
\begin{equation}\label{pequ:46}
\nbr{\begin{pmatrix}
c_t\\
g_{t_a}
\end{pmatrix}} \stackrel{\eqref{equ:SQP:direction:1}}{\leq} C_{13}\|\barDelta\bx_t\|, \quad\quad \|\blambda_{t_c}\| \stackrel{\eqref{pequ:30}}{\leq} C_{13}\nbr{\left(\begin{smallmatrix}
\barDelta\bx_t\\
J_t\bnabla_{\bx}\mL_t\\
G_t\bnabla_{\bx}\mL_t+ \Pi_c(\diag^2(g_t)\blambda_t)
\end{smallmatrix} \right)}.
\end{equation}
Combining \eqref{npequ:4}, \eqref{pequ:45}, and \eqref{pequ:46} together, we complete the proof.

\subsection{Proof of Lemma \ref{lem:10}}\label{pf:lem:10}

Analogous to the proof of Lemma \ref{lem:9}, we only track the constants $(\beta, \alpha_{max}, \kappa_{grad}, \kappa_{f}, p_{grad}, p_{f},\\ \chi_{grad}, \chi_{f})$. We use $\Upsilon_1, \Upsilon_2, \ldots$ to denote generic constants that are independent from $(\beta, \alpha_{max}, \\ \kappa_{grad}, \kappa_{f}, p_{grad}, p_{f}, \chi_{grad}, \chi_{f})$. Note that $\Upsilon_1$ in the proof may not~be~consistent with $\Upsilon_1$ in the statement, while the existence of $\Upsilon_1$ in the statement follows directly from our proof.

Let $\Upsilon_{\epsilon, \nu, \eta}$ be the upper bound of the generalized Hessian of $\mL_{\epsilon, \nu, \eta}$ in the compact set $(\mX\cap \mT_{\theta\nu})\times \mM\times \Lambda$ (see \cite{Pillo2002Augmented} for the definition of the generalized Hessian). In particular, $\Upsilon_{\epsilon, \nu, \eta} = \sup_{(\mX\cap \mT_{\theta\nu})\times \mM\times \Lambda}\|\partial^2\mL_{\epsilon, \nu, \eta}\|$. Without loss of generality, we suppose $\tepsilon$ in Theorem~\ref{thm:2a} satisfies $\tilde{\epsilon} = \barepsilon_0/\rho^{\tilde{i}}$ for some integer $\tilde{i}$. Then, with definition $\tilde{j}$ in \eqref{bound:nu}, we let 
\begin{equation*}
\Upsilon_{\tepsilon, \tnu, \eta} = \max\{\Upsilon_{\epsilon, \nu, \eta}: \epsilon = \barepsilon_0/\rho^i, \nu = \rho^j\barnu_0, 1\leq i\leq \tilde{i}, 1\leq j\leq \tilde{j}\}
\end{equation*}
and have $\Upsilon_{\barepsilon_{\bart}, \barnu_{\bart}, \eta} \leq \Upsilon_{\tepsilon, \tnu, \eta}$. Noting that $\bx_{s_t}, \bx_t \in \mT_{\barnu_{\bart}}$, we apply the Taylor expansion and~have
\begin{align}\label{pequ:33}
\mL_{\barepsilon_{\bart}, \barnu_{\bart}, \eta}^{s_t} & \leq \mL_{\barepsilon_{\bart}, \barnu_{\bart}, \eta}^t + \baralpha_t (\nabla\mL_{\barepsilon_{\bart}, \barnu_{\bart}, \eta}^t)^T\cDelta_t + \frac{\Upsilon_{\tepsilon, \tnu, \eta}\baralpha_t^2}{2}\|\cDelta_t\|^2 \nonumber\\
& = \mL_{\barepsilon_{\bart}, \barnu_{\bart}, \eta}^t + \baralpha_t(\bnabla\mL_{\barepsilon_{\bart}, \barnu_{\bart}, \eta}^t)^T\cDelta_t + \baralpha_t (\nabla\mL_{\barepsilon_{\bart}, \barnu_{\bart}, \eta}^t - \bnabla\mL_{\barepsilon_{\bart}, \barnu_{\bart}, \eta}^t)^T\cDelta_t + \frac{\Upsilon_{\tepsilon, \tnu, \eta}\baralpha_t^2}{2}\|\cDelta_t\|^2 \nonumber\\
& \leq \mL_{\barepsilon_{\bart}, \barnu_{\bart}, \eta}^t + \baralpha_t(\bnabla\mL_{\barepsilon_{\bart}, \barnu_{\bart}, \eta}^t)^T\cDelta_t + \baralpha_t\|\cDelta_t\|\cdot \barDelta(\nabla\mL_{\eta}^t) + \frac{\Upsilon_{\tepsilon, \tnu, \eta}\baralpha_t^2}{2}\|\cDelta_t\|^2 \nonumber\\
& \stackrel{\mathclap{\eqref{event:E1}}}{\leq} \mL_{\barepsilon_{\bart}, \barnu_{\bart}, \eta}^t + \baralpha_t(\bnabla\mL_{\barepsilon_{\bart}, \barnu_{\bart}, \eta}^t)^T\cDelta_t + \kappa_{grad}\baralpha_t^2\cdot\barR_t \|\cDelta_t\| + \frac{\Upsilon_{\tepsilon, \tnu, \eta}\baralpha_t^2}{2}\|\cDelta_t\|^2.
\end{align}
We consider the following two cases.

\noindent\textbf{Case 1}, $\cDelta_t = \barDelta_t$. Combining \eqref{cond:decrease:1} with \eqref{cond:decrease:2}, we have
\begin{equation}\label{pequ:34}
(\bnabla\mL_{\barepsilon_{\bart}, \barnu_{\bart}, \eta}^t)^T\barDelta_t \leq - \frac{\gamma_{B} \wedge \eta}{4}\nbr{\left(\begin{smallmatrix}
\barDelta\bx_t\\
J_t\bnabla_{\bx}\mL_t\\
G_t\bnabla_{\bx}\mL_t + \Pi_c(\diag^2(g_t)\blambda_t)
\end{smallmatrix} \right)}^2.
\end{equation}
By \eqref{pequ:29}, there exists $\Upsilon_1>0$ such that
\begin{equation}\label{pequ:49}
\|\barDelta_t\| \leq \Upsilon_1\nbr{\left(\begin{smallmatrix}
\barDelta\bx_t\\
J_t\bnabla_{\bx}\mL_t\\
G_t\bnabla_{\bx}\mL_t + \Pi_c(\diag^2(g_t)\blambda_t)
\end{smallmatrix} \right)}. 
\end{equation}
Furthermore, we have
\begin{equation}\label{pequ:16}
\barR_t \stackrel{\text{Lem. \ref{lem:3}}}{\leq} \frac{1}{\barepsilon_{\bart} q_{\barnu_{\bart}}^t\wedge 1}\nbr{\left(\begin{smallmatrix}
\bnabla_{\bx}\mL_t\\
c_t\\
\bw_{\barepsilon_{\bart}, \barnu_{\bart}}^t
\end{smallmatrix} \right)} \stackrel{\eqref{pequ:38}}{\leq} \frac{1}{\tepsilon\kappa_{\barnu_0} \wedge 1} \nbr{\left(\begin{smallmatrix}
\bnabla_{\bx}\mL_t\\
c_t\\
g_{t_a}\\
-\barepsilon_{\bart} q_{\barnu_{\bart}}^t\blambda_{t_c}
\end{smallmatrix} \right) } \leq \frac{\barepsilon_0\tnu \vee 1}{\tepsilon\kappa_{\barnu_0} \wedge 1} \nbr{\left(\begin{smallmatrix}
\bnabla_{\bx}\mL_t\\
c_t\\
g_{t_a}\\
\blambda_{t_c}
\end{smallmatrix}\right) },
\end{equation}
and thus, by \eqref{pequ:45}, \eqref{pequ:46}, there exists $\Upsilon_2>0$ such that
\begin{equation}\label{pequ:47}
\barR_t \leq \Upsilon_2\nbr{\left(\begin{smallmatrix}
\barDelta\bx_t\\
J_t\bnabla_{\bx}\mL_t\\
G_t\bnabla_{\bx}\mL_t + \Pi_c(\diag^2(g_t)\blambda_t)
\end{smallmatrix} \right)}.
\end{equation}
Plugging \eqref{pequ:49} and \eqref{pequ:47} into \eqref{pequ:33}, we have
\begin{align}\label{pequ:35}
&\mL_{\barepsilon_{\bart}, \barnu_{\bart}, \eta}^{s_t} \leq \mL_{\barepsilon_{\bart}, \barnu_{\bart}, \eta}^t + \baralpha_t(\bnabla\mL_{\barepsilon_{\bart}, \barnu_{\bart}, \eta}^t)^T\barDelta_t \nonumber\\
& \quad + \cbr{\Upsilon_1\Upsilon_2\kappa_{grad} + \frac{\Upsilon_{\tepsilon, \tnu, \eta}\Upsilon_1^2}{2}}\baralpha_t^2\nbr{\left(\begin{smallmatrix}
\barDelta\bx_t\\
J_t\bnabla_{\bx}\mL_t\\
G_t\bnabla_{\bx}\mL_t + \Pi_c(\diag^2(g_t)\blambda_t)
\end{smallmatrix}\right)}^2 \nonumber\\
& \stackrel{\mathclap{\eqref{pequ:34}}}{\leq} \;\; \mL_{\barepsilon_{\bart}, \barnu_{\bart}, \eta}^t + \baralpha_t(\bnabla\mL_{\barepsilon_{\bart}, \barnu_{\bart}, \eta}^t)^T\barDelta_t - \cbr{\Upsilon_1\Upsilon_2\kappa_{grad} + \frac{\Upsilon_{\tepsilon, \tnu, \eta}\Upsilon_1^2}{2}}\frac{4\baralpha_t^2}{\gamma_{B}\wedge\eta}(\bnabla\mL_{\barepsilon_{\bart}, \barnu_{\bart}, \eta}^t)^T\barDelta_t \nonumber\\
&\leq \mL_{\barepsilon_{\bart}, \barnu_{\bart}, \eta}^t  + \baralpha_t\cbr{1 - \Upsilon_3\rbr{\kappa_{grad} + 1}\baralpha_t }(\bnabla\mL_{\barepsilon_{\bart}, \barnu_{\bart}, \eta}^t)^T\barDelta_t,
\end{align}
where $\Upsilon_3 = 4\Upsilon_1\Upsilon_2/(\gamma_{B}\wedge\eta) \vee 2\Upsilon_1^2\Upsilon_{\tepsilon, \tnu, \eta}/(\gamma_{B}\wedge\eta)$. 

\noindent\textbf{Case 2}, $\cDelta_t = \hDelta_t$. By Lemma \ref{lem:9}(b), Lemma \ref{lem:3}, \eqref{cond:bound:fes:error}, and \eqref{pequ:38}, there exists $\Upsilon_4>0$ such~that
\begin{equation}\label{pequ:57}
\barR_t \leq \Upsilon_4\|\bnabla\mL_{\barepsilon_{\bart}, \barnu_{\bart}, \eta}^t\|.
\end{equation}
Plugging \eqref{cond:alter:dir} and \eqref{pequ:57} into \eqref{pequ:33}, we have
\begin{align}\label{pequ:42}
\mL_{\barepsilon_{\bart}, \barnu_{\bart}, \eta}^{s_t} \leq & \mL_{\barepsilon_{\bart}, \barnu_{\bart}, \eta}^t + \baralpha_t(\bnabla\mL_{\barepsilon_{\bart}, \barnu_{\bart}, \eta}^t)^T\hDelta_t + \Upsilon_4\chi_{u}\kappa_{grad} \baralpha_t^2\cdot\|\bnabla\mL_{\barepsilon_{\bart}, \barnu_{\bart}, \eta}^t\|^2 + \frac{\Upsilon_{\tepsilon, \tnu, \eta}\chi_{u}^2\baralpha_t^2}{2}\|\bnabla\mL_{\barepsilon_{\bart}, \barnu_{\bart}, \eta}^t\|^2 \nonumber\\
\stackrel{\mathclap{\eqref{cond:alter:dir}}}{\leq} &\;\; \mL_{\barepsilon_{\bart}, \barnu_{\bart}, \eta}^t + \baralpha_t(\bnabla\mL_{\barepsilon_{\bart}, \barnu_{\bart}, \eta}^t)^T\hDelta_t - \rbr{\Upsilon_4\chi_{u}^2\kappa_{grad} + \frac{\Upsilon_{\tepsilon, \tnu, \eta}\chi_{u}^3}{2}} \baralpha_t^2\cdot(\bnabla\mL_{\barepsilon_{\bart}, \barnu_{\bart}, \eta}^t)^T\hDelta_t \nonumber\\
\leq & \mL_{\barepsilon_{\bart}, \barnu_{\bart}, \eta}^t + \baralpha_t\cbr{1 - \Upsilon_5\rbr{\kappa_{grad} + 1 }\baralpha_t }(\bnabla\mL_{\barepsilon_{\bart}, \barnu_{\bart}, \eta}^t)^T\hDelta_t, 
\end{align}
where $\Upsilon_5 = \Upsilon_4\chi_{u}^2 \vee \Upsilon_{\tepsilon, \tnu, \eta}\chi_{u}^3/2$.

Combining \eqref{pequ:35} and \eqref{pequ:42}, and letting $\Upsilon_6 = \Upsilon_3\vee\Upsilon_5\vee 2$, we obtain
\begin{equation}\label{pequ:43}
\mL_{\barepsilon_{\bart}, \barnu_{\bart}, \eta}^{s_t} \leq \mL_{\barepsilon_{\bart}, \barnu_{\bart}, \eta}^{t} + \baralpha_t\cbr{1 - \Upsilon_6(\kappa_{grad} + 1)\baralpha_t} (\bnabla\mL_{\barepsilon_{\bart}, \barnu_{\bart}, \eta}^t)^T\cDelta_t.
\end{equation}
By the event $\E_2^t$, we have
\begin{align*}
\barL_{\barepsilon_{\bart}, \barnu_{\bart}, \eta}^{s_t} &\stackrel{\mathclap{\eqref{event:E2}}}{\leq}  \mL_{\barepsilon_{\bart}, \barnu_{\bart}, \eta}^{s_t} - \kappa_{f}\baralpha_t^2(\bnabla\mL_{\barepsilon_{\bart}, \barnu_{\bart}, \eta}^t)^T\cDelta_t\\
&\stackrel{\mathclap{\eqref{pequ:43}}}{\leq}  \;\; \mL_{\barepsilon_{\bart}, \barnu_{\bart}, \eta}^t + \baralpha_t\cbr{1 - \Upsilon_6\rbr{\kappa_{grad} + 1}\baralpha_t -\kappa_{f}\baralpha_t }(\bnabla\mL_{\barepsilon_{\bart}, \barnu_{\bart}, \eta}^t)^T\cDelta_t\\
&\stackrel{\mathclap{\eqref{event:E2}}}{\leq}  \barL_{\barepsilon_{\bart}, \barnu_{\bart}, \eta}^t + \baralpha_t\cbr{1 - \Upsilon_6\rbr{\kappa_{grad} + 1}\baralpha_t -2\kappa_{f}\baralpha_t }(\bnabla\mL_{\barepsilon_{\bart}, \barnu_{\bart}, \eta}^t)^T\cDelta_t\\
&\leq  \barL_{\barepsilon_{\bart}, \barnu_{\bart}, \eta}^t + \baralpha_t\cbr{1 - \Upsilon_6\rbr{\kappa_{grad} +\kappa_{f} + 1}\baralpha_t }(\bnabla\mL_{\barepsilon_{\bart}, \barnu_{\bart}, \eta}^t)^T\cDelta_t \quad (\text{since } \Upsilon_6\geq 2).
\end{align*}
Therefore, as long as
\begin{equation*}
1 - \Upsilon_6\rbr{\kappa_{grad} +\kappa_{f} + 1}\baralpha_t\geq \beta \Longleftrightarrow \baralpha_t\leq \frac{1-\beta}{\Upsilon_6(\kappa_{grad}+\kappa_{f} + 1) },
\end{equation*}
we have
\begin{equation*}
\barL_{\barepsilon_{\bart}, \barnu_{\bart}, \eta}^{s_t}\leq \barL_{\barepsilon_{\bart}, \barnu_{\bart}, \eta}^t + \baralpha_t\beta(\bnabla\mL_{\barepsilon_{\bart}, \barnu_{\bart}, \eta}^t)^T\cDelta_t.
\end{equation*}
This completes the proof.

\subsection{Proof of Lemma \ref{lem:12}}\label{pf:lem:12}

Algorithm \ref{alg:ASto} has three types of steps: a reliable step (Line 19), an unreliable step (Line~21), and an unsuccessful step (Line 24). For each type of step, $\cDelta_t = \barDelta_t$ or $\cDelta_t = \hDelta_t$. Thus, we analyze in the following six cases.
\vskip 5pt
\noindent\textbf{Case 1a, reliable step, $\cDelta_t = \barDelta_t$}. By Lemma \ref{lem:11}, we have
\begin{align}\label{pequ:44}
\mL_{\barepsilon_{\bart}, \barnu_{\bart}, \eta}^{t+1} - \mL_{\barepsilon_{\bart}, \barnu_{\bart}, \eta}^t & \leq  \frac{\baralpha_t\beta}{2}(\bnabla\mL_{\barepsilon_{\bart}, \barnu_{\bart}, \eta}^t)^T\barDelta_t
\stackrel{\eqref{cond:decrease:trust}}{\leq}  \frac{4\baralpha_t\beta}{9}(\bnabla\mL_{\barepsilon_{\bart}, \barnu_{\bart}, \eta}^t)^T\barDelta_t - \frac{\bardelta_t}{18} \nonumber\\
& \stackrel{\mathclap{\eqref{pequ:34}}}{\leq}  - \frac{\baralpha_t\beta(\gamma_{B}\wedge\eta)}{9}\nbr{\left(\begin{smallmatrix}
\barDelta\bx_t\\
J_t\bnabla_{\bx}\mL_t\\
G_t\bnabla_{\bx}\mL_t + \Pi_c(\diag^2(g_t)\blambda_t)
\end{smallmatrix} \right)}^2 - \frac{\bardelta_t}{18}.
\end{align}
Note that
\begin{equation*}
\|\nabla\mL_{\barepsilon_{\bart}, \barnu_{\bart}, \eta}^t\| \leq \barDelta(\nabla\mL_{\eta}^t) + \|\bnabla\mL_{\barepsilon_{\bart}, \barnu_{\bart}, \eta}^t\| \stackrel{\eqref{event:E1}}{\leq} \kappa_{grad}\baralpha_t\barR_t + \|\bnabla\mL_{\barepsilon_{\bart}, \barnu_{\bart}, \eta}^t\|.
\end{equation*}
Combining the above display with \eqref{pequ:47}, Lemma \ref{lem:9}(c), and using $\baralpha_t\leq \alpha_{max}$, there exists $\Upsilon_1>0$ such that 
\begin{equation}\label{pequ:51}
\|\nabla\mL_{\barepsilon_{\bart}, \barnu_{\bart}, \eta}^t\|  \leq \Upsilon_1(\kappa_{grad}\alpha_{max} + 1) \nbr{\left(\begin{smallmatrix}
\barDelta\bx_t\\
J_t\bnabla_{\bx}\mL_t\\
G_t\bnabla_{\bx}\mL_t + \Pi_c(\diag^2(g_t)\blambda_t)
\end{smallmatrix} \right)}.
\end{equation}
Combining the above inequality with \eqref{pequ:44}, we have
\begin{multline}\label{pequ:48}
\mL_{\barepsilon_{\bart}, \barnu_{\bart}, \eta}^{t+1} - \mL_{\barepsilon_{\bart}, \barnu_{\bart}, \eta}^t \leq  -\frac{\baralpha_t\beta(\gamma_{B}\wedge\eta)}{18}\nbr{\left(\begin{smallmatrix}
\barDelta\bx_t\\	
J_t\bnabla_{\bx}\mL_t\\
G_t\bnabla_{\bx}\mL_t + \Pi_c(\diag^2(g_t)\blambda_t)
\end{smallmatrix} \right)}^2\\
- \frac{\baralpha_t\beta(\gamma_{B}\wedge\eta)}{18\Upsilon_1^2(\kappa_{grad}\alpha_{max} + 1)^2}\|\nabla\mL_{\barepsilon_{\bart}, \barnu_{\bart}, \eta}^t\|^2 - \frac{\bardelta_t}{18}.
\end{multline}
By Line 20 of Algorithm \ref{alg:ASto}, $\bardelta_{t+1} - \bardelta_t = (\rho-1)\bardelta_t$. By the Taylor expansion and $\baralpha_{t+1} \leq \rho\baralpha_t$ (Line 18), there exists $\Upsilon_2>0$ such that
\begin{multline}\label{pequ:52}
\baralpha_{t+1}\|\nabla\mL_{\barepsilon_{\bart}, \barnu_{\bart}, \eta}^{t+1}\|^2 - \baralpha_t\|\nabla\mL_{\barepsilon_{\bart}, \barnu_{\bart}, \eta}^t\|^2\leq  2\rho\baralpha_t\cbr{\|\nabla\mL_{\barepsilon_{\bart}, \barnu_{\bart}, \eta}^t\|^2 + \Upsilon_{\tepsilon, \tnu, \eta}^2\baralpha_t^2\|\barDelta_t\|^2}\\
\stackrel{\eqref{pequ:49}}{\leq} 2\rho\baralpha_t\cbr{\|\nabla\mL_{\barepsilon_{\bart}, \barnu_{\bart}, \eta}^t\|^2 + \Upsilon_{\tepsilon, \tnu, \eta}^2\alpha_{max}^2\Upsilon_2\nbr{\left(\begin{smallmatrix}
\barDelta\bx_t\\
J_t\bnabla_{\bx}\mL_t\\
G_t\bnabla_{\bx}\mL_t + \Pi_c(\diag^2(g_t)\blambda_t)
\end{smallmatrix} \right)}^2}.
\end{multline}
Combining the above two displays with \eqref{equ:potential}, we obtain
\begin{align*}
\Theta_{\omega}^{t+1} - \Theta_{\omega}^t \leq& -\rbr{\frac{\omega\beta(\gamma_{B}\wedge\eta)}{18} - (1-\omega)\rho\Upsilon_{\tepsilon, \tnu, \eta}^2\alpha_{max}^2\Upsilon_2}\baralpha_t\nbr{\left(\begin{smallmatrix}
\barDelta\bx_t\\
J_t\bnabla_{\bx}\mL_t\\
G_t\bnabla_{\bx}\mL_t + \Pi_c(\diag^2(g_t)\blambda_t)
\end{smallmatrix} \right)}^2\\
& - \rbr{\frac{\omega\beta(\gamma_{B}\wedge\eta)}{18\Upsilon_1^2(\kappa_{grad}\alpha_{max} + 1)^2} - (1-\omega)\rho}\baralpha_t\|\nabla\mL_{\barepsilon_{\bart}, \barnu_{\bart}, \eta}^t\|^2\\
& - \rbr{\frac{\omega}{18} - \frac{(1-\omega)(\rho-1)}{2}}\bardelta_t.
\end{align*}
Let
\begin{align}\label{pequ:56}
\frac{\omega\beta(\gamma_{B}\wedge\eta)}{36}\geq (1-\omega)\rho\Upsilon_{\tepsilon, \tnu, \eta}^2\alpha_{max}^2\Upsilon_2 &\Longleftrightarrow  \frac{\omega}{1-\omega}\geq \frac{36\rho\Upsilon_{\tepsilon, \tnu, \eta}^2\alpha_{max}^2\Upsilon_2}{\beta(\gamma_{B}\wedge\eta)}, \nonumber\\
\frac{\omega\beta(\gamma_{B}\wedge\eta)}{36\Upsilon_1^2(\kappa_{grad}\alpha_{max} + 1)^2} \geq (1-\omega)\rho & \Longleftrightarrow\frac{\omega}{1-\omega} \geq \frac{36\rho\Upsilon_1^2(\kappa_{grad}\alpha_{max}+1)^2}{\beta(\gamma_{B}\wedge\eta)},\\
\frac{\omega}{36} \geq \frac{(1-\omega)(\rho-1)}{2} & \Longleftrightarrow \frac{\omega}{1-\omega} \geq 18(\rho-1), \nonumber
\end{align}
which is further implied by
\begin{equation}\label{pequ:53}
\frac{\omega}{1-\omega} \geq \frac{\Upsilon_3(\kappa_{grad}\alpha_{max}+\alpha_{max} + 1)^2 }{\beta} \vee 18(\rho-1)
\end{equation}
if we define $\Upsilon_3 = (36\rho\Upsilon_{\tepsilon, \tnu, \eta}^2\Upsilon_2\vee 36\rho\Upsilon_1^2)/(\gamma_{B}\wedge\eta)$. Then, we obtain
\begin{multline}\label{pequ:50}
\Theta_{\omega}^{t+1} - \Theta_{\omega}^t \leq -\frac{\omega\beta(\gamma_{B}\wedge\eta)}{36}\cdot\baralpha_t\nbr{\left(\begin{smallmatrix}
\barDelta\bx_t\\
J_t\bnabla_{\bx}\mL_t\\
G_t\bnabla_{\bx}\mL_t + \Pi_c(\diag^2(g_t)\blambda_t)
\end{smallmatrix} \right)}^2\\
- \frac{\omega\beta(\gamma_{B}\wedge\eta)}{36\Upsilon_1^2(\kappa_{grad}\alpha_{max} + 1)^2} \cdot\baralpha_t\|\nabla\mL_{\barepsilon_{\bart}, \barnu_{\bart}, \eta}^t\|^2 - \frac{\omega}{36}\bardelta_t.
\end{multline}
\noindent\textbf{Case 2a, unreliable step, $\cDelta_t = \barDelta_t$}. By Lemma \ref{lem:11}, we have
\begin{align*}
\mL_{\barepsilon_{\bart}, \barnu_{\bart}, \eta}^{t+1} &- \mL_{\barepsilon_{\bart}, \barnu_{\bart}, \eta}^t \leq  \frac{\baralpha_t\beta}{2}(\bnabla\mL_{\barepsilon_{\bart}, \barnu_{\bart}, \eta}^t)^T\barDelta_t\\
\stackrel{\eqref{pequ:34}}{\leq} & - \frac{\baralpha_t\beta(\gamma_{B}\wedge\eta)}{8}\nbr{\left(\begin{smallmatrix}
\barDelta\bx_t\\
J_t\bnabla_{\bx}\mL_t\\
G_t\bnabla_{\bx}\mL_t + \Pi_c(\diag^2(g_t)\blambda_t)
\end{smallmatrix} \right)}^2\\
\stackrel{\eqref{pequ:51}}{\leq} & - \frac{\baralpha_t\beta(\gamma_{B}\wedge\eta)}{16}\nbr{\left(\begin{smallmatrix}
\barDelta\bx_t\\
J_t\bnabla_{\bx}\mL_t\\
G_t\bnabla_{\bx}\mL_t + \Pi_c(\diag^2(g_t)\blambda_t)
\end{smallmatrix} \right)}^2 - \frac{\baralpha_t\beta(\gamma_{B}\wedge\eta)}{16\Upsilon_1^2(\kappa_{grad}\alpha_{max}+1)^2}\|\nabla\mL_{\barepsilon_{\bart}, \barnu_{\bart}, \eta}^t\|^2.
\end{align*}
By Line 22 of Algorithm \ref{alg:ASto}, $\bardelta_{t+1} - \bardelta_t = -(1-1/\rho)\bardelta_t$, while \eqref{pequ:52} still holds. Thus, under \eqref{pequ:53}, we have
\begin{multline}\label{pequ:54}
\Theta_{\omega}^{t+1} - \Theta_{\omega}^t \leq -\frac{\omega\beta(\gamma_{B}\wedge\eta)}{36}\cdot\baralpha_t\nbr{\left(\begin{smallmatrix}
\barDelta\bx_t\\
J_t\bnabla_{\bx}\mL_t\\
G_t\bnabla_{\bx}\mL_t + \Pi_c(\diag^2(g_t)\blambda_t)
\end{smallmatrix} \right)}^2\\
- \frac{\omega\beta(\gamma_{B}\wedge\eta)}{36\Upsilon_1^2(\kappa_{grad}\alpha_{max} + 1)^2} \cdot\baralpha_t\|\nabla\mL_{\barepsilon_{\bart}, \barnu_{\bart}, \eta}^t\|^2 - \frac{1}{2}\rbr{1-\omega}\rbr{1-\frac{1}{\rho}}\bardelta_t.
\end{multline}
\noindent\textbf{Case 3a, unsuccessful step, $\cDelta_t = \barDelta_t$}. In this case, $(\bx_{t+1},\bmu_{t+1}, \blambda_{t+1}) = (\bx_t, \bmu_t, \blambda_t)$, $\baralpha_{t+1} = \baralpha_t/\rho$ and $\bardelta_{t+1} = \bardelta_t/\rho$. Thus, we immediately have
\begin{equation}\label{pequ:55}
\Theta_{\omega}^{t+1} - \Theta_{\omega}^t \leq -\frac{1}{2}\rbr{1-\omega}\rbr{1 - \frac{1}{\rho}} \rbr{\baralpha_t\|\nabla\mL_{\barepsilon_{\bart}, \barnu_{\bart}, \eta}^t\|^2 + \bardelta_t}.
\end{equation}
Combining \eqref{pequ:50}, \eqref{pequ:54}, \eqref{pequ:55}, and noting that
\begin{align*}
\frac{\omega\beta(\gamma_{B}\wedge\eta)}{36\Upsilon_1^2(\kappa_{grad}\alpha_{max} + 1)^2} \geq \frac{1-\omega}{2}\rbr{1- \frac{1}{\rho}} &\Longleftarrow \frac{\omega}{1-\omega} \geq \frac{18\Upsilon_1^2(\kappa_{grad}\alpha_{max} + 1)^2}{\beta(\gamma_{B}\wedge \eta)},\\
\frac{\omega}{36}\geq \frac{1-\omega}{2}\rbr{1-\frac{1}{\rho}}&\Longleftarrow \frac{\omega}{1-\omega} \geq 18(\rho-1),
\end{align*}
with the right hand side being implied by \eqref{pequ:56} and further by \eqref{pequ:53}, we know \eqref{pequ:55}~holds for all three cases with $\cDelta_t = \barDelta_t$.

\noindent\textbf{Case 1b, reliable step, $\cDelta_t = \hatDelta_t$}. By Lemma \ref{lem:11}, we have
\begin{align*}
\mL_{\barepsilon_{\bart}, \barnu_{\bart}, \eta}^{t+1} - \mL_{\barepsilon_{\bart}, \barnu_{\bart}, \eta}^t & \leq \frac{\baralpha_t\beta}{2}(\bnabla\mL_{\barepsilon_{\bart}, \barnu_{\bart}, \eta}^t)^T\hDelta_t \\ 
& \stackrel{\mathclap{\eqref{cond:decrease:trust}}}{\leq}  \frac{\baralpha_t\beta}{3}(\bnabla\mL_{\barepsilon_{\bart}, \barnu_{\bart}, \eta}^t)^T\hDelta_t - \frac{\bardelta_t}{6}
\stackrel{\eqref{cond:alter:dir}}{\leq} -\frac{\baralpha_t\beta}{3\chi_{u}}\|\bnabla\mL_{\barepsilon_{\bart}, \barnu_{\bart}, \eta}^t\|^2 - \frac{\bardelta_t}{6}.
\end{align*}
Note that
\begin{equation*}
\|\nabla\mL_{\barepsilon_{\bart}, \barnu_{\bart}, \eta}^t\| \leq \barDelta(\nabla\mL_{\eta}^t) + \|\bnabla\mL_{\barepsilon_{\bart}, \barnu_{\bart}, \eta}^t\| \stackrel{\eqref{event:E1}}{\leq} \kappa_{grad}\baralpha_t\barR_t + \|\bnabla\mL_{\barepsilon_{\bart}, \barnu_{\bart}, \eta}^t\|.
\end{equation*}
Combining the above display with \eqref{pequ:57} and using $\baralpha_t\leq \alpha_{max}$, there exists $\Upsilon_4>0$ such~that
\begin{equation}\label{pequ:63}
\|\nabla\mL_{\barepsilon_{\bart}, \barnu_{\bart}, \eta}^t\| \leq \Upsilon_4(\kappa_{grad}\alpha_{max} + 1)\|\bnabla\mL_{\barepsilon_{\bart}, \barnu_{\bart}, \eta}^t\|.
\end{equation}
Combining the above three displays,
\begin{equation*}
\mL_{\barepsilon_{\bart}, \barnu_{\bart}, \eta}^{t+1} - \mL_{\barepsilon_{\bart}, \barnu_{\bart}, \eta}^t \leq -\frac{\baralpha_t\beta}{6\chi_{u}}\|\bnabla\mL_{\barepsilon_{\bart}, \barnu_{\bart}, \eta}^t\|^2 - \frac{\baralpha_t\beta}{6\Upsilon_4^2\chi_{u}(\kappa_{grad}\alpha_{max} + 1)^2}\|\nabla\mL_{\barepsilon_{\bart}, \barnu_{\bart}, \eta}^t\|^2 - \frac{\bardelta_t}{6}.
\end{equation*}
By Line 20 of Algorithm \ref{alg:ASto}, $\bardelta_{t+1} - \bardelta_t = (\rho-1)\bardelta_t$. By the Taylor expansion and $\baralpha_{t+1} \leq \rho\baralpha_t$ (Line 18),
\begin{multline}\label{pequ:60}
\baralpha_{t+1}\|\nabla\mL_{\barepsilon_{\bart}, \barnu_{\bart}, \eta}^{t+1}\|^2 - \baralpha_t\|\nabla\mL_{\barepsilon_{\bart}, \barnu_{\bart}, \eta}^t\|^2\leq  2\rho\baralpha_t\cbr{\|\nabla\mL_{\barepsilon_{\bart}, \barnu_{\bart}, \eta}^t\|^2 + \Upsilon_{\tepsilon, \tnu, \eta}^2\baralpha_t^2\|\hDelta_t\|^2}\\
\stackrel{\eqref{cond:alter:dir}}{\leq} 2\rho\baralpha_t\cbr{\|\nabla\mL_{\barepsilon_{\bart}, \barnu_{\bart}, \eta}^t\|^2 + \Upsilon_{\tepsilon, \tnu, \eta}^2\chi_{u}^2\alpha_{max}^2\|\bnabla\mL_{\barepsilon_{\bart}, \barnu_{\bart}, \eta}^t\|^2}.  
\end{multline}
Combining the above two displays,
\begin{align*}
\Theta_{\omega}^{t+1}-\Theta_{\omega}^t \leq & -\rbr{\frac{\omega\beta}{6\chi_{u}} -(1-\omega)\rho\Upsilon_{\tepsilon, \tnu, \eta}^2\chi_{u}^2\alpha_{max}^2}\baralpha_t\|\bnabla\mL_{\barepsilon_{\bart}, \barnu_{\bart}, \eta}^t\|^2\\
& - \rbr{\frac{\omega\beta}{6\Upsilon_4^2\chi_{u}(\kappa_{grad}\alpha_{max} + 1)^2} - (1-\omega)\rho }\baralpha_t\|\nabla\mL_{\barepsilon_{\bart}, \barnu_{\bart}, \eta}^t\|^2\\
& - \rbr{\frac{\omega}{6} - \frac{(1-\omega)(\rho-1)}{2}}\bardelta_t.
\end{align*}
Let
\begin{align}\label{pequ:62}
\frac{\omega\beta}{12\chi_{u}} \geq (1-\omega)\rho\Upsilon_{\tepsilon, \tnu, \eta}^2\chi_{u}^2\alpha_{max}^2 &\Longleftrightarrow \frac{\omega}{1-\omega}\geq \frac{12\rho\Upsilon_{\tepsilon, \tnu, \eta}^2\chi_{u}^3\alpha_{max}^2}{\beta}, \nonumber\\
\frac{\omega\beta}{12\Upsilon_4^2\chi_{u}(\kappa_{grad}\alpha_{max} + 1)^2} \geq (1-\omega)\rho&\Longleftrightarrow\frac{\omega}{1-\omega} \geq \frac{12\rho\Upsilon_4^2\chi_{u}(\kappa_{grad}\alpha_{max} + 1)^2}{\beta},\\
\frac{\omega}{12} \geq \frac{1-\omega}{2}(\rho-1) &\Longleftrightarrow \frac{\omega}{1-\omega} \geq 6(\rho-1), \nonumber
\end{align}
which is implied by \eqref{pequ:53} if we re-define $\Upsilon_3 \leftarrow \Upsilon_3\vee 12\rho\Upsilon_{\tepsilon, \tnu,\eta}^2\chi_{u}^3\vee12\rho\Upsilon_4^2\chi_{u}$. Then,
\begin{multline}\label{pequ:59}
\Theta_{\omega}^{t+1}-\Theta_{\omega}^t \leq -\frac{\omega\beta}{12\chi_{u}}\cdot\baralpha_t\|\bnabla\mL_{\barepsilon_{\bart}, \barnu_{\bart}, \eta}^t\|^2 \\- \frac{\omega\beta}{12\Upsilon_4^2\chi_{u}(\kappa_{grad}\alpha_{max} + 1)^2} \cdot\baralpha_t\|\nabla\mL_{\barepsilon_{\bart}, \barnu_{\bart}, \eta}^t\|^2 - \frac{\omega}{12}\bardelta_t.
\end{multline}
\noindent\textbf{Case 2b, unreliable step, $\cDelta_t = \hDelta_t$}. By Lemma \ref{lem:11}, we have
\begin{align*}
\mL_{\barepsilon_{\bart}, \barnu_{\bart}, \eta}^{t+1} - \mL_{\barepsilon_{\bart}, \barnu_{\bart}, \eta}^t \leq &  \frac{\baralpha_t\beta}{2}(\bnabla\mL_{\barepsilon_{\bart}, \barnu_{\bart}, \eta}^t)^T\hDelta_t 
\stackrel{\eqref{cond:alter:dir}}{\leq} -\frac{\baralpha_t\beta}{2\chi_{u}}\|\bnabla\mL_{\barepsilon_{\bart}, \barnu_{\bart}, \eta}^t\|^2\\
\stackrel{\eqref{pequ:63}}{\leq} & - \frac{\baralpha_t\beta}{4\chi_{u}}\|\bnabla\mL_{\barepsilon_{\bart}, \barnu_{\bart}, \eta}^t\|^2 - \frac{\baralpha_t\beta}{4\Upsilon_4^2\chi_{u}\rbr{\kappa_{grad}\alpha_{max} +1}^2 }\|\nabla\mL_{\barepsilon_{\bart}, \barnu_{\bart}, \eta}^t\|^2.
\end{align*}
By Line 22 of Algorithm \ref{alg:ASto}, $\bardelta_{t+1} - \bardelta_t = -(1-1/\rho)\bardelta_t$, while \eqref{pequ:60} still holds. Thus, under \eqref{pequ:53}, we have \begin{multline}\label{pequ:61}
\Theta_{\omega}^{t+1}-\Theta_{\omega}^t \leq  -\frac{\omega\beta}{12\chi_{u}}\cdot\baralpha_t\|\bnabla\mL_{\barepsilon_{\bart}, \barnu_{\bart}, \eta}^t\|^2 \\
- \frac{\omega\beta}{12\Upsilon_4^2\chi_{u}(\kappa_{grad}\alpha_{max} + 1)^2} \cdot\baralpha_t\|\nabla\mL_{\barepsilon_{\bart}, \barnu_{\bart}, \eta}^t\|^2 
- \frac{1}{2}\rbr{1-\omega}\rbr{1-\frac{1}{\rho}}\bardelta_t.
\end{multline}
\noindent\textbf{Case 3b, unsuccessful step, $\cDelta_t = \hDelta_t$}. In this case, \eqref{pequ:55} holds. Combining \eqref{pequ:59}, \eqref{pequ:61}, \eqref{pequ:55}, and noting that
\begin{align*}
\frac{\omega\beta}{12\Upsilon_4^2\chi_{u}(\kappa_{grad}\alpha_{max} + 1)^2} \geq \frac{1-\omega}{2}\rbr{1- \frac{1}{\rho}} &\Longleftarrow \frac{\omega}{1-\omega} \geq \frac{6\Upsilon_4^2\chi_{u}(\kappa_{grad}\alpha_{max} + 1)^2}{\beta},\\
\frac{\omega}{12}\geq \frac{1-\omega}{2}\rbr{1-\frac{1}{\rho}}&\Longleftarrow \frac{\omega}{1-\omega} \geq 6(\rho-1),
\end{align*}
as implied by \eqref{pequ:62} and further by \eqref{pequ:53}, we know \eqref{pequ:55} holds for all three cases with $\cDelta_t = \hDelta_t$. In summary, under \eqref{pequ:53}, \eqref{pequ:55} holds for all cases. This completes the proof.

\subsection{Proof of Lemma \ref{lem:13}}\label{pf:lem:13}

The proof follows the proof of Lemma \ref{lem:12}, except that \eqref{pequ:51} and \eqref{pequ:63} do not hold due to $(\E_1^t)^c$. We consider the following six cases.

\vskip5pt
\noindent\textbf{Case 1a, reliable step, $\cDelta_t = \barDelta_t$}. By Lemma \ref{lem:11}, we have
\begin{multline}\label{pequ:64}
\mL_{\barepsilon_{\bart}, \barnu_{\bart}, \eta}^{t+1} - \mL_{\barepsilon_{\bart}, \barnu_{\bart}, \eta}^t \leq  \frac{\baralpha_t\beta}{2}(\bnabla\mL_{\barepsilon_{\bart}, \barnu_{\bart}, \eta}^t)^T\barDelta_t
\stackrel{\eqref{cond:decrease:trust}}{\leq}  \frac{4\baralpha_t\beta}{9}(\bnabla\mL_{\barepsilon_{\bart}, \barnu_{\bart}, \eta}^t)^T\barDelta_t - \frac{\bardelta_t}{18}\\
\stackrel{\eqref{pequ:34}}{\leq}  - \frac{\baralpha_t\beta(\gamma_{B}\wedge\eta)}{9}\nbr{\left(\begin{smallmatrix}
\barDelta\bx_t\\
J_t\bnabla_{\bx}\mL_t\\
G_t\bnabla_{\bx}\mL_t + \Pi_c(\diag^2(g_t)\blambda_t)
\end{smallmatrix} \right)}^2 - \frac{\bardelta_t}{18}.
\end{multline}
By Line 20 of Algorithm \ref{alg:ASto}, $\bardelta_{t+1} - \bardelta_t = (\rho-1)\bardelta_t$, while \eqref{pequ:52} still holds. By the condition of $\omega$ in \eqref{pequ:56} and \eqref{pequ:53}, we know that under \eqref{cond:omega} (which implies \eqref{pequ:53}),
\begin{multline}\label{pequ:65}
\Theta_{\omega}^{t+1} - \Theta_{\omega}^t \leq -\frac{\omega\beta(\gamma_{B}\wedge\eta)}{36}\cdot\baralpha_t\nbr{\left(\begin{smallmatrix}
\barDelta\bx_t\\
J_t\bnabla_{\bx}\mL_t\\
G_t\bnabla_{\bx}\mL_t + \Pi_c(\diag^2(g_t)\blambda_t)
\end{smallmatrix} \right)}^2\\
+\rho(1-\omega)\baralpha_t\|\nabla\mL_{\barepsilon_{\bart}, \barnu_{\bart}, \eta}^t\|^2 - \frac{\omega}{36}\bardelta_t.
\end{multline}
\noindent\textbf{Case 2a, unreliable step, $\cDelta_t = \barDelta_t$}. By Lemma \ref{lem:11}, we have
\begin{equation*}
\mL_{\barepsilon_{\bart}, \barnu_{\bart}, \eta}^{t+1} - \mL_{\barepsilon_{\bart}, \barnu_{\bart}, \eta}^t \leq  \frac{\baralpha_t\beta}{2}(\bnabla\mL_{\barepsilon_{\bart}, \barnu_{\bart}, \eta}^t)^T\barDelta_t
\stackrel{\eqref{pequ:34}}{\leq}  - \frac{\baralpha_t\beta(\gamma_{B}\wedge\eta)}{8}\nbr{\left(\begin{smallmatrix}
\barDelta\bx_t\\
J_t\bnabla_{\bx}\mL_t\\
G_t\bnabla_{\bx}\mL_t + \Pi_c(\diag^2(g_t)\blambda_t)
\end{smallmatrix} \right)}^2.
\end{equation*}
By Line 22 of Algorithm \ref{alg:ASto}, $\bardelta_{t+1} - \bardelta_t = -(1-1/\rho)\bardelta_t$, while \eqref{pequ:52} still holds. Thus, under~\eqref{cond:omega}, 
\begin{multline}\label{pequ:66}
\Theta_{\omega}^{t+1} - \Theta_{\omega}^t \leq -\frac{\omega\beta(\gamma_{B}\wedge\eta)}{36}\cdot\baralpha_t\nbr{\left(\begin{smallmatrix}
\barDelta\bx_t\\
J_t\bnabla_{\bx}\mL_t\\
G_t\bnabla_{\bx}\mL_t + \Pi_c(\diag^2(g_t)\blambda_t)
\end{smallmatrix} \right)}^2\\
+\rho(1-\omega)\baralpha_t\|\nabla\mL_{\barepsilon_{\bart}, \barnu_{\bart}, \eta}^t\|^2 - \frac{1}{2}\rbr{1-\omega}\rbr{1-\frac{1}{\rho}}\bardelta_t.
\end{multline}
\noindent\textbf{Case 3a, unsuccessful step, $\cDelta_t = \barDelta_t$}. In this case, \eqref{pequ:55} holds. Combining \eqref{pequ:65}, \eqref{pequ:66}, and \eqref{pequ:55}, we have
\begin{equation}\label{pequ:67}
\Theta_{\omega}^{t+1} - \Theta_{\omega}^t \leq \rho(1-\omega)\baralpha_t\|\nabla\mL_{\barepsilon_{\bart}, \barnu_{\bart}, \eta}^t\|^2.
\end{equation}
\noindent\textbf{Case 1b, reliable step, $\cDelta_t = \hatDelta_t$}. By Lemma \ref{lem:11}, we have
\begin{align*}
\mL_{\barepsilon_{\bart}, \barnu_{\bart}, \eta}^{t+1} - \mL_{\barepsilon_{\bart}, \barnu_{\bart}, \eta}^t & \leq \frac{\baralpha_t\beta}{2}(\bnabla\mL_{\barepsilon_{\bart}, \barnu_{\bart}, \eta}^t)^T\hDelta_t \\
& \stackrel{\mathclap{\eqref{cond:decrease:trust}}}{\leq}  \frac{\baralpha_t\beta}{3}(\bnabla\mL_{\barepsilon_{\bart}, \barnu_{\bart}, \eta}^t)^T\hDelta_t - \frac{\bardelta_t}{6}
\stackrel{\eqref{cond:alter:dir}}{\leq} -\frac{\baralpha_t\beta}{3\chi_{u}}\|\bnabla\mL_{\barepsilon_{\bart}, \barnu_{\bart}, \eta}^t\|^2 - \frac{\bardelta_t}{6}.
\end{align*}
By Line 20 of Algorithm \ref{alg:ASto}, $\bardelta_{t+1} - \bardelta_t = (\rho-1)\bardelta_t$, while \eqref{pequ:60} still holds. By the condition of $\omega$ in \eqref{pequ:62}, we know that under \eqref{cond:omega} (which implies \eqref{pequ:62}), 
\begin{equation}\label{pequ:68}
\Theta_{\omega}^{t+1}-\Theta_{\omega}^t \leq -\frac{\omega\beta}{12\chi_{u}}\cdot\baralpha_t\|\bnabla\mL_{\barepsilon_{\bart}, \barnu_{\bart}, \eta}^t\|^2 +\rho(1-\omega)\baralpha_t\|\nabla\mL_{\barepsilon_{\bart}, \barnu_{\bart}, \eta}^t\|^2 - \frac{\omega}{12}\bardelta_t.
\end{equation}
\noindent\textbf{Case 2b, unreliable step, $\cDelta_t = \hDelta_t$}. By Lemma \ref{lem:11}, we have
\begin{equation*}
\mL_{\barepsilon_{\bart}, \barnu_{\bart}, \eta}^{t+1} - \mL_{\barepsilon_{\bart}, \barnu_{\bart}, \eta}^t \leq  \frac{\baralpha_t\beta}{2}(\bnabla\mL_{\barepsilon_{\bart}, \barnu_{\bart}, \eta}^t)^T\hDelta_t 
\stackrel{\eqref{cond:alter:dir}}{\leq} -\frac{\baralpha_t\beta}{2\chi_{u}}\|\bnabla\mL_{\barepsilon_{\bart}, \barnu_{\bart}, \eta}^t\|^2.
\end{equation*}
By Line 22 of Algorithm \ref{alg:ASto}, $\bardelta_{t+1} - \bardelta_t = -(1-1/\rho)\bardelta_t$, while \eqref{pequ:60} still holds. Thus, under~\eqref{cond:omega}, 
\begin{equation}\label{pequ:69}
\Theta_{\omega}^{t+1}-\Theta_{\omega}^t \leq  -\frac{\omega\beta}{12\chi_{u}}\cdot\baralpha_t\|\bnabla\mL_{\barepsilon_{\bart}, \barnu_{\bart}, \eta}^t\|^2 
+\rho(1-\omega) \baralpha_t\|\nabla\mL_{\barepsilon_{\bart}, \barnu_{\bart}, \eta}^t\|^2 
- \frac{1-\omega}{2}\rbr{1-\frac{1}{\rho}}\bardelta_t.
\end{equation}
\noindent\textbf{Case 3b, unsuccessful step, $\cDelta_t = \hDelta_t$}. In this case, \eqref{pequ:55} holds. Combining \eqref{pequ:68}, \eqref{pequ:69}, and \eqref{pequ:55}, we note that \eqref{pequ:67} holds as well. Thus, \eqref{pequ:67} holds for all six cases. This completes the proof.

\subsection{Proof of Lemma \ref{lem:14}}\label{pf:lem:14}

The proof follows the proof of Lemma \ref{lem:13}, except that Lemma \ref{lem:11} is not applicable. We consider the following six cases.

\vskip5pt
\noindent\textbf{Case 1a, reliable step, $\cDelta_t = \barDelta_t$}. We have
\begin{align*}
\mL_{\barepsilon_{\bart}, \barnu_{\bart}, \eta}^{t+1} - & \mL_{\barepsilon_{\bart}, \barnu_{\bart}, \eta}^t \\
\leq & \barL_{\barepsilon_{\bart}, \barnu_{\bart}, \eta}^{s_t} - \barL_{\barepsilon_{\bart}, \barnu_{\bart}, \eta}^t + \abr{\barL_{\barepsilon_{\bart}, \barnu_{\bart}, \eta}^{s_t} - \mL_{\barepsilon_{\bart}, \barnu_{\bart}, \eta}^{s_t}} + \abr{\barL_{\barepsilon_{\bart}, \barnu_{\bart}, \eta}^{t} - \mL_{\barepsilon_{\bart}, \barnu_{\bart}, \eta}^{t}}\\
\leq & \baralpha_t\beta(\bnabla\mL_{\barepsilon_{\bart}, \barnu_{\bart}, \eta}^t)^T\barDelta_t + \abr{\barL_{\barepsilon_{\bart}, \barnu_{\bart}, \eta}^{s_t} - \mL_{\barepsilon_{\bart}, \barnu_{\bart}, \eta}^{s_t}} + \abr{\barL_{\barepsilon_{\bart}, \barnu_{\bart}, \eta}^{t} - \mL_{\barepsilon_{\bart}, \barnu_{\bart}, \eta}^{t}}\\
\stackrel{\mathclap{\eqref{cond:decrease:trust}}}{\leq} &\; \frac{4\baralpha_t\beta}{5}(\bnabla\mL_{\barepsilon_{\bart}, \barnu_{\bart}, \eta}^t)^T\barDelta_t - \frac{\bardelta_t}{5} + \abr{\barL_{\barepsilon_{\bart}, \barnu_{\bart}, \eta}^{s_t} - \mL_{\barepsilon_{\bart}, \barnu_{\bart}, \eta}^{s_t}} + \abr{\barL_{\barepsilon_{\bart}, \barnu_{\bart}, \eta}^{t} - \mL_{\barepsilon_{\bart}, \barnu_{\bart}, \eta}^{t}}\\
\stackrel{\mathclap{\eqref{pequ:34}}}{\leq} &\; - \frac{\baralpha_t\beta(\gamma_{B}\wedge\eta)}{5}\nbr{\left(\begin{smallmatrix}
\barDelta\bx_t\\
J_t\bnabla_{\bx}\mL_t\\
G_t\bnabla_{\bx}\mL_t + \Pi_c(\diag^2(g_t)\blambda_t)
\end{smallmatrix}\right) }^2 - \frac{\bardelta_t}{5} \\
&\quad + \abr{\barL_{\barepsilon_{\bart}, \barnu_{\bart}, \eta}^{s_t} - \mL_{\barepsilon_{\bart}, \barnu_{\bart}, \eta}^{s_t}} + \abr{\barL_{\barepsilon_{\bart}, \barnu_{\bart}, \eta}^{t} - \mL_{\barepsilon_{\bart}, \barnu_{\bart}, \eta}^{t}}.
\end{align*}
By Line 20 of Algorithm \ref{alg:ASto}, $\bardelta_{t+1} - \bardelta_t = (\rho-1)\bardelta_t$, while \eqref{pequ:52} still holds. By the condition of $\omega$ in \eqref{pequ:56} and \eqref{pequ:53}, we know that under \eqref{cond:omega} (which implies \eqref{pequ:53}),
\begin{multline}\label{pequ:70}
\Theta_{\omega}^{t+1} - \Theta_{\omega}^t \leq -\frac{\omega\beta(\gamma_{B}\wedge\eta)}{36}\cdot\baralpha_t\nbr{\left(\begin{smallmatrix}
\barDelta\bx_t\\
J_t\bnabla_{\bx}\mL_t\\
G_t\bnabla_{\bx}\mL_t + \Pi_c(\diag^2(g_t)\blambda_t)
\end{smallmatrix} \right)}^2 +\rho(1-\omega)\baralpha_t\|\nabla\mL_{\barepsilon_{\bart}, \barnu_{\bart}, \eta}^t\|^2 \\
+ \omega\cbr{\abr{\barL_{\barepsilon_{\bart}, \barnu_{\bart}, \eta}^{s_t} - \mL_{\barepsilon_{\bart}, \barnu_{\bart}, \eta}^{s_t}} + \abr{\barL_{\barepsilon_{\bart}, \barnu_{\bart}, \eta}^{t} - \mL_{\barepsilon_{\bart}, \barnu_{\bart}, \eta}^{t}} }  - \frac{\omega}{36}\bardelta_t.
\end{multline}
\noindent\textbf{Case 2a, unreliable step, $\cDelta_t = \barDelta_t$}. We have
\begin{align*}
\mL_{\barepsilon_{\bart}, \barnu_{\bart}, \eta}^{t+1} - & \mL_{\barepsilon_{\bart}, \barnu_{\bart}, \eta}^t \\
\leq & \barL_{\barepsilon_{\bart}, \barnu_{\bart}, \eta}^{s_t} - \barL_{\barepsilon_{\bart}, \barnu_{\bart}, \eta}^t + \abr{\barL_{\barepsilon_{\bart}, \barnu_{\bart}, \eta}^{s_t} - \mL_{\barepsilon_{\bart}, \barnu_{\bart}, \eta}^{s_t}} + \abr{\barL_{\barepsilon_{\bart}, \barnu_{\bart}, \eta}^{t} - \mL_{\barepsilon_{\bart}, \barnu_{\bart}, \eta}^{t}}\\
\stackrel{\mathclap{\eqref{cond:armijo}}}{\leq} &\; \baralpha_t\beta(\bnabla\mL_{\barepsilon_{\bart}, \barnu_{\bart}, \eta}^t)^T\barDelta_t+ \abr{\barL_{\barepsilon_{\bart}, \barnu_{\bart}, \eta}^{s_t} - \mL_{\barepsilon_{\bart}, \barnu_{\bart}, \eta}^{s_t}} + \abr{\barL_{\barepsilon_{\bart}, \barnu_{\bart}, \eta}^{t} - \mL_{\barepsilon_{\bart}, \barnu_{\bart}, \eta}^{t}}\\
\stackrel{\mathclap{\eqref{pequ:34}}}{\leq} &\; - \frac{\baralpha_t\beta(\gamma_{B}\wedge\eta)}{4}\nbr{\left(\begin{smallmatrix}
\barDelta\bx_t\\
J_t\bnabla_{\bx}\mL_t\\
G_t\bnabla_{\bx}\mL_t + \Pi_c(\diag^2(g_t)\blambda_t)
\end{smallmatrix}\right) }^2 + \abr{\barL_{\barepsilon_{\bart}, \barnu_{\bart}, \eta}^{s_t} - \mL_{\barepsilon_{\bart}, \barnu_{\bart}, \eta}^{s_t}} \\
&\quad + \abr{\barL_{\barepsilon_{\bart}, \barnu_{\bart}, \eta}^{t} - \mL_{\barepsilon_{\bart}, \barnu_{\bart}, \eta}^{t}}.
\end{align*}
By Line 22 of Algorithm \ref{alg:ASto}, $\bardelta_{t+1} - \bardelta_t = -(1-1/\rho)\bardelta_t$, while \eqref{pequ:52} still holds. Thus, under~\eqref{cond:omega}, 
\begin{multline}\label{pequ:71}
\Theta_{\omega}^{t+1} - \Theta_{\omega}^t \leq -\frac{\omega\beta(\gamma_{B}\wedge\eta)}{36}\cdot\baralpha_t\nbr{\left(\begin{smallmatrix}
\barDelta\bx_t\\
J_t\bnabla_{\bx}\mL_t\\
G_t\bnabla_{\bx}\mL_t + \Pi_c(\diag^2(g_t)\blambda_t)
\end{smallmatrix} \right)}^2  - \frac{1}{2}\rbr{1-\omega}\rbr{1-\frac{1}{\rho}}\bardelta_t\\
+ \omega\cbr{\abr{\barL_{\barepsilon_{\bart}, \barnu_{\bart}, \eta}^{s_t} - \mL_{\barepsilon_{\bart}, \barnu_{\bart}, \eta}^{s_t}} + \abr{\barL_{\barepsilon_{\bart}, \barnu_{\bart}, \eta}^{t} - \mL_{\barepsilon_{\bart}, \barnu_{\bart}, \eta}^{t}} }
+\rho(1-\omega)\baralpha_t\|\nabla\mL_{\barepsilon_{\bart}, \barnu_{\bart}, \eta}^t\|^2.
\end{multline}
\noindent\textbf{Case 3a, unsuccessful step, $\cDelta_t = \barDelta_t$}. In this case, \eqref{pequ:55} holds. Combining \eqref{pequ:70}, \eqref{pequ:71}, and \eqref{pequ:55}, we obtain
\begin{equation}\label{pequ:72}
\Theta_{\omega}^{t+1} - \Theta_{\omega}^t \leq \omega\cbr{\abr{\barL_{\barepsilon_{\bart}, \barnu_{\bart}, \eta}^{s_t} - \mL_{\barepsilon_{\bart}, \barnu_{\bart}, \eta}^{s_t}} + \abr{\barL_{\barepsilon_{\bart}, \barnu_{\bart}, \eta}^{t} - \mL_{\barepsilon_{\bart}, \barnu_{\bart}, \eta}^{t}} }
+\rho(1-\omega)\baralpha_t\|\nabla\mL_{\barepsilon_{\bart}, \barnu_{\bart}, \eta}^t\|^2.
\end{equation}
\noindent\textbf{Case 1b, reliable step, $\cDelta_t = \hDelta_t$}. We have
\begin{align*}
\mL_{\barepsilon_{\bart}, \barnu_{\bart}, \eta}^{t+1} - & \mL_{\barepsilon_{\bart}, \barnu_{\bart}, \eta}^t \\
\leq & \barL_{\barepsilon_{\bart}, \barnu_{\bart}, \eta}^{s_t} - \barL_{\barepsilon_{\bart}, \barnu_{\bart}, \eta}^t + \abr{\barL_{\barepsilon_{\bart}, \barnu_{\bart}, \eta}^{s_t} - \mL_{\barepsilon_{\bart}, \barnu_{\bart}, \eta}^{s_t}} + \abr{\barL_{\barepsilon_{\bart}, \barnu_{\bart}, \eta}^{t} - \mL_{\barepsilon_{\bart}, \barnu_{\bart}, \eta}^{t}}\\
\leq & \baralpha_t\beta(\bnabla\mL_{\barepsilon_{\bart}, \barnu_{\bart}, \eta}^t)^T\hDelta_t+ \abr{\barL_{\barepsilon_{\bart}, \barnu_{\bart}, \eta}^{s_t} - \mL_{\barepsilon_{\bart}, \barnu_{\bart}, \eta}^{s_t}} + \abr{\barL_{\barepsilon_{\bart}, \barnu_{\bart}, \eta}^{t} - \mL_{\barepsilon_{\bart}, \barnu_{\bart}, \eta}^{t}}\\
\stackrel{\mathclap{\eqref{cond:decrease:trust}}}{\leq} &\; \frac{\baralpha_t\beta}{2}(\bnabla\mL_{\barepsilon_{\bart}, \barnu_{\bart}, \eta}^t)^T\hDelta_t - \frac{\bardelta_t}{2} + \abr{\barL_{\barepsilon_{\bart}, \barnu_{\bart}, \eta}^{s_t} - \mL_{\barepsilon_{\bart}, \barnu_{\bart}, \eta}^{s_t}} + \abr{\barL_{\barepsilon_{\bart}, \barnu_{\bart}, \eta}^{t} - \mL_{\barepsilon_{\bart}, \barnu_{\bart}, \eta}^{t}}\\
\stackrel{\mathclap{\eqref{cond:alter:dir}}}{\leq} &  - \frac{\baralpha_t\beta}{2\chi_{u}}\|\bnabla\mL_{\barepsilon_{\bart}, \barnu_{\bart}, \eta}^t\|^2   - \frac{\bardelta_t}{2} + \abr{\barL_{\barepsilon_{\bart}, \barnu_{\bart}, \eta}^{s_t} - \mL_{\barepsilon_{\bart}, \barnu_{\bart}, \eta}^{s_t}} + \abr{\barL_{\barepsilon_{\bart}, \barnu_{\bart}, \eta}^{t} - \mL_{\barepsilon_{\bart}, \barnu_{\bart}, \eta}^{t}}.
\end{align*}
By Line 20 of Algorithm \ref{alg:ASto}, $\bardelta_{t+1} - \bardelta_t = (\rho-1)\bardelta_t$, while \eqref{pequ:60} still holds. By the condition of $\omega$ in \eqref{pequ:62}, we know that under \eqref{cond:omega} (which implies \eqref{pequ:62}), 
\begin{multline}\label{pequ:73}
\Theta_{\omega}^{t+1}-\Theta_{\omega}^t \leq -\frac{\omega\beta}{12\chi_{u}}\cdot\baralpha_t\|\bnabla\mL_{\barepsilon_{\bart}, \barnu_{\bart}, \eta}^t\|^2 + \omega\cbr{\abr{\barL_{\barepsilon_{\bart}, \barnu_{\bart}, \eta}^{s_t} - \mL_{\barepsilon_{\bart}, \barnu_{\bart}, \eta}^{s_t}} + \abr{\barL_{\barepsilon_{\bart}, \barnu_{\bart}, \eta}^{t} - \mL_{\barepsilon_{\bart}, \barnu_{\bart}, \eta}^{t}} }\\ +\rho(1-\omega)\baralpha_t\|\nabla\mL_{\barepsilon_{\bart}, \barnu_{\bart}, \eta}^t\|^2 - \frac{\omega}{12}\bardelta_t.
\end{multline}
\noindent\textbf{Case 2b, unreliable step, $\cDelta_t = \hDelta_t$}. We have
\begin{align*}
\mL_{\barepsilon_{\bart}, \barnu_{\bart}, \eta}^{t+1} - & \mL_{\barepsilon_{\bart}, \barnu_{\bart}, \eta}^t \\
\leq & \barL_{\barepsilon_{\bart}, \barnu_{\bart}, \eta}^{s_t} - \barL_{\barepsilon_{\bart}, \barnu_{\bart}, \eta}^t + \abr{\barL_{\barepsilon_{\bart}, \barnu_{\bart}, \eta}^{s_t} - \mL_{\barepsilon_{\bart}, \barnu_{\bart}, \eta}^{s_t}} + \abr{\barL_{\barepsilon_{\bart}, \barnu_{\bart}, \eta}^{t} - \mL_{\barepsilon_{\bart}, \barnu_{\bart}, \eta}^{t}}\\
\leq & \baralpha_t\beta(\bnabla\mL_{\barepsilon_{\bart}, \barnu_{\bart}, \eta}^t)^T\hDelta_t+ \abr{\barL_{\barepsilon_{\bart}, \barnu_{\bart}, \eta}^{s_t} - \mL_{\barepsilon_{\bart}, \barnu_{\bart}, \eta}^{s_t}} + \abr{\barL_{\barepsilon_{\bart}, \barnu_{\bart}, \eta}^{t} - \mL_{\barepsilon_{\bart}, \barnu_{\bart}, \eta}^{t}}\\
\stackrel{\mathclap{\eqref{cond:alter:dir}}}{\leq} & -\frac{\baralpha_t\beta}{\chi_{u}} \|\bnabla\mL_{\barepsilon_{\bart}, \barnu_{\bart}, \eta}^t\|^2 + \abr{\barL_{\barepsilon_{\bart}, \barnu_{\bart}, \eta}^{s_t} - \mL_{\barepsilon_{\bart}, \barnu_{\bart}, \eta}^{s_t}} + \abr{\barL_{\barepsilon_{\bart}, \barnu_{\bart}, \eta}^{t} - \mL_{\barepsilon_{\bart}, \barnu_{\bart}, \eta}^{t}}.
\end{align*}
By Line 22 of Algorithm \ref{alg:ASto}, $\bardelta_{t+1} - \bardelta_t = -(1-1/\rho)\bardelta_t$, while \eqref{pequ:60} still holds. Thus, under~\eqref{cond:omega}, 
\begin{multline}\label{pequ:74}
\Theta_{\omega}^{t+1}-\Theta_{\omega}^t \leq  -\frac{\omega\beta}{12\chi_{u}}\cdot\baralpha_t\|\bnabla\mL_{\barepsilon_{\bart}, \barnu_{\bart}, \eta}^t\|^2 + \omega\cbr{\abr{\barL_{\barepsilon_{\bart}, \barnu_{\bart}, \eta}^{s_t} - \mL_{\barepsilon_{\bart}, \barnu_{\bart}, \eta}^{s_t}} + \abr{\barL_{\barepsilon_{\bart}, \barnu_{\bart}, \eta}^{t} - \mL_{\barepsilon_{\bart}, \barnu_{\bart}, \eta}^{t}} }\\
+\rho(1-\omega) \baralpha_t\|\nabla\mL_{\barepsilon_{\bart}, \barnu_{\bart}, \eta}^t\|^2 
- \frac{1}{2}\rbr{1-\omega}\rbr{1-\frac{1}{\rho}}\bardelta_t.
\end{multline}
\noindent\textbf{Case 3b, unsuccessful step, $\cDelta_t = \hDelta_t$}. In this case, \eqref{pequ:55} holds. Combining \eqref{pequ:73}, \eqref{pequ:74}, and \eqref{pequ:55}, we note that \eqref{pequ:72} holds as well. Thus, \eqref{pequ:72} holds for all six cases. This completes the proof.

\subsection{Proof of Theorem \ref{thm:3}}\label{pf:thm:3}

We suppose there are infinite many successful steps. Otherwise, $\baralpha_t$ decreases to zero (cf. Line 25 of Algorithm \ref{alg:ASto}) and the argument holds trivially. We use $\bart<t_1<t_2<\ldots$ to denote the subsequence with $t_i-1$, $\forall i\geq 1$, being a successful step. By Lemma \ref{lem:3}, Lemma \ref{lem:9}(b), and \eqref{pequ:38}, there exist $\Upsilon_1, \Upsilon_2>0$ such that for any $i\geq 1$,
\begin{equation}
R_{t_i} \stackrel{\text{Lem. }\ref{lem:3}}{\leq} \Upsilon_1\nbr{\begin{pmatrix}
\nabla_{\bx}\mL_{t_i}\\
c_{t_i}\\
\bw_{\barepsilon_{\bart}, \barnu_{\bart}}^{t_i}
\end{pmatrix}} \stackrel{\text{Lem. } \ref{lem:9}(b)}{\leq}\Upsilon_2\cbr{\|\nabla\mL_{\barepsilon_{\bart}, \barnu_{\bart}, \eta}^{t_i}\| + \nbr{\begin{pmatrix}
c_{t_i}\\
\bw_{\barepsilon_{\bart}, \barnu_{\bart}}^{t_i}
\end{pmatrix}}}.
\end{equation}
Since $t_i\geq \bart+1$, two parameters $\barepsilon_{\bart}, \barnu_{\bart}$ are fixed conditional on any $\sigma$-algebra $\mF\supseteq\mF_{\bart}$. Thus, for any $i\geq 1$,
\begin{align*}
& \nbr{\begin{pmatrix}
c_{t_i}\\
\bw_{\barepsilon_{\bart}, \barnu_{\bart}}^{t_i}
\end{pmatrix}} =  \mE\sbr{\nbr{\begin{pmatrix}
c_{t_i}\\
\bw_{\barepsilon_{\bart}, \barnu_{\bart}}^{t_i}
\end{pmatrix}} \mid \mF_{t_i-1}} \\
& = \mE\sbr{\nbr{\begin{pmatrix}
c_{t_i}\\
\bw_{\barepsilon_{\bart}, \barnu_{\bart}}^{t_i}
\end{pmatrix}}\b1_{\chi_{err}\|\bnabla\mL_{\barepsilon_{\bart}, \barnu_{\bart}, \eta}^{t_i}\| \leq \barR_{t_i}} \mid \mF_{t_i-1}} +\mE\sbr{\nbr{\begin{pmatrix}
c_{t_i}\\
\bw_{\barepsilon_{\bart}, \barnu_{\bart}}^{t_i}
\end{pmatrix}}\b1_{\barR_{t_i}< \chi_{err}\|\bnabla\mL_{\barepsilon_{\bart}, \barnu_{\bart}, \eta}^{t_i}\|} \mid \mF_{t_i-1}} \\
& \stackrel{\mathclap{\eqref{cond:bound:fes:error}}}{\leq }
\mE\sbr{\chi_{err}\|\bnabla\mL_{\barepsilon_{\bart}, \barnu_{\bart}, \eta}^{t_i}\|\cdot  \b1_{\chi_{err}\|\bnabla\mL_{\barepsilon_{\bart}, \barnu_{\bart}, \eta}^{t_i}\| \leq \barR_{t_i}} \mid \mF_{t_i-1}} \\
& \quad + (\barepsilon_{\bart}q_{\barnu_{\bart}}^{t_i}\vee 1) \mE\sbr{\nbr{\begin{pmatrix}
c_{t_i}\\
\max\{g_{t_i},-\blambda_{t_i}\}
\end{pmatrix}}\cdot\b1_{\barR_{t_i}< \chi_{err}\|\bnabla\mL_{\barepsilon_{\bart}, \barnu_{\bart}, \eta}^{t_i}\|} \mid \mF_{t_i-1} } \quad \text{(also use Lemma \ref{lem:3})}\\
& \stackrel{\mathclap{\eqref{pequ:38}}}{\leq}\chi_{err}\mE\sbr{\|\bnabla\mL_{\barepsilon_{\bart}, \barnu_{\bart}, \eta}^{t_i}\| \mid \mF_{t_i-1}} + (\barepsilon_0\tnu\vee 1)\mE\sbr{\barR_{t_i} \cdot\b1_{\barR_{t_i}< \chi_{err}\|\bnabla\mL_{\barepsilon_{\bart}, \barnu_{\bart}, \eta}^{t_i}\|} \mid \mF_{t_i-1}} \\
& \leq \cbr{1 + (\barepsilon_0\tnu\vee 1)}\chi_{err}\mE\sbr{\|\bnabla\mL_{\barepsilon_{\bart}, \barnu_{\bart}, \eta}^{t_i}\| \mid \mF_{t_i-1}}\\
& \leq \cbr{1 + (\barepsilon_0\tnu\vee 1)}\chi_{err}\cbr{\|\nabla\mL_{\barepsilon_{\bart}, \barnu_{\bart}, \eta}^{t_i}\| + \mE\sbr{\|\barDelta(\nabla\mL_{\eta}^{t_i})\| \mid \mF_{t_i-1} } }\\
& \stackrel{\mathclap{\eqref{cond:grad:var}}}{\leq}\cbr{1 + (\barepsilon_0\tnu\vee 1)}\chi_{err} \cbr{\|\nabla\mL_{\barepsilon_{\bart}, \barnu_{\bart}, \eta}^{t_i}\| + \chi_{grad} \sqrt{\bardelta_{t_i}/\baralpha_{t_i}}}.
\end{align*}
Combining the above two displays, we know there exists $\Upsilon_3>0$ such that
\begin{equation*}
R_{t_i} \leq \Upsilon_3(\chi_{grad} + 1)\cbr{\|\nabla\mL_{\barepsilon_{\bart}, \barnu_{\bart}, \eta}^{t_i}\| + \sqrt{\bardelta_{t_i}/\baralpha_{t_i} } },
\end{equation*}
which implies
\begin{equation}\label{bound:R:ti}
\baralpha_{t_i}R_{t_i}^2 \leq 2\Upsilon_3^2(\chi_{grad}+1)^2\cbr{\baralpha_{t_i}\|\nabla\mL_{\barepsilon_{\bart}, \barnu_{\bart}, \eta}^{t_i}\|^2 + \bardelta_{t_i}}.
\end{equation}
On the other hand, by Theorem \ref{thm:2}, we sum up the error recursion for $t\geq \bart +1$, take conditional expectation on $\mF_{\bart}$, and have
\begin{align}\label{pequ:100}
\sum_{t = \bart+1}^{\infty}&\mE[\baralpha_t\|\nabla\mL_{\barepsilon_{\bart}, \barnu_{\bart}, \eta}^t\|^2 + \bardelta_t \mid \mF_{\bart}] \nonumber\\
\leq& \frac{4\rho}{(1-p_{grad})(1-p_{f})(1-\omega)(\rho-1)}\sum_{t=\bart+1}^{\infty}\mE[\Theta_{\omega}^t\mid \mF_{\bart}] - \mE\sbr{\Theta_{\omega}^{t+1}\mid \mF_{\bart}} \nonumber\\
\leq &\frac{4\rho}{(1-p_{grad})(1-p_{f})(1-\omega)(\rho-1)}\rbr{\Theta_{\omega}^{\bart+1} - \min_{\mX\times \mM\times\Lambda} \omega\mL_{\barepsilon_{\bart}, \barnu_{\bart}, \eta}}<\infty.
\end{align}
Thus, applying the Fubini's theorem to exchange the summation and expectation, we~know~that $\mE[\limsup\limits_{t\rightarrow\infty} \baralpha_t\|\nabla\mL_{\barepsilon_{\bart}, \barnu_{\bart}, \eta}^t\|^2 + \bardelta_t\mid \mF_{\bart}]=0 $. Since $\baralpha_t\|\nabla\mL_{\barepsilon_{\bart}, \barnu_{\bart}, \eta}^t\|^2 + \bardelta_t$ is non-negative, we further obtain $\baralpha_t\|\nabla\mL_{\barepsilon_{\bart}, \barnu_{\bart}, \eta}^t\|^2 + \bardelta_t \rightarrow 0$ as $t\rightarrow \infty$ almost surely. By \eqref{bound:R:ti}, we have $\baralpha_{t_i}R_{t_i}^2\rightarrow 0$ as $i\rightarrow \infty$. Noting that $\baralpha_t R_t^2\leq \baralpha_{t_i}R_{t_i}^2$ for any $t_i\leq t<t_{i+1}$, we complete the proof.

\subsection{Proof of Theorem \ref{thm:4}}\label{pf:thm:4}

We adapt the proof of \cite[Theorem 4]{Na2022adaptive}. By Theorem \ref{thm:3}, it suffices to show that the ``limsup" of the random stepsize sequence $\{\baralpha_t\}_t$ is lower bounded away from zero. To show this, we define two stepsize sequences as follows. For any $t> \bart+1$, we let
\begin{align*}
\phi_t = &\log(\baralpha_t),\\
\varphi_t = &\min\{\log(c), \b1_{\E_1^{t-1}\cap\E_2^{t-1}}(\log (\rho) + \varphi_{t-1}) + (1 - \b1_{\E_1^{t-1}\cap\E_2^{t-1}})(\varphi_{t-1} - \log(\rho))\},
\end{align*}
and let $\phi_{\bart +1} = \varphi_{\bart+1} = \log(\baralpha_{\bart+1})$. Here, $c$ is a deterministic constant such that
\begin{align*}
c \leq \frac{1-\beta}{\Upsilon_1(\kappa_{grad}+\kappa_{f} + 1)}\wedge\alpha_{max}
\end{align*}
and $c = \rho^{-i}\alpha_{max}$ for some $i>0$. The first constant comes from Lemma \ref{lem:10}. We aim to show $\phi_t\geq \varphi_t$, $\forall t\geq \bart+1$.

First, we note that by the stepsize specification in Lines 18 and 25 of Algorithm \ref{alg:ASto} (Line 13 is not performed since $t\geq \bart + 1$), $\baralpha_t = \rho^{j_t}c$ for some integer $j_t$. Second, we note that $\phi_t$ and $\varphi_t$ are both $\mF_{t-1}$-measurable, that is, they are fixed conditional on $\mF_{t-1}$. Third, we show that $\phi_t\geq \varphi_t$ by induction. Note that $\phi_{\bart+1} = \varphi_{\bart+1}$. Suppose $\phi_t\geq \varphi_t$, we consider the following three cases.

\noindent\textbf{(a)}. If $\phi_t>\log(c)$, then $\phi_t\geq \log(c) + \log(\rho)$. Thus, $\phi_{t+1} \geq \phi_t - \log(\rho) \geq \log (c) \geq \varphi_{t+1}$.

\noindent\textbf{(b)}. If $\phi_t\leq \log(c)$ and $\b1_{\E_1^t\cap\E_2^t} = 1$, then Lemma \ref{lem:10} leads to
\begin{equation*}
\phi_{t+1} =\min\{\log(\alpha_{max}), \phi_t + \log(\rho)\} \geq \min\{\log(c), \varphi_t + \log(\rho)\} = \varphi_{t+1}. 
\end{equation*}

\noindent\textbf{(c)}. If $\phi_t\leq \log(c)$ and $\b1_{\E_1^t\cap\E_2^t} = 0$, then
\begin{equation*}
\phi_{t+1} \geq \phi_t - \log(\rho) \geq\varphi_t - \log(\rho) \geq \varphi_{t+1}.
\end{equation*}
Combining the above three cases, we have $\phi_t\geq \varphi_t$, $\forall t\geq \bart +1$. Note that, conditional on $\mF_{\bart}$, $\{\varphi_t\}_{t\geq \bart+1}$ is a random walk with a maximum and a drift upward (cf. \cite[Example 6.1.2]{Gallager2013Stochastic}). Thus, $\limsup_{t\rightarrow\infty}\varphi_t\geq \log(c)$ almost surely. In particular, we have
\begin{align*}
P\rbr{\limsup_{t\rightarrow\infty} \phi_t\geq \log(c)} & = \sum_{i = 0}^{\infty}\int_{\mF_i} P\rbr{\limsup_{t\rightarrow\infty} \phi_t\geq \log(c) \mid \mF_i, \bart = i} P\rbr{\mF_i, \bart = i}\\
&\hskip-0.2cm\stackrel{\phi_t\geq \varphi_t}{\geq} \sum_{i = 0}^{\infty}\int_{\mF_i} P\rbr{\limsup_{t\rightarrow\infty} \varphi_t\geq \log(c) \mid \mF_i, \bart = i} P\rbr{\mF_i, \bart = i}\\
& = \sum_{i = 0}^{\infty}\int_{\mF_i}P\rbr{\mF_i, \bart = i} \\
& = 1,
\end{align*}
which means that the ``limsup" of $\baralpha_t$ is lower bounded almost surely. Using Theorem \ref{thm:3}, we complete the proof.

\subsection{Proof of Theorem \ref{thm:5}}\label{pf:thm:5}

Suppose $\limsup_{t\rightarrow \infty}R_t = \epsilon>0$. By Theorem \ref{thm:4}, we know there exist two sequences $\{n_i\}_i$ and $\{m_i\}_i$ with $n_i<m_i<n_{i+1}$ for all $i$, such that
\begin{equation*}
R_{n_i}\geq \frac{2\epsilon}{3},\quad R_t\geq \frac{\epsilon}{3},\; t=n_i+1,\ldots,m_i-1,\quad R_{m_i}<\frac{\epsilon}{3}.
\end{equation*}
For each interval $[n_i, m_i]$, we use $\{t_{i,j}\}_{j=1}^{J_i}$ to denote a subsequence within the interval such that $n_i= t_{i,1} <\ldots< t_{i,j} < \ldots<t_{i,J_i} = m_i$ and $t_{i,j}-1$ is a successful step. In other words, $t_{i,j}$ is the first index that we arrive at the new point. Here, we suppose $n_i-1$ is a successful step; that is, the index $n_i$ is the first time we arrive at the point $(\bx_{n_i}, \bmu_{n_i}, \blambda_{n_i})$~(one can always choose $n_i$ to satisfy this condition). We also note that $t_{i,J_i} = m_i$ because $R_{m_i-1}\geq \epsilon/3$ while $R_{m_i}<\epsilon/3$. With these notation, there exist $\Upsilon_1, \Upsilon_2>0$ such that
\begin{align}\label{pequ:101}
\frac{\epsilon}{3} & \leq R_{n_i} - R_{m_i} \leq \sum_{t=n_i}^{m_i-1}\abr{R_{t+1} -  R_t} \leq \sum_{t=n_i}^{m_i-1}\nbr{\left(\begin{smallmatrix}
\nabla_{\bx}\mL_{t+1} - \nabla_{\bx}\mL_t\\
c_{t+1}-c_t\\
\max\{g_{t+1},-\blambda_{t+1}\} - \max\{g_t,-\blambda_t\}
\end{smallmatrix} \right)} \nonumber\\
& \leq \Upsilon_1\sum_{t=n_i}^{m_i-1}\|(\bx_{t+1} - \bx_t,\bmu_{t+1} - \bmu_t, \blambda_{t+1}-\blambda_{t})\| \quad \text{ (due to the Lip-continuity)} \nonumber\\
& = \Upsilon_1 \sum_{j=2}^{J_i}\|(\bx_{t_{i,j}} - \bx_{t_{i,j}-1}, \bmu_{t_{i,j}} - \bmu_{t_{i,j}-1}, \blambda_{t_{i,j}} - \blambda_{t_{i,j}-1})\| \nonumber\\
& = \Upsilon_1\sum_{j=2}^{J_i}\baralpha_{t_{i,j}-1}\|\cDelta_{t_{i,j}-1}\| \leq \Upsilon_2 \sum_{j=2}^{J_i}\baralpha_{t_{i,j}-1} \quad \text{(due to Assumption \ref{ass:3})} \nonumber\\
& \leq \Upsilon_2 \sum_{j=1}^{J_i-1}\baralpha_{t_{i,j}}\quad \text{(due to Line 25 of Algorithm \ref{alg:ASto})}.
\end{align}
Let us define the set $\T = \{t: t-1 \text{ is successful and } R_t\geq \epsilon/3\}$. We can see from \eqref{bound:R:ti} and \eqref{pequ:100} that $\sum_{t\in \T}\baralpha_t<\infty$. This contradicts \eqref{pequ:101} since $\sum_{t\in \T}\baralpha_t \geq \sum_i\sum_{j=1}^{J_i-1}\baralpha_{t_{i,j}}\stackrel{\eqref{pequ:101}}{=}\infty$. Thus, we know $\limsup_{t\rightarrow\infty}R_t=0$; and thus, we complete the proof.

\section{Auxiliary Lemmas}

\begin{lemma}\label{lem:3}
Let $\epsilon, \nu>0$ and $(\bx, \blambda)\in \mT_{\nu} \times \mR^r$. Then
\begin{equation*}
\frac{\|\bw_{\epsilon, \nu}(\bx, \blambda)\|}{\epsilon q_{\nu}(\bx, \blambda) \vee 1} \leq \|\max\{g(\bx), -\blambda\}\| \leq \frac{\|\bw_{\epsilon, \nu}(\bx, \blambda)\|}{\epsilon q_{\nu}(\bx, \blambda) \wedge 1}.
\end{equation*}
\end{lemma}

\begin{proof}
To prove Lemma \ref{lem:3}, we require the following lemma.

\begin{lemma}\label{lem:3.1}
For any two scalars $a, b$ and a scalar $c>0$, $|\max\{a, b\}| \leq \frac{1}{c\wedge 1}|\max\{a, cb\}|$.
\end{lemma}

\begin{proof}
Without loss of generality, we assume $b\neq 0$ and $c\neq 1$. We consider four cases.
	
\noindent{\bf Case 1:} $b>0$, $c<1$. If $a \leq cb< b$, then $|\max\{a, b\}| = b = \frac{1}{c}|\max\{a, cb\}|$. If $cb<a\leq b$, then $|\max\{a, b\}| = b \leq \frac{1}{c}a = \frac{1}{c}|\max\{a, cb\}|$. If $cb<b<a$, then $|\max\{a, b\}| = a\leq \frac{1}{c}|\max\{a, cb\}|$. Thus, the result holds.
	
\noindent{\bf Case 2: } $b>0$, $c>1$. If $a \leq b <cb$, then $|\max\{a, b\} |= b \leq cb = |\max\{a, cb\}|$. If $b<a\leq cb$, then $|\max\{a, b\}| = a \leq cb = |\max\{a, cb\}|$. If $b<cb<a$, then $|\max\{a, b\}| = a = |\max\{a, cb\}|$. Thus, the result holds.
	
\noindent{\bf Case 3:} $b<0$, $c<1$. If $a \leq b < cb$, then $|\max\{a, b\}| = |b| = \frac{1}{c}|\max\{a, cb\}|$. If $b<a \leq cb$, then $|\max\{a, b\}| = |a| \leq |b| = \frac{1}{c}|\max\{a, cb\}|$. If $b<cb< a$, then $|\max\{a, b\}| = |a| \leq \frac{|a|}{c} = \frac{1}{c}|\max\{a, cb\}|$. Thus, the result holds.
	
\noindent{\bf  Case 4:} $b<0$, $c>1$. If $a \leq cb<b$, then $|\max\{a, b\}| = |b| \leq c|b| = |\max\{a, cb\}|$. If $cb< a \leq b$,~then $|\max\{a, b\}| = |b| \leq |a| = |\max\{a, cb\}|$. If $cb<b<a$, then $|\max\{a, b\}| = |a| = |\max\{a, cb\}|$. Thus, the result holds.
	
\noindent Combining the above four cases, we complete the proof.
\end{proof}

Since $\epsilon, \nu>0$, $(\bx, \blambda)\in \mT_{\nu}\times \mR^r$, and $q_\nu(\bx, \blambda)>0$, we have for any $i\in\{1,2,\ldots, r\}$,
\begin{align*}
|(\bw_{\epsilon, \nu}&(\bx, \blambda))_i| =  |\max\{g_i(\bx), -\epsilon q_{\nu}(\bx, \blambda) \blambda_i\}| \leq \frac{1}{\frac{1}{\epsilon q_{\nu}(\bx, \blambda)} \wedge 1}|\max\{g_i(\bx), -\blambda_i\}|\\
= & \rbr{\epsilon q_\nu(\bx, \blambda) \vee 1}\cdot |\max\{g_i(\bx), -\blambda_i\}|
\leq  \frac{\epsilon q_\nu(\bx, \blambda) \vee 1}{\epsilon q_\nu(\bx, \blambda) \wedge 1}|\max\{g_i(\bx), -\epsilon q_{\nu}(\bx, \blambda) \blambda_i\}| \\
= & \frac{\epsilon q_\nu(\bx, \blambda) \vee 1}{\epsilon q_\nu(\bx, \blambda) \wedge 1}\abr{(\bw_{\epsilon, \nu}(\bx, \blambda))_i},
\end{align*}
where both inequalities are from Lemma \ref{lem:3.1}. Taking $\ell_2$ norm on both sides, we finish the~proof.

\end{proof}

\section{Auxiliary Experiments}\label{sec:aux:exp}

We follow the experiments in Section \ref{sec:5} and provide additional results. We first examine three proportions: (1) the proportion of the iterations with failed SQP steps, (2) the proportion of the iterations with unstabilized penalty parameters, (3) the proportion of the iterations with a triggered feasibility error condition. We then investigate a multiplicative noise, and apply the method on an inequality  constrained logistic regression problem.

\begin{figure}[!htp]
\centering     
\subfigure[$C =  2$]{\label{FSC1}\includegraphics[width=37mm]{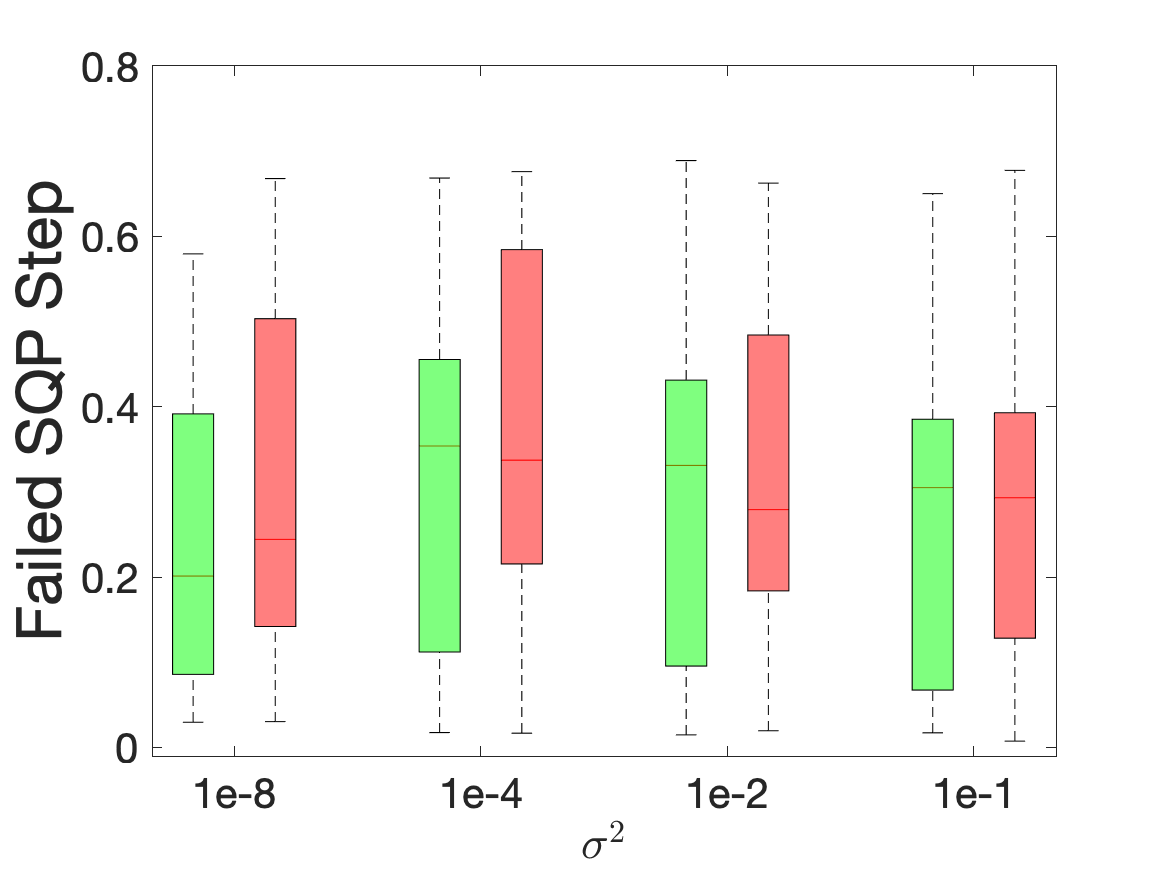}}
\subfigure[$C =  2^3$]{\label{FSC2}\includegraphics[width=37mm]{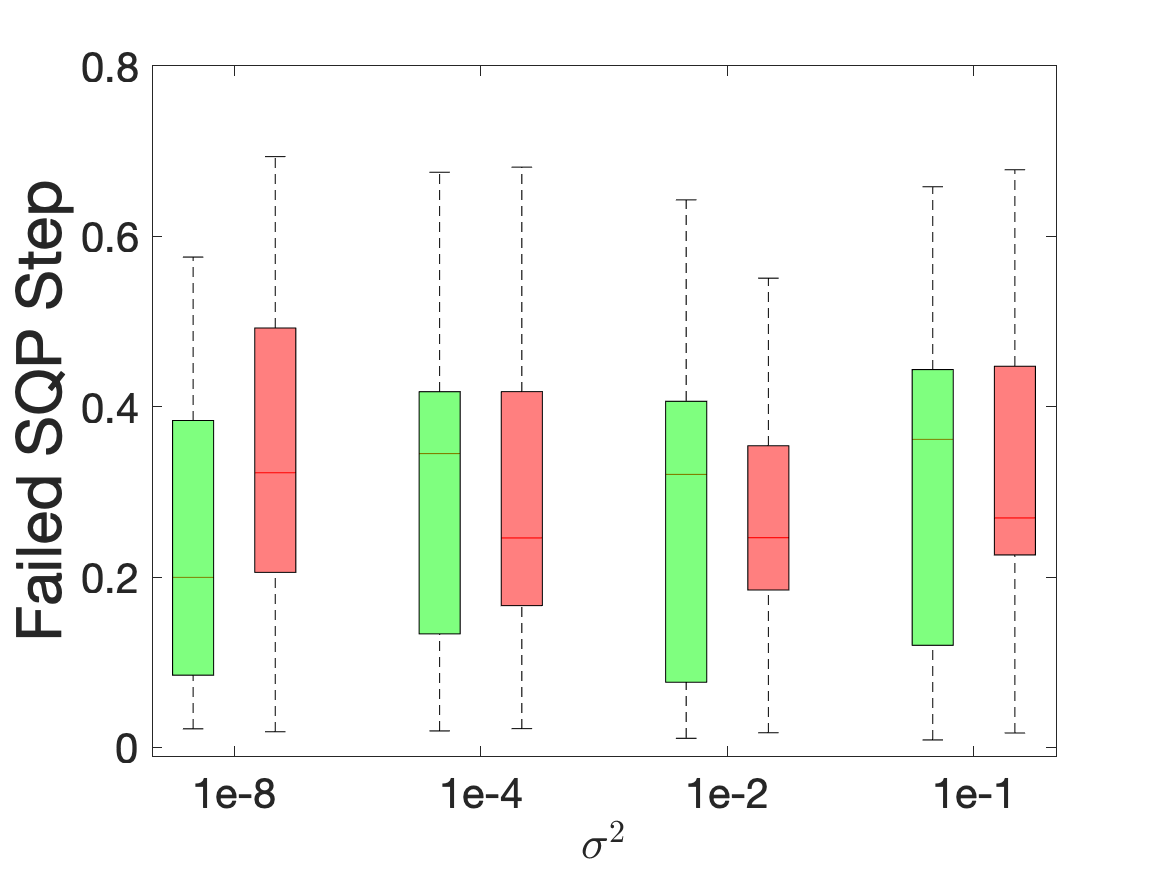}}
\subfigure[$C =  2^6$]{\label{FSC3}\includegraphics[width=37mm]{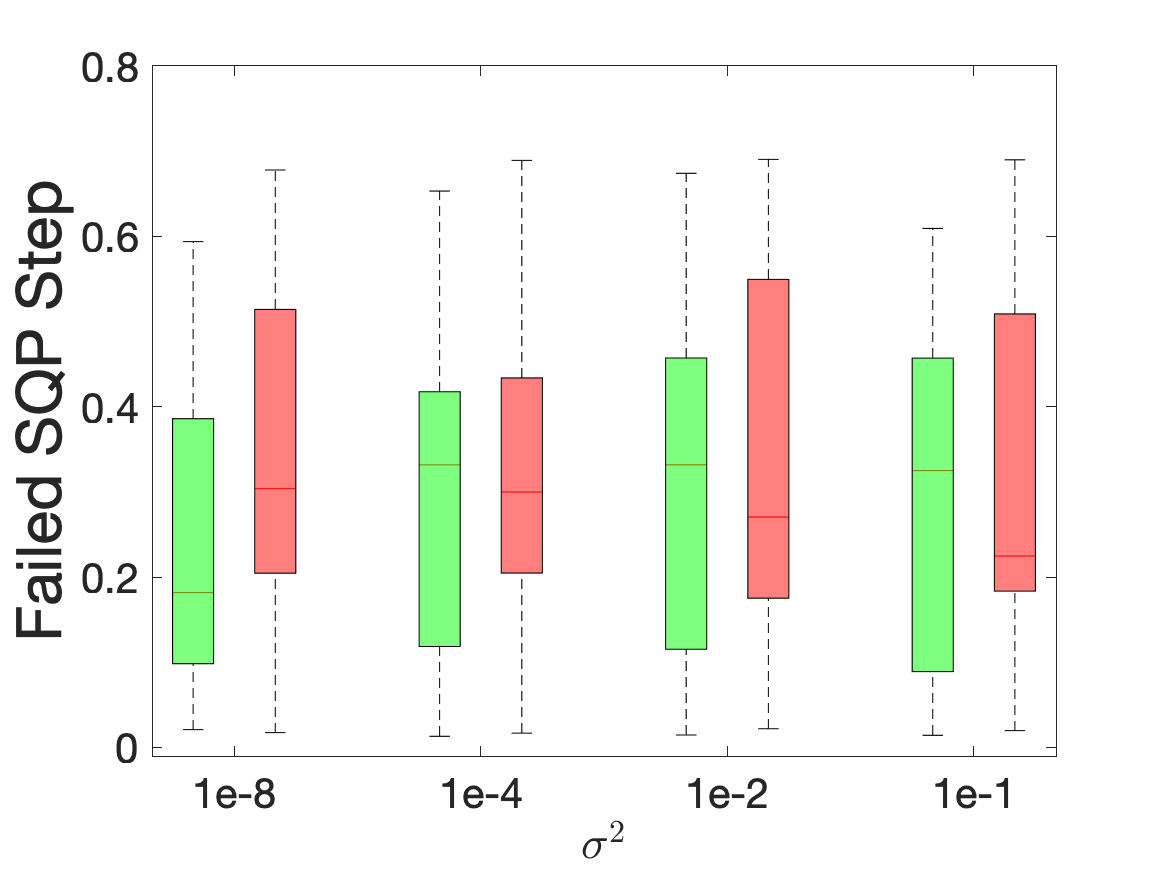}}
	
\subfigure[$\kappa =  2$]{\label{FSK1}\includegraphics[width=37mm]{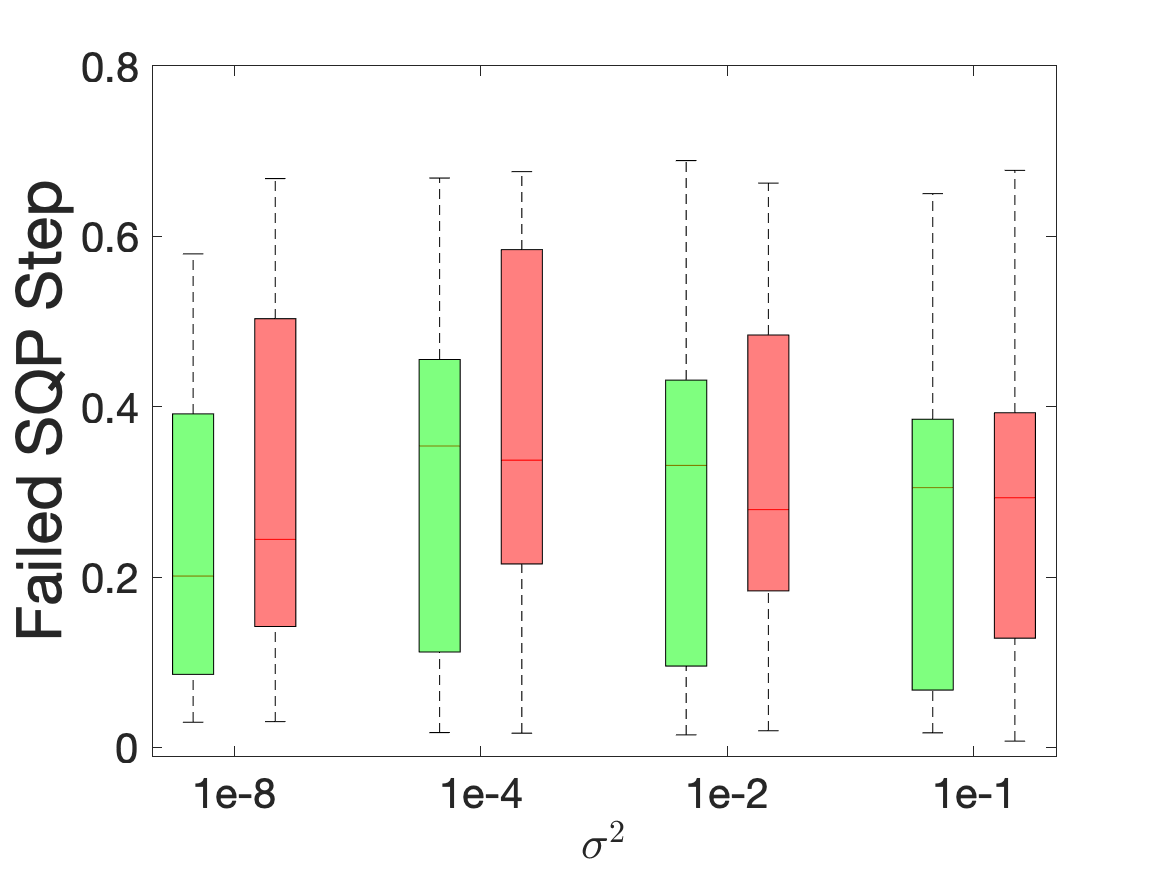}}
\subfigure[$\kappa =  2^3$]{\label{FSK2}\includegraphics[width=37mm]{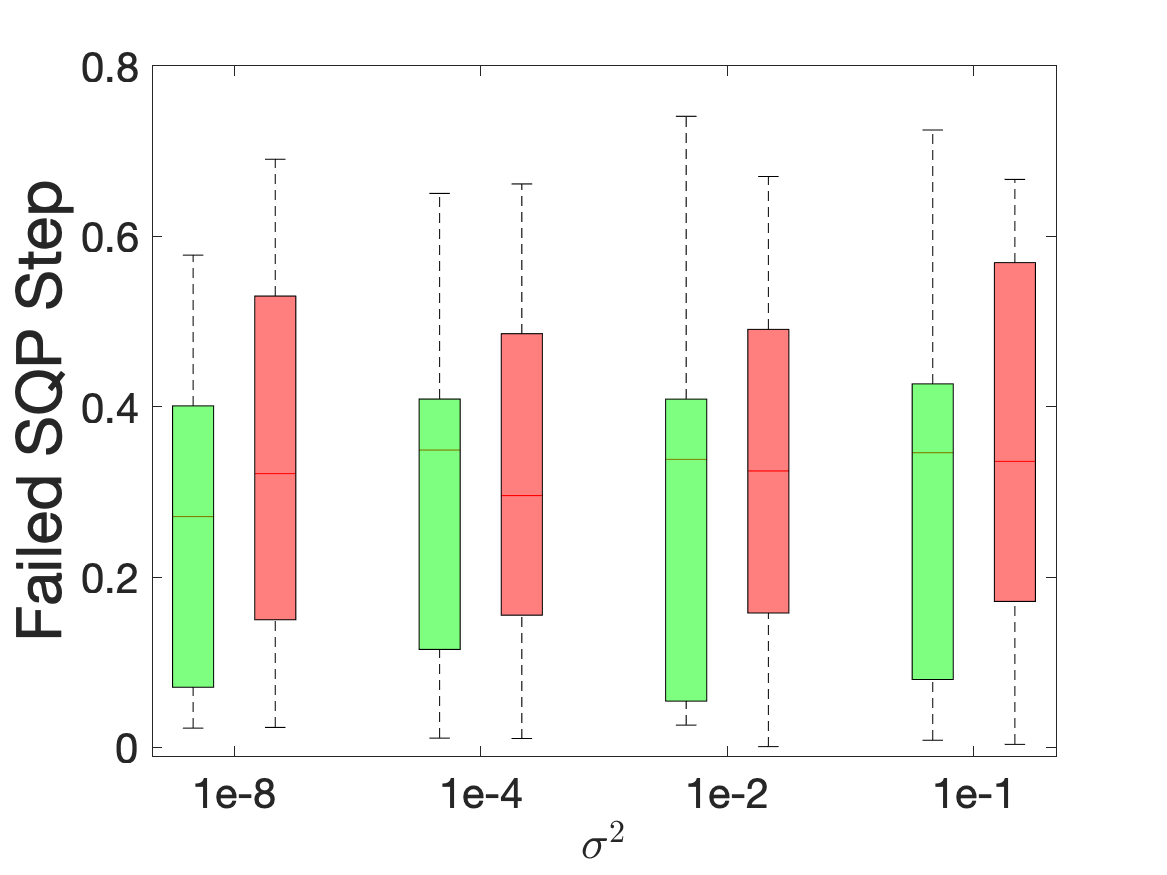}}
\subfigure[$\kappa =  2^6$]{\label{FSK3}\includegraphics[width=37mm]{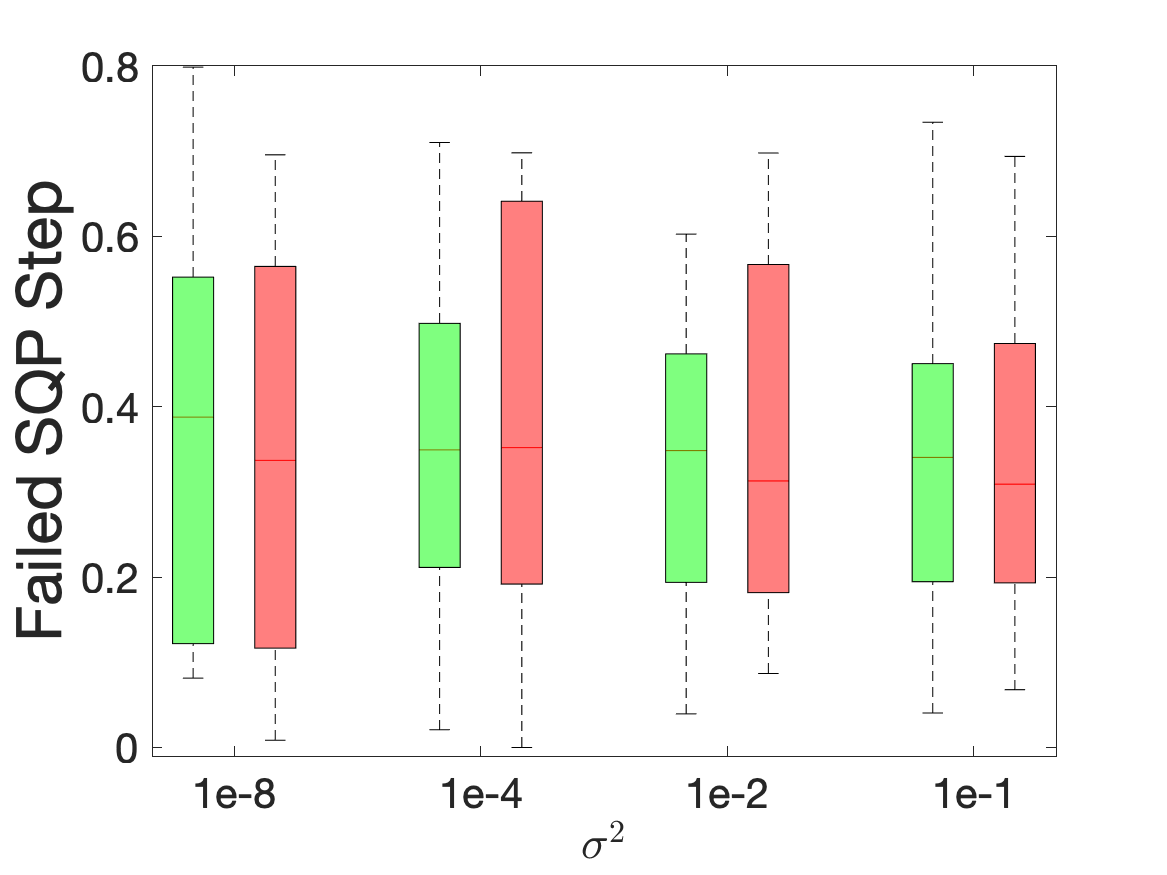}}
	
\subfigure[$\chi_{err} =  1$]{\label{FSChi1}\includegraphics[width=37mm]{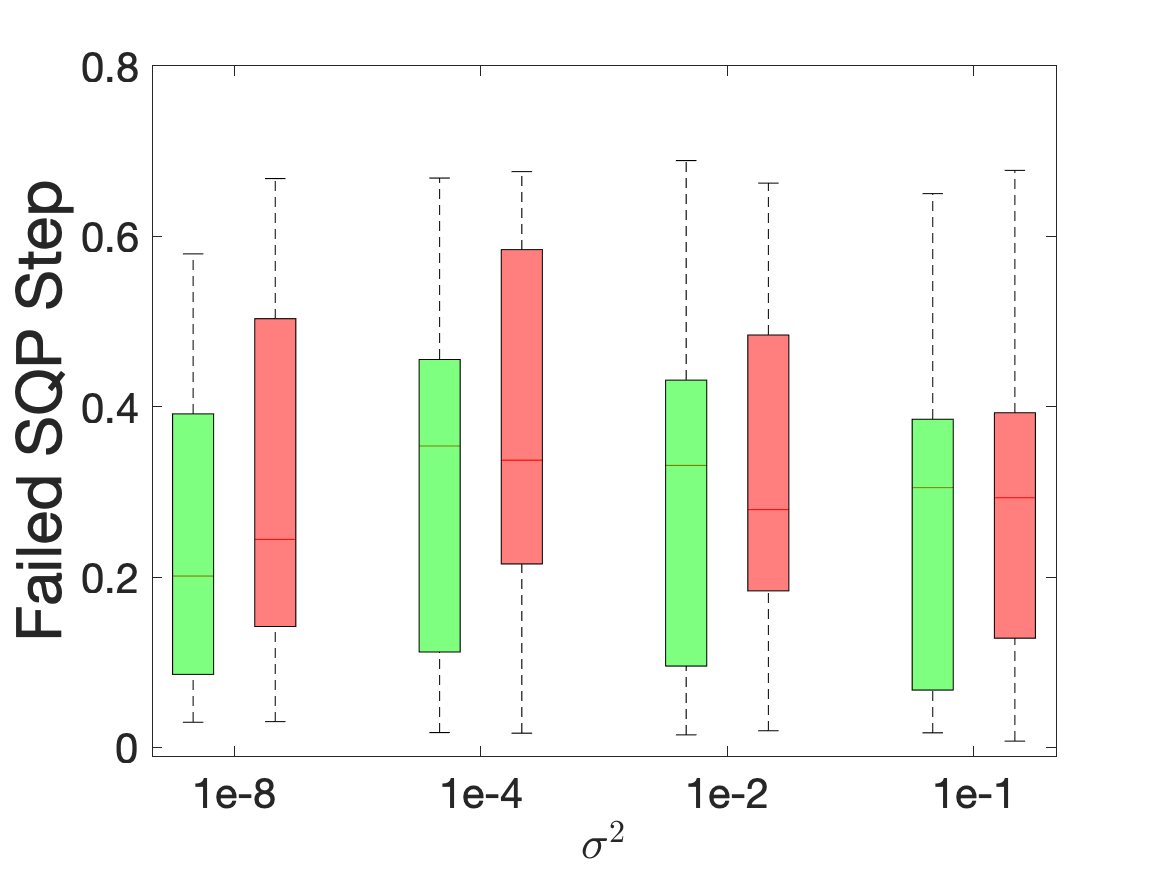}}
\subfigure[$\chi_{err} =  10$]{\label{FSChi2}\includegraphics[width=37mm]{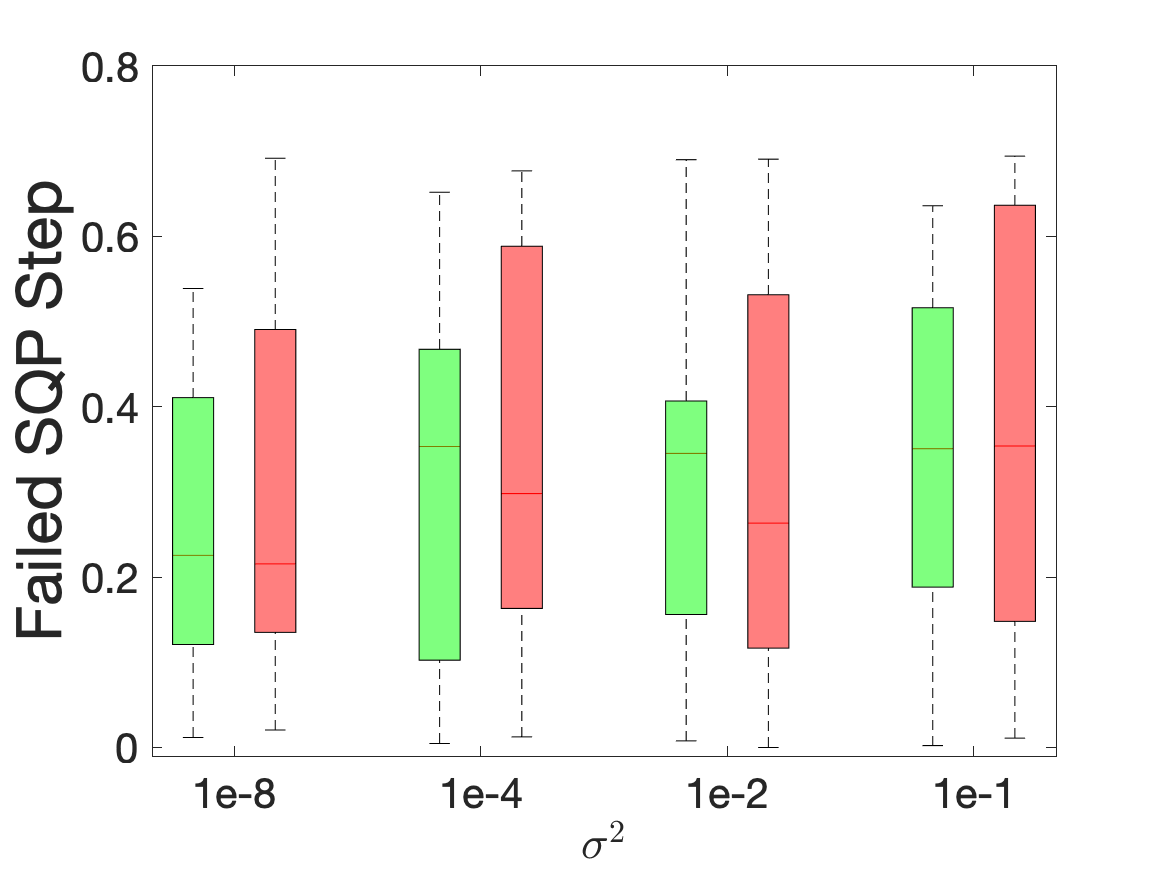}}
\subfigure[$\chi_{err} =  10^2$]{\label{FSChi3}\includegraphics[width=37mm]{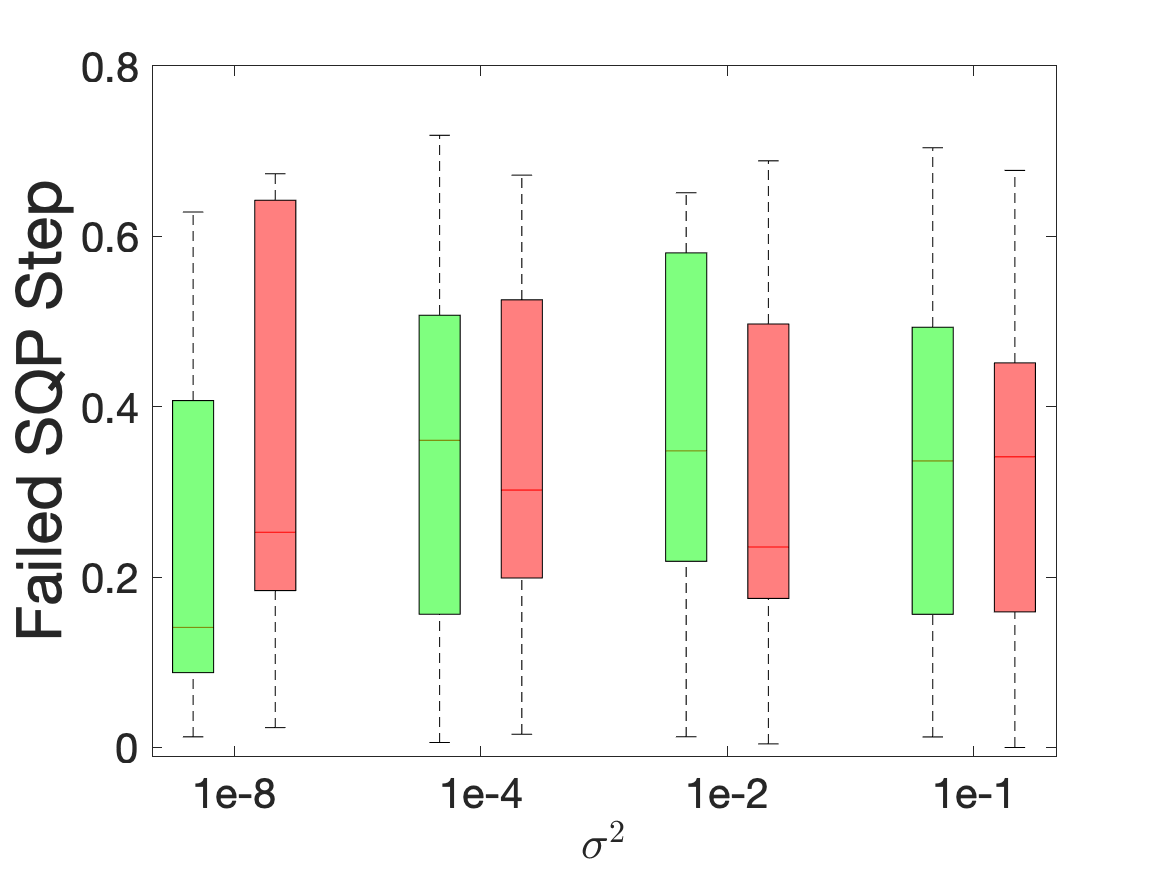}}
	
\includegraphics[width=0.3\textwidth]{Figure/Figure1/legend/KKTlegend}
\caption{Failed SQP step boxplots. Each panel corresponds to a setup of $(C, \kappa, \chi_{err})$. The default values are $C =\kappa= 2$ and $\chi_{err}=1$. When we vary one parameter, the other two are set as default. Thus, the three figures on the left column are the same.}\label{fig:5}
\end{figure}

\vskip4pt
\noindent\textbf{Failed SQP steps.} Figure \ref{fig:5} plots the proportion of the iterations with failed SQP steps. From the figure, we see that the proportion varies from $10\%$ to $60\%$ across the problems,~and AdapNewton tends to have a smaller proportion than AdapGD. Although the proportion~does not have a clear dependency on the variance $\sigma^2$, the noticeable proportion of failed SQP~steps illustrates the differences between equality and inequality constrained problems. As analyzed in Section \ref{sec:2}, the active-set SQP steps may not be informative if the identified active set is very distinct from the true active set. Due to the potential failure of the SQP steps, utilizing a safeguarding direction is critical in achieving the global convergence for the algorithm.

\begin{figure}[!htp]
\centering     
\subfigure[$C =  2$]{\label{NSC1}\includegraphics[width=37mm]{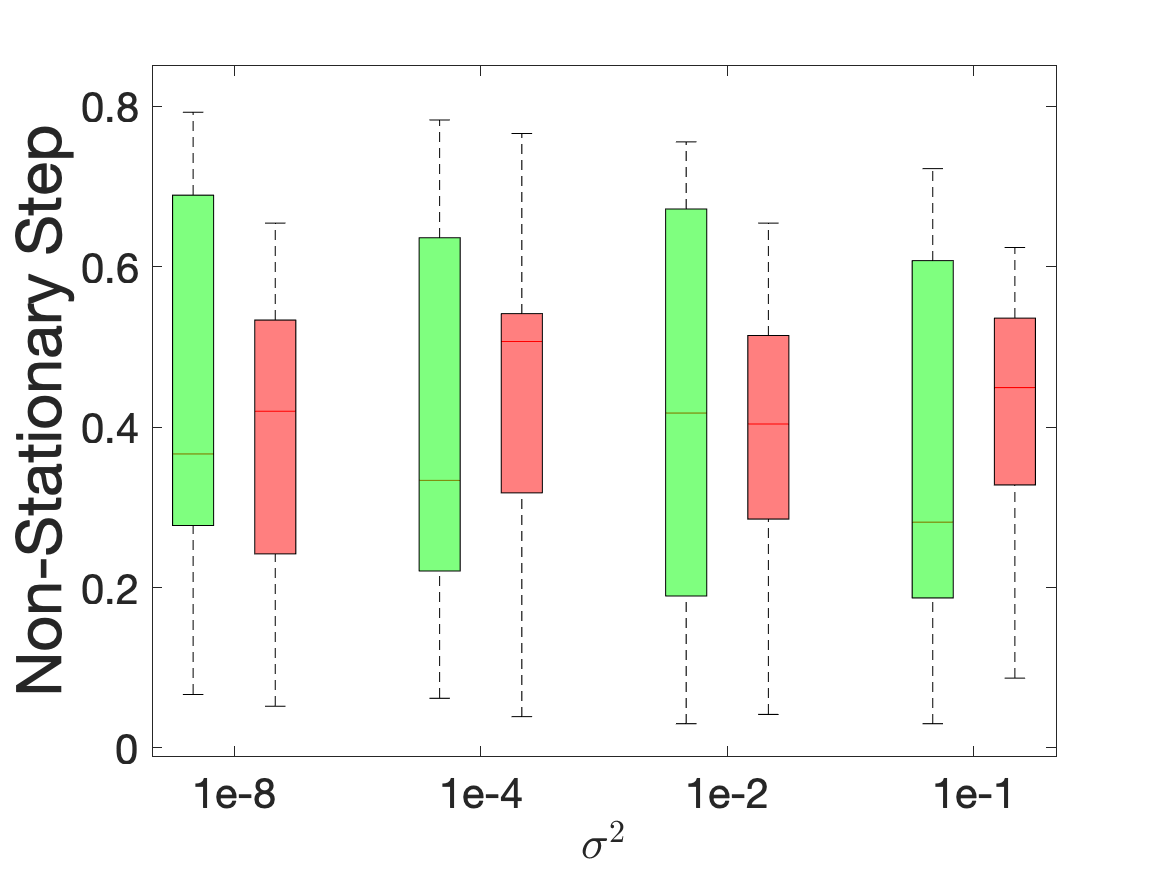}}
\subfigure[$C =  2^3$]{\label{NSC2}\includegraphics[width=37mm]{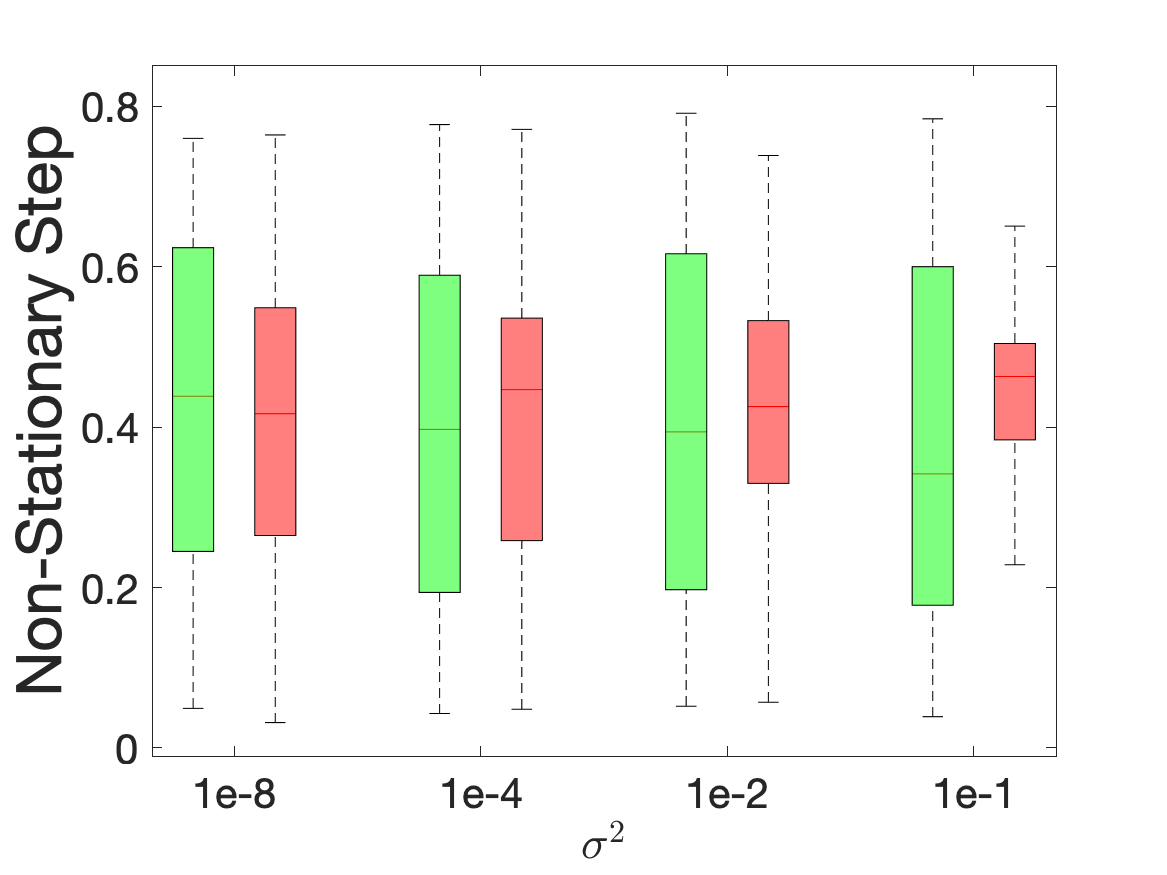}}
\subfigure[$C =  2^6$]{\label{NSC3}\includegraphics[width=37mm]{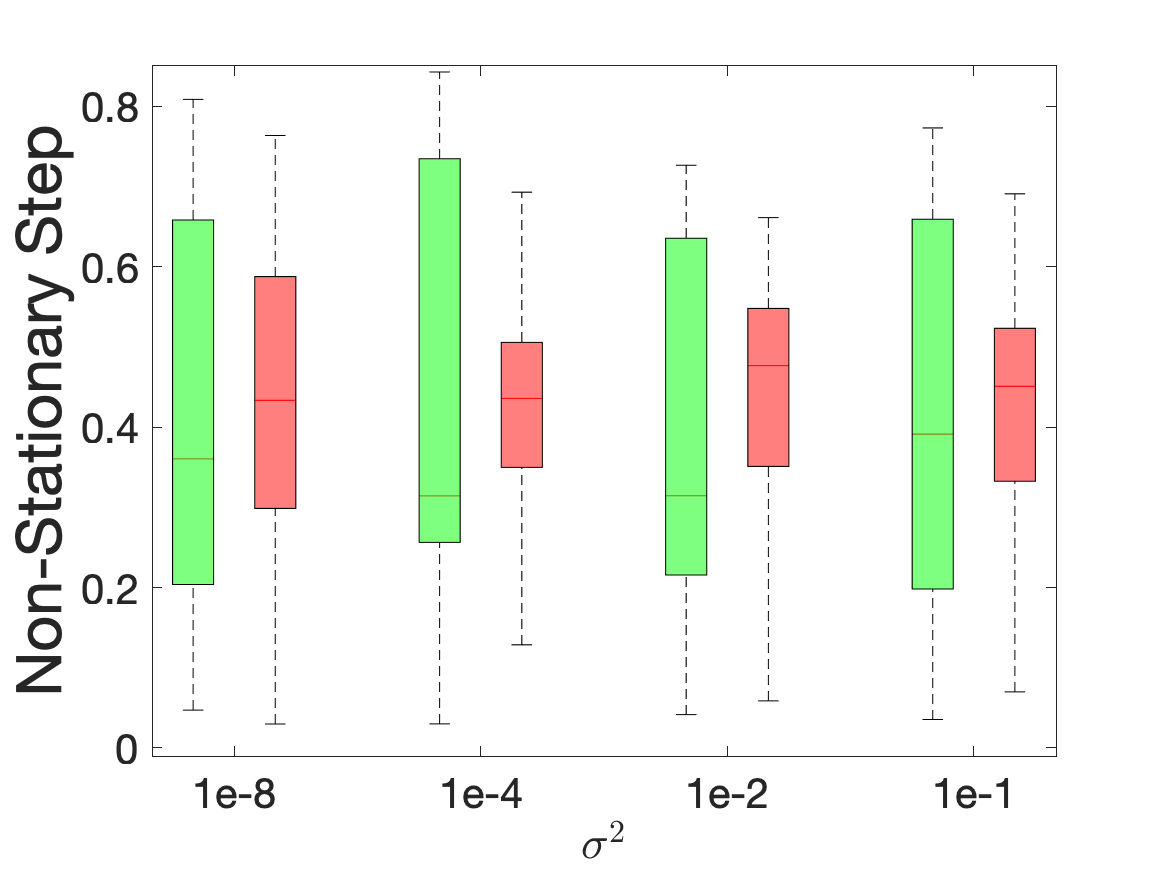}}
	
\subfigure[$\kappa =  2$]{\label{NSK1}\includegraphics[width=37mm]{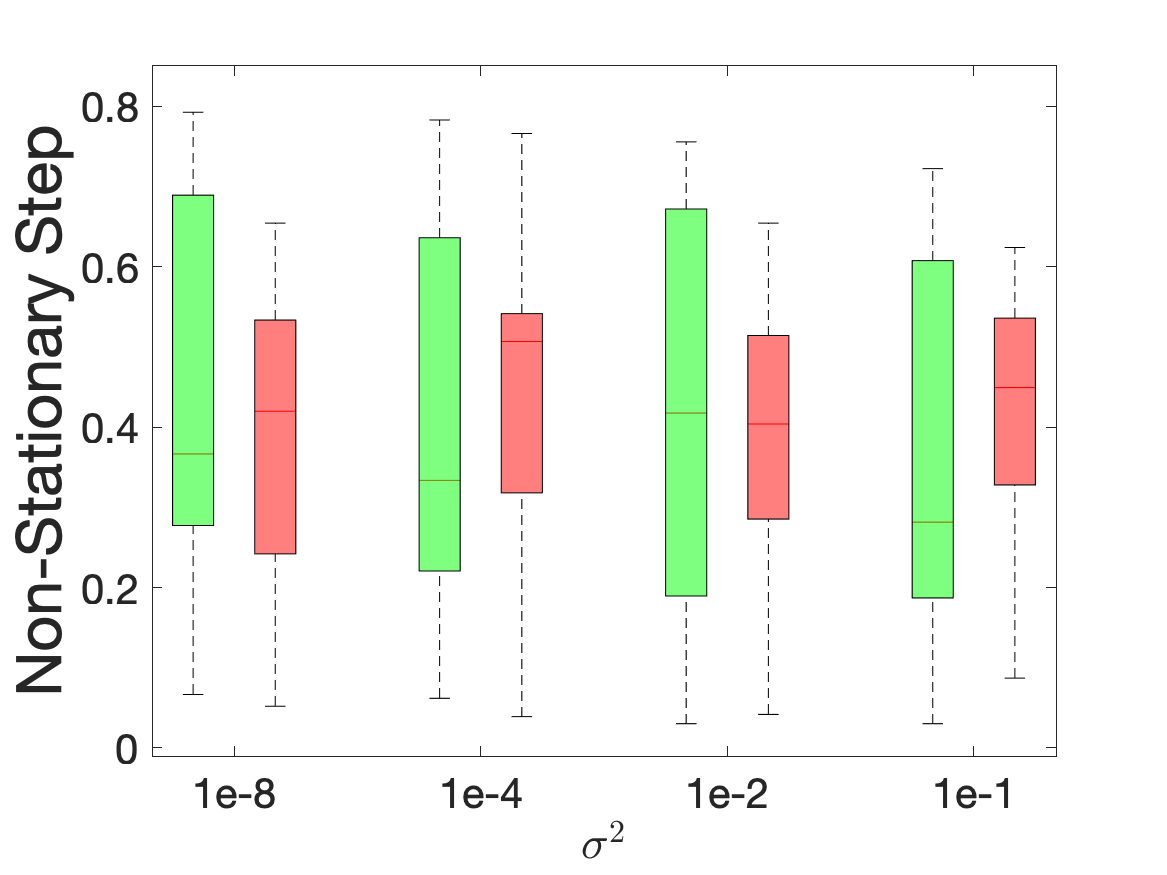}}
\subfigure[$\kappa =  2^3$]{\label{NSK2}\includegraphics[width=37mm]{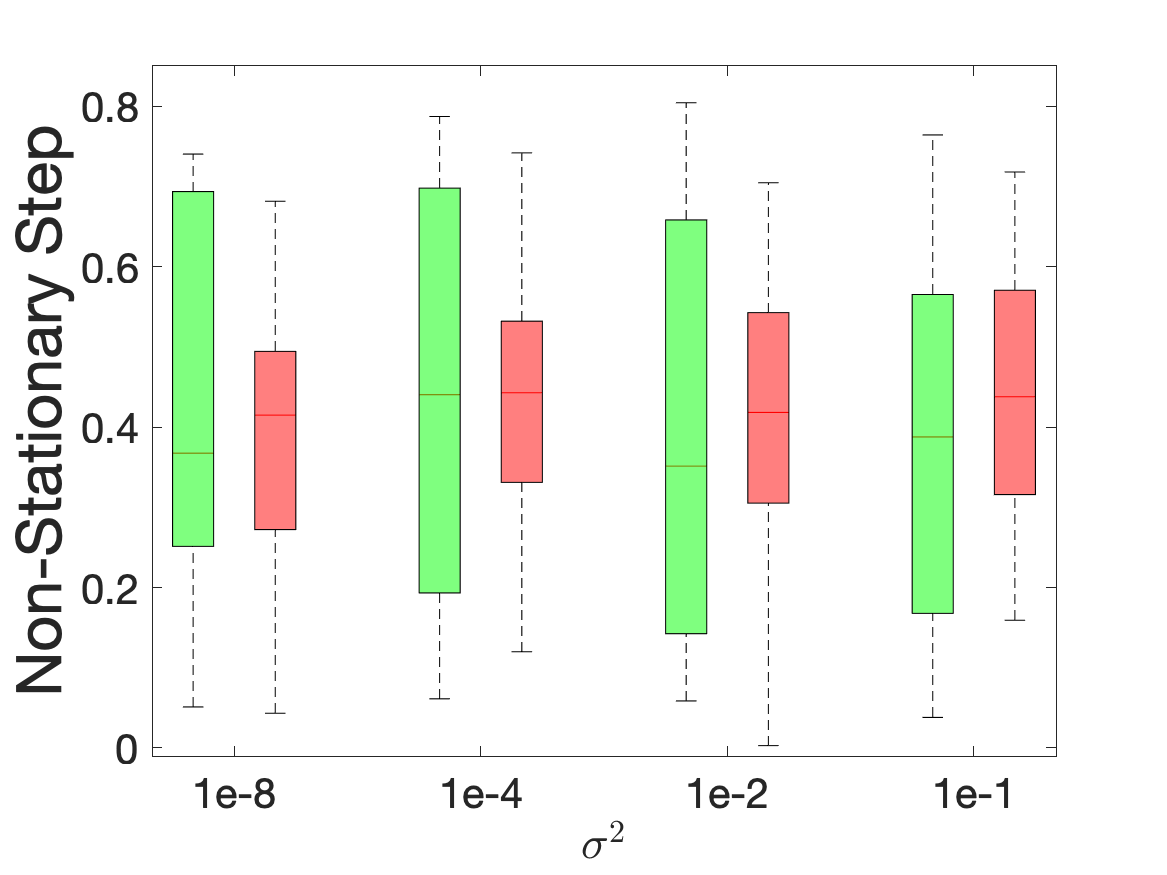}}
\subfigure[$\kappa =  2^6$]{\label{NSK3}\includegraphics[width=37mm]{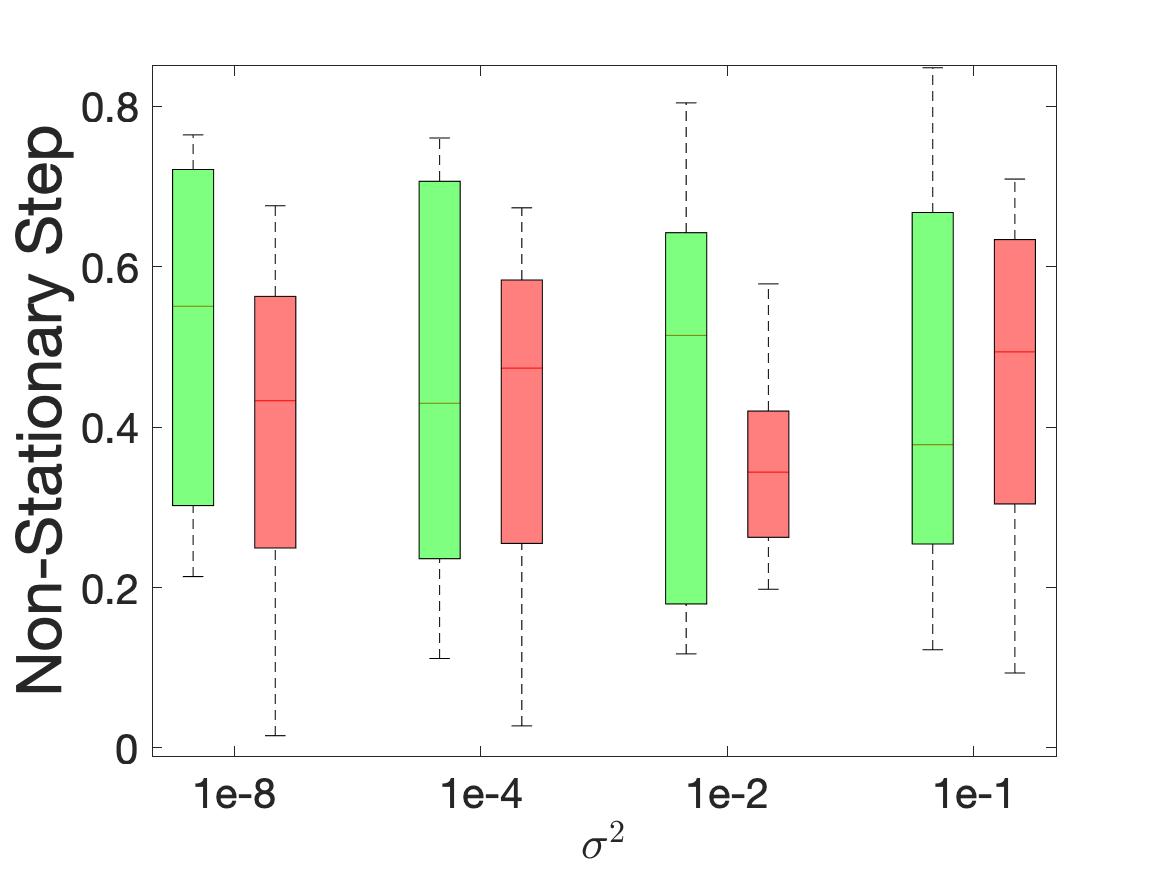}}
	
\subfigure[$\chi_{err} =  1$]{\label{NSChi1}\includegraphics[width=37mm]{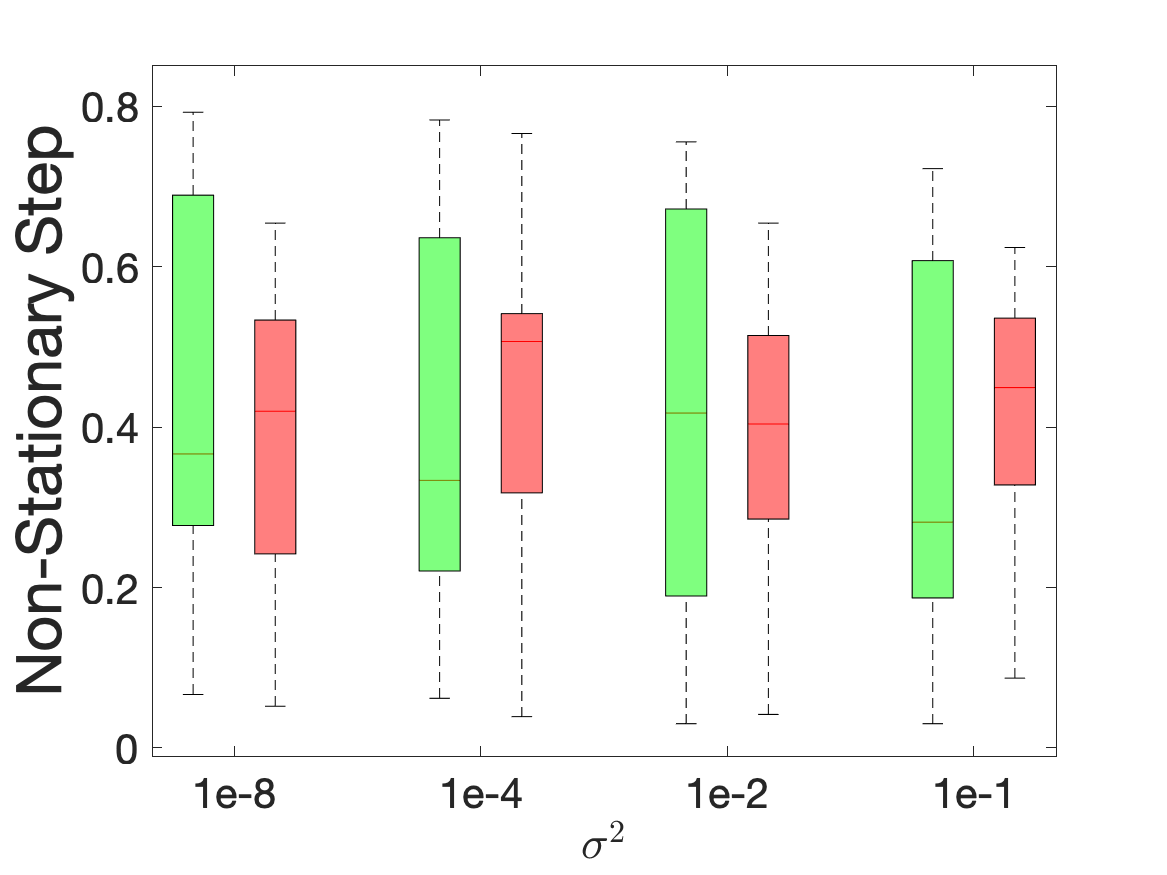}}
\subfigure[$\chi_{err} =  10$]{\label{NSChi2}\includegraphics[width=37mm]{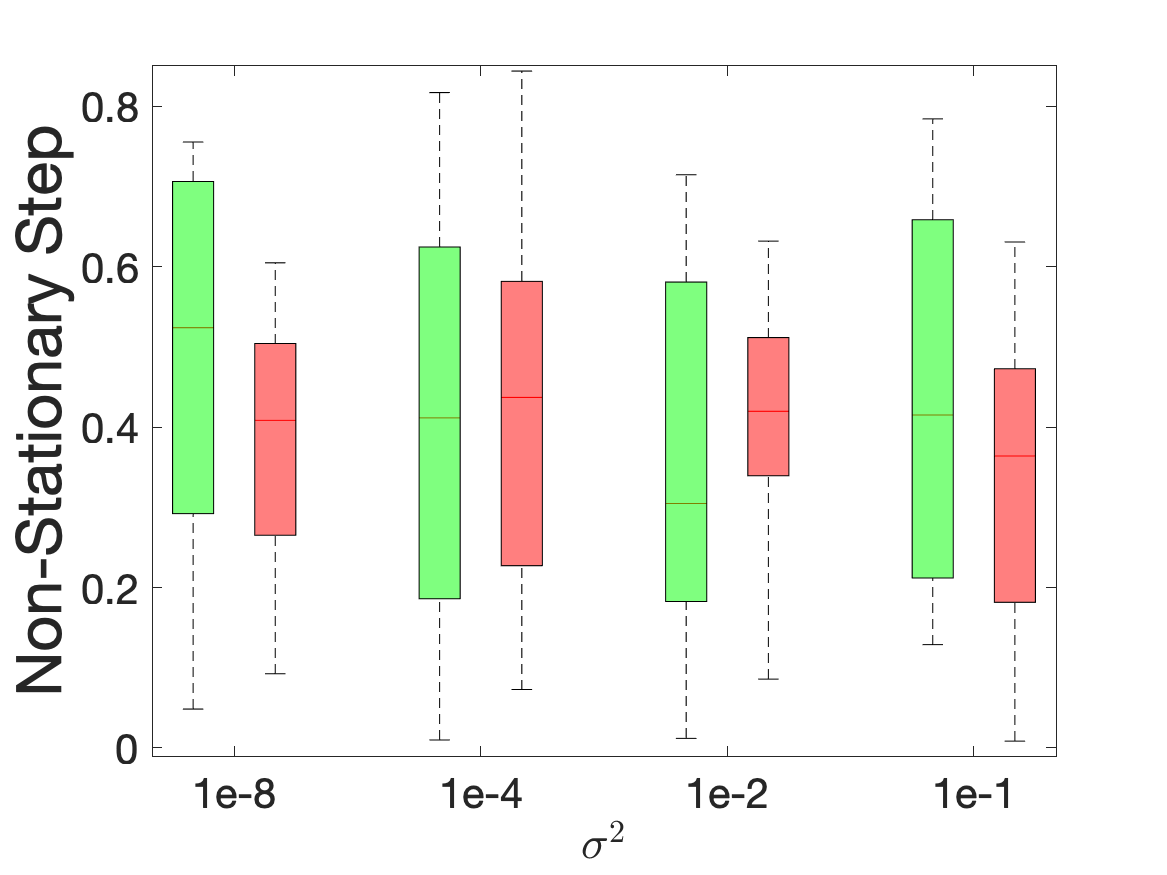}}
\subfigure[$\chi_{err} =  10^2$]{\label{NSChi3}\includegraphics[width=37mm]{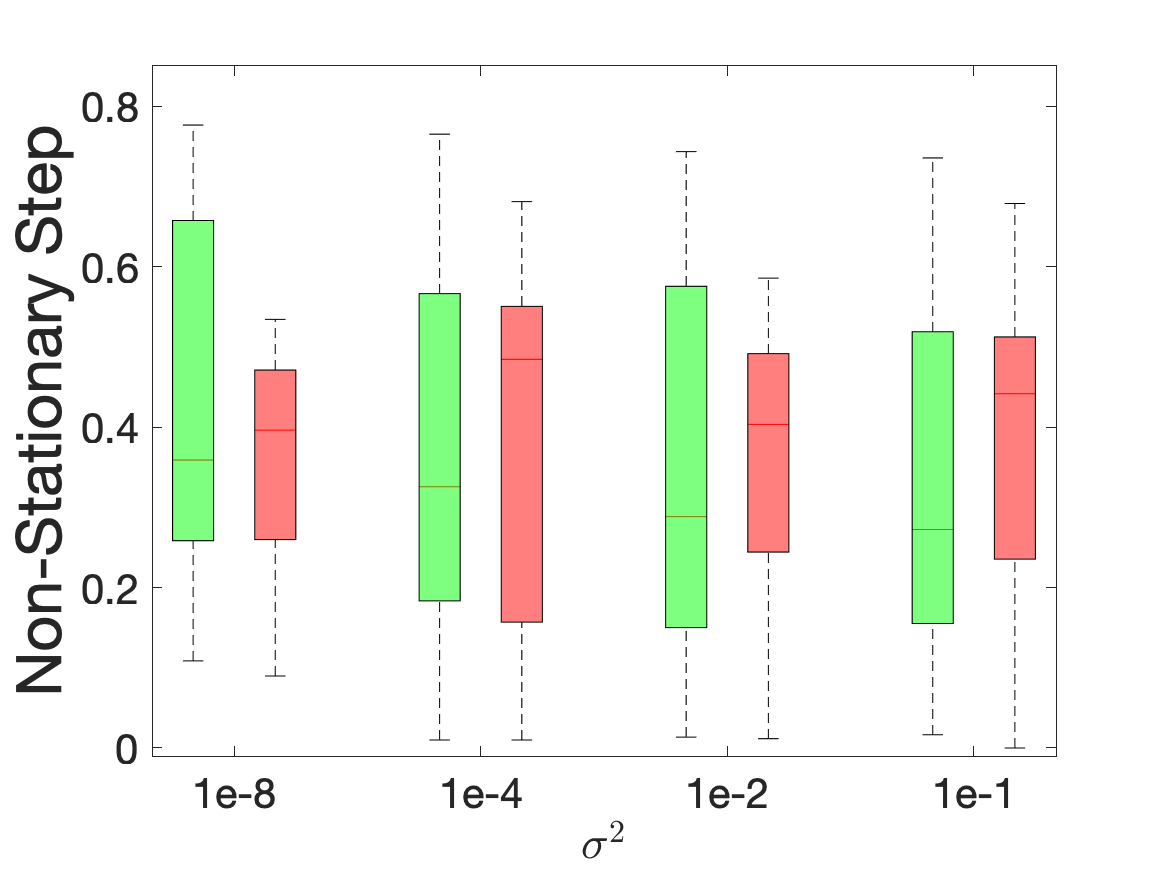}}
	
\includegraphics[width=0.3\textwidth]{Figure/Figure1/legend/KKTlegend}
\caption{Unstabilized penalty parameter boxplots. Each panel corresponds to a setup of $(C, \kappa, \chi_{err})$. The default values are $C =\kappa= 2$ and $\chi_{err}=1$. When we vary one parameter, the other two are set as default. Thus, the three figures on the left column are the same.}\label{fig:6}
\end{figure}

\vskip4pt
\noindent\textbf{Non-stationary penalty parameters.} Figure \ref{fig:6} plots the proportion of the iterations~with unstabilized penalty parameters; i.e., the last iteration that we update $\barepsilon_0$ over the total~number of the iterations. From the figure, we observe that the proportion varies from $20\%$ to $70\%$,~and AdapNewton and AdapGD have comparable results. In fact, the proportion highly depends on the adopted initial $\barepsilon_0$ and the updating rule of $\barepsilon_0$. For example, a large $\rho$ and a small $\barepsilon_0$ will reduce the proportion significantly; and the updating rules $\barepsilon_0\leftarrow\barepsilon_0/\rho$ and $\barepsilon_0\leftarrow \exp(-1/\barepsilon_0)$ will also lead to different proportions. The large variation in Figure \ref{fig:6} suggests that different problems stabilize $\barepsilon_0$ to different levels; thus, a problem-dependent tuning of $\barepsilon_0$ is desired in practice. We note in the experiments that the results on some problems can be improved if $\barepsilon_0=10^{-4}$, while such a setup may not be suitable for other problems. Thus, designing~a~robust scheme to select the penalty parameters deserves further studying.

\begin{figure}[!htp]
\centering     
\subfigure[$C =  2$]{\label{FEC1}\includegraphics[width=37mm]{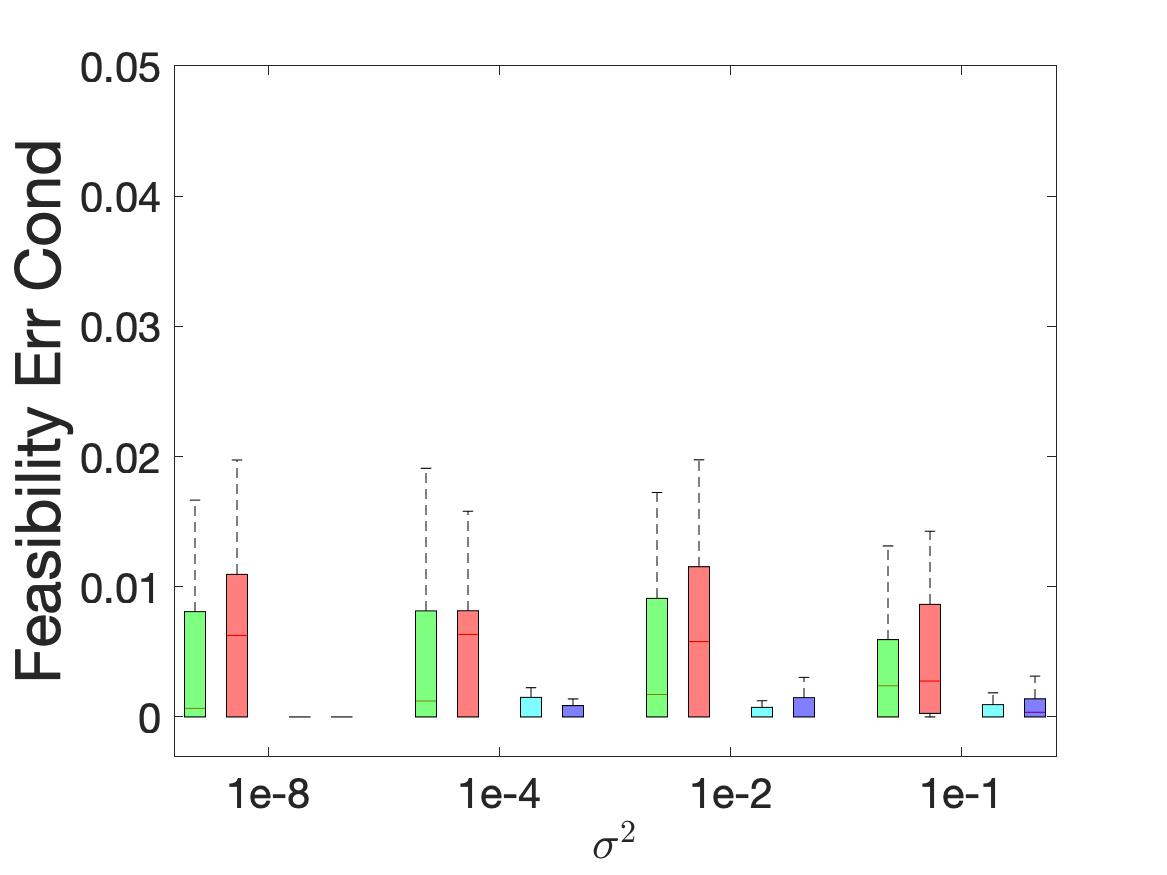}}
\subfigure[$C =  2^3$]{\label{FEC2}\includegraphics[width=37mm]{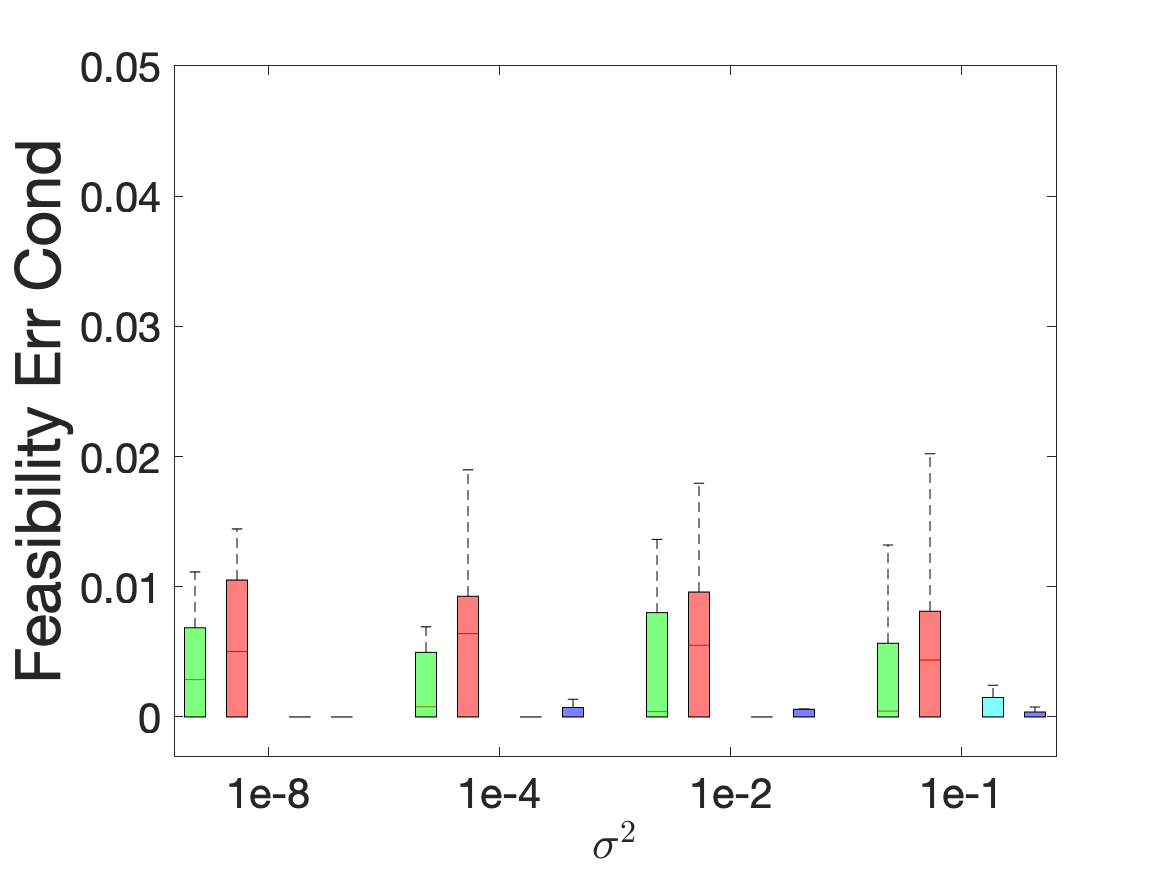}}
\subfigure[$C =  2^6$]{\label{FEC3}\includegraphics[width=37mm]{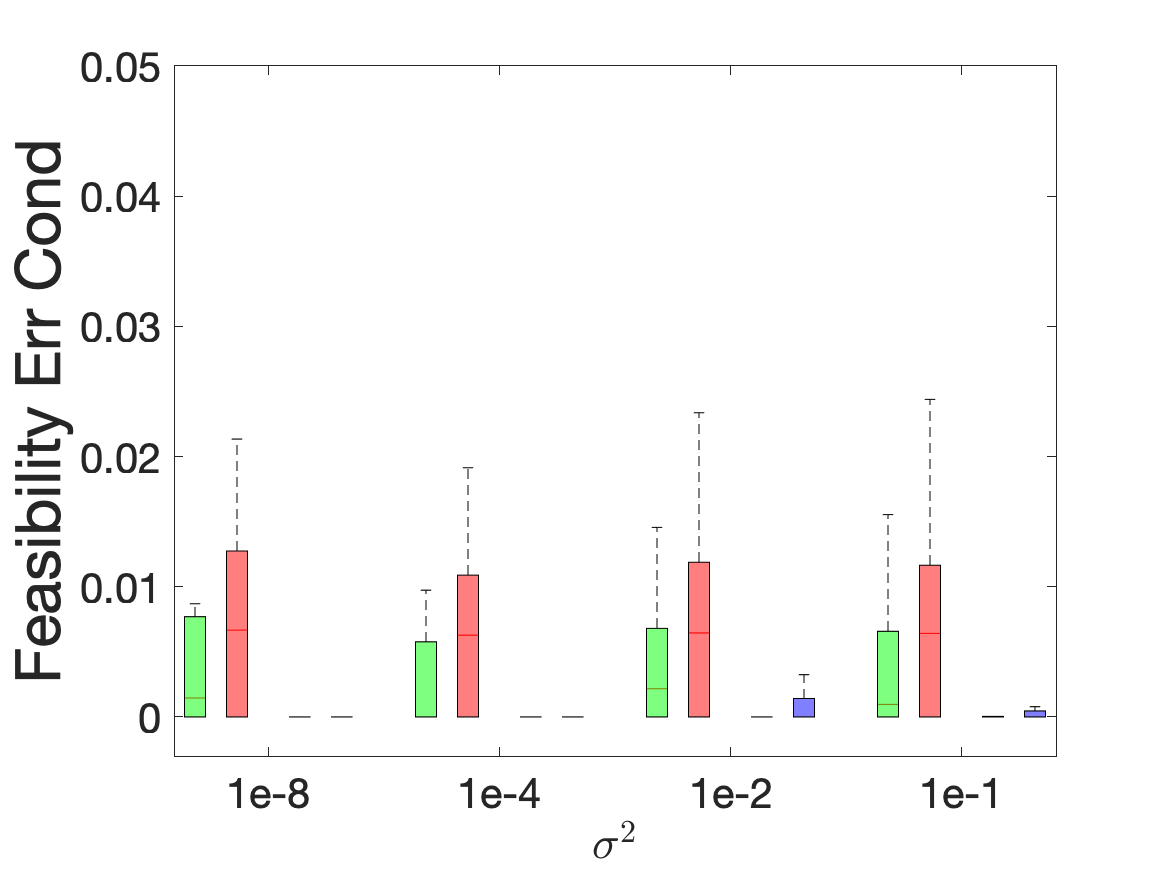}}
	
\subfigure[$\kappa =  2$]{\label{FEK1}\includegraphics[width=37mm]{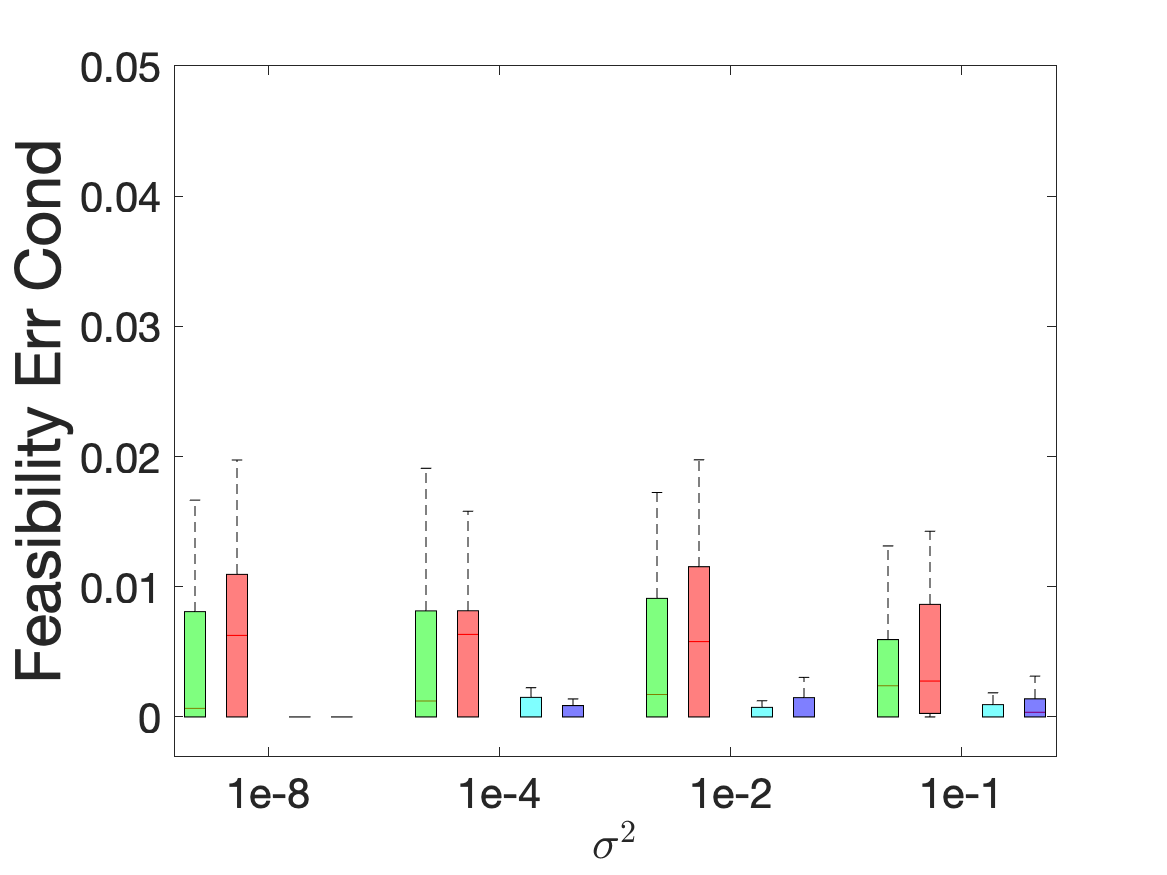}}
\subfigure[$\kappa =  2^3$]{\label{FEK2}\includegraphics[width=37mm]{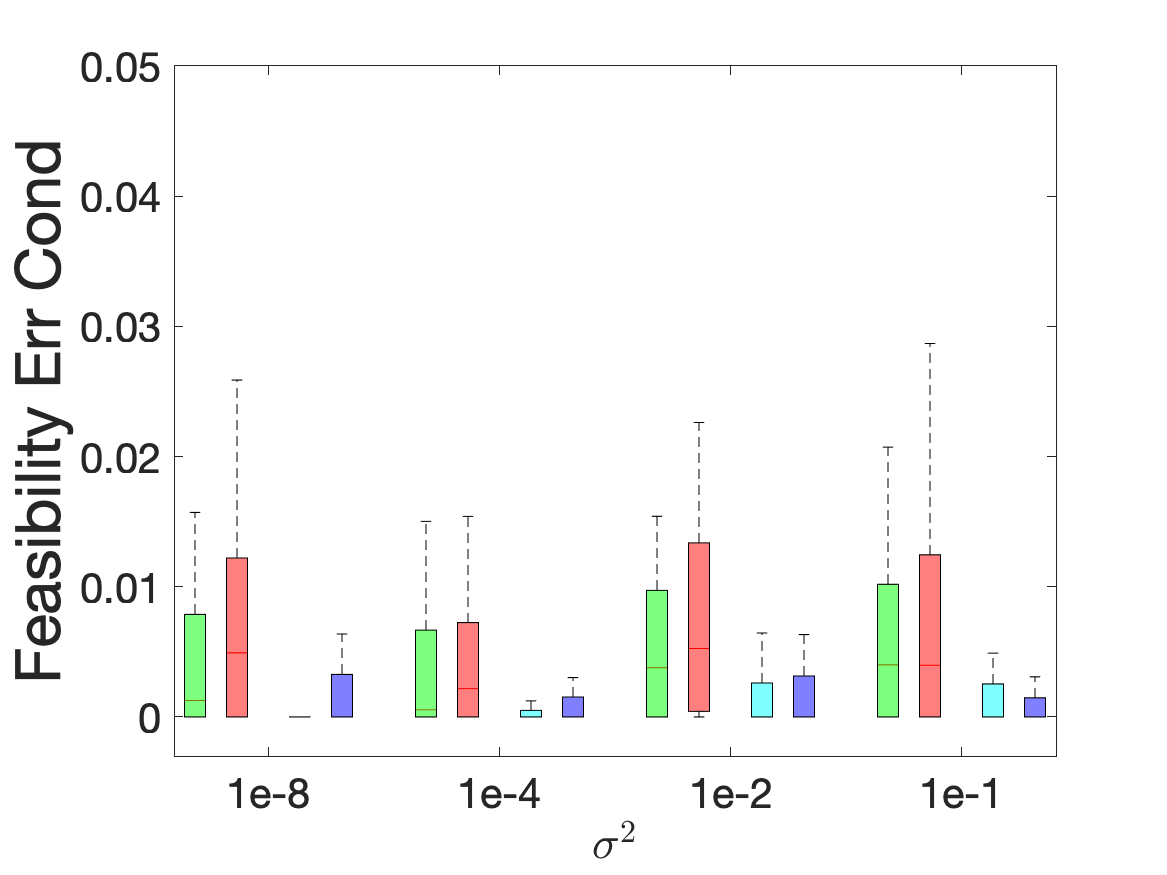}}
\subfigure[$\kappa =  2^6$]{\label{FEK3}\includegraphics[width=37mm]{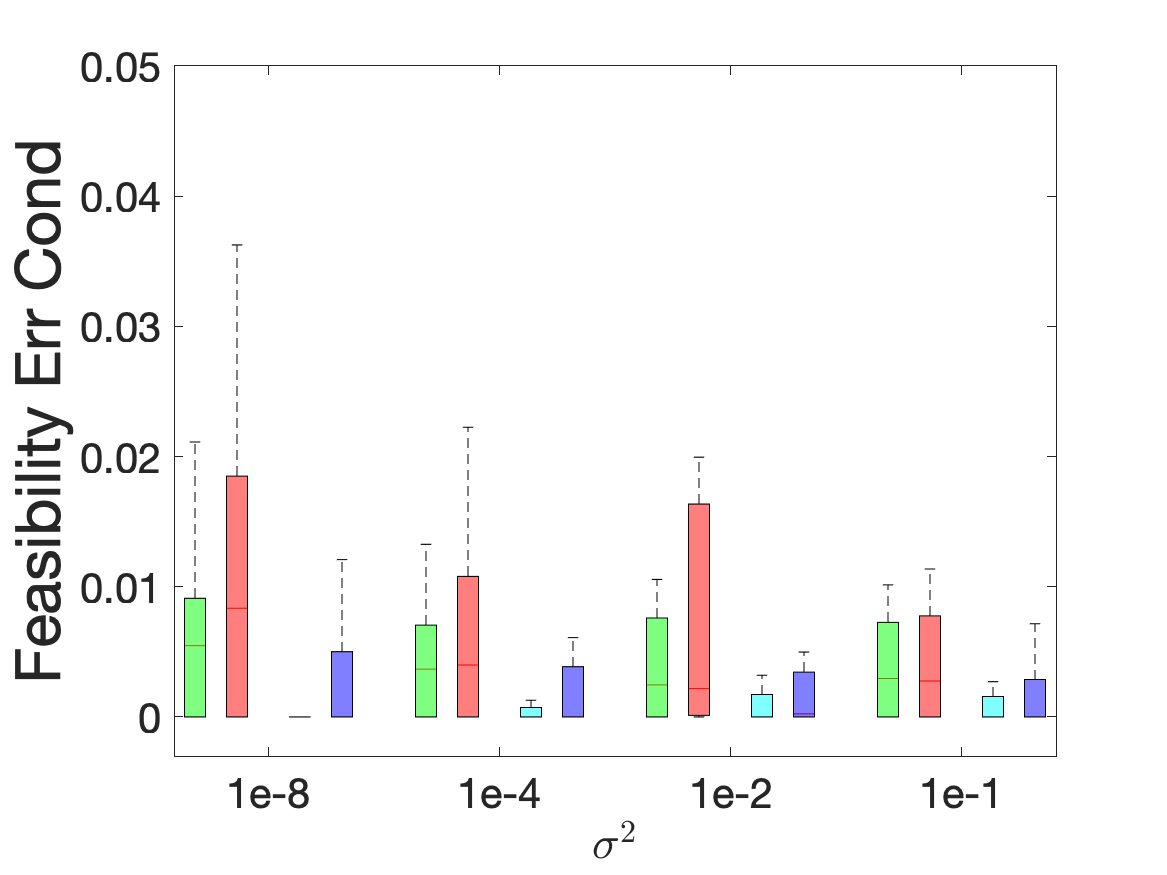}}
	
\includegraphics[width=0.7\textwidth]{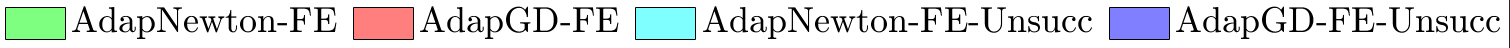}
\caption{Feasibility error condition boxplots. Each panel corresponds to a setup of $(C, \kappa)$. The default values are $C =\kappa= 2$. When we vary one parameter, the other parameter is set as default. Thus, the two figures on the left column are the same.}\label{fig:7}
\end{figure}

\vskip4pt
\noindent\textbf{Feasibility error condition.} Figure \ref{fig:7} plots the proportion of the iterations with a triggered feasibility error condition. We do not show the results for the different setups of $\chi_{err}$.~In~fact, when $\chi_{err}=1$, the results are identical to $C =2$ and $\kappa=2$ (see the left column of Figure~\ref{fig:7}). However, when $\chi_{err} = 10$ or $100$, the feasibility error condition is \textit{never} triggered.~From~Figure \ref{fig:7}, we see that the proportion is extremely small (e.g., as small as $1\%$). This suggests that the condition \eqref{cond:bound:fes:error} is hardly triggered in practice. Figure \ref{fig:7} also plots the iteration proportion that \eqref{cond:bound:fes:error} is triggered for an unsuccessful step. We see that such an proportion is even smaller (e.g., less than $0.5\%$). Given these negligible proportions, we can conclude that the condition \eqref{cond:bound:fes:error} does not negatively affect the performance of the designed StoSQP scheme.

\begin{figure}[!htp]
\centering     
\subfigure[KKT residual]{\label{MNKKT}\includegraphics[width=40.1mm]{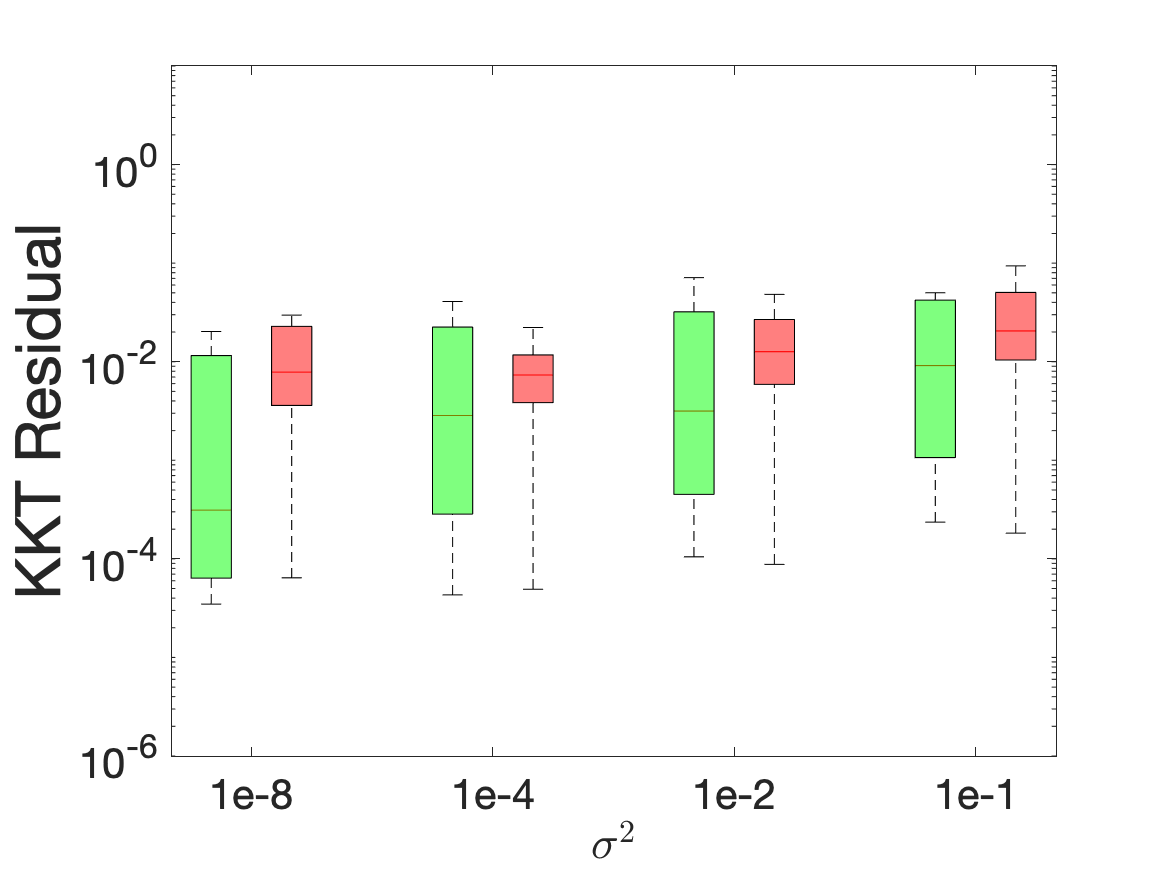}}
\subfigure[Sample size]{\label{MNSam}\includegraphics[width=40.1mm]{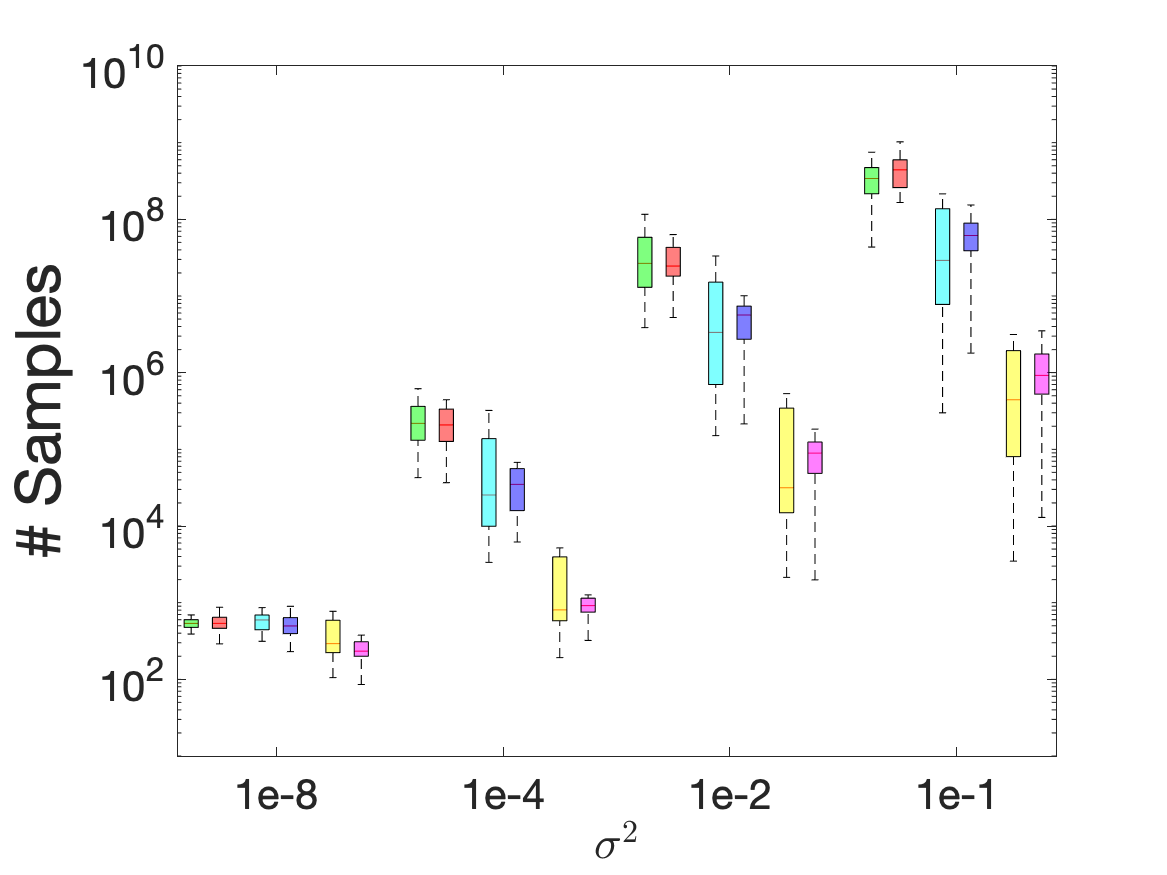}}
	
\includegraphics[width=0.8\textwidth]{Figure/Figure1/legend/Samlegend}
\caption{Multiplicative noise boxplots. The left figure shows the KKT residual boxplot and the right figure shows the sample size boxplot.}\label{fig:8}
\end{figure}

\vskip 4pt
\noindent\textbf{Multiplicative noise.} We also investigate a multiplicative noise in the experiments. In~particular, we employ the default setup $(C, \kappa, \chi_{err}) = (2,2,1)$ but replace the noise variance $\sigma^2$ by $(1 +\|\bx_t\|^2)\sigma^2$. Thus, the variance scales linearly with respect to the magnitude of the (primal) iterate. The KKT residual and sample size boxplots are shown in Figure \ref{fig:8}.~Compared to Figures \ref{fig:1} and \ref{fig:2}, we see that the algorithm achieves comparable results to additive noise. This observation is as expected because, regardless of the noise type, the algorithm enforces the same stochastic conditions on the model estimation accuracy in each~iteration,~and~adaptively selects the batch sizes that are mainly characterized by the current KKT residual.

\vskip4pt

\noindent\textbf{Logistic regression problem.} We study an inequality constrained logistic regression~problem, where we let
\vskip-5pt\begin{equation*}
F(\bx;(\xi_{\ba}, \xi_b)) = \log\{1+\exp(-\xi_b\cdot\xi_{\ba}^T\bx)\},\quad g(\bx) = C\bx + \bq.
\end{equation*}
We set $d = 10, r= 5$, and generate each entry of the matrix $C\in\mR^{5\times 10}$ and vector $q\in\mR^5$ from the standard Gaussian distribution. We let $\xi_b$ be a Rademacher variable (i.e., taking~$\{-1,1\}$ with equal probability), and consider different design distributions for $\xi_{\ba}$. In particular, we consider both a light tail design $(\xi_{\ba})_i\sim \N(\0, \sigma_{\ba}^2)$ and vary $\sigma_{\ba}^2\in\{10^{-8}, 10^{-4}, 10^{-2}\}$, and a heavy tail design $(\xi_{\ba})_i\sim \text{Exp}(\lambda_{\ba})$ and vary $\lambda_{\ba}\in\{10,10^2,10^4\}$. Note that $\text{Exp}(\lambda_{\ba})$ has the variance $1/\lambda_{\ba}^2$. For each design, we run AdapNewton and AdapGD for 20 times. The default algorithm setup is the same as in Section \ref{sec:5}. 

\begin{figure}[!htp]
\centering     
\subfigure[Gaussian design]{\label{KKT1}\includegraphics[width=40.1mm]{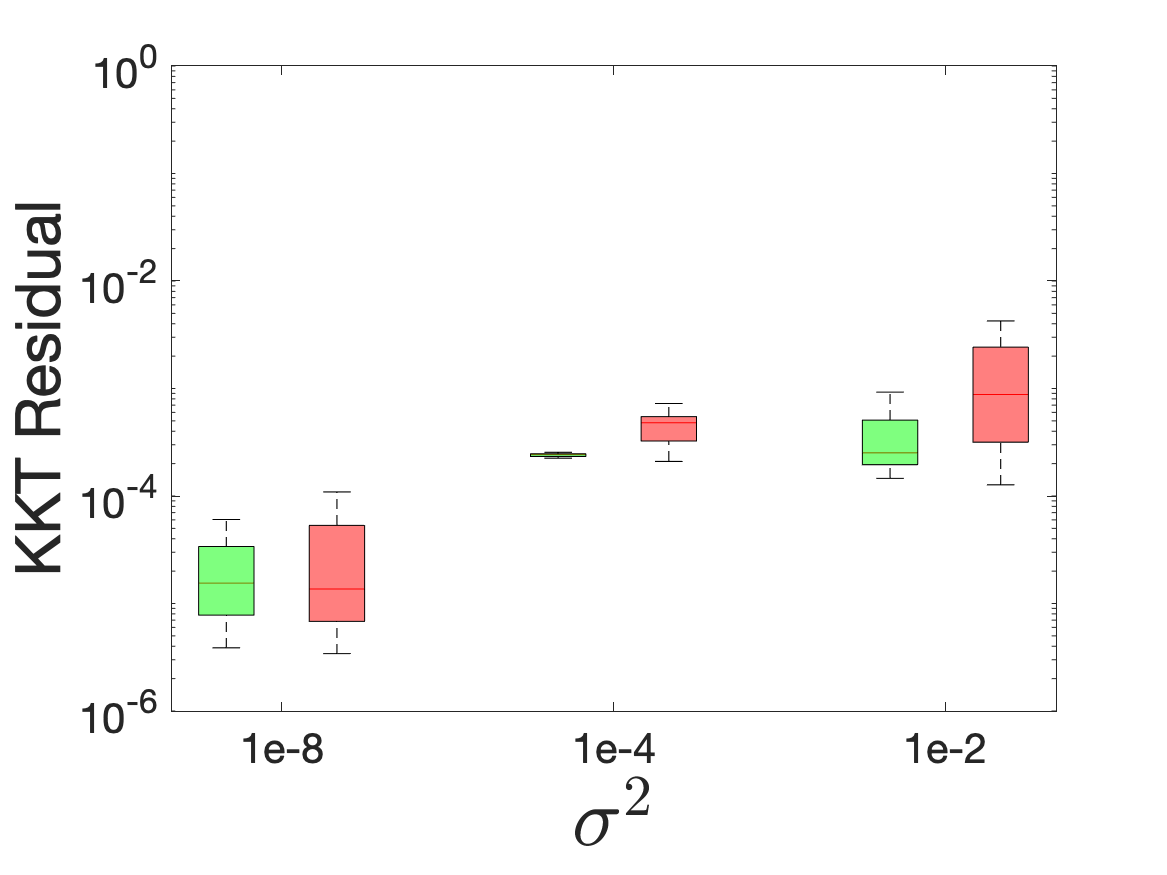}}
\subfigure[Exponential design]{\label{KKT2}\includegraphics[width=40.1mm]{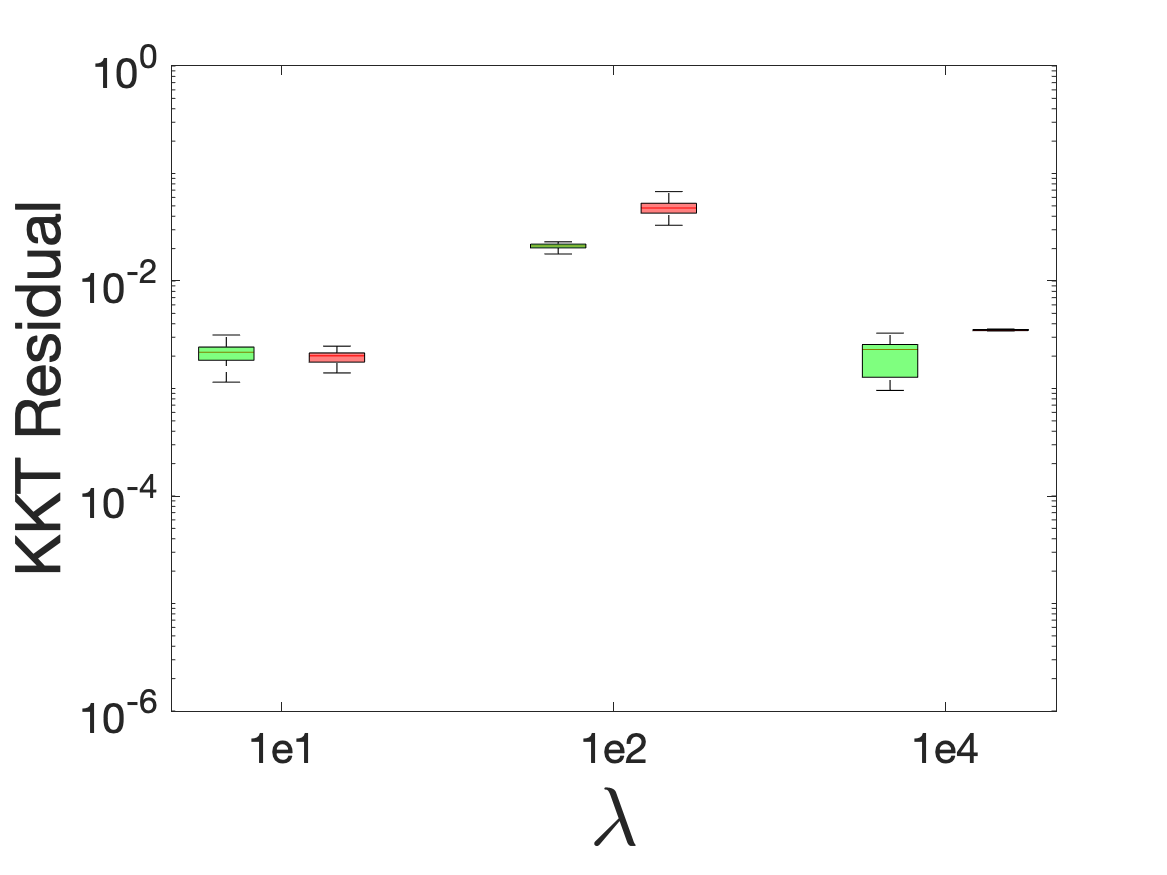}}

\includegraphics[width=0.3\textwidth]{Figure/Figure1/legend/KKTlegend}
\caption{KKT residual boxplots. The left figure shows the residual boxplot for the Gaussian design, and the right figure shows the residual boxplot for the Exponential design.}\label{fig:9}
\end{figure}

Figure \ref{fig:9} shows the KKT residual boxplots. From the figure, we observe that AdapNewton performs slightly better than AdapGD. Both methods achieve reasonable performance on all setups of the two designs, although the two methods perform better on the Gaussian design that has a lighter tail than the Exponential design. Overall, the experiments demonstrate the effectiveness of the proposed algorithm.

%

\bibliographystyle{spbasic}
\bibliography{ref}

\begin{flushright}
\scriptsize \framebox{\parbox{4.5in}{Government License: The submitted manuscript has been created by UChicago Argonne, LLC, Operator of Argonne National Laboratory (``Argonne"). Argonne, a U.S. Department of Energy Office of Science laboratory, is operated under Contract No. DE-AC02-06CH11357.  The U.S. Government retains for itself, and others acting on its behalf, a paid-up nonexclusive, irrevocable worldwide license in said article to reproduce, prepare derivative works, distribute copies to the public, and perform publicly and display publicly, by or on behalf of the Government. The Department of Energy will provide public access to these results of federally sponsored research in accordance with the DOE Public Access Plan. http://energy.gov/downloads/doe-public-access-plan. }}
\normalsize
\end{flushright}	

\end{document}